\documentclass [11pt]{article}
\usepackage{amsmath, amsthm, amssymb}
\usepackage{graphicx}
\newfont{\bbb} {msbm10}
\newcommand{\R}{\Bbb{R}}
\newcommand{\bS}{\Bbb{S}}
\newcommand{\bB}{\Bbb{B}}
\newcommand{\T}{\Bbb{T}}
\newcommand{\C}{\Bbb{C}}
\newcommand{\qH}{\Bbb{H}}
\newcommand{\Oc}{\Bbb{O}}

\newcommand{\Q}{\Bbb{Q}}
\newcommand{\D}{\Bbb{D}}

\newcommand{\sbs}{\subset}
\newcommand{\ra}{\rightarrow}

\newcommand{\bphi}{{\bar{\phi}}}
\newcommand{\br}{{\bar{r}}}
\newcommand{\bt}{{\bar{t}}}
\newcommand{\bg}{{\bar{g}}}
\newcommand{\bq}{{\bar{q}}}

\newcommand{\dsigma}{{\dot{\sigma}}}

\newcommand{\hg}{{\hat{g}}}
\newcommand{\hh}{{\hat{h}}}
\newcommand{\hf}{{\hat{f}}}
\newcommand{\hj}{{\hat{j}}}

\newcommand{\cF}{{\cal{F}}}

\newcommand{\cZ}{{\cal{Z}}}
\newcommand{\cY}{{\cal{Y}}}
\newcommand{\cX}{{\cal{X}}}
\newcommand{\ccP}{{\mbox{{\Large $\wp$}}}}

\newcommand{\p}{\partial}
\newcommand{\barp}{\bar{\partial}}
\newcommand{\pt}{\frac{\partial}{\partial t}}
\newcommand{\n}{\nabla}

\newcommand{\cL}{{\cal{L}}}

\newcommand{\cD}{{\cal{D}}}
\newcommand{\cS}{{\cal{S}}}
\newcommand{\cH}{{\cal{H}}}
\newcommand{\cC}{{\cal{C}}}
\newcommand{\cG}{{\cal{G}}}
\newcommand{\cA}{{\cal{A}}}
\newcommand{\cE}{{\cal{E}}}
\newcommand{\cB}{{\cal{B}}}
\newcommand{\cW}{{\cal{W}}}
\newcommand{\cN}{{\cal{N}}}

\newcommand{\B}{\Bbb{B}}

\newcommand{\HH}{\Bbb{H}}

\newcommand{\ua}{\underline{a}}
\newcommand{\ub}{\underline{b}}

\newcommand{\ssl}{{_\lambda}}
\newcommand{\sd}[1]{\Delta_{_{\bS^{#1}}}}

\newcommand{\s}[1]{{\sf{#1}}}

\newcommand{\0}[1]{_{_{#1}}}
\newcommand{\rC}{{\rm{C}}\,}

\newcommand{\sN}{{\sf{N}}}
\newcommand{\sL}{{\sf{Link}}}

\newcommand{\sB}{{\sf{B}}}
\newcommand{\sA}{{\sf{A}}}
\newcommand{\bs}{{\mbox{\bf s}}}
\newcommand{\dX}{\dot{X}}
\newcommand{\dDelta}{\dot{\Delta}}
\newcommand{\dsquare}{\dot{\square}}
\usepackage[all]{xy}

\begin{document}

\title{Pinched Smooth Hyperbolization}
\author{Pedro Ontaneda\thanks{The author was
partially supported by a NSF grant.}}
\date{}

\maketitle

\vspace{.3in}

\noindent {\bf \large  Section 0. Introduction.}\\

Classical flat geometry has formed part of basic human knowledge since ancient times. It is characterized by the almost universally known condition that the sum of the internal angles of a triangle $\triangle$ is equal to $\pi$. We write $\Sigma(\triangle)=\pi$.
Other fundamental geometries are defined by replacing the equality
$\Sigma(\triangle)=\pi$ by inequalities, thus positively curved geometries and
negatively curved geometries are determined by the inequalities
 $\Sigma(\triangle)>\pi$ and $\Sigma(\triangle)<\pi$, respectively, where
$\triangle$ runs over all small non-degenerate triangles in a space.
It is natural then to try to find spaces that admit such geometries, and this
task has been a driving force in Riemannian Geometry for many decades.
But surprisingly there are not too many examples of smooth closed manifolds
that support either a positively curved or a negatively curved metric.
For instance, besides spheres, in dimensions $\geq 17$ (and $\neq 24$) the only positively curved simply connected known
examples are complex and quaternionic projective spaces. In negative
curvature the situation is arguably more striking because negative curvature
has been studied extensively in many different areas in mathematics.
Indeed,  from the ergodicity of their geodesic flow in
Dynamical Systems to their topological rigidity in Geometric Topology;
from the existence of harmonic maps in Geometric Analysis to the
well-studied and greatly generalized algebraic properties of their fundamental groups,
negatively curved Riemannian manifolds are the main object in many 
important and well-known results in mathematics. Yet the fact remains that very
few examples of closed negatively curved Riemannian manifolds are known.
Besides the hyperbolic ones ($\R$, $\C$, $\qH$, $\Oc$), the other known examples are  the Mostow-Siu examples (complex dimension 2) which are local branched covers of complex hyperbolic space (1980, \cite{MS}), 
the Gromov-Thurston examples (1987, \cite{GT}) which are branched covers of  real hyperbolic
ones, the exotic Farrell-Jones examples (1989, \cite{FJ1}) which are homeomorphic but not diffeomorphic to real hyperbolic manifolds
(and there are other examples of exotic type), and the three
examples of Deraux (2005, \cite{Deraux}) which are of the Mostow-Siu
type in complex dimension 3. Hence, excluding the Mostow-Siu and
Deraux examples (in dimensions 4 and 6, respectively), all known examples
of closed negatively curved Riemannian manifolds are homeomorphic to
either to a hyperbolic one or a branched cover of a hyperbolic one.\\

\noindent {\bf Remark.} The Mostow-Siu and Deraux examples have
covers that have  maps
to complex hyperbolic space that look locally like branched covers.
It is not known whether all these examples are global branched covers.\\

This lack of examples in negative curvature changes dramatically if we allow singularities,
and a very rich and abundant class of negatively curved spaces (in the geodesic sense)
exists due to the strict hyperbolization process of Charney and Davis \cite{ChD}. The
hyperbolization process was originally introduced by M. Gromov \cite{G}, and later
studied by Davis and  Januszkiewicz \cite{DJ}, and Charney-Davis strict hyperbolization is built on these previous versions. The hyperbolization process is conceptually (but not technically) quite
simple since it has a lego type flavor: in the same way as simplicial complexes
and cubical complexes are built from a basic set of pieces,
basic ``hyperbolization pieces" are chosen, and anything that can be built or assembled 
with these pieces will be negatively curved. This conceptual simplicity could be in some sense a bit deceptive because hyperbolization produces an enormous class of examples
with a very fertile set of properties. But the richness and complexity of the hyperbolized
objects is matched by the richness and complexity of the singularities obtained,
and hyperbolized smooth manifolds are very far from being Riemannian.
Interestingly one can relax and lose even
more regularity and consider negative curvature from the algebraic point of
view, that is consider Gromov's hyperbolic groups, and it can be argued
\cite{Ol} that ``almost every group" is hyperbolic. So, negative curvature
is in some weak sense generic, but Riemannian negative curvature is very scarce.
It is natural then to inquire about the difference between the class of manifolds
with negatively curved metrics with singularities and its subclass of more
regular Riemannian counterparts. More specifically we can ask whether
the strict hyperbolization process can be brought into the Riemannian universe. 
In this paper we give a positive answer to this question, and we do this
by proving that all singularities of the Charney-Davis strict hyperbolization of
a closed smooth manifold can be smoothed, provided the ``hyperbolization
piece" is large enough (which can always be done). Moreover we prove that we can do this process in a $\epsilon$-pinched way. Here is the statement of our Main Theorem.\\

\noindent {\bf Main Theorem.} {\it Let $M^n$ be a closed smooth manifold
and let $\epsilon>0$. Then there is a closed Riemannian manifold $N^n$ and
a smooth map $f:N\ra M$ such that}
\begin{enumerate}
\item[{\it (i)}] {\it The Riemannian manifold $N$ has sectional curvatures in the interval
$[-1-\epsilon,-1]$.}
\item[{\it (ii)}] {\it The induced map  $f_*:H_*(N,R)\ra H_*(M,R)$
is surjective, for every (untwisted) $R$.}
\item[{\it (iii)}] {\it If $M$ is $R$-orientable then $N$ is $R$-orientable. In this case
(ii) implies that $f$ has degree one and
$f^*:H^*(M,R)\ra H^*(N,R)$ is injective.}
\item[{\it (iv)}] {\it The map $f^*$ sends the rational Pontryagin classes of $M$ to the
rational Pontryagin classes of $N$. }
\end{enumerate}
\vspace{.2in}

\noindent {\bf Addendum to Main Theorem.} {\it The manifold $N$ is
the Charney-Davis strict hyperbolization of $M$ but with a different smooth structure. The hyperbolization is done with a sufficiently
``large" hyperbolization piece $X$.}\vspace{.2in}

By ``large" above we mean that the width of the normal neighborhoods of the faces
of $X$ are very large. These large pieces always exist (see section 9.1).
Corollaries 1, 2 and 3 below are the $\epsilon$-pinched Riemannian versions of classical applications
of hyperbolization.\\

\noindent {\bf Corollary 1.} {\it Every closed smooth manifold is smoothly cobordant
to a closed Riemannian manifold with sectional curvatures in the interval
$[-1-\epsilon,-1]$, for every $\epsilon>0$.}\\

\noindent {\bf Corollary 2.} {\it The cohomology ring of any finite $CW$-complex 
embeds in the cohomology ring of closed  Riemannian manifold
with sectional curvatures in the interval
$[-1-\epsilon,-1]$, for every $\epsilon>0$.}\\

\noindent {\bf Proof.} Let $X$ be a finite $CW$-complex. Embed $X$ in some
$\R^n$ and let $P$ be a compact neighborhood of $X$ that retracts to $X$. Let $M$ be the double of $P$. Then there is a retraction $M\ra X$, and Corollary 2 follows from (iii) in the Main Theorem.\\

Since degree one maps between closed orientable manifolds are $\pi_1$-
surjective we obtain the following result.\\ 

\noindent {\bf Corollary 3.} {\it For every finite CW-complex $X$
there is a closed Riemannian manifold $N$ and a map $f:N\ra X$ such
that: (i) $N$ has sectional curvatures in the interval
$[-1-\epsilon,-1]$, (ii) $f$ is $\pi_1$-surjecive, (iii) $f$ is homology
surjective.}\\

All known examples of closed negatively curved Riemannian manifolds
with less than $\frac{1}{4}$ pinched curvature have zero rational Pontryagin classes
(for the Gromov-Thurston branched cover examples this was proved by
S. Ardanza \cite{A}). 
The next corollary gives examples of such manifolds with
nonzero rational Pontryagin classes. \\

\noindent {\bf Corollary 4.} {\it For every $\epsilon>0$ and $n\geq 4$ there
is a closed Riemannian $n$-manifold with sectional curvatures in the interval
$[-1-\epsilon,-1]$ and nonzero rational Pontryagin classes.}\\

\noindent {\bf Proof.} Take $M$ in the Main Theorem orientable with nonzero
Pontryagin classes.\\

All manifolds given in Corollary 4 are new examples of closed negatively
curved manifolds.
This follows from Novikov's topological invariance of the rational Pontryagin classes \cite{N}, and the $\frac{1}{4}$-pinched rigidity results given in
(or implied by)
the work
of Hern\'andez \cite{Her}, Yau and Zheng \cite{YZ}, Corlette \cite{Cor},
Gromov \cite{G3} and Mok-Siu-Yeung \cite{MSY}. We state this
in the next corollary.\\

\noindent {\bf Corollary 5.} {\it For any $\epsilon >0$ and $n\geq 4$ there are
closed Riemannian $n$-manifolds with sectional curvatures in the interval
$[-1-\epsilon,-1]$ that are neither homeomorphic to a hyperbolic manifold
($\R$, $\C$, $\qH$, $\Oc$) nor homeomorphic to the Gromov-Thurston branched cover of a real hyperbolic one, nor homeomorphic to one of the Mostow-Siu or Deraux examples.}\\

The next application was suggested to us by Stratos Prassidis some time ago and deals with cusps of negatively curved manifolds.
Recall that if $M$ is a complete finite volume noncompact real hyperbolic manifold
then $M$ has finitely many cusps isometric to manifolds of the form $Q\times [b,\infty)$
with metric $e^{-2t}h+dt^2$, where $(Q, h)$ is a closed flat manifold
and $b\in \R$. If $M$ has exactly one cusp diffeomorphic to
$Q\times [b,\infty)$
 we say that {\it the manifold
$Q$ bounds geometrically a hyperbolic manifold.} \\

More generally, in 1978 M. Gromov defined almost flat manifolds in \cite{G2} and
similar facts hold for them replacing hyperbolic manifolds by pinched
negatively curved manifolds. That is,
let $M$ be a complete finite volume noncompact manifold with pinched negative 
curvature (i.e all sectional curvatures lie in a fixed interval $[-a,-b]$, $0<b\leq a<\infty$).
Then $M$ has finitely many cusps diffeomorphic to manifolds of the form $Q\times [b,\infty)$, 
where $Q$ is an almost flat manifold. 
If $M$ has exactly one cusp diffeomorphic to 
 $Q\times [b,\infty)$
we say that {\it the manifold
$Q$ bounds geometrically a negatively curved manifold.}
Of course a necessary condition for $Q$ to bound geometrically as above is to
smoothly bound a compact manifold. \\

\noindent {\bf Remark.} Here we do not assume $Q$ to be connected.
Hence if $Q$ bounds geometrically the number of connected components of $Q$ is the same as the number of connected cusps.\\

It was proved by Hamrick and Royster \cite{HR} that every closed flat manifold bounds smoothly. This together with the work of Gromov in \cite{G1}, \cite{G2} motivated
Tom Farrell and Smilka Zdravkovska to make the following well-known 
conjectures in \cite{FZ}  30 years ago.\\

\noindent {\bf Conjecture 1.} Every closed almost flat manifold bounds smoothly. This conjecture was also proposed, independently, by S.-T. Yau in \cite{Y}.\\

\noindent {\bf Conjecture 2.} Every closed
flat manifold bounds geometrically a hyperbolic manifold.\\

\noindent {\bf Conjecture 3.} Every closed almost flat manifold bounds geometrically a negatively curved manifold.\\

It was showed by Long and Reid \cite{LR} that Conjecture 2 is false
by giving examples of three dimensional flat manifolds that do not bound.
The following result says that Conjecture 1 implies Conjecture 3.\\

\noindent{\bf Theorem A.} {\it Let $Q$ be a closed almost flat manifold. Assume that
$Q$ bounds smoothly. Then $Q$ bounds geometrically a negatively curved manifold $M$.}\\

Conjecture 1 has generated a lot of research in the last 30 years and it is
known to be true for an almost flat manifold in the following cases.
Let $Q$ be almost flat. Then $Q$ is covered by a nilmanifold, that is,
the quotient of a simply connected nilpotent Lie group $L$ by a uniform lattice. Denote by $G$ the holonomy of $Q$.\\

\begin{enumerate}
\item[(a)] The manifold $Q$ is a nilmanifold.
\item[(b)]  The holonomy $G$ has order $k$ or $2k$, where $k$ is odd,
due to Farrell-Zdravkovska \cite{FZ}.
\item[(c)] The holonomy $G$ of $Q$ acts effectively on the center of $L$, also due to
Farrell-Zdravkovska \cite{FZ}.
\item[(d)] The holonomy $G$ is cyclic, due to J. Davis and F. Fang.
Also  Upadhyay \cite{U} had proved that Conjecture 1 it true when
the following conditions hold: $G$ is cyclic, $G$ acts trivially on the center of $L$, and $L$ is 2-step nilpotent. \end{enumerate}

Hence in all of the above cases $Q$ bounds geometrically a pinched
negatively curved manifold. Note that for any closed $Q$ we have
$\p ( Q\times I)=Q\coprod Q$. Thus we get the following corollary of Theorem A.\\

\noindent {\bf Corollary 6.} {\it Let $Q$ be a closed almost flat manifold.
Then $Q\coprod Q$ bounds geometrically a pinched negatively curved
manifold.}\\

In other words, for every closed connected almost flat manifold there is
a complete finite volume pinched negatively curved manifold with
exactly {\sf two} connected cusps, each diffeomorphic to $Q\times [b,\infty)$.\\

A complete pinched negatively curved metric $g$ on $Q\times \R$ is called a {\it (pinched negatively curved) cusp metric} if the $g$-volume of $Q\times [0,\infty)$
is finite. And we say that a cusp
metric $g$ on $Q\times \R$ {\it is an eventually warped cusp metric} if $g=e^{-2t}h+dt^2$, for $t<c$, for some $c\in \R$ and a metric $h$ on $Q$.
I. Belegradek and V. Kapovich \cite{BK} (see also \cite{B}) show, based on earlier work by Z.M. Shen \cite{Shen}, that if $Q$ is almost flat then $Q\times\R$ admits an eventually warped cusp metric.\\

\noindent {\bf Addendum to Theorem A.} {\it 
Let $g$ be an eventually warped cusp metric on $Q\times\R$. If the sectional curvatures
of $g$ lie in $(a,b)$, with $a<-1<b$, then we can take $M$ in Theorem A with sectional curvatures also in $(a,b)$. Moreover the sectional curvatures of $M$ away from a cusp can be taken in $[-\epsilon-1,-1]$, for any $\epsilon>0$.}\\

Even though a flat manifold  may not necessarily bound geometrically a hyperbolic
manifold the next corollary says it does bound geometrically an $\epsilon$-pinched to -1
manifold, for any $\epsilon>0$. It follows from the Hamrick and Roster result \cite{HR},
Theorem A and its addendum.\\

\noindent {\bf Corollary 7.} {\it Every closed flat manifold
bounds geometrically a manifold with sectional curvatures in
$[-\epsilon-1,-1]$, for any $\epsilon>0$.}\\\\

We next give a rough and heuristic idea of some of the methods used in smoothing the singularities of
a Charney-Davis hyperbolized smooth manifold. We do this first in dimension two and then in dimension three
where we can visualize some of these methods. We have several items.\\

\noindent {\bf 0.1. Dimension two.}
A Charney-Davis hyperbolization piece $X^n$ of dimension $n$ is essentially a compact hyperbolic manifold with corners that has the symmetries of an $n$-cube, and all ``faces" intersect perpendicularly (see 9.1). We shall assume throughout this introduction
that $X$ is as ``large" as we need it to be (see 9.1).\\

Fix an $X^2$ and let $K$ be a cubical 2-complex.
Replace each cube by a copy of $X^2$ to obtain a piecewise hyperbolic space $K_X$. This is (essentially again) the Charney-Davis hyperbolization of $K$  (see figure 1).
We shall identify the vertices of $K$ with the vertices of $K_X$.

\begin{center}
\begin{figure}[ht]
\centering
\includegraphics[width=9cm,height=6cm]{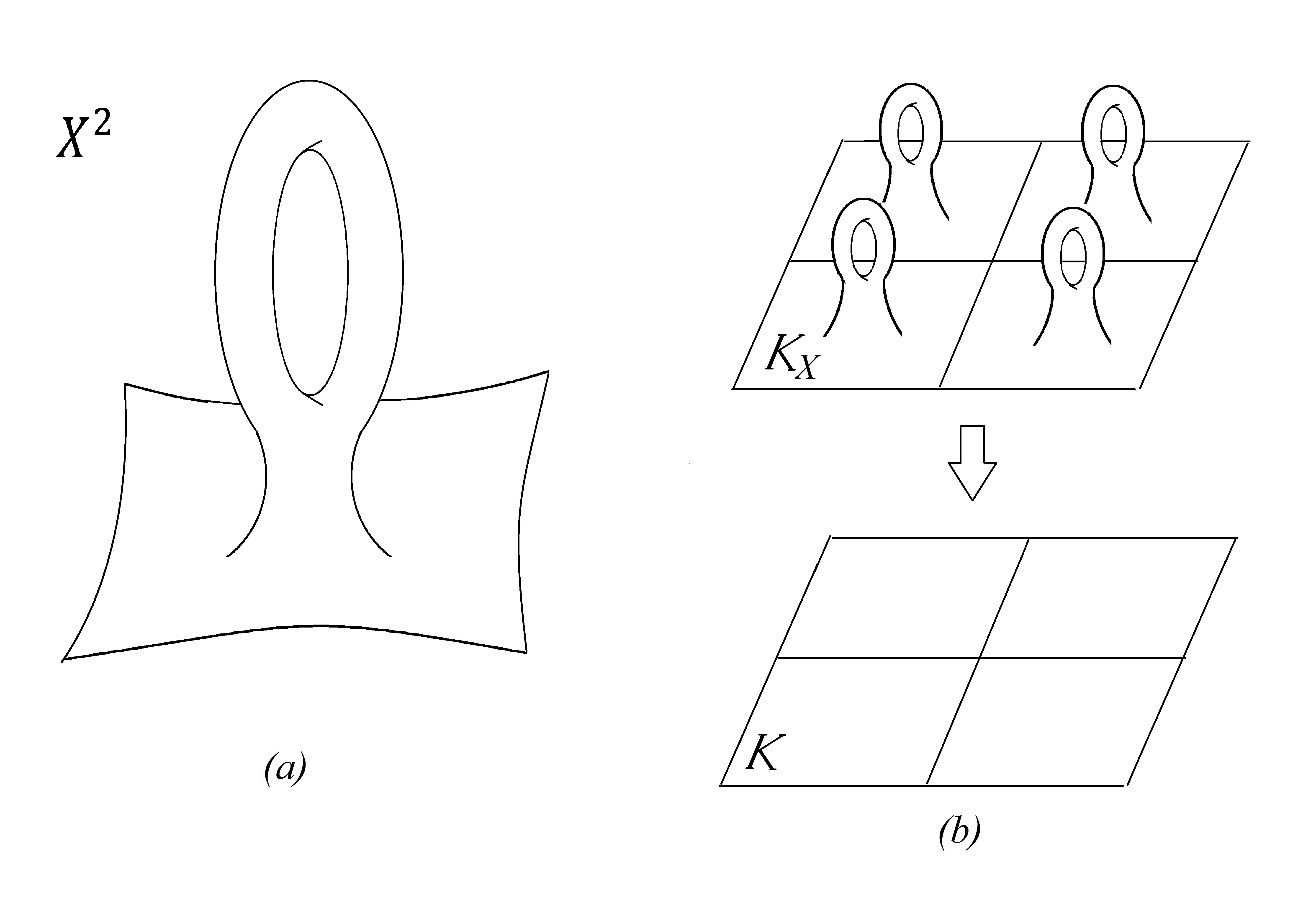} 
\caption{\small Picture (a) shows a 2 dimensional hyperbolization piece. Picture (b) 
shows the replacement of cubes in $K$ by copies of $X^2$.}
\end{figure}
\end{center}

Note that the edges (1-faces) match nicely and the piecewise hyperbolic metric (call it  $\sigma$)
is smooth away from the vertices. Near a vertex $o$ the metric is warped. Indeed
near $o$ it has the form $\sigma=\sigma\0{L}=sinh^2(t)\frac{L}{4}\sigma\0{\bS^1}+ dt^2$, where:\, (1) we are identifying
a punctured neighborhood of $o$ with $\bS^1\times (0,r+2)$, for some $r>0$, \,(2) $\sigma\0{\bS^1}$
is the canonical metric of the circle $\bS^1$, and \,(3) the number $L$ is the number of 2-cubes containing $o$ (which is equal to the number of copies of $X$ containing $o$).
Of course if $L=4$ the metric $\sigma\0{L}$ is already hyperbolic and smooth near $o$.
The problem arises when $L\neq 4$. In this case the solution is given by the {\it Gromov-Thurston trick}. Choose $d>0$ with $d<r$ and subdivide $(0,r+2)$ in three pieces $I_1=(0,r-d]$, $I_2=[r-d,r]$, $I_3=[r,r+2)$ and
let $\rho=\rho\0{L,r,d}$ be a smooth function on $(0,r+2)$ such that $\rho\equiv 1$ on $I_1$, $\rho\equiv \frac{L}{4}$ on
$I_3$ ($\rho$ is obtained by rescaling a fixed function). Consider now $h\0{L}=h\0{L,r,d}=sinh^2(t)\rho(t)\sigma\0{\bS^1}+ dt^2$. Then $h\0{L}$ and $\sigma\0{L}$ coincide on
$\bS^1\times I_3$, hence we can define the {\it smoothed metric} $\cG_L=\cG(L,r,d)$ near $o$ to be equal to $\sigma\0{L}$
outside $\bS^1\times (I_1\cup I_2)$ and equal to $h$ on $\bS^1\times (0,r)$ (see fig. 2). Using the
Bishop-O'Neill formula in \cite{BisOn} it can be shown that by choosing $r$ and
$d$ large enough
(depending on $L$) the metric $\cG_L$ will have curvature very close to -1. 
Furthermore, since  $\cG_L$ is canonically hyperbolic on $\bS^1\times I_1$
we can extend the metric $\cG_L$ to a smooth metric on the whole ($r+2$)-ball centered at $o$ which is hyperbolic on the ($r-d$)-ball.
We do this for every vertex and we are done. 
Note that we do need the pieces $X^2$ to be very large, that is the injectivity
radius of the vertices has to be very large.\\

\begin{center}
\begin{figure}[ht]
\centering
\includegraphics[width=9cm,height=6cm]{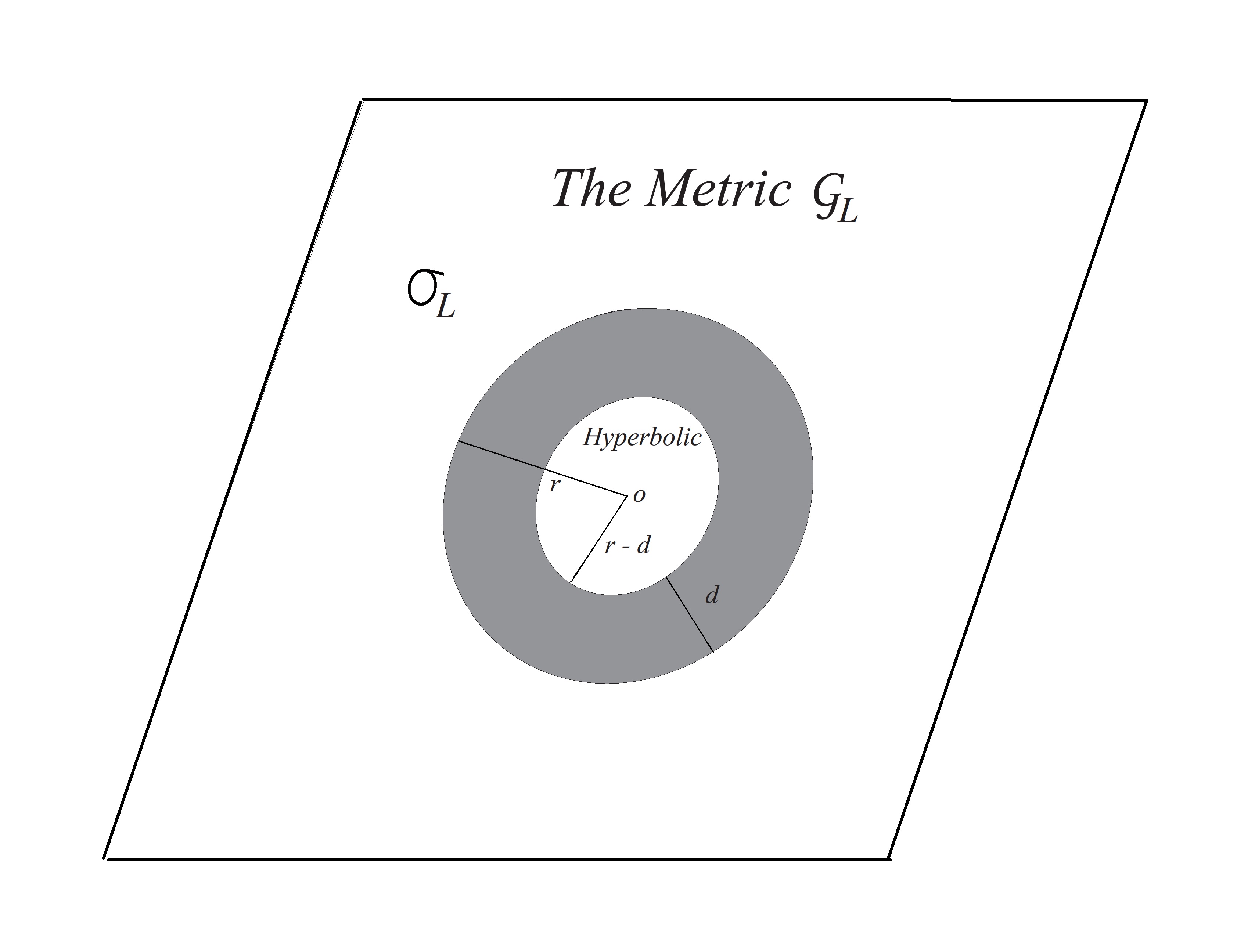}
\caption{\small The metric $\cG_L$ is canonically hyperbolic on the ball of radius $r-d$ and equal to
$\sigma_L$ outside the ball of radius $r$. }
\end{figure}
\end{center}

\newpage

\noindent {\bf 0.2. Codimension two and the Gromov-Thurston trick.}
If $N$ is a closed codimension two totally geodesic submanifold of a hyperbolic
manifold $(M,g)$, with trivial normal bundle, then $N$ has a neighborhood $\cN_{r+2}$
isometric to $N\times \B_{r+2}$ (where $\B_{r+2}\sbs\HH^2$ is the ball,
centered at $0\in\HH^2$, of radius $r+2$) with metric
$cosh^2(t)\,h+\sigma\0{\HH^2}$, where $h=g|\0{N}$, $\sigma\0{\HH^2}$ is the canonical metric on $\HH^2$ and $t$ is the distance to $0\in\HH^2$. We call this metric
a {\it hyperbolic extension of} $\sigma\0{\HH^2}$. 
Suppose now that we have a singular metric on $M$, which is smooth outside
$N$, and on $\cN_r-N$ is isometric to $N\times (\B_{r+2}-\{0\})$, with metric
$cosh^2(t)\,h+\sigma\0{L}$. Then we can smooth the metric $g$ 
to obtain a smooth metric $\cG_N=\cG(N,L,r,d)$ by changing $g$
using the smooth metric $cosh^2(t)\,h+\cG_L$ (where $\cG_L$ is as in 0.1) instead of the
singular metric $cosh^2(t)\,h+\sigma\0{L}$. This
method was used by Gromov and Thurston \cite{GT} to smooth singular metrics
obtained using branched covers.
The smoothed metric $cosh^2(t)\,h+\cG_L$
is a {\it hyperbolic extension of} $\cG_L$. Note also that $\cG_N$ is hyperbolic
on $N\times\B_{r-d}$ and equal to $g$ outside $\cN_r$ (see fig. 3).
\\

\begin{center}
\begin{figure}[ht]
\centering
\includegraphics[width=9cm,height=6cm]{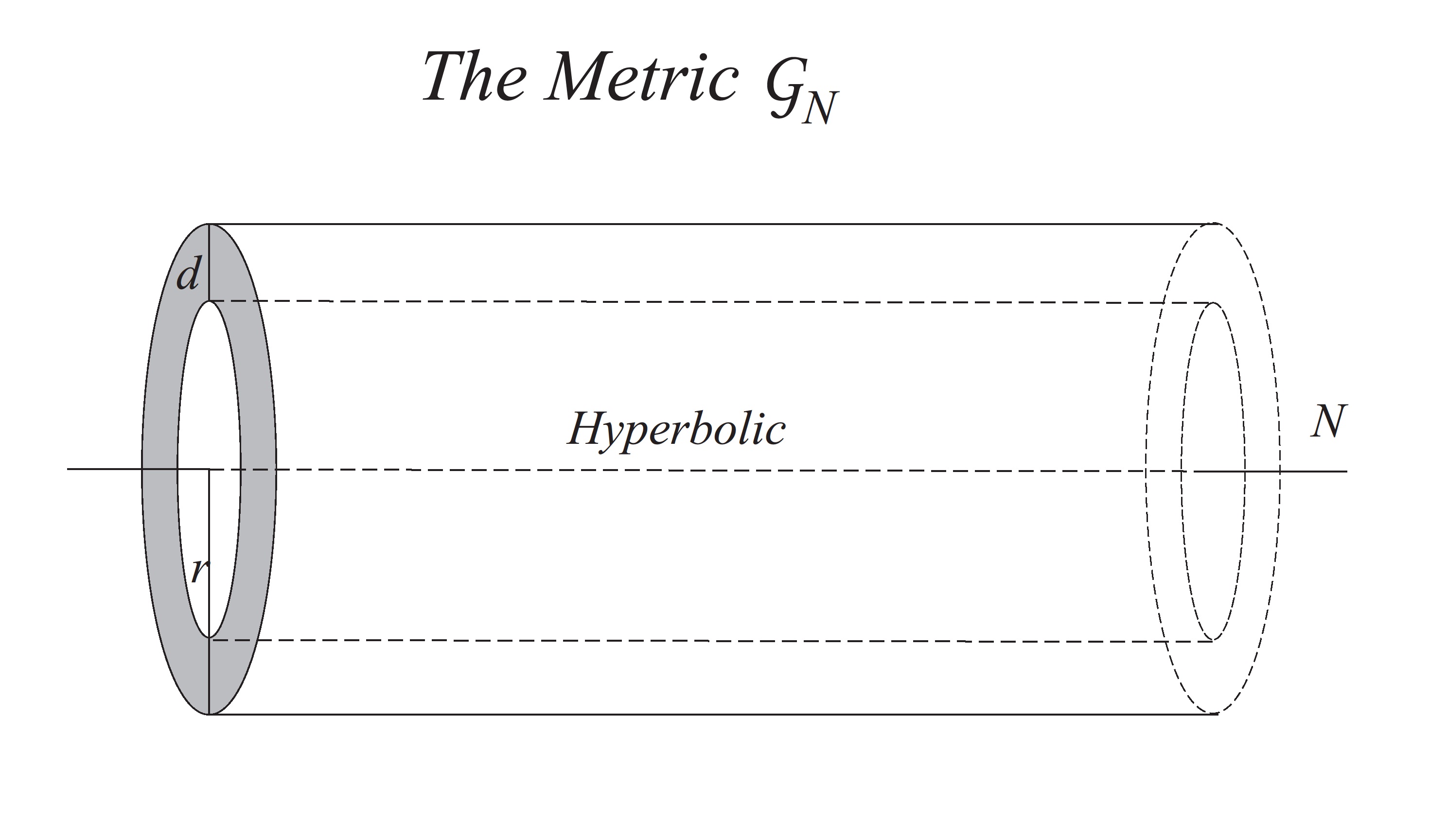}
\caption{\small The metric $\cG_N$ is canonically hyperbolic on the 
$(r-d)$ normal neighborhood of $N$ and equal to $g$ outside the
$r$ normal neighborhood. }
\end{figure}
\end{center}

\noindent {\bf 0.3. The Farrell-Jones warping trick.}
Before we deal with the dimension three case we have to discuss the
{\it Farrell-Jones warping trick}. This trick is in some sense a generalization of
the Gromov-Thurston trick in dimension 2.\\

Suppose we have a metric $h$ on the sphere
$\bS^n$. Consider the warped metric $g=sinh^2(t)\, h +dt^2$ on $\R^{n+1}-\{0\}=\bS^n\times (0,\infty)$. If $h=\sigma\0{\bS^n}$, the canonical metric on $\bS^n$, then $g$ is hyperbolic and, in particular, smooth everywhere. But for general
$h$ the metric $g$ is singular at 0.  Before we continue here is an important observation
that can easily be deduced from Bishop-O'Neill curvature formula in \cite{BisOn}.\\

\noindent {\bf (0.3.1.)}\,\, {\it Given $\epsilon>0$ there is $t_0$ such that the sectional curvatures of $g$ at $(x,t)$ are within \linebreak
\hspace*{.63in} $\epsilon$ of -1, provided $t>t_0 $.}\\

\noindent To smooth the metric $g$ consider the family of metrics (see \cite{FJ1}): $$g\0{\alpha}(x,t)=
sinh^2(t)\,\Big( \big(1-\rho\0{\alpha}(t)\big)\, \sigma\0{\bS^n}(x)\,\,+\,\, \rho\0{\alpha}(t)\,h(x) \Big)\,\,\,+\,\,\,dt^2$$
\noindent where $\rho\0{\alpha}(t)=\rho(\frac{t}{\alpha})$, and $\rho:\R\ra [0,1]$ is
a fixed smooth function with $\rho(t)=0$ for $t\leq 1$ and $\rho(t)=1$ for $t\geq 2$.
Hence, for $t\leq \alpha$, the metric $g\0{\alpha}$ is hyperbolic, for $t\geq 2\alpha$
we have $g\0{\alpha}=g$ and in between $t=\alpha$ and $t=2\alpha$ the metric
$\sigma\0{\bS^n}$ deforms to $h$. The metrics $g\0{\alpha}$ have two important
properties: (1) they are all hyperbolic for $t\leq\alpha$, hence smooth everywhere,
and (2) given $\epsilon>0$ there is $\alpha\0{0}$ such that all sectional curvatures of 
$g\0{\alpha}$ lie within $\epsilon$ of -1, provided $\alpha>\alpha\0{0}$.
Here is a heuristic explanation why (2) holds. If $\alpha$ is very large the deformation
between $\sigma\0{\bS^n}$ and $h$ happens very slowly (on the 
``stretched interval" $[\alpha,2\alpha]$), so $g\0{\alpha}$
is ``almost warped", hence the Bishop-O'Neill formula should give a good approximation
of the curvatures of $g\0{\alpha}$. Therefore, by 0.3.1 the curvatures of
$g\0{\alpha}$ should be close to -1, provided we are far away from 0. But since
we are assuming $\alpha$ large, we are in fact far away from 0. Interestingly
the actual proof of (2) given in \cite{FJ1} does not follow exactly this intuitive
explanation because there is a more direct (though more technical) proof.\\

In this paper we need a more elaborate version of the Farrell-Jones trick, which
we call {\it hyperbolic forcing}.
We need to know to what extend the ``stretching" (which in the Farrell-Jones trick \cite{FJ1} is given by a variable $\alpha$)
to be independent of the ``far-away constant" (given in the Farrell-Jones trick also
by $\alpha$).  Moreover, we need a more quantitative version also. Interestingly
again our proof of this more elaborate version does follow more closely the
intuitive approach explained above. Before we finish this subsection we make
an important remark.\\

\noindent {\bf (0.3.2.)} {\it Given $\epsilon>0$, the stretching
and the far-away constants needed in hyperbolic forcing \\
\hspace*{.556in} (to obtain an $\epsilon$-
pinched to -1 metric) do depend on the metric $h$.}
\vspace{.2in}

\noindent {\bf 0.4. Dimension three.} Suppose we have a cubical complex $K$ of
dimension 3. As in 0.1 choose $X^3$ and construct $K_{X}$. 
Call the piecewise hyperbolic metric on $K_X$ by $\sigma=\sigma\0{K_X}$. Again as in 0.1 the codimension one faces (the 2-faces) of $X$ match nicely and there are singularities only on the
``1-skeleton" of $K_X$, that is, along the edges (i.e 1-faces) and vertices (i.e 0-faces).
The singularities along the 1-faces can be smoothed using the Gromov-Thurston
trick as in 0.1 and 0.2 (i.e. using smoothing in dimension two plus hyperbolic
extension). In this way we obtain a metric $\sigma'$ which is smooth
near (part of) the edges.
Let $\cN_{r+2}(e)$ be the normal neighborhood of width $r+2$ of the edge $e$.
Then $\sigma=\sigma'$ outside the union $\bigcup_{e\, {\mbox{\tiny edge}}}\cN_{r+2}(e)$. Notice that
there is some ambiguity in the definition of the metric $\sigma'$ because
for different edges $e$, $e'$ with a common vertex the neighborhoods
$\cN_{r+2}(e)$, $\cN_{r+2}(e')$ have nonempty intersection. So $\sigma'$ is only well-defined outside the $s$-neighborhoods (i.e $s$-balls ) $\cN_s(o)$ of the vertices $o$,
where $s$ is large enough (see fig. 4). Let $L(e)$ be the number of copies of $X$ that
contain the edge $e$, and write $\cG_e=\cG_{L(e)}$ and $\sigma_e=\sigma_{L(e)}$ (see 0.1).
Therefore on each $\cN_{r+2}(e)$, and outside
$\bigcup_{o\, {\mbox{\tiny vertex}}}\cN_s(o)$, the metric $\sigma'$ is equal to
the metric $cosh^2(t)\, \sigma\0{\R}+\cG_e$\, which is the hyperbolic
extension of the metric  $\cG_e$  (see 0.2).\\

\begin{center}
\begin{figure}[ht]
\centering
\includegraphics[width=9cm,height=7.5cm]{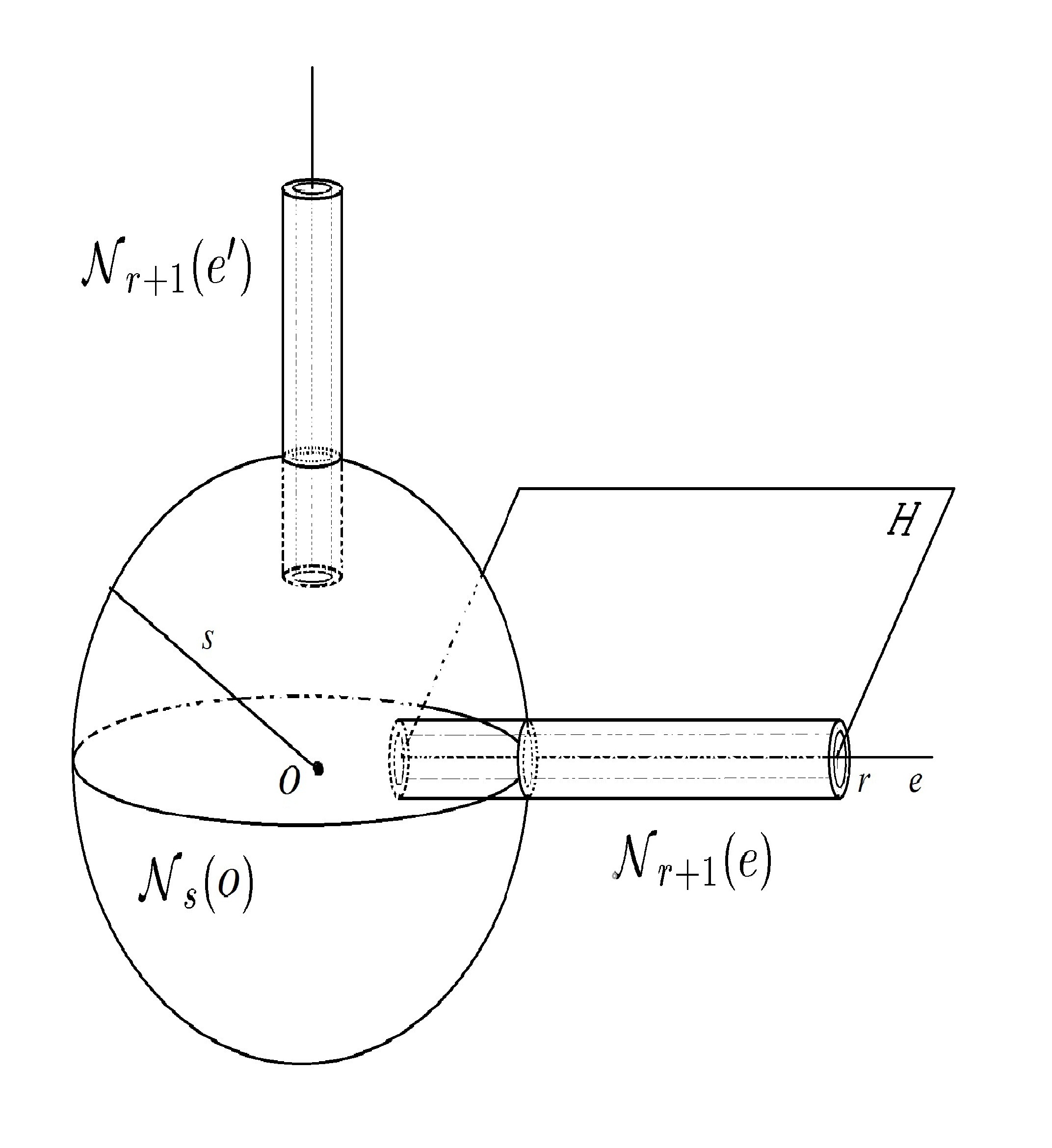}
\caption{\small The piecewise hyperbolic metric $\sigma$ can be smoothed 
near each edge $e$. In this way we obtain the metric $\sigma'$. The plane $H$ is a plane that
contains the edge $e$. In this picture we only show a ``half" of $H$.}
\end{figure}
\end{center}

We are left to smooth the metric near the vertices.   
Fix a vertex $o$. Let $\{e_i\}$ be the edges containing $o$ and write $e=e_1$.
Let $P$ be the link of $o$. Then $P$ is a $PL$-sphere of dimension 2
and it has a natural all-right spherical metric $\sigma\0{P}$, that is
$\sigma\0{P}$ is the piecewise spherical metric with all edges in $P$ having length $\pi/2$. Note that the metric $\sigma$ near $o$ is the warped (piecewise hyperbolic) metric
$sinh^2(s) \,\sigma\0{P}+ds^2$, where $s$ is the distance to $o$.
Write $L_i=L(e_i)$ and $L=L_1$. Then the metric $\sigma$ on $\cN_{r+2}(e_i)$ is
equal to the hyperbolic extension metric $cosh^2(t_i)\, \sigma\0{\R}+\sigma\0{e_i}$,
where $t_i$ is the distance to $e_i$. Hence $\sigma$ is $sinh$-warped from $o$
and $cosh$-warped from each $e_i$ (near $e_i$). What about the metric $\sigma'$?
It is also $cosh$-warped from $e_i$ because, by definition, the metric $\sigma'$
is equal to $cosh^2(t_i)\, \sigma\0{\R}+\cG_{e_i}$ near $e_i$ (i.e on $\cN_{r+2}(e_i)-\bigcup_{o'\,{\mbox{\tiny vertex}}}\cN_s(o')$). But, and this is a key observation, the metric $\sigma'$ is {\sf not }
(in general) warped from $o$. Here is a heuristic idea why this is so.
Write $\cE_i=cosh^2(t_i)\, \sigma\0{\R}+\cG_{e_i}$ and $\cE=\cE_1$.
Note that $\cE$ has rotational symmetry, that is it is invariant by rotations
fixing $e=e_1$. Let $H$ be a plane containing $e$.
Then $\cN_r(e)$, $\cN_{r-d}(e)$ intersect $H$ in two lines each.
Let $p\in H\cap\cN_{r+2}(e)$, and $x\in e\sbs H$ be the closest point in $e$ to $p$. Also
let $v$ a vector at $p$ perpendicular to $H$,
and denote the circle centered at $x$, perpendicular to $H$ and passing through
$p$ (hence tangent to $v$) by $S(p)$. Let $t=t(p)$ be the distance from
$p$ to $x$ and $s=s(p)$ the distance from $o$ to $p$ (see fig. 5).\\

\begin{center}
\begin{figure}[ht]
\centering
\includegraphics[width=9cm,height=6cm]{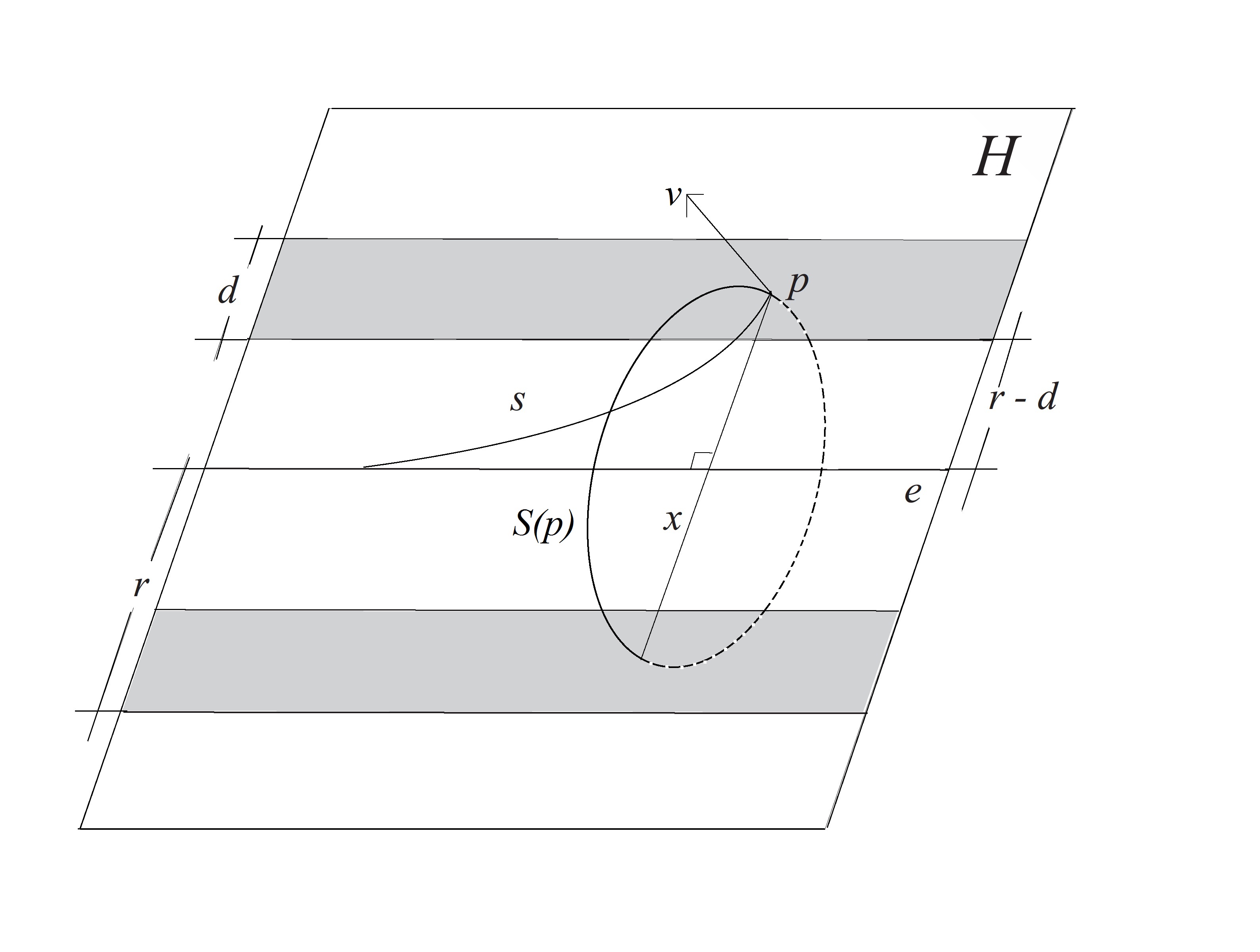}
\caption{\small The horizontal lines are the intersections of $\cN_r(e)$ and $\cN_{r-d}(e)$ with the plane $H$. The circle $S(p)$ is perpendicular to $H$. The vector $v$ is tangent to $S(p)$ at $p$.}
\end{figure}
\end{center}

Now, it can be checked from the definitions that the metrics $\sigma$ and $\sigma'$ coincide on vectors tangent to $H$. They differ on their values on the vectors $v$
as above. These values are directly proportional to the lengths $\ell (S(p))$,
$\ell'(S(p))$ of the circle $S(p)$ with respect to the metrics $\sigma$ and $\sigma'$,
respectively. Hence these metrics can be understood by looking at these lengths.
We have $\ell (S(p))=2\pi\,\frac{L}{4}\,  sinh(t)$ and $\ell' (S(p))=2\pi\,\rho(t)\,sinh(t)$ (see 0.1).
Let $\theta$ be the angle at $o$ between $e$ and the geodesic segment $\beta=
[o,p]$ (which lies on $H$). Let $p(s)$ be the point in $\beta$ at distance $s$
from $o$. Now if $\sigma'$ were $sinh$-warped from $o$  the lengths of the circles $S(s)=S(p(s))$ would have the form
$c\, sinh (s)$ for some constant $c$. But from the hyperbolic law of sines
we have $sin\,\theta=\frac{sinh\,t}{sinh\,s}$, hence $t(s)=sinh^{-1}(sin\,\theta\,
sinh\, s)$ and we get\\

\noindent {\bf (0.4.1.)} \hspace{.5in} $\ell'(S(s))=2\,\pi\,\rho(t(s))\, sinh(t(s))\,\,=\,\,2\,\pi\, \rho(t(s))\,sin\,\theta\, sinh\, s $\\

\noindent Note that if $L=4$, then $\rho\equiv 1$ and the formula above shows
why hyperbolic three space $\HH^3$ is at the same time $sinh$-warped from 
a point $o$ and $cosh$-warped from a line containing $o$.
But in general $\rho(t(s))$ is not a constant, hence $\sigma'$ is not, in general,
$sinh$-warped from $o$, as we wanted to show. Note that $\rho(t)$ is constant
for $t\notin [r-d,r]$.\\

Why do we want $\sigma'$ to be $sinh$-warped from $o$?
Because  in this case we could apply hyperbolic forcing (see 0.3) and force (or extend) the metric
$\sigma'$ to be hyperbolic near $o$, hence smooth near $o$. (Even in this case
there would be a problem in using hyperbolic forcing because of 0.3.2, but
more on this in a moment.) Now, even though $\sigma'$ is not warped from $o$
it is ``very close to being warped", provided $r$ and $d$ are large. Here is an idea why this is true. Since $sin\,\theta=\frac{sinh\, t}{sinh\, s}$, if $t$ is large, so is $s$ and
we get $s\approx t-ln\,sin\, \theta$, and if $r$ and $d$ are large then
the function $\rho(t(s))$ in 0.4.1 even though is not constant, it does, in this case,
change very slowly, hence behaves (locally) almost like a constant. Therefore
$\sigma'$ is ``almost warped" in this case, and we can ``deform" $\sigma'$ to
a $sinh$-warped from $o$ metric. We call this process {\it warp forcing}.
Therefore the idea is to first warp force the metric $\sigma'$ near $o$ and then
hyperbolic force it to make it hyperbolic near $o$, hence smooth. In our particular
case the $sinh$-warped metric to which we deform $\sigma'$ is obtained in the
following way. Let $\bS_{s_0}=\bS_{s_0}(o)$ be the sphere of radius $s_0$ in $K_X$
centered  at $o$ (recall $X$ is as large as needed). The metric $h_{s_0}=\sigma'|_{\bS_{s_0}}$, ie, the restriction of $\sigma'$ to $\bS_{s_0}$ is called
the {\it (warped) spherical cut of $\sigma'$ at $s_0$} and the metric
$\hat{h}_{s_0}=\frac{1}{sinh\, s_0} h_{s_0}$ is called the {\it unwarped
spherical cut of $\sigma'$ at $s_0$}. The metric $\sigma'$ is deformed, 
by warp forcing, to the warped metric $sinh^2(s)\, \hat{h}_{s_0}+ds^2$.
\vspace{.3in}

Now we deal with the problem mentioned above. Suppose we succeeded
in warp forcing the metric $\sigma'$ and obtained the $sinh$-warped from $o$
metric $sinh^2(s)\, \hat{h}_{s_0}+ds^2$. Recall that we needed to assume $r$ and $d$ large. By 0.3.2 the constants needed for hyperbolic forcing  (call them $\alpha_1,\, \alpha_2$) depend on
$\hat{h}_{s_0}$, which in turn depend on $r,d$. It may happen that
the $\alpha_i=\alpha_i(r,d)$ are too big for $s_0$ and we have no space for
hyperbolic forcing. And in fact this may happen if we don't do things in a 
precise way. To solve this problem we proceed in the following way.
Fix an angle $\theta_0>0$ and $d$ large as needed but fixed.
Consider the plane $H$ as before (see fig.5 and 6). Let $q$ be
a point in $H$ at distance $r$ from $e$ such that the geodesic (in $H$)
$[o,q]$ makes an angle $\theta_0$ at $o$. Let $s=s(r)$ be the distance from
$o$ to $q$. Now let $p=p(r)$ be the point in $H$ such that the distance
from $p$ to $o$ is also $s$, and the distance to $e$ is $r-d$. Let $\theta_1(r)$
be the angle at $o$ between $e$ and $[o,p]$. \\

\begin{center}
\begin{figure}[ht]
\centering
\includegraphics[width=9cm,height=6cm]{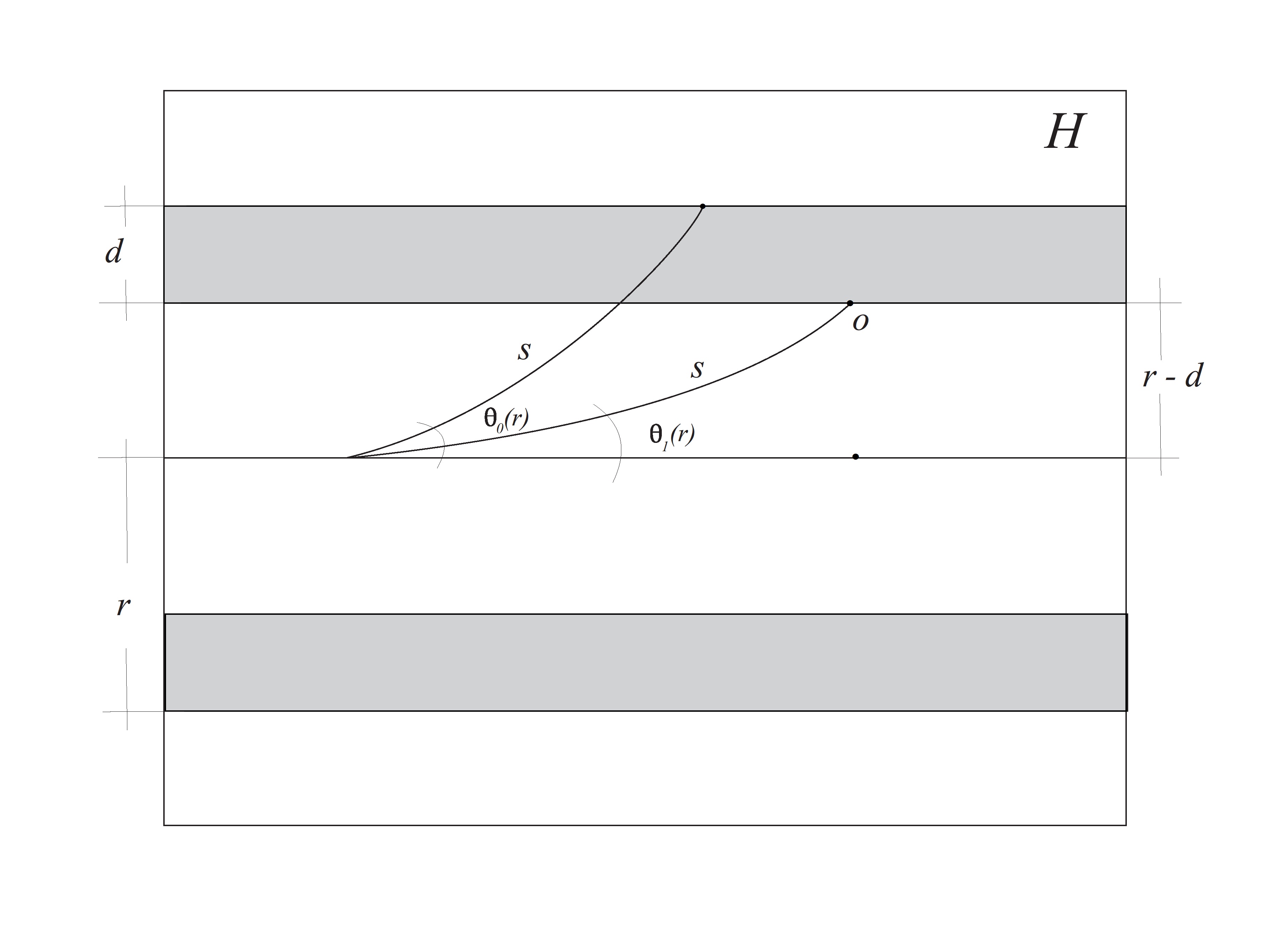}
\caption{\small Both the number $d$ and the angle $\theta_0$ are fixed.
The number $s=s(r)$ depends on $r$. Likewise the angle $\theta_1(r)$ and the {\it spherical cut}
${\hat{h}}_{s(r)}$ depend on $r$ and have a limit, as $r\ra\infty$.}
\end{figure}
\end{center}

It can be shown with a 
straightforward calculation that in this particular case the corresponding
metrics $\hat{h}_{s(r)}$ $C^2$-converge to a smooth metric $\hat{h}$.
(And the angles $\theta_1(r)$ also converge. Note $\theta_0$ is fixed.) 
And it can be shown that in this case 0.3.2 does not pose a problem any more
because all metrics obtained are in fact very close, so the corresponding
constants $\alpha_i$ are close. In particular they don't grow indefinitely.
And this is what we needed. Consequently, we first take
$d$ very large (for curvature and other considerations) and then take $r$ large
so that everything else works, and during all this process we have to make
sure to be far enough away from $e$, and this is given by the constant $\theta_0>0$.\\

Of course dimensions 2 and 3 are very special cases. For instance, in general dimensions the metrics $\cG_e$ do not have rotational symmetry and the
simple calculations shown above don't work in general. Also there
are all sorts of matching problems and many constants and widths of
neighborhoods have to be chosen in a very precise way. 
For these and other
reasons we warn the reader that the actual proofs presented in this paper
do not follow exactly the pattern given in this introduction, though some key ideas
have in fact been introduced. There are other important
issues that have been omitted from the introduction. For instance the
smoothing problem of finding a good differentiable structure to work
with is quite considerable  (the ``natural" one of Davis-Januszkiewicz
and Charney-Davis is not useful in our case).\\

Here is a brief description of the paper. In section 1 we recall some generalities about
warped metrics and introduce notation for weaker versions of it. In section 2 we define and study hyperbolic extensions.  Section 3 deals with approximation of metrics and hyperbolic forcing. We  treat spherical cuts and warp forcing in section 4. We introduce
limit metrics and study how to ``continue" (i.e extend) families of metrics
 in section 5. In section 6 we study neighborhoods of cubic, spherical and
piecewise hyperbolic complexes. In this section we introduce a technical 
device that we called {\it sets of widths}. These are sets of positive real numbers
 that are used as widths for normal neighborhoods of faces of a complex.
We prove that there are sets of widths, independent of the complex, that
satisfy very useful properties. These objects make all matching processes
work. In section 7 we deal with the smoothing issue for cubical and
all-right spherical complexes. We prove that smooth manifolds with a smooth
cubification admit a ``normal differentiable structure" that is diffeomorphic
to the original one. The proof is by decreasing induction on the dimension
of the skeleta. Since heavy duty classical smoothing techniques are
used (which do not work in dimension 4) a special analysis had to be done
for dimension four, and this is done in appendix C. We put everything together
in section 8 to smooth cones. Section nine is dedicated to the Charney-Davis
strict hyperbolization process. In particular we treat also the smoothing issue
for the strictly hyperbolized complexes and use the results in section 7
to show that strictly hyperbolized smooth manifolds also have nice
``normal differentiable structures". This is quite a delicate issue and three proofs are
relegated to the appendices. Finally we prove the Main Theorem in
section 10 and Theorem A in section 11. Some subsections have been added
at the end of a few sections that deal with generalizations to the case
of manifolds with codimension zero singularities. These subsections
are used in section 11. There are seven appendices. Several numerical
estimates appear throughout the paper. Most of these estimates are far from being sharp.\\

We are grateful to Tom Farrell and Ross Geoghegan for their comments
and suggestions. We are also thankful to C.S. Aravinda, M. Deraux and
J.-F. Lafont for the useful information provided to us.

\vspace{1in}

\noindent {\bf \large  Section 1. Constructing Metrics.}\\

We recall some classical constructions of metrics and introduce some notation for
other constructions. In this paper by ``geodesic" we mean both the function
and its image depending on the context. We will also use the same notation for both.\\

\noindent {\bf A.  Product metrics.}\\
Let $M_1$, $M_2$ be smooth Riemannian manifolds with Riemannian metrics $g_1$, $g_2$. Recall that the product metric
$g$ on the product $M_1\times M_2$ is given by $g(v_1+v_2 , w_1+w_2)_{(p,q)}=g_1(v_1,w_1)_p+g_2(v_2, w_2)_q$, for $v_1+v_2,\, w_1+w_2\in T_pM_1\oplus 
T_qM_2= T_{(p,q)}(M_1\times M_2)$. We write $g=g_1+g_2$. In particular, if $M$ has metric $h$, the product metric on $M\times \R$
is $h+dt^2$, where $dt^2=dt\otimes dt$ is the canonical metric on $\R$.\\

\noindent {\bf B.  Warped metrics.}\\
Let $f:M_2\ra \R^+=(0,\infty)$ be smooth. The metric $g=f^2g_1+g_2$ is called a {\it warped metric} on $M_1\times M_2$, and $f$ is the 
{\it warping function}. $M_1$ is called the {\it fiber} and $M_2$ is the {\it base}.
With this warped metric the fiber $M_1\times \{q\}$ in $M_1\times M_2$ over $q\in M_2$ is scaled 
by the positive number $f(q)$.
(For a study of warped metrics see \cite{BisOn}.) In the one-dimensional
case, that is, when $M_2=I\sbs\R$ is an interval, the warped metric on $M\times I$ is written $f^2h+dt^2$, where
$f:I\ra \R^+$ and $h$ is a metric on $M$. In this case, the vector field $\pt$ is called the {\it parameter
vector field}.\\

For some classical examples of warped metrics see next section (section 2).\\

\noindent {\bf C.  Doubly warped metrics.}\\
More generally, let $N,\, M_1,\, M_2$ be Riemannian manifolds with metrics
$h,\, g_1,\, g_2$, respectively.
Let $f_1, \,f_2:N\ra \R^+$ be smooth. The metric $g=f_1^2g_1+f^2_2g_2+h$ is a {\it doubly warped metric} on $M_1\times M_2\times N$. Here $f_1$
and $f_2$ are the warping functions. Analogously, triply warped metrics,..\, etc, can be defined.\\

\noindent {\bf D.  Variable metrics.}\\
Let $g_t$, $t\in I\sbs\R$, be a one-parameter family of metrics on the manifold $M$.
(All one-parameter families in this paper will be assumed to be smooth, were ``smooth" means $(x,t)\mapsto g(x,t)=g_t\mid_x$ is smooth.)
 We call the metric $g=g_t+dt^2$ on $M\times I$ a {\it variable metric
with metrics $g_t$}. Note that the $t$-lines $t\mapsto (x,t)$, $x\in M$, are the integral curves of the vector field $\pt$ on $M\times I$.
Also $\pt$ has $g$-length one, and it is perpendicular to the (integrable) distribution $(x,t)\mapsto T_xM\sbs T_xM\times\R=T_{(x,t)}M\times I$.\\

\noindent {\bf Lemma 1.1.} {\it The t-lines are geodesics.}\\

\noindent {\bf Proof.} By Koszul's formula $2\,g(\nabla_\pt\pt,\pt)=\frac{d}{dt}g(\pt,\pt)+\frac{d}{dt}g(\pt,\pt)-\frac{d}{dt}g(\pt,\pt)=0$.
Also, for a $t$-independent vector field $V$ tangent to the fibers $M\times\{t\}$ we have
$g((\nabla_\pt\pt,V)=\frac{d}{dt}g(\pt,V)+\frac{d}{dt}g(V,\pt)-Vg(\pt,\pt)=0$. Hence 
$\nabla_\pt\pt=0$. This proves the Lemma.\\

Conversely, let $g$ be a metric on $M\times I$ and assume that $\pt$ has $g$-length one. Assume also that $T_xM\sbs T_{(x,t)}M\times I$
is $g$-perpendicular to $\pt$, for all $(x,t)$. Then $g$ is a variable metric on $M\times I$. For instance, on a Riemannian manifold $(N^n,g)$ we
can identify, using the exponential map, a geodesic ball $B(p,\epsilon)$, with its center $p$  deleted, with the cylinder $\bS^{n-1}\times (0,\epsilon)$.
Then, by Gauss lemma, $g$ is a variable metric on this cylinder, and we write $g=g_r+dr^2$, where $r$ is the distance to $p$.
For a point $z\in B(p, \epsilon)-\{ p\}$, $z=\exp_p(rv)$, $v\in T_pM$, $|v|=1$, the pair $(v, r)$ are called the {\it polar coordinates expression }
of $z$ and sometimes we will write $z=rv$.\\

Note that every warped metric $f^2h+dt^2$ on $M\times I$ is a variable metric on $M\times I$, with $h_t=f^2(t)h$.\\

Let $N$ be a smooth manifold with Riemannian metric $g$.
Suppose $N$ has a decomposition $N=M\times I$ and $g$ is a variable metric
with respect to this decomposition. We say that the set of lines $t\mapsto (x,t)$
(which are $g$-geodesics) is the {\it ray structure of $g$ corresponding to the
decomposition $N=M\times I$.}\\


\noindent {\bf E.  Warped variable metrics.} \\
Let $g_t$, $t\in I\sbs\R$, be a one-parameter family of metrics on the manifold $M$ and let $f:I\ra\R^+$ be smooth. 
We call the metric $g=f^2g_t+dt^2$ $M\times I$ a {\it warped variable metric with metrics $g_t$
and warping function $f$}. Note that every variable metric $g=g_t+dt^2$ is a warped variable metric by any $f$ because we can write
$g=f^2(\frac{1}{f^2} g_t)+dt^2$, but the one-parameter family of metrics does depend on the chosen $f$.
Note also that the $t$-lines are geodesics.\\

Also, doubly, triply. etc, warped variable metrics are defined in the obvious way.\\

\noindent {\bf F. Warped metrics on piecewise Riemannian manifolds.}\\
We can assume the fibers in items B, C, D, E to be
piecewise Riemannian manifolds, by warping each piece and gluing the pieces (which will match because the ``unwarped"
pieces match). It is possible to get the piecewise Riemannian case for the base space under certain conditions, 
but we will not need this case.


\vspace{1in}

\noindent {\bf \large  Section 2. Classical Warped Metrics and Hyperbolic Extensions.}\\

In the first part of this section we give some constant curvature examples of warped (and doubly warped) metrics. The spherical and hyperbolic examples
will be used later. The examples also motivate the main construction in section 2, the ``hyperbolic extension", which is given in the second part 
of this section.\\

\noindent {\bf 2.1. Classical warped metrics.}\\
We will denote by  $\sigma\0{M}=\sigma$ the canonical Riemannian metric on $M$, where $M$ is either euclidean space $\R^k$, the unit sphere $\bS^k$
or hyperbolic space $\HH^k$.\\

\noindent {\bf a. Punctured euclidean space.} $\R^n-\{ 0\}$ is canonically isometric to $\bS^{n-1}\times (0,\infty)$ with warped metric $r^2\sigma_{\bS^{n-1}}+dr^2$. The warping function is $r\mapsto r$, i.e. the identity function. Here $r$ corresponds to the distance to the origin in $\R^n$.\\

\noindent {\bf b. Normal neighborhoods of euclidean subspaces.} More generally, let $0\leq k<n$. Then $\R^n-\R^k$ is 
isometric to $(\R^k\times\bS^{n-k-1})\times(0,\infty)$ with
doubly warped metric $\sigma_{\R^k}+r^2\sigma_{\bS^{n-k-1}}+dr^2$. The warping functions are $r\mapsto 1$ and $r\mapsto r$, and $r$ corresponds to
the distance to $\R^k\sbs\R^n$. \\

Also note that $\R^n$ is the normal neighborhood of $\R^k$ of width $\infty$. $\R^n$ isometric to $\R^{k}\times \R^{n-k}$ with
warped (product) metric $\sigma_{\R^{k}}+\sigma_{\R^{n-k}}$. The warping function is $r\mapsto 1$, where we can consider both  $\R^{n-k}$ and $\R^{k}$
as the base space. \\

\noindent {\bf c. Punctured  unit sphere.} Fix a point $u\in\bS^n$. Then $\bS^n-\{ u,-u\}$ is isometric to $\bS^{n-1}\times (0,\pi)$ with metric
$(sin^2r)\,\sigma_{\bS^{n-1}}+dr^2$. The warping function is $r\mapsto sin\, r$ and $r$ corresponds to the distance to $u$ (or $-u$).\\

\noindent {\bf d. Normal neighborhoods of spheres.} Let $0\leq k<n$. Let $\bS^k,\,\,\bS^{n-k-1}\sbs\bS^n$ be induced by the
embeddings $\R^{k+1}\hookrightarrow \R^{k+1}\times\{0\}$ and $\R^{n-k}\hookrightarrow \{0\}\times\R^{n-k}$. Hence $\bS^k,\,\,\bS^{n-k-1}$ are 
disjoint and linked.
Then $\bS^n-(\bS^k\cup\bS^{n-k-1})$ is 
isometric to $(\bS^k\times\bS^{n-k-1})\times(0,\pi/2)$ with
doubly warped metric $(cos^2\,r)\,\sigma_{\bS^k}+(sin^2\,r)\,\sigma_{\bS^{n-k-1}}+dr^2$. The warping functions are $r\mapsto cos\, r$ and $r\mapsto sin\,r$, 
and $r$ corresponds to the distance to $\bS^k\sbs\bS^n$. \\

Also, note that $\bS^n-\bS^{n-k-1}$ is the open geodesic normal neighborhood of width $\pi/2$ of $\bS^k$ in $\bS^n$.
This neighborhood is isometric to $\bS^{k}\times B^{n-k}$ with
warped metric $(cos^2\, r)\,\sigma_{\bS^{k}}+\sigma_{\bS^{n-k}}$, where $B^{n-k}\sbs\bS^{n-k}$ is an open hemisphere, i.e. it is one of the components
of $\bS^{n-k}-\bS^{n-k-1}$, and $r:B^{n-k}\ra [0,\pi/2)$ is the distance to the center of $B^{n-k}$.
The warping functions are $r\mapsto cos\,r$ and $r\mapsto 1$, and $r$ corresponds to the distance to $\bS^k\sbs\bS^n$.\\

\noindent {\bf e. Punctured hyperbolic space.} Fix a point $u\in\HH^n$. Then $\HH^n-\{ u\}$ is isometric to $\bS^{n-1}\times (0,\infty)$ with metric
$(sinh^2r)\,\sigma_{\bS^{n-1}}+dr^2$. The warping function is $r\mapsto sinh\, r$ and $r$ corresponds to the distance to $u$.\\

\noindent {\bf f. Normal neighborhoods of hyperbolic subspaces.} Let $0\leq k<n$. 
Then $\HH^n-\HH^k$ is  isometric to $(\HH^k\times\bS^{n-k-1})\times(0,\infty)$ with
doubly warped metric $(cosh^2\,r)\,\sigma\0{\HH^k}+(sinh^2\,r)\,\sigma_{\bS^{n-k-1}}+dr^2$. The warping functions are 
$r\mapsto cosh\, r$ and $r\mapsto sinh\,r$,  and $r$ corresponds to the distance to $\HH^k\sbs\HH^n$. \\

Also note that $\HH^n$ is the normal neighborhood of $\HH^k$ in $\HH^n$ of width $\infty$. And $\HH^n$ is isometric to $\HH^{k}\times \HH^{n-k}$ with
warped metric $(cosh^2\, r)\,\sigma_{\HH^{k}}+\sigma_{\HH^{n-k}}$, where $r:\HH^{n-k}\ra[0,\infty)$ is the distance to a fixed point in $\HH^{n-k}$. 
Notice that if we interchange $k$ with $n-k$ we get a different decomposition
of $\HH^n$.\\

For instance, in the case $n=2$, since $\HH^1=\R^1$ we have that $\HH^2$ is isometric to $\R^2=\{(u,v)\}$ with warped metric $cosh^2v\, du^2+ dv^2$.\\ 

\noindent {\bf g. Horospheres in hyperbolic space.} Hyperbolic space is also isometric to $\R^{n-1}\times\R$ with warped metric
$e^t\sigma_{\R^{n-1}}+dt^2$. The sets $\R^{n-1}\times\{ x\}$ correspond to horospheres, and all $t$-lines $t\mapsto (x,t)$ tend to the
same point at infinity, as $t\ra\-\infty$.
\vspace{.8in}

\noindent {\bf 2.2.  Hyperbolic Extensions.}  \\
Let $M^n$ be a complete Riemannian manifold with {\it center} $o=o_{_M}\in M$, that is, the exponential map $exp_o:T_oM\ra M$ is a diffeomorphism.
In particular $M$ is diffeomorphic to $\R^n$. 
For instance if $M$ is Hadamard manifold every point is a center point.
Denote the metric on $M$ by $h$.\\

\noindent {\bf Remarks.} 

\noindent {\bf 1.} In this paper we will use the same symbol ``$o$" to denote a center of a Riemannian manifold unless it is necessary to
specify the manifold, in which case we will write $o_{_M}$ if $o$ is a center of $M$.

\noindent {\bf 2.} In section 6.5 we will also consider piecewise Riemannian
spaces with a center.\vspace{.2in}

Let $r:M\ra [0,\infty)$ be the distance to $o$. Then 
\begin{enumerate}
\item[{\bf i.}] we have that $r(exp_o v)=h_o(v,v)^{1/2}$ hence $r$ is continuous and smooth on $M-\{ o\}$. Also $r^2$ is smooth on $M$.

\item[{\bf ii.}]  The (images of the) geodesic rays $exp_o(\R^+ v)$ are convex sets in $M$, and the geodesics lines $exp_o(\R v)$ are totally geodesic in $M$.
Here $\R^+= (0,\infty)$.

\item[{\bf iii.}]  the function $dr$ is strictly distance decreasing on non-radial vectors. That is, for $v\in TM-T_oM$ we have $|dr(v)|\leq h(v,v)^{1/2}$ and $|dr(v)|= h(v,v)^{1/2}$ if and only if $v$ is radial, i.e. tangent to a geodesic passing through $o$.
\end{enumerate}\vspace{.2in}

Using the diffeomorphism $exp_o$ onto $M$ and an identification
of $T_oM$ with
$\R^n$ via some fixed choice of an orthonormal basis in $T_oM$,
we can identify $M$ with $\R^n$ and $M-\{o\}$ with $\bS^{n-1}\times\R^+$.
 As mentioned in section 1D
 the metric $h|_{M-\{ o\}}$ can be written 
as a variable metric $h_r+dr^2$ on $\bS^{n-1}\times \R^+$.
Also, as in 1D, we shall call the set of ($h$-geodesic) rays
$t\mapsto (x,t)\in \bS^{n-1}\times \R^+$ the {\it ray structure of
$h$ with respect to o}.\\

The warped metric $g=(cosh^2 r)\, \sigma\0{\HH^k}+h$ on $\HH^k\times M$
is a {\it hyperbolic extension} 
of the metric $h$ on $M$. We will also say that the space $\HH^k\times M$ with metric $g$ 
is a hyperbolic extension of the Riemannian manifold $M$  and we will denote this extension by $\cE_k(M)$, and sometimes we will also write
$g=\cE_k(h)$.
Note that, even though $r$ is not smooth at $o$, the warping function $cosh\,r$ is smooth on
$M$ because $cosh$ is a smooth  even function.
Since $M$ is complete we have that $\cE_k(M)$ is also complete (see \cite{BisOn}, p.23).\\

For instance, if  $M=\HH^{l}$ then the hyperbolic extension $\cE_k(\HH^{l})$ is 
hyperbolic $(k+l)$-space $\HH^{k+l}=\HH^k\times\HH^{l}$, with metric $(cosh^2\, r)\,\sigma\0{\HH^{k}}+\sigma\0{\HH^{l}}$ (see item f above).\\

For a subset $A\sbs\HH^k$ we shall write $\cE_A(M)= A\times M\sbs\cE_k(M)$, with
the metric $\cE_k(h)$ restricted to the set $\HH^k\times A$.\\

We will write $\HH^k=\HH^k\times \{ o\}\sbs\cE_k(M)$. The hyperbolic extension, away from 
$\HH^k\sbs \cE_k(M)$ can be described in an alternative way: it is isometric to
$(\HH^k\times \bS^{n-1})\times (0,\infty)$ with doubly warped variable metric $(cosh^2r)\,\sigma\0{\HH^k}+h_r+dr^2$.\\

Note that $\HH^k$ and every $\{y\}\times M$ are convex in $\cE_k(M)$ (see \cite{BisOn}, p.23).
Let $\eta$ be a complete geodesic line in $M$ passing though $o$
and let $\eta^+$ be one of its two geodesic rays (beginning at $o$) . Then $\eta$ is 
a totally geodesic subspace of $M$ and $\eta^+$ is convex (see item (ii) above). Also, let $\gamma$ be a complete geodesic line in $\HH^k$.\\

\noindent {\bf Lemma 2.2.1.} {\it 
We have that\,\, $\gamma\times \eta^+$ is a convex subspace of $\cE_k(M)$
and $\gamma\times \eta$ is totally geodesic in $\cE_k(M)$.}\\

\noindent {\bf Proof.} Let $\pi_\gamma:\HH^k\ra\gamma\sbs\HH^k$ denote the orthogonal projection, and note that
$d\pi_\gamma$ is distance non-increasing, i.e. $\sigma\0{\HH^k}(v,v)\geq \sigma\0{\HH^k}
(d\pi_\gamma(v),d\pi_\gamma(v))$, for $v\in T\HH^k$. Moreover, the equality holds if and only if $v\in T\gamma$.\\

We assume $\eta:[0,\infty)\ra\eta\sbs M$ to be parametrized by the arc-length, that is, it is a speed-one geodesic ray.
Let $\pi_\eta:M\ra\eta^+$ denote the proper map $\pi_\eta(p)=\eta(r(p))$ . Note that $\pi_\eta$ is
smooth on $M-\{ o\}$ and  item (iii) above implies that $\pi_\eta$ is
strictly distance decreasing on non-radial tangent vectors on $M-\{ o\}$.\\

 Let $\alpha:[0,1]\ra\HH^k\times M$ and write $\alpha(u)=(a(u),b(u))\in\HH^k\times M$. Assume $b(u)\neq o$
for all $u\in [0,1]$. Let $\beta=(\pi_\gamma\, a,\pi_\eta\, b)$. Since $r(b(u))=r(\pi_\eta\, b (u))$ we have that
the length of $\beta$ is greater than the length of $\alpha$, unless $\alpha=\beta$. And, by continuity the same holds
without the assumption that $b(u)\neq o$. Therefore $\gamma\times\eta^+$ is convex because $\beta$
is a path in $\gamma\times\eta^+$. \\

We prove that $\gamma\times \eta$ is totally geodesic. Let $\eta^-=\overline{\eta-\eta^+}$ be the
``other" geodesic ray of $\eta$. Then $\gamma\times\eta^-$ is also convex.
Therefore $\gamma\times \eta-\gamma$ is totally geodesic hence the second fundamental of $\gamma\times\eta$ vanishes there. By continuity this form vanishes on the whole of $\gamma\times\eta$. This proves the lemma.\\

\noindent {\bf Corollary 2.2.2.} {\it We have that \,$\HH^k\times\eta^+$  and $\gamma\times M$ are convex in $\cE_k(M)$. Also  \,$\HH^k\times\eta$
is totally geodesic in $\cE_k(M)$.}\\

\noindent {\bf Proof.} For $\HH^k\times\eta$ just replace $\gamma$ by $\HH^k$ and $\pi_\gamma$ by the identity in the proof of
Proposition 2.2.1. For $\gamma\times M$ replace $\beta$ in the proof of Proposition 2.2.1 by  $\beta=(\pi_\gamma\, a,  b)$. 
This proves the Corollary.\\

As before, we have that $\HH^k\times\eta$
is isometric to $\HH^k\times \R$ with warped metric $cosh^2 v\, \sigma\0{\HH^k}+dv^2$, which is just hyperbolic $(k+1)$-space $\HH^{k+1}$.\\

\noindent {\bf Remarks 2.2.3.}  \\
\noindent {\bf 1.}  Note that $\gamma\times\eta$
is isometric to $\R\times \R$ with warped metric $cosh^2 v\, du^2+dv^2$, which is just hyperbolic 2-space $\HH^2$.
In particular $\cE_k(M)$ is also complete and every point in $\HH^k=\HH^k\times\{ o\} \sbs\cE_k(M)$ is a center point.\\
\noindent {\bf 2.} It follows from lemma 2.2.1 and remark 1 that
the ray structure of $\cE_k(h)$ with respect to any center $o\in\HH^k
\sbs\cE_k(M)$ only depends on the ray structure of $M$ and the center $o$.\\
\noindent {\bf 3.} Denote by $\B_r(M)$ the ball of radius $r$ of $M$.
Note that if $h$ and $h'$ on $M$ have the same ray structures
then the balls $\B_r(M)$ coincide.\\
\noindent {\bf 4.} Recall that $\HH^k$ is convex in $\cE_k(M)$.
Moreover, for $l\leq k$, we also have $\HH^l\sbs\HH^k\sbs\cE_k(M)$ is
convex. If $h$ and $h'$ on $M$ have the
same ray structures then the $r$-neighborhoods of the convex subsets
$\HH^l$ coincide.

\vspace{.8in}

\noindent {\bf 2.3. Coordinates on $\cE_k(M)$.} \\
Recall that we are identifying $M-\{ o\}$ with $\bS^{n-1}\times \R^+$ and we are denoting a point
$v=(u,r)\in \bS^{n-1}\times\R^+=M-\{ o\}$ as $v=ru$.
Fix a center $o\in \HH^k\in \cE_k(M)$.  Then, for $y\in\HH^k-\{ o\}$ we can also write $y=t\,w$, $(w,t)\in \bS^{k-1}\times\R^+$. 
Similarly, using the exponential map we can identify $\cE_k(M)-\{ o\}$
with $\bS^{k+n-1}\times \R^+$, and for $p\in\cE_k(M)-\{ o\}$
we can write $p=s\,x$, $(x,s)\in\bS^{k+n-1}\times\R^+$.\\

As before denote the metric on $\cE_k(M)$ by $g$ and we can write $g=g_s+ds^2$.
Since $\HH^k$ is convex in $\cE_k(M)$ we can write $\HH^k-\{ o\}=\bS^{k-1}\times \R^+\sbs\bS^{k+n-1}\times \R^+$
and $\bS^{k-1}\sbs \bS^{k+n-1}$.
\\

A point $p\in\cE_k(M)\, -\, \HH^k$ has two sets of coordinates: the {\it polar coordinates}
$(x,s)=(x(p),s(p))\in \bS^{k+n-1}\times \R^+$ and the {\it hyperbolic extension coordinates} $(y,v)=(y(p), v(p))\in \HH^k\times M$. Write $M_o=\{o\}\times M$.
Therefore we have the following functions:
$$
\begin{array}{lll}
{\mbox{the distance to {\it o} function:}}  & s:\cE_k(M)\ra [0,\infty), & s(p)=d\0{\cE_k(M)}(p,o)\\\\
{\mbox{the direction of {\it p} function:}}  & x:\cE_k(M)-\{o\}\ra \bS^{n+k-1}, & p=s(p)x(p)\\\\
{\mbox{the distance to {\it $\HH^k$} function:}}  & r:\cE_k(M)\ra [0,\infty), & r(p)=d\0{\cE_k(M)}(p,\HH^k)\\\\
{\mbox{the projection on $\HH^k$ function:}}  & y:\cE_k(M)\ra \HH^k, &\\\\
{\mbox{the projection on $M$ function:}}  & v:\cE_k(M)\ra M, & \\\\
{\mbox{the projection on $\bS^{n-1}$ function:}}  & u:\cE_k(M)-\HH^k\ra \bS^{n-1}, & v(p)=r(p)u(p)\\\\
{\mbox{the length of $y$ function:}}  & t:\cE_k(M)\ra [0,\infty), & t(w)=d_{\HH^k}(y,o)\\\\
{\mbox{the direction of $y$ function:}}  & w:\cE_k(M)-M_o\ra \bS^{k-1}, & y(p)=t(p) w(p)
\end{array}
$$

Note that $r=d_M(v, o)$. Note also that, by 2.2.1, the functions $w$ and $u$ are constant on geodesics emanating from $o\in\cE_k(M)$, that is
$w(sx)=w(x)$ and $u(sx)=u(x)$.\\

Let $\p_r$ and $\p_s$ be the gradient vector fields of $r$ and $s$, respectively. Since the $M$-fibers $M_y=\{ y\}\times M$ are convex
the vectors $\p_r$ are the velocity vectors of the speed one geodesics of the form $a\mapsto (y, a\, u)$, $u\in\bS^{n-1}\sbs M$. These geodesics
emanate orthogonally from $\HH^k\sbs \cE_k(M)$.
Also the vectors  $\p_s$ are the velocity vectors of the speed one geodesics 
emanating from $o\in\cE_k(M)$. For $p\in\cE_k(M)$, denote by $\bigtriangleup =\bigtriangleup (p)$ the right triangle with vertices $o$, $y=y(p)$, $p$
and sides the geodesic segments $[o,p]\in\cE_k(M)$, $[o,y]\in\HH^k$, $[p,y]\in\{ y\}\times M\sbs\cE_k(M)$.
(These geodesic segments are unique and well-defined because:\, (1) $\HH^k$ is
convex in $\cE_k(M)$,\, (2) $(y,o)=o_{_{\{ y\}\times M}}$ and $o$ are centers in $\{ y\}\times M$ and $\HH^k\sbs\cE_k(M)$, respectively.)\\

\noindent {\bf Lemma 2.3.1} {\it Let $\eta^+$  (or $\eta$) be a geodesic ray (line) in $M$ through $o$ containing
$v=v(p)$ and $\gamma$ a geodesic line in $\HH^k$ through $o$ containing $y=y(p)$. Then $\bigtriangleup (p)\sbs \gamma\times \eta^+\sbs \gamma\times \eta$.}\\

\noindent {\bf Proof.} We have that $[o,v]\sbs\eta$ and $[o, y]\in\gamma$. By lemma 2.2.1 we have $[o,p]\in \gamma\times \eta^+$.
This proves the lemma.\\

Let $\alpha:\cE_k(M)-\HH^k\ra \R$ be the angle between 
$\p_s$ and $\p_r$ (in that order), thus  $cos\, \alpha=g(\p_r,\p_s)$, $\alpha\in [0,\pi]$. 
Then $\alpha=\alpha(p)$ is the interior angle, at $p=(y,v)$, of the right triangle $\bigtriangleup =\bigtriangleup (p)$.
We call  $\beta(p)$ the interior angle of this triangle at $o$, that is $\beta(p)=\beta(x)$ is the spherical distance 
between $x\in \bS^{k+n-1}$ and the totally geodesic sub-sphere $\bS^{k-1}$. Alternatively, $\beta$ is the angle between the geodesic segment
$[o,p]\sbs\cE_k(M)$ and the convex submanifold $\HH^k$.
Therefore $\beta$ is constant on geodesics emanating from $o\in\cE_k(M)$, that is
$\beta(sx)=\beta(x)$.\\ 

Note that the right geodesic triangle $\bigtriangleup (p)$ has sides of length $r=r(p)$, $t=t(p)$ and $s=s(p)$. By lemma 2.3.1 and remark 2.2.3
we can consider $\bigtriangleup$ as contained in hyperbolic 2-space.
Hence using hyperbolic trigonometric identities
we can find relations between $r$, $t$, $s$, $\alpha$ and $\beta$. For instance, using the hyperbolic law of cosines we get:\\

\noindent {\bf (2.3.2.)}\hspace{1.5in}
$cosh\, (s)\, =\, cosh\, (r)\,\, cosh\, (t)$\\

\noindent Note that this implies $t\leq s$.
Here is an application of this equation.\\

\noindent {\bf Proposition 2.3.3 (Iterated hyperbolic extensions)} {\it We have that $$\cE_l\big( \cE_k(M)   \big)=\cE_{l+k}(M)$$
where we are identifying $\HH^{l+k}$ with $\HH^l\times\HH^k$ with warped metric $(\cosh^2t)\,\sigma_{\HH^l}+\sigma\0{\HH^k}$.}\\

\noindent {\bf Remarks.} 

\noindent {\bf 1.} Note that the identification depends on the order of the dimensions $l$ and $k$, that is, on the order in which
the hyperbolic extensions are taken.

\noindent {\bf 2.} As before, here the function $t:\HH^k\ra [0,\infty) $ is the distance in $\HH^k$ to the  point $o\in\HH^k$. \\

\noindent {\bf Proof.} As above let $s:\HH^k\times M\ra [0,\infty)$ be the distance in $\cE_k(M)$ to $o$,
$r(p)=d_M(v(p),o)$, and $t$ as in the statement of the Proposition.
Then $\cE_l\big( \cE_k(M)   \big)$ is $\HH^l\times (\HH^k\times M)$ with metric $$(cosh^2 s)\,\sigma_{\HH^l}+\big[  (cosh^2r)\,\sigma\0{\HH^k} +h \big]$$
On the other hand $\cE_{l+k}(M)$ is  $(\HH^l\times \HH^k)\times M$ with metric $$(cosh^2 r)\,\big[(cosh^2t)\, \sigma_{\HH^l}+  \sigma\0{\HH^k}\big] +h $$

\noindent Hence the Proposition follows from identity (2.3.2) above. This proves the Proposition.\\

\noindent {\bf Proposition 2.3.4.} {\it  We have the following identity defined outside $\HH^k\cup \big(\{o \}\times M$\big)}

$$
\big( sinh^2 (s) \big)\, d\beta\,^2\,\,\, +\,\,\, ds^2\,\,\, =\,\,\, cosh^2 (r) \, dt^2\,\,\,+\,\,\, dr^2
$$

\noindent {\bf Proof.} First a particular case. Take $M=\R$ and $k=1$, hence $\cE_k(M)=\cE_1(\R)=\HH^2$. In this case the left-hand side of
the identity above is the expression of the metric of $\HH^2$ in polar coordinates $(\beta, s)$, and right hand side of the equation is
the expression of the same metric in the hyperbolic extension coordinates $(r,t)=(v,y)$. (Here $r$ and $t$ are ``signed" distances.)
Hence the equation holds in this particular case.\\

Now, the general case can be reduced to this particular case using lemma 2.2.1 and remark 2.2.3. This proves the proposition.\\

A direct (and longer) proof of the lemma above can be given using hyperbolic trigonometric identities.

\vspace{.8in}

\noindent {\bf 2.4. Coordinates on the spheres $\bS^\circ_s\big( \cE_k(M)\big)$.} \\
The open geodesic ball of radius $r$ on a complete Riemannian manifold
$N^n$ centered at $o$ is denoted by $\B_r=\B_r(N)$ and we can identify $\B_r-\{ o\}$ with $\bS^{n-1}\times (0,r)$.
The geodesic sphere of radius $r$ centered at $o$ is denoted by $\bS_r=\bS_r(N)$ and we can also identify $\bS_r$ with $\bS^{n-1}\times\{ r\}$.\\

Let $M$ have center $o$ and metric $h$ as before.  Consider the hyperbolic extension $\cE_k(M)$ of $M$
with center $o\in\HH^k=\HH^k\times \{ o\}\sbs \cE_k(M)$ and metric $g$.
Since $\HH^k\sbs \cE_k(M)$ is convex, we can write $\bS_s(\cE_k(M))\cap\HH^k=\bS_s(\HH^k)$. Equivalently 
$\big(\bS^{k+n-1}\times\{s\}\big)\cap \HH^k=\bS^{k-1}\times\{s\}$.\\

Recall $M_o=\{ o\}\times M$. Write

$$
\cE^\circ_k(M)\, =\, \cE_k(M)\, -\, \big( \HH^k\, \coprod \, M_o  \big)
$$

and

$$
\bS^\circ_s\big(\, \cE_k(M)\, \big)\, =\,\bS_s\big( \cE_k(M) \big)\, \bigcap\, \cE^\circ_k(M)\,=
\, \bS_s\big(\, \cE_k(M)\, \big)\, -\, \big( \HH^k\, \coprod \, M_o  \big)
$$\

Note that the functions $\alpha$ and $\beta$ are well-defined and smooth on $\cE^\circ_k(M)$. Also we have that $0<\beta(p)<\pi/2$.
Moreover, for $p\in \cE^\circ_k(M)$ there is a unique well-defined 2-dimensional totally geodesic subspace of $\cE_k(M)$ of the form
considered in lemma 2.2.1, namely $\gamma_w\times\eta_u$, where $\gamma_w=\R w(p)\sbs \HH^k$ and $\eta_u=\R u(p)\sbs M$.
We write $\HH^2(p)=\gamma_w\times\eta_u$. Lemma 2.3.1 implies $\bigtriangleup (p)\sbs\HH^2(p)$. By the identification between $\bS^{n+k-1}\times\{s\}$ with $\bS_s(\cE_k(M))$
and lemma 2.2.1 we have that $\HH(p)\cap\bS_s(\cE_k(M))$ gets identified with a
geodesic circle $\bS^1(p)\sbs\bS^{n+k-1}$. 
Moreover, since $\HH^2(p)$ and $\HH^k$ intersect orthogonally on $\gamma_w$, we have that the spherical geodesic segment 
$\big[ x(p), w(p)\big]\0{\bS^{n+k-1}}$
intersects $\bS^{k-1}\sbs\bS^{n+k-1}$ orthogonally at $w$. This together with the fact that $\beta<\pi/2$ imply that $\big[ x(p), w(p)\big]_{\bS^{n+k-1}}$ is a length minimizing spherical geodesic in $\bS^{k+n-1}$
joining $x$ to $w$. Consequently $\beta$ is the length of $\big[ x(p), w(p)\big]\0{\bS^{n+k-1}}$. We will denote this spherical geodesic segment also
by $\beta=\beta(p)$.\\

We now give a set of coordinates on $\bS^\circ_s\big(\cE_k(M)\big)$. For 
$p\in \bS^\circ_s\big(\cE_k(M)\big)$  define 

$$
\Xi\, (p)\,=\, \Xi_s(p)\,=\,\big(\,  w\, ,\,  u  \, ,\,  \beta  \,   \big)\,\in\, \bS^{k-1}\times\, \bS^{n-1}\times \big( 0\, ,\, \pi/2  \big)
$$ \

\noindent where $w=w(p)$, $u=u(p)$, $\beta=\beta(p)$.
Note that $\Xi$ is constant on geodesics emanating from $o\in\cE_k(M)$, that is $\Xi(sx)=\Xi(x)$.\\

Using hyperbolic trigonometric identities we can find well-defined and smooth functions $r=r(s,\beta)$ and $t=t(s,\beta)$ such that
$r$, $s$, $t$ are the lengths of the sides of a right geodesic triangle on $\HH^2$ with angle $\beta$ opposite the the side with length $r$.
With these functions we can construct explicitly a smooth inverse to $\Xi$.\\

\noindent {\bf 2.4.1 Remark. } We have $r=sinh^{-1}\big( sin\, (\beta)\,\,sinh\, (s) \big)$ and 
$t=cosh^{-1}\bigg( \frac{cosh\, (s)}{\sqrt{1+sin^2(\beta)\,sinh^2(s)}}  \bigg)$.\\

We have the following facts.

\begin{enumerate}
\item[{\bf 1.}] For $(w,u)\in \bS^{k-1}\times\bS^{n-1}$ we have 
$$\Xi\bigg(  \big(\gamma_w\times\eta_u\big) \,\cap\, \bS^\circ_s\big(\cE_k(M)\big)\bigg)\,=\,\{ \pm w\}\times\{\pm u\}\times (0,\pi/2)$$
\noindent which are the images of the paths $a\mapsto (\pm w,\pm u, a)$. By lemma 2.2.1 we have that this set is a spherical geodesic segment emanating orthogonally
from $\bS^{k-1}$. 

\item[{\bf 2.}] For $w\in \bS^{k-1}$ we have 
$$\Xi\bigg(  \big(\gamma_w\times M\big)\,\cap\, \bS^\circ_s\big(\cE_k(M)\big) \bigg)\,=\,\{\pm w\}\times\bS^{n-1}\times (0,\pi/2)$$ 
\noindent By lemma 2.2.2 we have that this set is a spherical geodesic ball of radius $\pi/2$ and of dimension $n$ (with its center deleted)
intersecting  $\bS^{k-1}$ orthogonally at $w$. Note that the geodesic segments on this ball emanating from $w$ are the spherical
geodesic segments of item 1, for all $u\in \bS^{n-1}$.

\item[{\bf 3.}] For $w\in \bS^{k-1}$ and $r$ with $0<r< s$ we have 
$$\Xi\big(  \gamma_w\times \bS_r(M) \big)\, =\, \{ w\}\times\bS^{n-1}\times \beta(r)$$
\noindent where $\beta(r)$ is the angle of the right geodesic hyperbolic triangle with sides of length $s$ (opposite to the right angle)
and $r$, opposite to $\beta$. The angle $\beta$ can be computed using hyperbolic trigonometric identities, namely
$\beta=sin^{-1}\big(\frac{sinh (r)}{sinh(s)}  \big)$.

\item[{\bf 4.}]Since the $M$-fibers $\{y\}\times M$ are orthogonal in $\cE_k(M)$ to the $\HH^k$-fibers $\HH^k\times\{ v\}$, items 1,2, and 3 above
imply that the $\bS^{k-1}$-fibers, the $\bS^{n-1}$-fibers and $(0,\pi/2)$-fibers are mutually orthogonal 
in $\bS^{k-1}\times\, \bS^{n-1}\times \big( 0\, ,\, \pi/2  \big)$ with the metric $\Xi_* g$.

\item[{\bf 5.}]
The map $$\Xi'=(\Xi, s):\cE^\circ_k(M)\ra \bS^{k-1}\times\, \bS^{n-1}\times \big( 0\, ,\, \pi/2  \big)\times\R^+$$ 
\noindent gives coordinates on $\cE^\circ_k(M)$.

\end{enumerate}

\vspace{1in}

\noindent {\bf \large  Section 3. Approximating Metrics.}\\

\noindent {\bf 3.1. The basic local hyperbolic model and $\epsilon$-hyperbolic metrics.}\\ Let $\bB=\bB^n\sbs\R^n $ be the unit ball, 
with the flat metric $\sigma\0{\R^n}$.  Write  $I_\xi=(-(1+\xi),1+\xi)$, $\xi>0$.
Our basic models are $\T^{n+1}_\xi=\T_\xi=\bB\times I_\xi$,  with hyperbolic metric $\sigma=e^{2t}\sigma\0{\R^n}+dt^2$. 
In what follows we may sometimes suppress the sub index $\xi$, if the context is clear.
The number $\xi$ is called the {\it excess} of $\T_\xi$.\\

\noindent {\bf Remarks.}

\noindent {\bf 1.}  There are two main reasons to introduce the model spaces $\T_\xi$. First we want the ``size" of the model
to be sufficiently big so that we have enough space to glue metrics locally while keeping the gluing functions bounded.
On the other hand, a couple of our constructions reduce the size of the
charts in the $t$-direction (which depends on $\xi$).
The first construction is hyperbolic extension and
we will prove in section 3.5 that the hyperbolic extension of an $\epsilon$-hyperbolic space is $\epsilon'$-hyperbolic, but the excess of
the charts is reduced by a controllable small amount.
The second construction is hyperbolic forcing (see section 4) and
it reduces the excess of the charts by one.

\noindent {\bf 2.} In the applications we will actually need warped metrics with warping functions that are multiples of hyperbolic
functions. All these functions are close to the exponential $e^t$ (for $t$ large), so instead of introducing one model for each hyperbolic function
we introduced only the exponential model. The estimates given in section 3.6 will be used to change the exponential by the hyperbolic functions.\\

Let $|.|_{C^k}$ denote the uniform $C^k$-norm of $\R^l$-valued functions on $\T_\xi=\bB\times I_\xi\sbs\R^{n+1}$, and we will write $|.|=|.|_{C^2}$.
Given a metric $g$ on $\T$, $|g|_{C^k}$ is computed considering $g$ as the $\R^{(n+1)^2}$-valued function $(x,t)\mapsto (g_{ij}(x,t))$ where, as usual,
$g_{ij}=g(e_i,e_j)$, and the $e_i$'s are the canonical vectors in $\R^{n+1}$. We will say that a metric $g$ on $\T$ is $(\epsilon, C^k)$-{\it hyperbolic}
if $|g-\sigma|_{C^k}<\epsilon$.\\

\noindent {\bf Lemma 3.1.1.} {\it Let $g_i$ be $(\epsilon_i, C^2)$-hyperbolic
on $\T_\xi$, for $i=1,\,2$. Let $\lambda:\T_\xi\ra [0,1]$ be smooth with $|\lambda|$
finite. Then
the metric $\lambda\, g_1\,+\, (1-\lambda)\,g_2$ is $\Big(4\,|\lambda|\,(\epsilon_1+\epsilon_2),\, C^2\Big)$-hyperbolic.}\\

\noindent {\bf Proof.} The proof follows from the triangular inequality,
Leibniz rule and the equality $\Big(\lambda\, g_1\,+\, (1-\lambda)\,g_2\Big)-\sigma=
\lambda\,(g_1-\sigma)+(1-\lambda)(g_2-\sigma)$.
This proves the lemma.\\

A Riemannian manifold $(M,g)$ is $(\epsilon, C^k)$-{\it hyperbolic} if there is $\xi>0$ such that for every $p\in M$ there is an $(\epsilon, C^k)$-hyperbolic 
chart with center $p$, that is, there is a chart
$\phi :\T_\xi\ra M$, $\phi(0,0)=p$,  such that $\phi^*g$ is $(\epsilon, C^k)$-hyperbolic. Note that all charts are defined on the same model space
$\T_\xi$. The number $\xi$ is called the {\it excess } of the charts
(which  is fixed).
More generally, a subset $S\sbs M$ is $(\epsilon, C^k)$-hyperbolic if every $p\in S$ is the center of an $(\epsilon, C^k)$-hyperbolic chart in $M$ with fixed
excess $\xi$. \\

\noindent {\bf Remark.} 
Note that $\HH^n$ is $(\epsilon, C^k)$-hyperbolic for every $\epsilon >0$, but if a
hyperbolic manifold is not complete or have very small injectivity radius then it is not  $(\epsilon, C^k)$-hyperbolic.\\

Consider $M\times I$ with variable metric $g=g_t+dt^2$. We say that a subset $S$ of $(M\times I ,g)$ is {\it warped} $(\epsilon, C^k)$-{\it hyperbolic}
if every $p\in S$ is the center of a {\it warped $(\epsilon, C^k)$-hyperbolic  chart of fixed excess $\xi$}, that is, by charts
$\phi :\T_\xi\ra M$, $\phi(0,0)=p$,
 such that $\phi^*g$ is $(\epsilon, C^k)$-hyperbolic and $\phi$ respects the ``structure" of $g$, i.e. $\phi$ satisfies the following two
conditions:

\begin{enumerate}
\item[{\bf (i)}] the map $\phi$ respects the $\B$-direction, i.e. $\phi \,(\bB\times\{ t\})\sbs
M\times\{ t'\}$

\item[{\bf (ii)}] the map $\phi$ preserves the $t$-direction, i.e. $\phi\, (\{ x\}\times I_\xi)\sbs
\{ x'\}\times I$ and $\phi |_{(\{ x\}\times I_\xi)}$ is an isometry for all $x\in\B$.
\end{enumerate}

Equivalently, the chart has the form $\phi(x,t)=(\phi_1(x),\, t+a)$, for some fixed number $a\in\R$ and chart $\phi_1$ on $M$. 
Note that, in this case, the pullback metric has the form $\phi^*g=\phi^*g_t+dt^2$, ie. it is a variable metric on $\T$.
Note also that every warped $\epsilon$-hyperbolic chart $\phi$ can be extended to $\phi:\B\times (I-\{ a\})$ by the same formula $\phi (x,t)=(\phi_1(x),t+a)$,
but this extension may fail to be $\epsilon$-hyperbolic. (Here $I-\{a\}=\{t-a,\, t\in I\}$.)
Of course a  warped $(\epsilon, C^k)$-hyperbolic manifold is $(\epsilon, C^k)$-hyperbolic.\\

To simplify our notation, in what follows we will always consider $k= 2$, and drop the $C^k$ symbol, unless it is necessary in the discussion.
Almost all of the results that will be proved hold for $k\geq2$.\\

\noindent {\bf Example.} Consider $M^{n+1}=\bS^n\times\R^+$ with warped metric $sinh^2t\,\sigma\0{\bS^n}+dt^2$. Hence $M$ is isometric to
a punctured hyperbolic space (see item (e) in section 2). As stated in the remark above $M$ is not $\epsilon$-hyperbolic,
but $S_L=\bS^n\times (L,\infty)\sbs M$ is $\epsilon$-hyperbolic, provided $L>2$. On the other hand $S_L$ is not warped $\epsilon$-hyperbolic 
when $\epsilon$ is small.
But the following is true: given $\epsilon>0$ there is $L>0$ such that $S_l$  is $\epsilon$-hyperbolic, for $l\geq L$. Actually, a more general 
statement is proven below (see 3.4.4).\\

The next lemma says that the ``radius" of the chart is bounded by
$(2+\xi)+n^2\,\epsilon$.\\

\noindent {\bf Lemma 3.1.2.} {\it Let $\phi:\T_\xi\ra M$ be a warped 
$\epsilon$-hyperbolic chart centered at $p\in M$. Then, for every $q\in \T_\xi$ we have}

$$
d\0{M} (\phi(q),p)\,\leq\,  (2+\xi)+n^2\,\epsilon$$

\noindent {\bf Proof.} Write $q=(x_0,t_0)\in\B\times I_\xi$. Consider the path
$\alpha(t)=(t\,x_0,0)$, $t\in [0,1]$, $\beta(t)=(x_0,t\,t_0)$, $t\in [0,1]$, and
$\gamma=\alpha\ast\beta$.
Write $g'=\phi^*g$ and we have $g'=\sigma+h$, with $|h|<\epsilon$. Then the $g$-length $\ell_g(\phi\circ\gamma)$ of $\phi\circ\gamma$ is
$$\ell_{g'}(\gamma)=\ell_{g'}(\alpha)+\ell_{g'}(\beta)\leq
\ell_\sigma(\alpha)+\ell_h(\alpha)+(1+\xi)\leq
1+\epsilon n^2+ (1+\xi)$$
\noindent hence $d\0{M} (\phi(q),p)\,\leq\,\ell_{g}(\phi\circ\gamma)\leq\,
(2+\xi)+n^2\epsilon$. This proves he lemma.

\vspace{.8in}

\noindent {\bf 3.2. Sectional curvatures.} \\
For every $n$ there is a function $\epsilon'=\epsilon'(\epsilon,\xi)$ with the following  property:  if a Riemannian metric $g$ on a 
manifold $M^{n+1}$ is $\epsilon'$-hyperbolic, with charts of excess $\xi$,  then the sectional curvatures of $g$ all lie $\epsilon$-close to -1.
This choice is possible, and depends only on $n$ and $\xi$, because the curvature depends only of the derivatives up to order 2 of $\phi^*g$ on $\T_\xi$,
where $\phi$ is an $\epsilon$-hyperbolic chart with excess $\xi$.\vspace{.6in}

\noindent {\bf 3.3. Slow variable metrics.}\\
Consider the family of metrics $g\0{t}$, $t\in I\sbs \R$, on $M$, $M$ closed. We say that 
$\{ g\0{t}\}$ is {\it  $\epsilon$-slow} if $g\0{t}$ and its derivatives in the $M$-direction change $\epsilon$-slowly with
$t$ (up to order 2 and 1, respectively). That is:

\begin{enumerate}
\item[(i)] for every $t_0\in I$, \,$k=1,\,2$ and \,$u\in TM$ we have
 $\Big|\frac{d^k}{d\,t^k}\,g\0{t}(u,u)|\0{t_0}\,\Big|\,<\,\epsilon\,\, g\0{t_0}(u,u)$
\item[(ii)] for every $t_0\in I$, $v\in TM$ and   $u$ vector field on $M$ we have 
\end{enumerate}
$$\Big|\,\frac{d}{d\,t}\,v\,g\0{t}(u,u)|\0{t_0}\,\Big|\,<\,\epsilon\,\Big( \,g\0{t_0}(u,u)\,
g\0{t_0}^{{\mbox{\tiny 1/2}}}(v,v)\,+\,
g\0{t_0}^{{\mbox{\tiny 1/2}}}(u,u)\,\,g\0{t_0}^{{\mbox{\tiny 1/2}}}(\n _vu,\n_vu)\,\Big)$$\\

For instance, if the family $\{ g\0{t}\}$ is constant, then it is $\epsilon$-slow, for every $\epsilon>0$. We will  need the following lemma later, which follows from the chain rule.\\

\noindent {\bf Lemma 3.3.1.} {\it Consider the family of metrics $\{ g\0{t}\}_{t\in I}$, $I\sbs \R$ an interval, and assume
$\{ g\0{t}\}$ is $\epsilon$-slow. Let $\varphi: J\ra I$ be a diffeomorphism with $|\varphi'|,\,|\varphi''|<a$. Then
the family $\{ g_{\varphi(s)}\}_{s\in J}$ is $\big(\,(a+a^2)\epsilon)\,\big)$-slow.}\\


A set of metrics $\{ g\0{\lambda}\}$ on the compact manifold $M$ is {\it c-bounded} if there is $c >1$ such 
that $|g\0{\lambda}|< c$ and $|\, det \,g\0{\lambda}\,|_{C^0} > 1/c$, for all $\lambda$. The set $\{g\0{\lambda}\}$ is bounded if it is $c$-bounded for some $c >1$ \\

\noindent {\bf Remarks.}

\noindent {\bf 1.} Here the uniform $C^2$-norm $| .|$ is taken with respect to a fixed locally finite atlas $\cA$ with ``extendable" charts, i.e.
charts that can be extended to the (compact) closure of their domains.

\noindent {\bf 2.} For a metric $g$ on $M$ we will use the same symbol $g$ to denote a matrix representation in a chart of the
fixed finite atlas.

\noindent {\bf 3.} We are taking $c >1$ just to simplify  the statement of lemma 3.3.3.

\noindent {\bf 4.} If $\{ g\0{t}\}_{t\in I}$ is a (smooth) family and $I\sbs \R$ is compact then, by compacity and continuity, the family
 $\{ g\0{t}\}$ is bounded.

\noindent {\bf 5.} If $\{ g\0{t}\}_{t\in I}$ is $c$-bounded then clearly
$\{ g\0{t(s)}\}_{s\in J}$ is also $c$-bounded, for any reindexation (or reparametrization)
$t=t(s)$.\\

The variable metric $g=g\0{t}+dt^2$ on $M\times I$ is $\epsilon$-slow, $c$-bounded or bounded if the family $\{ g\0{t}\}$ is $\epsilon$-slow, $c$-bounded
or bounded, respectively. 
For $I\sbs\R$ we write $I(\xi)=\Big\{\,  t\in I,\,  \big(\,t-(1+\xi), t+(1+\xi)\,\big)\sbs I\,\Big\}$.  We also consider $inf\, I$ to be positive. We write 
$T=inf \,I(\xi)=1+\xi+inf\, I$.
Hence $t\geq T$, for all $t\in I(\xi)$.\\

\noindent {\bf Proposition 3.3.2.} {\it Let $M^n$ be a closed smooth manifold
and $\xi>0$.
If the variable metric $g=g\0{t}+dt^2$ on $M\times I$ is $\epsilon$-slow and $c$-bounded then
the variable warped metric $h=e^{2t}g\0{t}+dt^2$ is warped $\eta$-hyperbolic on $M\times I(\xi)$
with charts of excess $\xi$, provided $$C\,(e^{-T}+\epsilon)\leq\eta$$ 
\noindent where $C=C(c,n, \xi)$.}\\

\noindent {\bf Remark.} The constant $C$ depends solely on the dimension $n$
of $M$, the constant $c$ and the desired excess $\xi$. An explicit formula for $C$ is given
at the end of the proof of proposition 3.3.2. Note that if $T$ is large and $\epsilon$ small
then we can take $\eta$ small.\\

\noindent {\bf Corollary 3.3.3.} {\it Under the same conditions as in 3.3.2 we have that 
the variable warped metric $h=sinh^2t\,g\0{t}+dt^2$ is warped $\eta$-hyperbolic on $M\times I(\xi)$
with charts of excess $\xi$, 
provided $$\,C_1\,(e^{-T}+\epsilon)\leq\eta$$ 
\noindent where $C_1=C_1(c,n,\xi)$.}\\

\noindent {\bf Remark.} We can take $C_1=128\, C(2^{4n}c,n,\xi)$, where
$C$ is as in proposition 3.3.2.\\

We will use the following lemmas.\\

\noindent {\bf Lemma 3.3.4.} {\it Let $\{ g\0{\lambda}\}$ be a $c$-bounded family of metrics on the unit ball $\bB^n\sbs\R^n$.
Then we can write $g\0{\lambda}= F_\lambda^TF_\lambda$, where $|F_\lambda|_{C^0}, |F_\lambda^{-1}|_{C^0}<\, \sqrt{n\,n!\,c^{n+1}}$.}\\

\noindent {\bf Proof.} Since $g\0{\lambda}$ is symmetric we can write $g\0{\lambda}=O^T D O$, where the columns of $O^{T}$ form an orthonormal
basis (in $\R^n$) of eigenvectors of $g\0{\lambda}$ and $D$ is a diagonal matrix whose diagonal entries are the corresponding eigenvalues $\mu$ of $g\0{\lambda}$. 
Since $|g\0{\lambda}|_{C^0}<c$ and $det \,g\0{\lambda} > 1/c$ we have that 
$\frac{1}{n!\,c^{n+1}} \, <\, \mu \, <\, n c$. Take $F=\sqrt{D}O$. This proves the lemma.\\

\noindent {\bf Remark.} Let $u\in\R^n$, and let $|u|$ be the euclidean norm.
The following estimate follows from the proof of 3.3.4:
$$
\frac{1}{n!\, \,c^{n+1}}\,\,|u|\,<\, g\0{\lambda}(u,u)\,<n\,c\, |u|
$$\\

\noindent {\bf Lemma 3.3.5.} {\it  Let $\{ g\0{t}\}$ be an $\epsilon$-slow $c$-bounded family of metrics on the unit ball $\bB^n\sbs\R^n$.
Then (see remark 2 above)}
\begin{enumerate}
\item[]  $|\frac{\p^k}{\p t^k}\big( g\0{t}\big)\0{ij}|< 3\,\epsilon \, c $,\,\, $k=1,\,2$.
\item[]  $|\frac{\p^2}{\p t\,\p x\0{l}}\big( g\0{t}\big)\0{ij}|< \epsilon \, c\0{3} $
\end{enumerate}
\noindent {\it where $c\0{3}=c\0{3}(c,n)=4\,c^{3/2}+\frac{27}{4}\big(n!\big)^2
c^{^{{\mbox{\tiny $\frac{4n+11}{2}$}}}}$.}\\

\noindent {\bf Proof.} We prove the first inequality for $k=1$. The proof for $k=2$
is similar. To simplify our notation write $g'\0{ij}=\frac{\p}{\p t}(g\0{t})\0{ij}$.
We will also sometimes omit the variable $t_0$.
From the definition we have that 
$|\frac{\p}{\p t} g\0{t}(u,u)(t_0)|< \epsilon \,  g\0{t_0}(u,u)$. Taking $u=\p_i$ 
and using the fact that the family is $c$-bounded we get that $|g'\0{ii}|<\epsilon\, c$.
And taking $u=\p_i+\p_j$ we get $|g'\0{ii}+2g'\0{ij}+g'\0{jj}|<
\epsilon |g\0{ii}+2g\0{ij}+g\0{jj}|\leq 4\epsilon\,c$. These two inequalities imply
$$
|2g'\0{ij}|\leq |g'\0{ii}+2g'\0{ij}+g'\0{jj}|+|g'\0{ii}|+|g'\0{jj}|<\,6\,\epsilon\,c
$$
\noindent This proves the first inequality.  We prove the second inequality. Since
the family of metrics is $c$-bounded we have that $|g^{ij}|<(n-1)!\,c^{n+1}$, where $(g^{ij})$ is the inverse of $g=(g\0{ij})$. Hence we obtain the following estimate for
the Christoffel symbols

\begin{equation*} \Big|\Gamma^k\0{ij}\Big|\,\,<\,\,\frac{3}{2} \,(n-1)!\, c^{n+2} 
\tag{a}
\end{equation*}

\noindent In what follows of this proof we use the summation notation. It
follows from (a) above that

{\small \begin{equation*}
g\Big( \n_{\p_l}\p_i\,,\, \n_{\p_l}\p_j \Big)\,\,=\,\, g\Big( \Gamma\0{li}^k\p_k,
\Gamma\0{lj}^m\p_m \Big)\,\,=\,\, g\0{km}\Gamma\0{li}^k\Gamma\0{lj}^m
\,\,<\,\, \frac{9}{4}\,\big(n!\big)^2\,c^{2n+5}\,=\, c\0{1}
\tag{b}
\end{equation*}}

\noindent Write $g'\0{ij,l}=\frac{\p^2}{\p t\,\p x\0{l}} g\0{ij}$.
Now take $v=\p_l$ and $u=\p_i$ in (ii) of the definition of $\epsilon$-slow
metrics, and use (b) to obtain
\begin{equation*}
\Big| g'\0{ii,l}\Big|\,<\,\epsilon\, \Big(\, c\, c^{1/2}\, +\,
c^{1/2}\, c\0{1}\,  \Big)\,=\, \epsilon\, c\0{2}
\tag{c}
\end{equation*}
\noindent where $c\0{2}=c\0{2}(c,n)= c^{3/2}+c\,c\0{1}=c^{3/2}+
\frac{9}{4}\,\big(n\big)^2\,c^{^{{\mbox{\tiny $\frac{4n+11}{2}$}}}}$.
Take $v=\p_l$,  $u=\p_i+\p_j$ now in (ii) in the definition of $\epsilon$-slow
metrics and use (b), (c) to obtain
$$\begin{array}{lll}
\Big| 2g'\0{ij,l} \Big|&\leq&
\Big| g'\0{ii,l}+2g'\0{ij,l}+g'\0{jj,l} \Big|\,+\,\Big| g'\0{ii,l} \Big|\,+\,
\Big| g'\0{jj,l} \Big|\\  \\&<& \epsilon\,\Big(\,c^{1/2}\,(4\,c)\,+\,2\, c^{1/2}\,c\0{1} 
\,+\,2\,c\0{2}\Big)
\,\,=\,\, 2\,\epsilon\, c\0{3}
\end{array}
$$

\noindent This proves the lemma.\\

We shall prove a sort of a converse to this lemma at the end of section 3.3.\\ 

\noindent {\bf Proof of Proposition 3.3.2.} First we reduce the problem to $\R^n$ using the fixed finite atlas: 
for each $p\in M$ choose a chart $\psi=\psi_p=(U,\varphi)$ such that (1) $\psi$ is the restriction of one of the charts in the
fixed finite atlas (2) $\varphi(U)=\bB(0, \mu_0)$, for some fixed $\mu_0>0$,  (3) $\psi (p)=0$. For simplicity we assume $\mu_0=1$,
that is $\varphi(U)=\bB$.\\

\noindent \Big[ The charts of the atlas $\cA$  (see remark afer 3.3.1)
can be composed with a dilation to obtain $\mu_0=1$. In general 
the constant $C$ in proposition 3.3.2 would depend on 
$\mu_0$.\Big]\\

Fix $\xi>0$, $p\in M$ and $t_0\in I(\xi)$. Thus $(t_0-(1+\xi),t_0+(1+\xi))\sbs I$. Write $g\0{t_0}=F^TF$ with $F$ as in lemma 3.3.4.
(Here $g\0{t_0}$ denotes also the matrix representation of $g\0{t_0}$ in the chart $\xi_p$.) Let $A=F^{-1}$.
Then $|A|_{C^0}< c\0{4}$, where $c\0{4}=\sqrt{n\,n!\, c^{n+1}}$. In what follows we identify via $\psi$ the neighborhood $U$ of $p$ with the unit ball $\bB$.
Define the chart $\phi :\T=\bB\times I_\xi\ra \bB\times I\sbs M\times I$ by
$$ \phi (x,t)=(e^{-t_0}A x, t+t_0)$$

\noindent Write $f=\phi^*h$ and  $f=e^{2t}f_t+dt^2$, where $f_t=e^{2t_0}\,\phi^*g\0{t+t_0}$. Then 
\begin{equation*} f_t(x)=A^T\, g\0{t+t_0}(e^{-t_0}Ax)\, A\tag{1}\end{equation*}

\noindent Hence we have that
\begin{equation*} f(0,0)=\phi^*h(0,t)=e^{2t} I=\sigma (0,0)\tag{2}
\end{equation*}

\noindent where $I$ is the identity matrix. Differentiating equation (1) and using lemma 3.3.5 and the fact that
$|A|_{C^0}< c\0{4} $ we get the following estimates:

\begin{equation*}  |\p_J f|=e^{2t}|\p_J f_t|< e^{2t}\, n^2\, c^2\0{4}\, \Big[ n\, c\,
e^{-t_0}c\0{4}\Big]^{|J|} \tag{3}\end{equation*}

\noindent where $J$ is a multi-index of order $|J|=1,\,2$ in the $\bB$-direction, i.e. no $t$-derivatives are considered.
Also, from lemma 3.3.5  and (1) we get 

\begin{equation*} |\frac{\p^k}{\p t^k}\big( f_t\big)_{ij}|< 3\,n^2\,c^2\0{4}\,c\,  \epsilon\,=\,c\0{5}\tag{4}\end{equation*}

\begin{equation*} |\frac{\p^{2}}{ \p x\0{l}\p t}\big( f_t\big)_{ij}|< \epsilon \,\bigg[ n^3\,
c\0{3}\,c\0{4}^3\, e^{-t_0}\bigg]\,=\,c\0{6} 
\tag{5}\end{equation*}

\noindent From (2) and (4), for every $t$  we get
\begin{equation*} |f_t(0)-I|<\,6\,n^2\,c\0{4}\, c \,\epsilon\,=\,c\0{7}
\tag{6}\end{equation*}

\noindent Write $c\0{8}=n^3\,c\0{4}^3\,c\,e^{-t_0}$. Then from (3) we have
$|\frac{\p}{\p x\0{l}}f_t|<c\0{8}$.
This together with (6) and the mean value Theorem imply
$$
|f(x,t)-\sigma(x,t)|\,\leq\,e^{2t}\bigg( |f_t(x)-f_t(0)|+|f_t(0)-I| \bigg)\,<\,
e^{2(1+\xi)}\Big( n^{1/2}\,c\0{8}\,+c\0{7} \Big)\,=\, c\0{9}
$$
\noindent Hence
\begin{equation*} |f-\sigma|_{C^0}\,<\, c\0{9}
\tag{7}\end{equation*}

\noindent Now, since the derivatives in the $\bB$-direction of $\sigma$ all vanish,  equation (3)
holds replacing $f$ by $f-\sigma$ and we get
\begin{equation*} \Big| \,\p^J(f-\sigma)\,\Big|\,<\,e^{2(1+\xi)}\, n^4\, c^4\0{4}\, c^2\,
e^{-t_0}\,=\, c\0{10}
\tag{8}\end{equation*}
\noindent where we are assuming $t_0>0$ so that $e^{-2t_0}<e^{-t_0}$.
Now, note also that $\pt (f_t-I)=\pt f_t$. This together with
the definitions of $f$ and $\sigma$ imply
$\pt (f-\sigma)=2(f-\sigma)+ e^{2t}\pt f_t$, which together with (4) and (7)
imply 
\begin{equation*} |\pt (f-\sigma)|\,<\, 2\,c\0{9}\,+\, e^{2(1+\xi)}\, c\0{5}\,=\, c\0{11}
\tag{9}
\end{equation*}

\noindent Analogously, differentiating with respect to $t$ again and using (4) we get
\begin{equation*} |\frac{\p^2}{\p t^2} (f-\sigma)|\,<\, 4\,c\0{9}\,+\, 4\, e^{2(1+\xi)}\,c\0{5}\,+\, e^{2(1+\xi)}\,c\0{5}\,=\, 4\,c\0{9}\,+\, 5\,e^{2(1+\xi)}\,c\0{5}\,=\, c\0{12}
\tag{10}
\end{equation*}
\noindent And from (5) and (8) we finally get

\begin{equation*} |\frac{\p^2}{\p x\0{l} \p t} (f-\sigma)|\, <\, 2\,c\0{10}\,+\, e^{2(1+\xi)}\, c\0{6} \,=\,c\0{13}
\tag{11}
\end{equation*}

\noindent Hence $|f-\sigma|<c\0{14}$, where $c\0{14}=$max$\{
c\0{9},\,c\0{10},\,c\0{11},\,c\0{12},\,c\0{13} \}\leq C\big(\,e^{-t_0}\,+\,\epsilon\big)$, and it can be verified that we can take
(assuming $n\geq 2$ and $t_0>0$)

$$C\,=\,C(c,n,\xi)\,=\,2\,e^{2(1+\xi)}\, n^3\, c\0{4}^3\, \Big( 4c^{3/2}+\frac{27}{4} (n!)^2 
c^{^{{\mbox{\tiny $\frac{4n+11}{2}$}}}}\Big)$$

\noindent This completes the proof of Proposition 3.3.2.\\\\

\noindent {\bf Proof corollary 3.3.3.}
Let $k(t)=\frac{sinh^2(t)}{e^{2t}}=\frac{(1-e^{-2t})^2}{4}$.
A quick calculation shows that for $t>1$ we have 
(note that we indeed have $t>1$ because $\xi>0$)

\begin{equation*}\begin{array}{l}
\frac{1}{16}\,<k\,<\frac{1}{4}\\\\
 0<\,k'\,<\,e^{-2t}\,<\,1/2\\\\
0<\,|k''|\,<\,2\,e^{-2t}\,<\,1
\end{array}
\tag{12}
\end{equation*} \\

We can write $h=sinh^2(t)\,g\0{t}+dt^2=e^{2t}\,k\,g\0{t}+dt^2$. Then the corollary
follows from the following claim.\\

\noindent {\bf Claim.} {\it The family of metrics $\{k\,g_t\}$ is
$128\,(e^{-2T}+\epsilon)$-slow and $(2^{4n}\,c)$-bounded.}\\

\noindent {\bf Proof of claim.} Using (12) we have 
$$|k\,g\0{t}|=k\,|g\0{t}|<\frac{1}{4}\,c\,<\, c$$
\noindent Also, $|\, det\,(k\,g\0{t})\,|\,=\,k^n\,|\,det\,g\0{t}\,|\,>\,\frac{1}{16^n}\,|\,det\,g\0{t}\,|\,>\,\frac{1}{16^n\,c}$.
Hence $\{g\0{t}\}$ is $(2^{4n}\,c)$-bounded. We prove now the first statement of the definition of $\epsilon$-slow metrics.
Using (12) and the fact that $\{g\0{t}\}$ is $\epsilon$-slow we get
\begin{equation*} \begin{array}{lllll}\Big( (k\, g)(u,u)\Big)'&=&k'\,g(u,u)\,+\, k\,g'(u,u)
\,\,\,\,<\,\,\,\,
\frac{1}{k}\,\big(k'\,+\, k\,\epsilon\big)\,(k\,g)(u,u)\\\\
&<&16\,\big(\,e^{-2T}+\epsilon\,\big)\,(k\,g)(u,u)\,\,\,\,<\,\,\,\,
128\,\big(e^{-2T} +\epsilon\big)\,(k\,g)(u,u)\end{array}\tag{13}\end{equation*}
\noindent and
\begin{equation*} {\small
\begin{array}{lllll}\Big( k\, g(u,u)\Big)''&=&k''\,g(u,u)\,+\, 2\,k'\,g'(u,u)\,+\,k\,g''(u,u)
\,\,\,\,<\,\,\,\,
\frac{1}{k}\,\big(k''\,+\,2\, k'\,\epsilon\,+\, k\,\epsilon\big)\,(k\,g)(u,u) \\ \\
&<&16\,\big(\,2\,e^{-2T}\,+\,\epsilon\,+\, \epsilon\,\big)\,(k\,g)(u,u)\,\,\,\,<\,\,\,\,128\,\big(\,e^{-2T} +\epsilon\,\big)\,(k\,g)(u,u)
\end{array}}
\tag{14}\end{equation*}
\noindent Also
\begin{equation*}{\small \begin{array}{lll}
\Big(v\,k\,g(u,u)\Big)'&=&k'\Big(v\,g(u,u)\Big)\,+\, k\,\Big(v\,g(u,u)\Big)'\\\\
&<&k'\Big(v\,g(u,u)\Big)\,+\,\epsilon\,k\,\Big(\,g(u,u)
\,g^{^{\mbox{\tiny $1/2$}}}(v,v)\,+\,
\,g^{^{\mbox{\tiny $1/2$}}}(u,u)\, g^{^{\mbox{\tiny $1/2$}}}
(\n_vu,\n_vu)   \Big)
\end{array}}
\tag{15}
\end{equation*}
\noindent But
\begin{equation*}
\Big|\,v\,g(u,u)\,\Big|\,=\,2\,\Big|\,g(u,\n_vu)\,\Big|\,\leq\,
2\,g^{^{\mbox{\tiny $1/2$}}}(u,u)\,\,g^{^{\mbox{\tiny $1/2$}}}(\n_vu,\n_vu)
\tag{16}
\end{equation*}
\noindent And from (15),(16) and (12) we get
\begin{equation*}{\small\begin{array}{lllllll}
\Big(v\,k\,g(u,u)\Big)'&<&\Big(\,2\,k'\,+\,k\,\epsilon\Big)
\Big(\,g(u,u)
\,g^{^{\mbox{\tiny $1/2$}}}(v,v)\,+\,
\,g^{^{\mbox{\tiny $1/2$}}}(u,u)\, g^{^{\mbox{\tiny $1/2$}}}
(\n_vu,\n_vu)   \Big)&&&\\\\&<&
\Big(\,2\,k'\,+\,k\,\epsilon\Big)
\Big(\,g(u,u)
\,g^{^{\mbox{\tiny $1/2$}}}(v,v)\,+\, k^{-1/2}
\,g^{^{\mbox{\tiny $1/2$}}}(u,u)\, g^{^{\mbox{\tiny $1/2$}}}
(\n_vu,\n_vu)   \Big)&&&\\\\
&<&k^{-3/2}\,
\Big(\,2\,k'\,+\,k\,\epsilon\Big)
\Big(\,k^{3/2}\,g(u,u)
\,g^{^{\mbox{\tiny $1/2$}}}(v,v)\,+\, k
\,g^{^{\mbox{\tiny $1/2$}}}(u,u)\, g^{^{\mbox{\tiny $1/2$}}}
(\n_vu,\n_vu)   \Big)&&&\\\\
&=&64\,
\Big(2\,k'+k\,\epsilon\Big)
\Big((k\,g)(u,u)
\,(k\,g)^{^{\mbox{\tiny $1/2$}}}(v,v)+ 
\,(k\,g)^{^{\mbox{\tiny $1/2$}}}(u,u)\, (k\,g)^{^{\mbox{\tiny $1/2$}}}
(\n_vu,\n_vu)   \Big)&&&
\\\\&<&128\,\big(\,e^{-2T} +\epsilon\,\big)\,
\Big((k\,g)(u,u)
\,(k\,g)^{^{\mbox{\tiny $1/2$}}}(v,v)+ 
\,(k\,g)^{^{\mbox{\tiny $1/2$}}}(u,u)\, (k\,g)^{^{\mbox{\tiny $1/2$}}}
(\n_vu,\n_vu)   \Big)&(17)&&
\end{array}}
\end{equation*}

\noindent And the claim follows from (13), (14) and (17). This proves the claim.\\

The corollary now follows from proposition 3.3.2, the claim and the fact that
(recall $T>1$) $$
e^{-T}\,+\,128\,(\,e^{-2T}\,+\, \epsilon\,)\, <\,128\,(\,e^{-T}\,+\,\epsilon\,)
$$
\noindent This proves corollary 3.3.3.
\vspace{.5in}

Here is a simple particular case of Proposition 3.3.2.  Assume that the family $\{ g\0{t}\}$ is constant, i.e. $g\0{t}=g\0{0}$, for all $t\in I$, 
then the variable metric
$g=g\0{t}+dt^2=g\0{0}+dt^2$ is just a product metric, and the warped variable metric $h=e^{2t}g\0{0}+dt^2$ is just a simply warped metric.
In this case Proposition 3.3.2 says that, given $\epsilon >0$, the warped metric $h$ is $\epsilon$-hyperbolic, provided $t$ is large
enough (how large depending on $\epsilon$, the dimension $n$ and the metric $g\0{0}$).\\

Here is a sort of a converse to lemma 3.3.5.\\

\noindent {\bf Lemma 3.3.6.} {\it Suppose that $\{g\0{t}\}$ is $c$-bounded and that}
\begin{enumerate}
\item[]  $|\frac{\p^k}{\p t^k}\big( g\0{t}\big)\0{ij}|< \epsilon\0{1} $,\,\, $k=1,\,2$.
\item[]  $|\frac{\p^2}{\p t\,\p x\0{l}}\big( g\0{t}\big)\0{ij}|< \epsilon\0{2} $
\end{enumerate}
\noindent {\it Then $\{ g\0{t}\}$ is $\epsilon\0{3}$-slow, where
$\epsilon\0{3}=\Big(\epsilon\0{2}\,n\,+\,2\,\epsilon\0{1}\, c\0{14}\,    \Big)\,n^2
\Big(n!\,c^{n+1}\Big)^3$, and $c\0{14}=\frac{3}{2}n^{3/2}n!c^{n+2}$.}\\

\noindent {\bf Remark.} Note that $\epsilon\0{3}=a\,\epsilon\0{1}\,+\,b\,\epsilon\0{2}$
where $a=a(n,c)$ and $b=b(n,c)$ are constants that depend solely on $c$ and $n$.\\

\noindent {\bf Proof.} 
First we prove (i) of the definition of $\epsilon$-slow metrics. We do this just for $k=1$
because the proof for $k=2$ is similar. We use the summation notation and denote
the derivative with respect to $t$ by a prime. Also the euclidean norm on $\R^n$
will be denoted by bars $|.|$. \\

Write $u=u^i\p_i$. Then using the remark after 3.3.3 we get 
$$g(u,u)'\,=\,(g\0{ij}u^iu^j)'\,=\,g\0{ij}'u^iu^i\,<\,\epsilon\0{1}n^2\,|u|^2\,<\,
\epsilon\0{1}\,n^2\,n!\, c^{n+1}\, g(u,u)$$
\noindent This proves (i). To prove (ii) we need the following estimate.
Let $u=u^k\p_k$ and $v=v^l\p_l$. We have
\begin{equation*}
\big| v(u^k) \big|\,\leq\,\big| v(u^k)\p_k \big|\,\leq\, \big|  \n_vu \big|\,+\,c\0{14} \big|v \big|\,
\big| u \big|
\tag{d}
\end{equation*}

\noindent where $c\0{14}=\frac{3}{2}n^{3/2}n!c^{n+2}$. To prove this note that 
$$\n_vu=\n_vu^k\p_k=v(u^k)\p_k+u^kv^l\n_{\p_l}\p_k=v(u^k)\p_k+
u^kv^l\Gamma^s_{kl}\p_s$$
\noindent and equation (d) follows from (a) in the proof of 3.3.5.\\

Now, to prove (ii) we use (d)  to compute:
$$\begin{array}{lllll}
\Big| \frac{d}{dt}v\,\big(g\0{ij}  \big)(u,u) \Big|&=&
\bigg[ v\,\Big( g\0{jk}u^ju^k  \Big)   \bigg]'& = &v^l\Big(\p\0{l}\,g\0{jk}\Big)'u^ju^k\,+
\,2\,g'\0{ij}\,u^j\,v(u^k)\\ &&&<&\epsilon\0{2}\, n^3\, |v|\,|u|^2\,+\,2\,
\epsilon\0{1}n^2\,|u|\,\Big(|\n_vu|+c\0{4}|u|\,|v|\Big)\\
&&&=&\Big(\epsilon\0{2}\,n+2\,\epsilon\0{1}\,c\0{14}\Big)\,n^2\,|v|\,|u|^2\,+\,
2\,\epsilon\0{1}\,n^2\,|u|\,|\n_vu|\\
\end{array}
$$

\noindent This together with the remark after 3.3.3 gives:
$$\Big| \frac{d}{dt}v\,\big(g\0{ij}  \big)(u,u) \Big|\,\,<\,\,
 A\,\,g\0{t_0}(u,u)\,
g\0{t_0}^{{\mbox{\tiny 1/2}}}(v,v)\,\,+\,\, B\,\,
g\0{t_0}^{{\mbox{\tiny 1/2}}}(u,u)\,\,g\0{t_0}^{{\mbox{\tiny 1/2}}}(\n _vu,\n_vu)
$$ 

\noindent where $A=\Big(\epsilon\0{2}\,n\,+\,2\,\epsilon\0{1}\, c\0{14}\,    \Big)\,n^2
\Big(n!\,c^{n+1}\Big)^3$ and $B=\Big(2\,\epsilon\0{1}\, n^2  \Big)\,
\Big( n!\, c^{n+1} \Big)^2$. Since $A>B$ the lemma follows.
This proves the lemma.\\

Recall we are assuming $M$ closed. If $I$ is compact then we can find
$\epsilon\0{1}$ and $\epsilon\0{2}$ as in  the statement of lemma 3.3.6.
Hence we obtain the following corollary.\\

\noindent {\bf Corollary 3.3.7.} {\it If $I$ is compact
then $\{g\0{t}\}$ is $c$-bounded and $\epsilon$-slow, for some $c>1$ and
$\epsilon>0$.}\\

Here is a sort of a local converse of proposition 3.3.2.\\

\noindent {\bf Corollary 3.3.8.} {\it Let $g=e^{2t}g\0{t}+dt^2$ be a variable metric on $\T_\xi$. Assume $g$ is $\epsilon$-hyperbolic. Then we have}

\begin{enumerate} 
\item[{\bf 1.}] {\it the family $\{ g\0{t}\}$ is $\epsilon'$-slow, where
$\epsilon'=a'\epsilon$, with $a'=a'(n,\xi)$ }

 \item[{\bf 2.}] {\it if $\epsilon<\frac{1}{2^{n+1}\,e^{2(1+\xi)}\,n!}$, then the family $\{ g\0{t} \}$ is $2$-bounded.}
\end{enumerate}

\noindent {\bf Remark.} We can take $a'=9e^{2(1+\xi)}\,(n+2\,c\0{4}(2))\,n^2\big(n!\,2^{n+1}  \big)^3$, where $c\0{4}(2)=\sqrt{n\,n!\,2^{n+1}}$.\\

\noindent {\bf Proof.} The metric $g=e^{2t}g\0{t}+dt^2$ is $\epsilon$-hyperbolic means 
$$e^{2t}\big|g\0{t}-\sigma\0{\R^n} \big|\,=\,\big|(e^{2t}g\0{t}+dt^2)-(e^{2t}\sigma\0{\R^n}+dt^2)\big|\, =\, \big|g-\sigma\big|\, <\,\epsilon$$
\noindent Thus the symmetric bilinear forms $g\0{t}-\sigma\0{\R^n}$ have $C^2$-norms bounded by $\epsilon\, e^{2(1+\xi)}$. Hence
$$ \big|g\0{t}\big|\, <\, \big|g\0{t}-\sigma\0{\R^n}\big|\, +\, \big|\sigma\0{\R^n}\big|\, \leq \,  (\epsilon \, e^{2(1+\xi)})\, +\, 1 \, < \,2 $$
\noindent Also, since we have that the matrices $g\0{t}$ can be written as $I+(\epsilon \, e^{2(1+\xi)})\,A$, with $|A|<1$ (here $| A|$ is the $C^0$-uniform norm of $A$).
We get that $det \,g\0{t}> 1-\epsilon\, 2^ne^{2(1+\xi)} n!>1/2$. \,\, This proves part 2.\\

Now, we can write $g\0{t}=\sigma\0{\R^n}+e^{-2t}\tau_t$, with $|\tau|<\epsilon$. 
Hence for $u\in \R^n$ we have
$$ g\0{t} (u, u)\, =\, |u|^2 + e^{-2t}\tau_t (u,u)
$$
\noindent where $|u|^2=\sigma\0{\R^n}(u,u)$. Therefore
$ g'\0{t} (u, u)\, =\,  e^{-2t} \tau'\0{t} (u,u)\, -\, 2\,e^{-2t}\tau\0{t} (u,u)$
which implies $|g'\0{t}|<3e^{2(1+\xi)}\epsilon$. Similar calculations yield
$|g''\0{t}|<9e^{2(1+\xi)}\epsilon$ and $|\frac{\p}{\p x\0{l}}g'\0{t}|<3e^{2(1+\xi)}\epsilon$.
The corollary now follows from proposition 3.3.6 by taking $c=2$ and
$\epsilon_1=\epsilon_2=9e^{2(1+\xi)}\epsilon$.
This completes the proof of proposition 3.3.6.\\

As we mentioned in 3.1 our definition of warped $\epsilon$-hypebolicity is not perfect.
For instance hyperbolic space is not warped $\epsilon$-hyperbolic near the
(chosen) center. The next result tells us how far we have to be from the center to
get $\epsilon$-hyperbolicity of hyperbolic space. It follows from corollary 3.3.3
and the fact that we can take $\epsilon=0$ in this particular case.\\

\noindent {\bf Corollary 3.3.9.} {\it Let $o\in \HH^{n+1}$. Then hyperbolic $(n+1)$-space
$\HH^{n+1}$ is warped $\epsilon$-hyperbolic (with respect to $o$) outside
$\B_a(\HH^{n+1})$, with charts o excess $\xi$, provided}
$$
C_1' \, e^{-a}\,\leq \, \epsilon
$$ 

\noindent {\it where $C_1'=C_1'(n,\xi)$.}\\

\noindent We can take $C_1'=C_1(c\0{\bS^n},n,\xi)$, with $C_1$ is as in 3.3.3.\\

We can rephrase this result in the following way.\\

\noindent  {\bf Corollary 3.3.10.} {\it Let $o\in \HH^{n+1}$. Then hyperbolic $(n+1)$-space
$\HH^{n+1}$ is warped $\epsilon$-hyperbolic (with respect to $o$) outside
$\B_a(\HH^{n+1})$, with charts o excess $\xi$, provided
$a\geq a\0{n+1}(\epsilon,\xi)$,
where  $a_{n+1}(\epsilon,\xi)=ln\Big(\frac{C_1'}{\epsilon}  \Big)$. Here $C_1'=C_1(c\0{\bS^n},n,\xi)$, with $C_1$ is as in 3.3.3.}\\

\vspace{.8in}

\noindent {\bf 3.4.  Hyperbolic forcing.} \\
Let $g\0{*}$ be a metric on the $n$-sphere $\bS^n$ and consider the warped metric
$g=sinh^2t\, g\0{*} +dt^2$ on $\bS^n\times \R^+$. We will use corollary 3.3.3. to change the metric $g$ to a new metric $h$ with
the following four properties 
\begin{enumerate}
\item[(1)] the metric $h$ is a warped (by $sinh\,t $) variable metric, 
\item[(2)] the metric $h$  coincides with $g$ for $t$ greater than
some number $b>0$, 
\item[(3)] the metric $h$ is $\epsilon$-hyperbolic, 
\item[(4)] the metric $h=sinh^2t\, \sigma\0{\bS^n}+dt^2$ (hence hyperbolic) for $t$ less than some positive number  $a<b$.
\item[(5)] the ray structure of $h$ (see section 1?) coincides with the that of $g$, that is
the rays $t\mapsto t\,u$, $u\in\bS^n$ are $h$-geodesic rays and these rays are
perpendicular to the spheres $\bS^n(r)=\{x\,,\, |x|=r\}$. 
\end{enumerate}

\noindent We will show this is possible provided $a$ and $d=b-a$ are large (how large depending on $n$, $\epsilon$ and  $g\0{*}$). That is, we are forcing the
space to be hyperbolic near the ``cone vertex".\\

We fix a function $\rho:\R\ra[0,1]$ with $\rho (t)=0$ for $t\leq 0$ and $\rho(t)=1$ for $t\geq 1$.  It can be checked that we can find
such a $\rho$ with the following properties: (i) $|\rho'|< 3$, (ii) $|\rho''|<12$.
Given positive numbers $a$ and $d$ define $\rho\0{a,d}(t)=\rho(2\,\frac{t-a}{d})$. Then
$|\rho'\0{a,d}|< 6/d$ \, and \, $|\rho''\0{a,d}|<48/d^2$.\\

Also fix and atlas $\cA\0{\bS^n}$ on $\bS^n$ as in 3.3.
All norms and boundedness constants will be taken with respect to
this atlas.\\

For simplicity write $\sigma\0{0}=\sigma\0{\bS^n}$.
From corollary 3.3.7 we have that there are $\epsilon\0{0}=\epsilon\0{0}(g\0{*})$ and $c=c(g\0{*})$\, (or $\epsilon\0{0}=\epsilon\0{0}(g)$ and $c=c(g)$)\,
such that $\big\{\sigma\0{0}+s(g\0{*}-\sigma\0{0})\big\}\0{s\in[0,1]}$ is $\epsilon\0{0}$-slow
and $c$-bounded. By lemma 3.3.1, remark 5 after 3.3.1, and the fact that 
$|\rho'\0{a,d}|,\,|\rho''\0{a,d}|<6/d$ (assuming $d\geq 8$) we  have that\\

\noindent {\bf Remark.} Let $c\0{\bS^n}$ be a fixed constant such that $\sigma\0{0}$ is
$c\0{\bS^n}$-bounded. Also let $c\0{*}$ be such that $g\0{*}$ is
$c\0{*}$-bounded. Then we can take $c=c\0{*}+c\0{\bS^n}$.\\

\noindent {\bf (3.4.1.)}\,\,\,\, The family of metrics $\big\{ \sigma\0{0} +\rho\0{a,d}(t) (g\0{*}-\sigma\0{0})\big\}\0{t\in\R^+}$ is $c$-bounded and  $\Big(\frac{12}{d}\,\epsilon\0{0}\Big)$-slow.\\

\noindent Write $g\0{t}=\sigma\0{0} +\rho\0{a,d}(t) (g\0{0}-\sigma\0{0})$ and define the {\it forced partially hyperbolic} metric
$$
\cH_{_{a,d}}\, g\, =\, sinh^2\, t\,\,g\0{t}+dt^2
$$
\noindent (Sometimes we may write $\cH\0{a,d}\,g\0{*}$ instead of $\cH\0{a,d}\,g$). Hence, by 3.4.1, corollary 3.3.3 and the definition of
$\cH_{_{a,d}}\, g$ we get the following corollary.\\

\noindent {\bf Corollary 3.4.2.} {\it Let $\xi>0$, $r>1+\xi$ and $a>d$. Outside the ball $\B_{r}(\R^{n+1})$
the warped variable metric $\cH_{_{a,d}}\, g\,$ is $\epsilon$-hyperbolic
with charts of excess $\xi$, provided} $$\,C_1\,\Big(e^{-r}+\frac{12}{d}\,\epsilon\0{0}\Big)\leq\epsilon$$
\noindent {\it where $C_1=C_1(c,n,\xi)$ (see 3.3.3), $\epsilon\0{0}=\epsilon\0{0}(g)$,
$c=c(g)$.
Moreover the metric $\cH_{_{a,d}}\, g\,$  satisfies the five properties 
(1)-(5) mentioned at the beginning of this section. }\\

Actually the constant $\epsilon\0{0}=\epsilon\0{0}(g\0{*})$ depends only on the constant
$c=c(g\0{*})$. This is stated in the next result.\\

\noindent {\bf Lemma 3.4.3.}  {\it If $g\0{*}$ is $c\0{*}$-bounded then we can take}
$$\epsilon\0{0}=\epsilon\0{0}(c)\,=\,C_2\,\big(\,n\,, \,c\,+\,c\0{\bS^n} \big)$$\\
\noindent {\it Here $c\0{\bS^n}$ is (a fixed constant) such that $\bS^n$ is 
$c\0{\bS^n}$-bounded, and the function $C_2(n,x)$ is given by
$$C_2(n,x)\,=\, x\, \Big( a(n,x')+b(n,x')\Big)$$
\noindent where $a$ and $b$ are as in 3.3.6 (see also the remark after the statement of 3.3.6), and $x'=\big[n!\,x^{n+1}\big]^n$.}\\

\noindent {\bf Proof.} Write $\sigma=\sigma\0{\bS^n}$ and $c'=c\0{\bS^n}$.
Also, to simplify the notation write $c=c\0{*}$ in this proof.
We want to apply proposition 3.3.6 to the family
$\{g\0{s}\}\0{s\in [0,1]}$, where $g\0{s}=\sigma+s(g\0{*}-\sigma)=
(1-s)\sigma+sg\0{*}$.  First we need the following claim.\\

\noindent {\bf Claim.} {\it The family $\{g\0{s}\}\0{s\in [0,1]}$ is
$c''$-bounded, where $c''=\big[n!\,(c+c')^{n+1}\big]^n$.}\\

\noindent First note that (recall $c>1$)
$$|g\0{s}|_{C^2}\,\leq\, |g\0{*}|_{C^2}\,+\,|\sigma|_{C^2}\,<\, 
c+c'\,<\,c''$$
\noindent It remains to prove that $det\, g\0{s}>\frac{1}{c''}$.
We will use the following fact:\\

\noindent {\it Fact: Let $A$ be symmetric positive definite. Then
$u^TAu>d\,|u|^2$, for all $u$, if and only if all eigenvalues of $A$ are $>d$.}\\
  
From the proof of lemma 3.3.4 we get that all eigenvalues of
$g\0{*}$ are $>\frac{1}{n!c^{n+1}}$ and 
all eigenvalues of $\sigma$ are $>\frac{1}{n!(c')^{n+1}}$.
Hence all eigenvalues of either $g\0{*}$ or $\sigma$ are $>e=\frac{1}{n!(c+c')^{n+1}}$. 
Therefore (using the fact above) we get $h(u,u)> e\, |u|^2$, for every $u\in\R^n$, where $h$ is either
$g\0{*}$ or $\sigma$. Then
$$
g\0{s}(u,u)\,=\,(1-s)\,\sigma(u,u)\,+\,s\,g\0{*}(u,u)\,>\,e\,|u|^2
$$
\noindent which implies, by the fact above, that all eigenvalues of $g\0{s}$
are $>e$. Thus $det \,g\0{s}>e^n=c''$. This proves the claim.\\

To finish the proof of the lemma apply lemma 3.6.6 the the family
$\{g\0{s}\}\0{s\in [0,1]}$, which, by the claim, is $c''$-bounded.
A simple calculation shows that in this case we can take $\epsilon\0{1}$
and  $\epsilon\0{2}$ in lemma 3.3.6 satisfying
 $\epsilon\0{1}=\epsilon\0{2}=|\sigma|_{C^2}+|g\0{*}|_{C^2}<c+c'$.
This proves the lemma.\\

In the next corollary we use the constants $a\0{n}$ introduced in 3.3.10.
We also use the following constants: \\

\noindent {\bf (3.4.4.)}\hspace{.7in}$\begin{array}{lll}d_{n+1}(c\0{*},\epsilon,\xi)&=&
\frac{12\,\epsilon\0{0}(c\0{*})\,C_1(c\0{*},n,\xi)}{\epsilon}+(2+4\xi)\\\\
a'\0{n+1}(c\0{*},\epsilon,\xi)&=&ln\big(\frac{C_1(c\0{\bS^n}+c\0{*},n,\xi)}{\epsilon}\big)+(2+2\xi)\end{array}$\\

\noindent where $C_1$ is as 3.3.3.\\

\noindent {\bf Corollary 3.4.5.} {\it  Let $\xi>0$. Then
the warped variable metric $\cH_{_{a,d}}\, g\,$ is $(2\epsilon)$-hyperbolic
with charts of excess $\xi$, outside the ball of radius $a\0{n+1}(\epsilon,\xi)$,
provided we take} 

$$a\,>\, a'\0{n+1}(c\0{*},\epsilon,\xi)\hspace{.4in} and\hspace{.4in}  d\,>\, d\0{n+1}
(c\0{*},\epsilon,\xi)$$
\noindent {\it Moreover the metric $\cH_{_{a,d}}\, g\,$  satisfies the five properties 
(1)-(5) mentioned at the beginning of this section. }\\

\noindent {\bf Proof.} Write $r=a-(1+\xi)$. Then $a\,>\, a'\0{n+1}(\epsilon,c\0{*},\xi)$
implies:
(1) $e^{-r}<\frac{\epsilon}{C_1(c\0{\bS^n}+c\0{*},n,\xi)}
$,\, (2) $r>1+\xi$. Note that $r>a\0{n+1}(\epsilon,\xi)+(1+\xi)$, because
(it can easily be checked that) $C_1$ is increasing in the first variable.
Hence we can apply 3.4.2 to get that outside the ball of radius $r$
the metric is warped $(2\epsilon)$-hyperbolic. Inside the ball of radius $r$
(and outside the ball of radius $a\0{n+1}(\epsilon,\xi)$) we can apply
3.3.10, because the metric is hyperbolic on the ball of radius
$r+(1+\xi)=a$. This proves the corollary.\\

\noindent {\bf Remark.} Hyperbolic forcing is a generalization and a more detailed version of (and was motivated by) the
``Farrell-Jones warping trick" in \cite{FJ1}, which has the ``stretching of the cylinder $\bS^n\times[0,1]$" (given here by $d$) and the
far-away constant (given here by $a$)  parametrized by just one variable
$\alpha$. Proposition 3.3.2 shows to what extend the stretching and the far-away constant are independent. In fact, in next sections we will keep the stretching $d$ fixed but sufficiently large and then take the far-away constant $a$ as large as we will need
it to be for some other reasons. Another detail is worth mentioning. In \cite{FJ1} it is showed that the sectional curvatures
of the metrics obtained are $\delta$-closed to -1, while Proposition 3.3.1 shows that the spaces actually approximate, locally,
hyperbolic space, in the sense of being $\epsilon$-hyperbolic, which implies (see 3.2.)
that the sectional curvatures
of the metrics obtained are $\delta$-closed to -1.

\vspace{.8in}

\noindent {\bf 3.5. The hyperbolic extension of an $\epsilon$-hyperbolic metric.} \\
Let $(M^n,h)$ have center $o$ as above and consider the hyperbolic extension $\cE_k(M)$ (see section 2). As before the metric $\cE_k(h)$ on $\cE_k(M)$ is denoted by $g$.
The ball of radius $r$ on $M$ centered at $o$ will be denoted by $B_r$,
that is $B_r=\B_r(M)$.
Choose $o\in\HH^k\sbs \cE_k(M)$.
Recall that $o$ is a center of $\cE_k(M)$ (see 2.2.3), hence we can express $\cE_k(M)-\{ o\}$ as $\bS^{n+k-1}\times\R^+$ with
variable metric $g=g_u+du^2$, $u>0$. Recall we are assuming $\xi>0$.\\

\noindent {\bf Proposition 3.5.1.} {\it Let $\xi>1$ and
assume $r>7+3\xi$. Let $h$ be warped $\epsilon$-hyperbolic
on $S\sbs M^n-\{ o\}$ with charts of excess $\xi$. Then $g$
is warped $\eta$-hyperbolic on $\HH^k\times \big(S- B_r\big)$, provided}
$$
C_3\,\Big(\,\epsilon\,+\, e^{-r}\,\Big)\,\,\leq\,\, \eta
$$

\noindent {\it where $C_3=C_3(k,\xi)$.  Moreover,  $g$
is warped $\eta$-hyperbolic with charts of excess $\xi'$, provided} $$0\,<\,\xi'\,<\,\xi\,-\,e^{-\big(r-(7+3\xi)\big)}$$\\

\noindent {\bf Remarks.} 

\noindent {\bf 1.} Note that in the statement above 
$h$ is warped $\epsilon$-hyperbolic with respect to the decomposition $M-\{ o\}=\bS^{n-1}\times\R^+$ and
$g$ is warped $\eta$-hyperbolic  with respect to the decomposition $\cE_k(M)-\{ z\}=\bS^{n+k-1}\times\R^+$.

\noindent {\bf 2.} We can take $C_3(k,\xi)=e^{^{11+4\xi}}\, +\,C_4(k,\xi)$, where $C_4(k,\xi)=
\,C_1(c\0{\bS^k},k,\xi)$, with $C_1$ is as in corollary 3.3.3, and $c\0{\bS^k}$
is as in lemma 3.4.3 and is taken with respect to the atlas obtained by using the exponential map.\\

Before we prove the Proposition we shall give a particular type of charts on $\HH^{k+1}=\cE_k(\R)$.
For this we consider the following two ways of describing $\HH^{k+1}$. First we can express $\HH^{k+1}$ as the
the hyperbolic extension $\cE_k(\R)$, which is $\HH^k\times\R$ with metric $(cosh^2r)\, \sigma\0{\HH^k}+dr^2$. 
Hence a point $z\in\HH^{k+1}$ has coordinates $(y,r)
\in\HH^k\times\R $ and we call $(y,r)$ the $\cE_k(\R)$-{\it coordinates of} $z$ (see 2.3).\\

Also, as mentioned in item (e) of section 2.
1, we can express  the space 
$\HH^{k+1}-\{ o\}$ as $\bS^k\times \R^+$ with warped metric $(sinh^2\, t)\,\sigma_{\bS^k}+dt^2$.
Thus a point $z\in\HH^{k+1}-\{o\}$ has coordinates $(x,t)\in\bS^{k}\times\R^+$ and we call $(x,t)$ the {\it polar} coordinates of $z$.
Using the law of sines and the laws of cosines for right triangles in $\HH^2$ we can find transformation rules between the polar coordinates $(x,t)$
and the $\cE_k(\R)$-coordinates $(y,r)$. We are only interested in the explicit expression for $r=r(x,t)$. In this case we have\\

\noindent {\bf (3.5.2)}\hspace{1.3in}$ r(x,t)\, =\ sinh^{-1}\bigg(   sinh\, ( t)\, sin\, \beta (x) \bigg)$\\

\noindent where $\beta (x)$ is the spherical distance from $x\in \bS^k$ to the equator $\bS^k\cap \HH^k\sbs\HH^{k+1}$.\\

Let $z_0=(x_0,t_0)\in\bS^k\times (2,\infty)$ and let $(y_0,r\0{0})$ be the $\cE_k(\R)$-coordinates of $z_0$.
We define the chart $\psi=\psi_{z_0}:\T_\xi
\ra\bS^k\times\R^+=\HH^{k+1}-\{ o\}$ by
$$ \psi (x,t)\,=\,\bigg( exp_{x_0}\big[\frac{1}{sinh\, t_0}\, x\big] \, ,\,  t_0+t \bigg)
$$

\noindent where we are identifying the euclidean unit ball $\B^k$ with the unit ball in the tangent space $T_{x_0}\bS^k$,
and $exp_{x_0}:T_{x_0}\bS^k\ra\bS^k$ is the exponential map. Note that in the formula above the output of the map $\psi$ is given in polar coordinates.\\

\noindent {\bf Lemma 3.5.3.} {\it The chart $\psi$ is a warped  
$\epsilon$-hyperbolic chart, provided}$$C_4\,e^{-\Big[t_0-(1+\xi)\Big]}\,\leq\,\epsilon$$
\noindent {\it where $C_4=C_4(k,\xi)$ is as in remark 2 after proposition 3.5.1.}\\

\noindent {\bf Proof.} We claim that $\psi$ is exactly the chart used
in corollary 3.3.3, which is the chart $\phi$ in the proof of proposition
3.3.2 for the metric $k\,\sigma\0{\bS^k}$, where $k=\frac{sinh^2t}{e^{2t}}$ (see the proof of 3.3.3). To see this just recall that the derivative of the
exponential map $exp\0{x\0{0}}$ at $x\0{0}$ is the identity
hence the matrix $A$ in the proof of 3.3.2 is $\frac{e^{2t\0{0}}}{sinh\,t\0{0}}\,I$.  Also note that the formula above follows
because the family of metrics $\{\sigma\0{\bS^k}\}$ is constant
hence ``zero"-slow (i.e $\epsilon$-slow for every $\epsilon>0$).
Therefore, by the proof of 3.3.3 the family $\{k\,\sigma\0{\bS^k}\}$
is $(e^{-2T}+\epsilon)$-slow, for every $\epsilon>0$,
where, in our case, we have $T=t\0{0}-(\xi+1)$. This proves the lemma.\\

Denote the $\cE_k(\R)$-coordinates functions of \, $\psi$ by \,$y=y_{z_0}:\T_1\ra \HH^k$ and \,$\br=\br_{z_0}:\T_1\ra \R$. That is
$$\psi (x,t)\,\,=\,\,\big( y(x,t)\, ,\, \br(x,t)\big)\, \in\, \HH^k\times \R\, =\, \cE_k(\R)
$$

\noindent Using equation (3.5.2) above we can write\\

\noindent {\bf (3.5.4.)}\hspace{1.1in} $\br(x,t)\, =\ sinh^{-1}\bigg(   sinh\, ( t_0+t)\, sin\, \beta (x') \bigg)$\\

\noindent where $x'=exp_{x_0}\big[\frac{1}{sinh\, t_0}\, x\big] $. \\ \\

\noindent {\bf Lemma 3.5.5.} {\it We have that}\\

$$\bigg|\, \br (x,t)\,\, -\,\, \big(  t\, +\, r\0{0} \big)\,\bigg|_{C^2}\,\, \leq \,\, \frac{4\,e^{^{2(1+\xi)}}}{cosh \big(\,r\0{0}-(2+\xi)\,\big)}$$\\

\noindent {\it provided $r\0{0}\geq 5+2\,\xi$ (recall that $(y_0,r\0{0})$ are the $\cE_k(\R)$-coordinates of $z_0$).}\\

The proof of lemma 3.5.5 is given in Appendix A.\\

The next result is the reason why we introduced the variable $\xi$ in the definition of the models $\T_\xi$.\\

\noindent {\bf Corollary 3.5.6.} {\it We have that}

$$ \psi (\T_{\xi'})\,\sbs\, \HH^k\times \big[r\0{0}-(1+\xi)\, ,\, r\0{0}+(1+\xi)\big ] 
$$

\noindent {\it provided $0<\xi' <\xi-
\frac{4\,e^{^{2(1+\xi)}}}{cosh \big(\,r\0{0}-(2+\xi)\,\big)}$
and $r\0{0}>5+2\,\xi$}.\\

\noindent {\bf Proof.} Write $\kappa =
\frac{4\,e^{^{2(1+\xi)}}}{cosh \big(\,r\0{0}-(2+\xi)\,\big)}$. By lemma 3.5.5 we have $(t+r\0{0})-\kappa\leq \br(x,t)\leq (t+r\0{0})+\kappa$.
Since $t\in (-1-\xi', 1+\xi')$ we have $r\0{0}-(1+\xi'+\kappa)\leq \br(x,t)\leq r\0{0}+(1+\xi'+\kappa)$. Hence to prove the corollary we
need $\xi'+\kappa\leq \xi$. This proves the corollary.\\

\noindent {\bf Proof of proposition 3.5.1.} First some notation. 
Recall that we are denoting the metric on $M-\{ o\}$ by $h=h_r+dr^2$ and the one on $\cE_k(M)$ by $g$.
We write $I_\xi =[-1-\xi,1+\xi]$. For $u\in \bS^{n-1}$ we denote by
$Ru$ the complete geodesic line passing through $o$ with direction $u$, i.e. $Ru=exp_o (\R u)=\{ p\in M$ such that $p=ru$, i.e. 
$p$ has polar coordinates $(u,r),\, r\in\R \}$. Also write $R^+u=\exp_o(\R^+u)$. Then $Ru=R(-u)$ but $R^+u\cap R^+(-u)=\emptyset$.
Hence we get  $M-\{ o\}=\coprod _{u\in\bS^{n-1}} R^+u$. Therefore

\begin{equation*} 
\cE_k(M)-\HH^k\,\,\, =\,\,\, \HH^k\, \times \, \big(M-\{ o\}\big)\, \,\,=\, \,\,\coprod_{u\in\bS^{n-1}}\bigg( \HH^k\times R^+u\bigg)
\tag{1}\end{equation*}

\noindent Specifically, if $w=(y,p)\in \HH^k\times M=\cE_k(M)$, $p\neq o$, then $w\in \HH^k\times R^+u$, provided $p$ has polar coordinates $(u,r)$, for some $r>0$.\\

Since we are taking $u$ with length one,
we have an obvious identification of $Ru$ with $\R$, given by $r\mapsto ru$ (this identification does depend on the "sign" of $u$).
This identification gives a canonical (metric) identification of $\cE_k(Ru)=\HH^k\times Ru$ (with metric $g|_{\cE_k(Ru)}$) 
\, with\, $\cE_k(\R)=\HH^k\times\R=\HH^{k+1}$ (with the canonical warped metric).\\

Write $\cE_k(R^+u)=\HH^k\times R^+u\sbs\cE_k(Ru)$ and
note that we can canonically identify $\cE_k( R^+u)$\, with half hyperbolic $(k+1)$-space \, 
 $\HH^{k+1}_+=\HH^k\times \R^+\sbs \HH^k\times \R=\HH^{k+1}$.\\

For $r>0$ and $y\in \HH^k$ denote by $\bS_{r,y}$ the set $\{ (y,ru)\in \HH^k\times M,\,\,\, u\in\bS^{n-1}\,  \}$. Then $\bS_{r,y}$ is the geodesic
sphere of radius $r$ of the convex submanifold $\{ y\}\times M\sbs\cE_k(M)$. Note that every $\cE_k(R^+u)$ intersects every $\bS_{r,y}$ orthogonally in the single point $(y, ru)$.\\

Let $w_0\in \HH^k\times S \sbs\cE_k(M)$. Write $w_0=(y_0,p_0)\in \HH^k\times M$ and let $(u_0,r\0{0})$ be the polar coordinates of $p\in S\sbs M$.
Also let $t_0$ be the distance in $\cE_k(M)$ from $o$ to $w_0$.
Since $p_0\in S$ and $S$ is warped $\epsilon$-hyperbolic, there is an warped $\epsilon$-hyperbolic chart  $\phi : \T^{n+1}_\xi\ra M$ with center $p_0$.
For $(x,r)\in \T_\xi=\B^{n-1}\times I_\xi$ we can write 

\begin{equation*}
\phi\,(x,r)\,\,=\,\,\big(u(x), r+r\0{0}\big)
\tag{2}
\end{equation*}

Write $z_0=( y_0, r\0{0})\in\cE_k(\R^+u)$ and let $\psi=\psi_{z_0}=(y_{z_0}, \br_{z_0})=(y,\br)$ as defined before the beginning of the  proof of this proposition. 
We now define a chart $\bphi:\B^k\times\B^{n-1}\times \R^+\ra \cE_k(M)$ by

\begin{equation*}
\bphi \,(x_1, x_2, t)\,\,=\,\,\bigg(y(x_1,t)\, ,\,  \br(x_1,t)\,u(x_2)\bigg)\, \in \HH^k\times M
\tag{3}
\end{equation*}

\noindent Note that, by lemma 3.5.6 we have that

\begin{equation*}
\bphi\,\big(    \B^k\times\B^{n-1}\times I_{\xi'}     \big)\, \sbs\, \HH^k\times \phi (\T_\xi)
\tag{4}
\end{equation*}

\noindent provided $\xi'<\xi-\frac{4\,e^{2(1+\xi)}}{cosh(r\0{0}-(2+\xi))}$. By the definition of $\bphi$ (see equation (3)) we have

\begin{equation*}
\bphi \bigg(\{ x_1\}\times\B^{n-1}\times \{ t\}\bigg)\,\,=\,\, \phi \bigg( \B^{n-1}\times \{\br(x_1,t) \}  \bigg)\,\,\sbs\,\,
\bS_{\br(x_1,t), y(x_1,t)}
\tag{5}
\end{equation*}

\noindent and

\begin{equation*}
\bphi \bigg(\B^{k}\times \{ x_2\}\times I_{\xi'}\bigg)\,\,\sbs\,\,\cE_k\big(    R^+u(x_2) \big)
\tag{6}
\end{equation*}\\

\noindent Moreover, using (6), the metric canonical identification of $\cE_k(Ru)$ with $\cE_k(\R)=\HH^{k+1}$,
and the obvious identification of $\B^{k}\times \{ x_2\}\times I_{\xi'}$ with $\T_{\xi'}= \B^{k}\times I_{\xi'}$,
we can say that the chart $\bphi$ satisfies

\begin{equation*}
\bphi \big| _{\B^{k}\times \{ x_2\}\times I_{\xi'}}\,\,=\,\, \psi 
\tag{7}
\end{equation*}\\

\noindent and using the obvious identifications of $\{ x_1\}\times \B^{n-1}\times\{ t\}$\, with $\B^{n-1}\times\{ \br (x_1,t)\}$
\, and $\{ y(x_1,t)\}\times M$\, with \, $M$, we can write

\begin{equation*}
\bphi \big| _{\{ x_1\}\times \B^{n-1}\times\{ t\}}\,\,=\,\, \phi\big|_{\B^{n-1}\times\{ \br (x,t)\}} 
\tag{8}
\end{equation*}\\

Since, as mentioned above, every $\cE_k(R^+u)$ intersects 
every $\bS_{r,y}$\, $g$-orthogonally in a single point, we have that the $\B^{n-1}$-fibers \,\, $\{ x_1\}\times\B^{n-1}\times \{ t\}$\,,  and the \,$(\B^k\times I)$-fibers
\,\,$\B^{k}\times \{ x_2\} \times I_{\xi'}$\,\, are $\bphi^*(g)$-orthogonal. Also, by (7),  $\bphi^*(g)$ restricted to a $\B^k\times I$-fiber is canonically
hyperbolic, hence, by lemma 3.5.3 and the fact that $r\0{0}\leq t_0$, we have\\

\noindent {\bf (3.5.7)} \,\,\,{\it
The metric $\bphi^*(g)$, restricted to a $\B^k\times I$-fiber, 
 is $\epsilon$-hyperbolic, provided $C_2\,e^{-r\0{0}}\leq \epsilon$}\\

Therefore $\bphi$ is a warped hyperbolic chart  and has the form $\bphi^*g=f_1+f_2+dt^2$
where $f_1$ is the restriction of $\bphi^*g$ to the $\B^{n-1}$-fibers and $f_2$ is the restriction of $\bphi^*g$ to the $\B^k$-fibers.
Also, $f_2+dt^2$ is the restriction of $\bphi^*g$ to the $(\B^k\times I)$-fibers, and $f_2+dt^2$ is  a hyperbolic metric.
Furthermore, again by (7), we have that  $f_2+dt^2$ (hence also $f_2$) is independent of the variable $x_2$. This together
with (3.5.7) imply that it is enough to consider $f_2$, that is, to prove the following claim.\\

\noindent {\bf Claim.} {\it We have that}
$$\big|\, f_2(x_1,x_2,t)\,\, -\,\,e^t\sigma\0{\R^{n-1}} \big|_{C^2}\,\, \leq\,\, e^{11+4\xi}\,\Big( \epsilon\,+\, e^{-r\0{0}}\Big)$$
\noindent {\it provided $r\0{0}>7+3\xi$.}\\

\noindent {\bf Proof of claim.} Let $a_{ij}$ be the entries of the matrix $f_2$. We have to prove that
$\big|\, a_{ij}(x_1,x_2,t)\,\, -\,\,e^t\delta_{ij} \big|_{C^2}\,\, <\,\, 
e^{11+4\xi}\,\big( \epsilon\,+\, e^{-r\0{0}}\big)$, assuming $r\0{0}>7+3\xi$.
Let $b_{ij}(x,r)$ be the entries of the matrix $\phi^*h_r$. Since, by hypothesis, $\phi$ is $\epsilon$-hyperbolic, we
have that

\begin{equation*}
\big|\, b_{ij}(x\,,\,r)\,\, -\,\,e^r\delta_{ij} \big|_{C^2}\,\, <\,\, \epsilon 
\tag{9}
\end{equation*}\

\noindent On the other hand, equation (8) implies:

\begin{equation*}
a_{ij}(x_1\,,\,x_2\,,\,t)\,=\, b_{ij}\big(x_2\, ,\,\br(x_1\, ,\,t)\,-\, r\0{0}\,\big)
\tag{10}
\end{equation*}\

The proof of the claim is obtained by calculating the derivatives of $a_{ij}(x_1,x_2,t)\,\, -\,\,e^t\delta_{ij}$ 
up to order 2 and finding estimates of these derivatives using (9), (10) and lemma 3.5.5.
This is done in appendix B. (The idea here is that, by (10), $a_{ij}(x_1,x_2,t)\,\, -\,\,e^t\delta_{ij}$ 
is equal to $b_{ij}(x_2,\br (x_1,t)-r\0{0})\,\, -\,\,e^t\delta_{ij}$, which, by lemma 3.5.5 is $C^2$-close to
$b_{ij}(x_2,t)\,\, -\,\,e^t\delta_{ij}$ which, by (9), is small.)
This completes the proof of proposition 3.5.1.\\\\

\noindent {\bf Proposition 3.5.7.}  {\it Let $\xi>1$ and $a>7+3\xi$. Let $(M^n,h)$ have center $o$. Suppose:}
\begin{enumerate}
\item[{\it (1)}] {\it the metric $h$ is warped $\epsilon$-hyperbolic, with charts of excess $\xi$, outside $\B_{a}(M)$.}
\item[{\it (2)}] {\it The metric $h$ is hyperbolic
on $\B_{a+2+\xi}(M)$.}
\end{enumerate}

\noindent {\it Then }

\begin{enumerate}
\item[{\it (1')}] {\it the metric $\cE_k(h)$ is warped $\eta$-hyperbolic, with charts of excess $\xi'$, outside $\B_a(\cE_k(M))$, provided $0<\xi'<\xi-e^{-\big(a-(7+3\xi)\big)}$.}
\item[{\it (2')}] {\it The metric $\cE_k(h)$ is hyperbolic
on $\B_{a+2+\xi}(\cE_k(M))$}
\end{enumerate}

\noindent {\it provided}

$$
C_1'\,\, e^{-a}\,\,+\,\,C_3\,\, \epsilon\leq\eta
$$\\

\noindent {\bf Remarks.} \\
{\bf 1.} Here $C_3=C_3(k,\xi)$ is as in proposition 3.5.1,
and $C_1'=C_1'(n+k, \xi)$ is as in 3.4.5.
The metric $\cE_k(h)$ in the proposition above is warped $\eta$-hyperbolic
with respect to any center $o\in\HH^k\sbs\cE_k(M)$.\\
{\bf 2.} By ``$h$ is hyperbolic on $\B_s$" we mean that the metric $h$ has
all sectional curvatures equal to -1 on $\B_s$. In particular $\B_s$ is isometric
to $\B_s(\HH^n)$.\\

\noindent {\bf Proof.} Denote the center of $\cE_k(M)$ by $o=(o\0{\HH^k},o\0{M})$.
Note that condition (2) in the proposition imply

\begin{equation*}
{\mbox{\it  the metric $\cE_k(h)$ is hyperbolic on 
$\HH^k\times\B_{a+2+\xi}(M)$}}
\tag{a}
\end{equation*}

\noindent Note that (a) implies (2').
Let $p=(y,v)\in\cE_k(M)$. We use the functions (coordinates) in section 2.3. 
In particular $s=d\0{\cE_k(M)}(o,p)$, $t=d\0{\cE_k(M)}(p,\HH^k)=d\0{M}(o\0{M},v)$.
We can write $\B_{a}(\cE_k(M))=\{s<a\}$. Note that $t\geq a$ implies $s\geq a$.
We have two cases.\\

\noindent {\bf First case.} $t\geq a$.\\
It can be checked that $C_1'(n+k,\xi)>C_3(k,\xi)$,
hence the hypothesis
$C_1'\,\, e^{-a}\,\,+\,\,C_3\,\, \epsilon\leq\eta$ implies
$C_3\,\big( e^{-a}\,\,+\,\, \epsilon\big)\leq\eta$. This together with
proposition 3.5.1 imply that $\cE_k(h)$ is warped $\eta$-hyperbolic on $p$
(i.e there is a warped $\eta$-hyperbolic chart centered at $p$).\\

\noindent {\bf Second case.} $t < a$, $s\geq a$.\\
The hypothesis
$C_1'\,\, e^{-a}\,\,+\,\,C_3\,\, \epsilon\leq\eta$ implies
$C_1'\, e^{-a}\leq\eta$.  Since $s\geq a$ we can use (a) and apply
proposition 3.3.9 to obtain that $\cE_k(h)$ is warped $\eta$-hyperbolic on $p$.
Therefore there is a warped $\eta$-hyperbolic chart $\phi$ centered at $p$,
with image (a priori) contained in $(n+k)$-hyperbolic space.
Since $\phi$ is centered at $p$ and we are assuming $t\leq a$ we get that
$\phi(\T_\xi)\sbs \{ t\leq a+1+\xi\}$ (because $\phi$ is warped $\eta$-hyperbolic).
Therefore, again by (a), the chart $\phi$ is also a chart for $\cE_k(h)$.
This proves the proposition.\\

In the next corollary we use the constants $a\0{n}(\epsilon,\xi)$ introduced
in 3.3.10.\\

\noindent {\bf Theorem 3.5.8.}  {\it Let $\xi>1$. Let $(M^n,h)$ have center $o$. Suppose:}
\begin{enumerate}
\item[{\it (1)}] {\it the metric $h$ is warped $\epsilon$-hyperbolic, with charts of excess $\xi$, outside the ball of radius $a\0{n}(\epsilon,\xi)$.}
\item[{\it (2)}] {\it The metric $h$ is hyperbolic
on $\B\0{R}(M)$.}
\end{enumerate}

\noindent {\it Then, if  $R\geq R\0{n,k}(\epsilon,\xi)$, we have}

\begin{enumerate}
\item[{\it (1')}] {\it the metric $\cE_k(h)$ is warped $(C_3'\epsilon)$-hyperbolic, with charts of excess $\xi'$, outside the ball of radius $a\0{n+k}(C_3'\epsilon)$. Here \,$\xi'=\xi-e^{-R/2}>0$.}
\item[{\it (2')}] {\it The metric $\cE_k(h)$ is hyperbolic
on $\B\0{R}(\cE_k(M))$}
\end{enumerate}

\noindent {\bf Remarks.} \\
{\bf 1.} The constant  $R\0{n,k}$ is defined by
$R\0{n,k}(\epsilon,\xi)=ln\big(\frac{1}{\epsilon}\big)+
a\0{n}(\epsilon,\xi)+a\0{n+k}(C_3'\epsilon)+18+8\xi$.
Here $C_3'=C_3'(n,k,\xi)=C_1'+C_3$, where $C_1'$, $C_3$ are as in 3.5.7.

\noindent {\bf 2.} Note that being $\epsilon$-hyperbolic (or $(C_3'\epsilon)$-hyperbolic)
outside the ball of radius $a\0{n}(\epsilon,\xi)$ (or $a\0{n+k}(C_3'\epsilon)$)
is the best we can hope for, because this is what happens to hyperbolic space
(see 3.3.10).\\

\noindent {\bf Proof.} Just take $a=R-(2+\xi)$.
Then we get:
(1) $e^{-a}<\epsilon$,\, (2) $R\geq a+2+\xi$, \,(3) $a>7+3\xi$,\, (4) $a>a\0{n+k}(C_3'\epsilon,\xi)$,\, (5) $a>a\0{n}(\epsilon,\xi)$,\,  (6) $a-(7+3\xi)>R/2$. 
We can now apply 3.5.7. Note that it can be verified that
$a\0{n} $ is increasing on the variable $\xi$, hence 
$a\0{n+k}(C_3'\epsilon,\xi)\geq a\0{n+k}(C_3'\epsilon,\xi')$.
This proves the corollary.

\vspace{.8in}

\noindent {\bf 3.6.  Warping with  $sinh\, t$.}\\
The metric of our basic hyperbolic model $\T$ is an exponentially warped metric. Here we show that we can easily
change the exponential by multiples of $sinh(t)$, for $t$ large.\\

\noindent {\bf Lemma 3.6.1.} 
{\it  For $t_0> 2$ \, we have }
{\scriptsize \begin{enumerate}
\item[{\it (1)}] 
$\big|e^{-t}\bigg(\frac{sinh\, (t+t_0)}{sinh\, (t_0)}\bigg)\, -\, 1   \big|_{C^0(\R^+)}\, <  \,(\frac{49}{48})\,e^{-2t_0} $
\item[{\it (2)}] 
$\big|  \frac{d}{dt}e^{-t}\bigg(\frac{sinh\, (t+t_0)}{sinh\, (t_0)}\bigg)\,  \big|_{C^0(\R^+)}\, <  \,(\frac{1}{24})\,e^{-2t_0} $
\item[{\it (3)}] 
$\big| \frac{d^2}{dt^2} e^{-t}\bigg(\frac{sinh\, (t+t_0)}{sinh\, (t_0)}\bigg)\,    \big|_{C^0(\R^+)}\, <  \,(\frac{1}{12})\,e^{-2t_0} $
\item[{\it (4)}]
$\big|e^{-t}\bigg(\frac{sinh\, (t+t_0)}{sinh\, (t_0)}\bigg)\, -\, 1   \big|_{C^2(\R^+)}\, <  \,(\frac{49}{48})\,e^{-2t_0} $
\item[{\it (5)}] 
$\big| e^{-t}\bigg(\frac{sinh\, (t+t_0)}{sinh\, (t_0)}\bigg) \big|_{C^2(\R^+)}\, <\, 1.03$
\item[{\it (6)}] 
$\big| e^{-2t}\bigg(\frac{sinh\, (t+t_0)}{sinh\, (t_0)}\bigg)^2\, -\, 1   \big|_{C^2(\R^+)}\, <\,16\,e^{-2t_0} $.
\end{enumerate}}

\noindent {\bf Proof.} For (1), (2) and (3)
just write $e^{-t}\big(\frac{sinh\, (t+t_0)}{sinh\, (t_0)}\big)\, -\, 1 \,=\,\big(\frac{e^{-2t_0}}{1-e^{-2t_0}}\big)(1-e^{-2t})$
and differentiate twice. Items (4) and (5) follow from (1), (2), (3). 
For (5) write\ 

{\scriptsize $$e^{-2t}\big(\frac{sinh\, (t+t_0)}{sinh\, (t_0)}\big)^2\, -\, 1 =
\bigg(e^{-t}\big(\frac{sinh\, (t+t_0)}{sinh\, (t_0)}\big)\, -\, 1\bigg)\bigg( e^{-t}\big(\frac{sinh\, (t+t_0)}{sinh\, (t_0)}\big)\, +\, 1\bigg) $$} 
\noindent and use the previous items together with the Leibniz rule.
This proves the lemma.\\

Recall $I_\xi=\big(-(1+\xi),\, 1+\xi\big)$.
 Let $\nu:I_\xi\ra \R^+$ be smooth. For a metric $g=g\0{t} + dt^2$ on $\T_\xi=\B^k\times I_\xi$, we write
$g\0{\nu}=\nu g\0{t}+dt^2$.\\

\noindent {\bf Lemma 3.6.2.} {\it We have } 
\begin{enumerate}
\item[{\it (1)}] 
$\big| g-g\0{\nu}  \big|_{C^2}\, <\, 4\, \big| 1-\nu   \big|_{C^2}\, \big| g\big|_{C^2}$.  
\item[{\it (2)}] 
$\big| h\0{\nu}-g\0{\nu}  \big|_{C^2}\, <\, 4\, \big| \nu   \big|_{C^2}\, \big| h- g\big|_{C^2}$.  \\
\end{enumerate}

\noindent {\bf Proof.} For (1) just note that $g-g\0{\nu}=(1-\nu)g$ and differentiate twice. Item (2) is similar.
This proves the lemma.\\

Recall that the metric on our model $\T_\xi$ is $\sigma=e^{2t}\sigma_{\R^k}+dt^2$. 
Consider now the metric $\sigma\0{t_0}=\Big(\frac{sinh\, (t+t_0)}{sinh\, (t_0)}\Big)^2\sigma\0{\R^k}+dt^2$, $t_0>2$.\\

\noindent {\bf Lemma 3.6.3.} {\it Let $g=g\0{t}+dt^2$ be $\epsilon$-hyperbolic on $\T_\xi$, $\xi>0$. 
Let $h=e^{-2t}\big(\frac{sinh\, (t+t_0)}{sinh\, (t_0)}\big)^2g\0{t}+dt^2$. 
Assume $t_0>2$. Then }
\begin{enumerate}
\item[{\it (1)}] $ \big | h-g \big |\, <\,2^6\, (\epsilon +e^{2(1+\xi)})\, e^{-2t_0}$.
\item[{\it (2)}] $ \big | \sigma-\sigma\0{t_0} \big |\, <\,2^6\, e^{2(1+\xi)}\, e^{-2t_0}$.
\item[{\it (3)}] $ \big | h-\sigma\0{t_0} \big |\, <\,  8\, \epsilon$
\item[{\it (4)}] $ \big | g-\sigma\0{t_0} \big |\, <\, (2^6\, e^{2(1+\xi)}\, e^{-2t_0})\,+\, \epsilon$
\item[{\it (5)}] $ \big | h-\sigma \big |\, <\,\epsilon\,+ \,2^6\, (\epsilon +e^{2(1+\xi)})\, e^{-2t_0}\,
\leq \,2^6\, e^{2(1+\xi)}\,\Big( e^{-2t_0}\,+\,\epsilon\Big)$
\item[{\it (6)}] $h $ is $\eta$-hyperbolic, with $\eta\geq \, 2^6\,  e^{2(1+\xi)}\, \Big( e^{-2t_0}\, +\,\epsilon\Big)$.
\end{enumerate}

\noindent {\bf Proof.} Note that $h=g\0{\nu}$, with $\nu=e^{-2t}\big(\frac{sinh\, (t+t_0)}{sinh\, (t_0)}\big)^2$. Therefore (1)
follows from 3.6.1 (6),  3.6.2 (1) and the fact that $|g|_{C^2}\leq |g-\sigma|_{C^2}+|\sigma|_{C^2}<\epsilon+e^2$.
Also $\sigma\0{t_0}=\sigma_\nu$, hence (2) is similar to (1). Item (3) follows from 3.6.1 (4) and 3.6.2 (2) (we use $t_0>2$ to have $16 e^{-2t_0}<1$).
Item (4) follows from the previous items together with triangular inequality. 
Similarly item (6) follows from item (1), the triangular inequality and the fact that
$g$ is $\epsilon$-hyperbolic. Item (6) is the same as item (5).
This proves the lemma.\\

Hence $h$ is \,\,$\big(8\, \epsilon \big)$-hyperbolic with respect to the warped $sinh$ model $(\B^k\times I_\xi, \sigma\0{t_0})$.
Also $g$ is \,\,$\big(2^6\,   e^{2(1+\xi)}\, e^{-t_0}\,+\, \epsilon \big) $-hyperbolic with respect to this same model. Results for $cosh$ are similar.\\

\noindent {\bf 3.6.4. Remarks. }\\
\noindent {\bf 1.} Note that $h$ in lemma 3.6.3 is just $g\0{\nu}$, with
$\nu(t)=e^{-2t}\big(\frac{sinh\, (t+t_0)}{sinh\, (t_0)}\big)^2$.\\
\noindent {\bf 2.} It is straightforward to 
verify that lemma 3.6.1 holds if we replace the variable $t$ by $t-s$, for some fixed $s$. Hence items (1) and (6) in lemma 3.6.3 hold if we replace $h=g\0{\nu}$ (with
$\nu$ as in remark 1), by $g\0{\nu\0{s}}$, where $\nu\0{s}(t)=\nu(t-s)$, for fixed $s$.

\vspace{1in}

\noindent {\bf \large  Section 4. Spherical Cuts and Warp Forcing.}\\

In this section we will show how  to force an $\epsilon$-hyperbolic metric to be  warped. 
We also define spherical cuts and describe the spherical cut of a hyperbolic extension.
As in section 2, let $(M^n, h)$ be a complete Riemannian manifold with center $o\in M$. Recall that we
can write the metric on $M-\{ o\}=\bS^{n-1}\times\R^+$ as \,\,$h=h_t+dt^2$. Let $S\sbs\bS^{n-1}$ and denote by $\rC S$ the open cone
$S\times\R^+\sbs\bS^{n-1}\times\R^+\sbs M$. We write $S_r=\rC S\cap\bS_r(M)=S\times\{r\}$.
We will use the notation of sections 2.2, 2.3 and 2.4.
\vspace{.8in}

\noindent {\bf 4.1.  Spherical cuts.}\\
The metric $h_r$ on $\bS_r$ is called the {\it (warped by sinh) spherical cut of $h$  at $r$}, and the metric
$$\hh\0r\, =\, \bigg(\frac{1}{sinh^2(r)}\bigg)\, h_r$$
\noindent is called the {\it unwarped (by sinh) spherical cut of $h$ at $r$ }. \\

If $h=h_t+dt^2$ is a variable metric defined only on $\rC S$ the above definitions still make sense
and we say that $h_r$ is the {\it spherical cut of $h$  on $S$, at $r$}, and
$\hh\0r$  is the {\it unwarped spherical cut of $h$  on $S$, at $r$}.

\vspace{.8in}

\noindent {\bf 4.2.  Spherical cuts of hyperbolic extensions.} \\
As before we denote by $g$ the metric on the hyperbolic extension $\cE_k(M)$, and we write
$g=g_s+ds^2$ on $\cE_k(M)-\{ o\}$.
We use the map $\Xi=\Xi_s$ of section 2.4 that gives coordinates on $\bS^\circ_s\big( \cE_k(M) \big)$.
Let $g_s'=\Xi_* g_s$, that is the metric $g'_s$ on $\bS^{k-1}\times\bS^{n-1}\times (0,\pi/2)$
is the expression of $g_s$ in the $\Xi$-coordinates.\\

\noindent {\bf Proposition 4.2.1} {\it We have that} $$g'_s\,=\,\Big( sinh^2(s)  \, cos^2\, (\beta)\,\Big)\, \sigma\0{\bS^{k-1}}\, \, +\,\,
h_r\,\,+\,\,\Big( \,sinh^2(s)  \,\Big)\,d\beta^2$$
\noindent {\it where $r=\bs(s,\beta)=sinh^{-1}\big( sin\, (\beta)\,\,sinh\, (s) \big)$  (see  remark 2.4.1).}\\

\noindent {\bf Proof.} By item 4 in section 2.4 we have that $g'_s$ has the form $A+B+C$, where $A(u,\beta)$ is a metric on
$\bS^{k-1}\times\{ u\}\times \{ \beta\}$, $B(w,\beta)$ is a metric on
$\{ w\}\times\bS^{n-1}\times \{ \beta\}$ and $C(u,\beta)$ is a metric on
$\{w\}\times \{ u\}\times (0,\pi/2)$, i.e $C=f(w,u,\beta)\,d\beta^2$, for some positive function $f$.\\ 

Now, by definition we have 
$$g\,=\,cosh^2(r)\, \sigma\0{\HH^k}\, +\, h_r+dr^2\, =\,cosh^2(r)\bigg( sinh^2(t)\sigma\0{\bS^{k-1}}+dt^2 \bigg)\, +\, h_r\, +\, + dr^2 $$
\noindent By lemma 2.3.4? and the identity $cosh(r)\, sinh(t)= sinh(s)\, cos\,(\beta)$
(which follows from the law of sines and the second law of cosines) we can write

$$g\,=\,\big( sinh^2(s)  \, cos^2\, (\beta)\,\big) \sigma\0{\bS^{k-1}}\, \, +\,\,
h_r\,\,+\,\,\big( sinh^2(s)  \big)\,d\beta^2\,\,+\,\, ds^2$$
\noindent 
If we fix $s$ we get $r$ as a function of $\beta$, and also $ds^2\equiv 0$. This proves the proposition.\\





\noindent {\bf Corollary 4.2.2} {\it We have that} $$\Xi_*\big(\,\, \hg_s\,\,\big)\,=\, cos^2\, (\beta)\, \sigma\0{\bS^{k-1}}\, \, +\,\,
sin^2\, (\beta)\,\hh\0r\,\,+\,\,\,d\beta^2$$
\noindent {\it where $r=\bs(s,\beta)$ as in proposition 4.2.1.}\\

\noindent {\bf Proof.} Since $sinh^2(r)\,\hh\0r=h_r$,
the corollary follows from proposition 4.2.1 and
the identity $\frac{sinh (r)}{ sinh (s)}=sin\, (\beta)$, which is deduced from the hyperbolic law of sines.

\vspace{.8in}

\noindent {\bf 4.3.  Local warp forcing.}  \\
Here we give a kind of a local version the the method of warp forcing, which is given
in section 4.4.\\

Let $a$ be a metric on $\B^k$. For fixed $s\in I_\xi=(-1-\xi,1+\xi)$, 
$\xi>0$, we denote by $\ua=\ua_s$ 
the warped metric $e^{2(t-s)}a+dt^2$ on $\B^k\times I_\xi$.\\

\noindent {\bf Lemma 4.3.1.} {\it Fix $s$. Let $a$, $b$ be metrics on $\B^k$ with \,\,$|\, a\, -\, b\, |_{C^2(\B^k)}\, < \,\epsilon$. Then 
\,$|\, \ua\, -\, \ub\, |_{C^2}\, <\,4 \, e^{4(1+\xi)}\,\epsilon$.}\\

\noindent {\bf Proof.}
Just compute the derivatives of $\ua-\ub=e^{2(t-s)}(a-b)$. This proves the lemma.\\

\noindent {\bf Lemma 4.3.2.} {\it Let $g=g\0{t}+dt^2$ be an $\epsilon$-hyperbolic metric on $\T_\xi=\B^k\times I_\xi$. 
Fix $s\in I_\xi$ and
consider ${\underline{g_s}}=e^{2(t-s)}g_s+dt^2$. Then  we have that }
\begin{enumerate}
\item[{\it (i)}] {\it  $|\,g-{\underline{g_s}}\,|_{C^2}\,<\, (1+4\,e^{4(1+\xi)})\,\epsilon$.}
\item[{\it (ii)}] {\it the metric ${\underline{g_s}}$ is $(4\,e^{4(1+\xi)}\,\epsilon)$-hyperbolic.} 
\end{enumerate}  

\noindent {\bf Proof.} By hypothesis
we have $|\, (g\0{t}+dt^2)\,-\, (e^{2t}\sigma\0{\R^k}+dt^2)\, |_{C^2}\,<\,\epsilon$. Therefore, taking $t=s$ we get
$|\, g_s\, -\, e^{2s}\sigma\0{\R^k}\, |_{C^2(\B^k)}\, < \,\epsilon$.
Lemma 4.3.1 implies then that $|{\underline{g_s}}-\sigma |<4\,e^{4(1+\xi)}\epsilon$, and this
together with the triangular inequality imply (i).
This completes the proof of lemma 4.3.2.\\

Let $\rho :\R\ra [0,1]$ be the function of section 3.4, that is $\rho$ is
smooth and such that: {\bf (i)}\, $\rho|_{(-\infty, 0]}\equiv 0$, and\, \,{\bf (ii)}\, $\rho|_{[1\, ,\infty)}\equiv 1$.
Let ${\bar{\rho}}(t)=\rho(1-2t)$. Then we have
{\bf (i)}\, ${\bar{\rho}}|_{(-\infty, 0]}\equiv 1$, and\, \,{\bf (ii)}\, ${\bar{\rho}}|_{[1/2\, ,\infty)}\equiv 0$.\, For $s\in \R$ write $\rho\0{s}(t)={\bar{\rho}}(t-s)$.\\

\noindent {\bf Lemma 4.3.3} {\it Let $g=g\0{t}+dt^2$ be an $\epsilon$-hyperbolic metric on $\T_\xi=\B^k\times I_\xi$. Fix $s\in I_\xi$.
Write $h=\rho\0{s}\,{\underline{g_s}}+(1-\rho\0{s})\,g$. Then we have that }
\begin{enumerate}
\item[{\it (i)}] {\it \,\, $|\,h-{\underline{g_s}}\,|_{C^2}\,<\,\bigg( |\rho\0{s}|_{C^2}+2|\rho\0{s}|_{C^1}+1  \bigg)\, (1+4\,e^{4(1+\xi)})\,\epsilon$.}
\item[{\it (ii)}] {\it \,\, $|\,h-g\,|_{C^2}\,<\, \bigg( |\rho\0{s}|_{C^2}+2|\rho\0{s}|_{C^1}+1  \bigg)\, (1+4\,e^{4(1+\xi)})\,\epsilon$.}
\item[{\it (iii)}] {\it the metric $h$ is $(e^{6+4\xi}\,\epsilon)$-hyperbolic.} 
\end{enumerate}  

\noindent {\bf Proof.} Note that $h-{\underline{g_s}}=(1-\rho\0{s})(g-{\underline{g_s}})$. Using lemma 4.3.2 (i) and the Leibniz rule to differentiate
$(1-\rho\0{s})(g-{\underline{g_s}})$ we obtain (i). Since $|\rho\0{s}|_{C^k}=|1-\rho\0{s}|_{C^k}$, item  (ii) is proved similarly. It follows from (ii),
the fact that $g$ is $\epsilon$-hyperbolic and the triangular inequality that
 the metric $h$ is $(A\,\epsilon)$-hyperbolic, with $A=1+\bigg( |\rho\0{s}|_{C^2}+2|\rho\0{s}|_{C^1}+1  \bigg)\, (1+4\,e^{4(1+\xi)})$. But a calculation shows that we can take $\rho$ with $|\rho|_{C^1}<6$ and $|\rho|_{C^2}<48$ (see 3.4).
Again a calculation shows we can take $A=e^{6+4\xi}$. This proves the lemma.\\
\vspace{.8in}

\noindent {\bf 4.4.  Warp forcing.} \\
Let $M$ have center $o$ and metric $g=g\0{r}+dr^2$. Fix $r\0{0}>0$.
We define the warped by $sinh$ metric $\bg\0{r\0{0}}$ by:

$$\bg \0{r\0{0}}\,=\,sinh^2(t)\hg\0{r\0{0}}\,+\, dr^2     \,=\, sinh^2\, (t)\big(  \frac{1}{sinh^2(r\0{0})}\big) g\0{r\0{0}}\, +\, dr^2
$$\

We now force the metric $g$ to be equal to $\bg\0{r\0{0}}$ on  $B_{r\0{0}}=\B_{r\0{0}}(M)$ and stay equal to $g$ outside $B_{r\0{0} +\frac{1}{2}}$.
For this we define the {\it warped forced } (on $B_{r\0{0}}$) metric as:\\

$$
\cW\0{r\0{0}}\, g\, =\, \rho\0{r\0{0}}\, \bg\0{r\0{0}}\, +\, (1-\rho\0{r\0{0}})\, g
$$\\

\noindent where $\rho\0{r\0{0}}$ is as in section 4.3. Hence we have\\\\

\noindent{\bf (4.4.1)}\hspace{1.2in}
 $\cW\0{r\0{0}} g\,=\,\left\{ \begin{array}{lllll}
\bg\0{r\0{0}}&& {\mbox{on}}& & B_{r\0{0}}\\\\
g&&{\mbox{outside}}& & B_{r\0{0}+\frac{1}{2}}
\end{array}\right.$\\\\

The next result shows that if $g$ is $\epsilon$-hyperbolic then
the warp forced metric  $\cW\0{r\0{0}} g$ is $\eta$-hyperbolic,
where $\eta$ depends on $\epsilon$ and how far we are from
the center $o$. Also, the excess of the $\eta$-hyperbolic charts
of  $\cW\0{r\0{0}} g$ is one less than the excess of the
$\epsilon$-hyperbolic charts of $g$. Therefore some excess is
lost in the warp forcing process.\\

\noindent {\bf Proposition 4.4.2.} {\it  Let $\xi>1$ and let the metric $g$ be warped $\epsilon$-hyperbolic on the subset $S\sbs M-\{ o\}$
with charts of excess $\xi$.
Then $\cW\0{r\0{0}} g$ is warped $\eta$-hyperbolic on $S-B\0{r\0{0}-(1+\xi)}$ with charts of excess $\xi-1$, provided $\eta\geq e^{18+6\xi}\big(e^{-2r\0{0}}+\epsilon\big)$.}\\

\noindent {\bf Proof.} At some points in the proof we will use the notation of section 3.6. Let $p\in S$ and outside
$B\0{r\0{0}-(1+\xi)}$. We have  three cases.\\

\noindent {\bf First  case.}  $p\notin B_{r\0{0}+\frac{1}{2}+(1+\xi)}$\\
 In this case we can completely fit a warped $\epsilon$-hyperbolic chart of $g$ of excess $\xi$ outside 
$B_{r\0{0}+\frac{1}{2}}$. But, by (4.4.1), this chart is also
a warped $\epsilon$-hyperbolic chart for $\cW\0{r\0{0}} g$.
This shows  the metric 
$\cW\0{r\0{0}} g$ is $\epsilon$-hyperbolic outside $B_{r\0{0}+\frac{1}{2}+(1+\xi)}$, with charts of excess $\xi$.\\

\noindent {\bf Second case.}  $p\in B_{r\0{0}+\frac{1}{2}+(1+\xi)}-
B_{r\0{0}+\frac{1}{2}+\xi}$\\
Let $\phi:\T_\xi\ra M $ be an $\epsilon$-hyperbolic chart of $g$ centered at $p=(x,t_0)$. Then the image of the restriction $\phi|\0{\T_{\xi-1}}$
of $\phi$ to $\T_{\xi-1}$ does not intersect $B_{r\0{0}+\frac{1}{2}}$,
hence as in the first case, by 4.4.1, the chart $\T\0{\xi-1}$ is an
$\epsilon$-hyperbolic chart for $\cW\0{r\0{0}} g$, but with excess
$\xi-1$.\\

\noindent  {\bf Third case.}  $p\in B_{r\0{0}+\frac{1}{2}+\xi}$\\
In this case the interval $I_{t\0{0},\xi}=(t\0{0}-(1+\xi), \, t\0{0}+(1+\xi) )$  
contains $r\0{0}$ (recall we are assuming $p\notin B\0{r\0{0}-(1+\xi)}$). 
Write $g'=g'\0{t}+dt^2=\phi^*g$. Then $g'\0{t}=\phi^*g\0{t+t_0}$.  Write $ s=r\0{0}-t_0$, thus $-\frac{1}{2}-\xi<s<1+\xi$. 
In particular $s\in I_\xi$.
Therefore $g'\0{s}=\phi^*g\0{r\0{0}}$.\\

Since $\frac{sinh^2 (t+t_0)}{sinh^2(r\0{0})}=\frac{sinh^2 ((t-s)+r\0{0})}{sinh^2 (r\0{0})}$ we have that
$e^{2(t-s)}\nu_s(t)=\frac{sinh^2 ((t-s)+r\0{0})}{sinh^2 (r\0{0})}$, where $\nu(t)=e^{-2t}\frac{sinh^2(t+r\0{0}))}{sinh^2(r\0{0})}$ and $\nu_s(t)=\nu(t-s)$.
Hence 
(see notation in 3.6)

\begin{equation*}
\phi^*\big(\, {\bar{g}}\0{r\0{0}}\,\big)\,\,=\,\,({\underline{g'\0{s}}})_{\nu\0{s}}
\tag{1}
\end{equation*}

\noindent And since $g'$ is $\epsilon$-hyperbolic, 4.3.2 (ii) implies
that ${\underline{g'\0{s}}}$ is $(4\,e^{4(1+\xi)}\epsilon)$-hyperbolic.
This together with 3.6.3 (6) (see also remark 3.6.4)  imply 

\begin{equation*} ({\underline{g'\0{s}}})\0{\nu\0{s}} 
\,\,\,\,\,\,\,\,\, is\,\,\,\,\,\,\,\,\,
\Big(2^6\,e^{2(1+\xi)}\,\big(e^{-2t_0}+4\,e^{4(1+\xi)}\epsilon\big)\,\Big)-hyperbolic
\tag{2}
\end{equation*}

\noindent Note that

\begin{equation*}
\phi^*\big( \cW\0{r\0{0}}\, g \big )\,=\,\rho\0{r\0{0}}\phi^*\big(\bg\0{r\0{0}}\big)
\,+\,\big(1-\rho\0{r\0{0}}\big) \,\phi^* g\,=\,
\rho\0{r\0{0}}\big(\underline{g\0{s}}\big)\0{\nu\0{s}}
\,+\,\big(1-\rho\0{r\0{0}}\big) \, g'
\tag{3}
\end{equation*}

\noindent From (1), (2), (3) and lemma 3.1.1 we get that $\phi^*\big( \cW\0{r\0{0}}\, g \big )$ is $\epsilon'$-hyperbolic with

$$
\epsilon'\,\,=\,\,4\,\, |\,\rho \0{r\0{0}}\,|\,\, \bigg(\Big(2^6\,e^{2(1+\xi)}\,\big(\,e^{-2t_0}+4\,e^{4(1+\xi)}\epsilon\,\big)\,\Big)\,\,+\,\, \epsilon   \bigg)
$$

\noindent Since we can take $|\rho \0{r\0{0}}|<48$ (see 3.4) and $t_0>r\0{0}-(1+\xi)$ we get
$\epsilon'<4\, 48\, (1+2^8e^{6(1+\xi)})(e^{-2r\0{0}}+\epsilon)<e^{18+
6\xi}(e^{-2r\0{0}}
+\epsilon)$. Note that the excess of the charts in this third case
is also $\xi$. This proves the proposition.\\

\noindent {\bf 4.4.3. Remarks.}\\
\noindent {\bf 1.} Note that the definitions, methods and discussion above
still make sense if the metric $g=g\0{r}+dt^2$ is defined only outside
some $B_R$, provided $r\0{0}>R-(1+\xi)$. 

\noindent{\bf 2.} Note also that, by construction, $\cW\0{r\0{0}}\, g$
is a variable metric, hence it has the same ray structure a the original
metric $g$ (see section 1?).

\vspace{1in}

\noindent {\bf \large  Section 5. Limit Metrics and Continuation.}\\

 We will define a family of partially defined metrics on $\R^{n+1}$ and show how to extend these metrics provided they have a ``limit".\\

\noindent {\bf 5.1. Families of metrics and limits.}\\
As before we identify $\R^{n+1}-\{0\}$ with $\bS^n\times \R^+$, and let $r$ denote the euclidean distance to the origin. All metrics and families 
of metrics here are assumed to be smooth. We write $B_r=\B_r(\R^{n+1})$.\\\

Fix $\xi>0$.
We say that $\{ g_\ssl\}_{\lambda> 2+\xi}$ is a $\odot${\it-family of metrics} if each $g_\ssl$ is 
a variable metric defined (at least) on $\R^n-{\bar{B}}_{_{\lambda-(2+\xi)}}$, that is $g_\ssl=\big(g_\ssl\big)_r+dr^2$, with 
$\big(g_\ssl\big)_r$ defined for (at least) $r>\lambda-(2+\xi)$.
Of course every family of variable metrics defined on the whole of $\R^n-\{ 0\}$ is a $\odot$-family of metrics.
We say that the  $\odot$-family $\{g_\ssl\}$ has {\it cut limit at $b$} if there is a 
$C^2$ metric $\hg_{_{\infty+b}}$ on $\bS^{n}$ such that 
$$
{\widehat{\big(g_\ssl\big)}}_{_{\lambda+b}} \,\,\,\,\stackrel{C^2}{\longrightarrow}\,\,\,\, \hg_{_{\infty+b}}
$$

\noindent {\bf Remarks.}

\noindent {\bf 1.} The metric ${\widehat{\big(g_\ssl\big)}}_{_{\lambda+b}}$ is the unwarped spherical cut of $g_\ssl$ at
$\lambda+b$. See section 4.1.

\noindent {\bf 2.} The arrow above means convergence in the $C^2$-norm on $\bS^n$.

\noindent {\bf 3.} The definition above implies that $\big(g_\ssl\big)_{_{r+b}}$ is defined for large $\lambda$,
even if $b<-(2+\xi)$.

\noindent {\bf 4.} Note that the concept of cut limit at $b$ depends strongly on the indexation of the family.\\

Let $I\sbs\R$ be an interval (compact or noncompact). We say the  $\odot$-family $\{g_\ssl\}$ has {\it cut limits on $I$} if 
$$
{\widehat{\big(g_\ssl\big)}}_{_{\lambda+b}} \,\,\,\,\stackrel{C^2(\bS^n\times I )}{\longrightarrow}\,\,\,\, \hg_{_{\infty+b}}
$$

\noindent where we consider $b\in I$. Here the convergence
on $C^2(\bS^n\times I)$ is $C^2$-convergence on compact supports.
In particular $\{g_\ssl\}$ has a cut limit at $b$, for every $b\in I$ and 
the family $\big\{\hg_{_{\infty+b}}\big\}_{b\in I}$ of metrics on $\bS^{n}$ is smooth.\\

We also have ``partial" versions of the definitions above. Let $S\sbs\bS^n$.
We say that $\{ g_\ssl\}_{\lambda> (2+\xi)}$ is a $\odot${\it-family of metrics on $S$} if each $g_\ssl=g_\ssl=\big(g_\ssl\big)_r+dr^2$ is 
a variable metric defined (at least) on $\rC S-{\bar{B}}_{_{\lambda-(2+\xi)}}$,
where $\rC S$ is the cone on $S$. The definitions of {\it cut limit over $S$, at $b$}  and {\it cut limits  over $S$, on $I$} are similar.
In this case the $C^2$ convergence is uniform $C^2$ convergence
with compact supports
in the $I$-direction and uniform $C^2$ convergence in the $S$-direction.\\

If $M$ is a Riemannian manifold with center $o$, them by identifying $M$ with $\R^n$ (using $exp_o$) all definitions above make
sense for $M$.

\vspace{.8in}

\noindent {\bf 5.2. Continuation.}\\
Consider the  $\odot$-family $\{g_\ssl\}$ and assume 

\begin{enumerate}
\item[{\bf (a)}] $\{g_\ssl\}$ has a cut limit metric $\hg_{_\infty}=\hg_{_{\infty+0}}$ at $b=0$. 
\item[{\bf (b)}] each $g_\ssl$ is warped $\epsilon$-hyperbolic (with respect to the decomposition $\R^{n+1}-\{0\}=\bS^n\times\R^+$).
\end{enumerate}

We give a process to extend the metrics $g_\ssl$, for large $\lambda$, to $\eta$-hyperbolic metrics (with charts of excess $\xi$) that have
 ``good" properties" on the whole $\R^n$.
Here $\eta=\eta(\lambda,\epsilon, \xi)$.\\

Given $g_\ssl$ and  $d>0$, with $\lambda>d$,  we define
the {\it continuation metric} $\cC_{_d}\big( g_\ssl \big)$ in the following way.
First warp-force the metric $g_\ssl$, i.e take $\cW_\ssl \big(g_\ssl \big)$ (see section 4.4 and remark 4.4.3 (1)). More specifically, outside $B_{_{\lambda-1}}$ take
$$ \cC_{_d}\big( g_\ssl \big)\,=\,\cW_\ssl \big(g_\ssl\big)\,$$ 
and recall $\cW_\ssl \big(g_\ssl\big)$ is a warped metric on $B_\ssl-B_{_{\lambda-1}}$. Hence $\cW_\ssl \big(g_\ssl\big)$ is defined on the whole of $\R^{n+1}-\{0\}$ (i.e on $B\0{\lambda}$ it is just the $sinh$-warped metric $ \big({\overline{g_\ssl}}\big)_\ssl $).
Next, use  hyperbolic forcing (see section 3.4) and on $B\0{\lambda}$ take
$$\cC_{_d}\big( g_\ssl \big)\,=\, \cH_{_{(\lambda-d),d}}\bigg( \,\, \cW_\ssl \big(g_\ssl\big)\,\,   \bigg)\,=\, \cH_{_{(\lambda-d),d}}\bigg( \,\,\,{\overline{(g_\ssl )}}_\ssl\,\,\bigg) $$
\noindent for $\lambda>d$.
Write $h_\ssl=\cC_{_d}\big( g_\ssl \big)$. We say that the family $\{h_\ssl\}$ is the {\it continued family} 
corresponding to the $\odot$-family $\{g_\ssl\}$.
Since $h_\ssl$ is still a variable metric we can write $h_\ssl=\big( h_\ssl \big)_r+dr^2$.
We can explicitly describe $\big( h_\ssl \big)_r$:\\\\

\noindent {\bf (5.2.1.)}\,\,\,\,\,\,\,\,{\small $
\big( h_\ssl \big)_r\,\, =\,\,\left\{
\begin{array}{lll}
\big( g_\ssl \big)_r&&   \lambda +\frac{1}{2}\leq r\\\\
\rho_\ssl(r) \,sinh^2 (r)\, {\widehat{\big(g_\ssl \big)}}_\ssl \, +\, \big(1-\rho_\ssl (r)\big)\,\big( g_\ssl \big)_r&
&\lambda \leq r\leq \lambda +\frac{1}{2}\\\\
sinh\,^2 (r) \bigg(  \rho_{_{(\lambda -d), d}}(r)  \,{\widehat{\big(g_\ssl \big)}}_\ssl\,\,+\,\,
\big( 1-  \rho_{_{(\lambda -d), d}}(r) \big)\, \sigma\0{\bS^n}  \bigg)&&  \lambda -d \leq r\leq \lambda  \\\\
sinh\,^2(r)\, \sigma\0{\bS^n}&&r\leq \lambda-d
\end{array}
\right.
$}\\\\

\noindent where the gluing functions $\rho_\ssl$ and $ \rho_{_{(\lambda -d), d}} $ are defined in sections 4.4 and 3.4, respectively.
(Both these functions are just reparametrizations of the function $\rho$ of section 3.4.)\\\\

\noindent {\bf Proposition 5.2.2.}  {\it The metrics $h_\ssl$ have the following properties.}

{\it \noindent \begin{enumerate}
\item[(i)]  The metrics $h_\ssl$ are canonically hyperbolic on $B_{_{\lambda-d}}$, i.e equal to $sinh^2(r)\sigma\0{\bS^n}+dr^2$ on
$B\0{\lambda-d}$, provided $\lambda >d$.
\item[(ii)] We have that $g_\ssl=h_\ssl$ outside $B_{_{\lambda+\frac{1}{2}}}.$
\item[(iii)] The metric $h$ coincides with $\cW\0{\lambda} \big(g\0{\lambda}\big)$ outside $B\0{\lambda-\frac{d}{2}}$.
\item[(iv)] The metric $h$ coincides with 
$\cH_{_{(\lambda-d),d}}\Big( {\overline{(g_\ssl )}}_\ssl\,\Big)$
 on $B\0{\lambda}$.
\item[(v)] The $\odot$-family $\big\{ h_\ssl\big\}$ has  cut limits on $(-\infty,0]$. In fact we have
$$
\hh\0{_{\infty+b}} \,\, =\,\,\left\{
\begin{array}{lll}
\hg_{_{\infty}}&&  b=0\\\\
 \rho(1+\frac{b}{d})  \,\hg_{_{\infty}}\,\,+\,\,
\big( 1-  \rho(1+\frac{b}{d})  \big)\, \sigma\0{\bS^n}  &&   -d \leq b\leq 0  \\\\
 \sigma\0{\bS^n}&&b\leq -d
\end{array}
\right.
$$
\noindent where $\rho$ is as in 3.4.
\item[(vi)]  If we additionally assume that $\{g_\ssl\}$ has  cut limits on $[0,\frac{1}{2}]$, then
$\big\{ h_\ssl\big\}$ has also  cut limits on $[0,\frac{1}{2}]$. In fact, for $b\in [0,\frac{1}{2}]$ we have
$$
\hh\0{_{\infty+b}} \,\, =\,\, \rho(1-2b)\,\hg_{_{\infty}}\,\,+\,\, \big( 1-\rho(1-2b) \big)\,\hg_{_{\infty+b}}
$$
\noindent where $\rho$ is as in 3.4. Of course if $\{g_\ssl\}$ has a cut limit at $b> \frac{1}{2}$ then $\{h_\ssl\}$ 
has the same cut limit at $b$ (see item 2).
\item[(vii)] If $\{g_\ssl\}$ has cut limits on $I$, with $0\in I$, then the continued family $\{ h_\ssl\}$ has
cut limits on $I\cup (-\infty, 0]$.
\item[(viii)] All the rays $r\mapsto ru$, $u\in \bS^n$, emanating from the origin are geodesics of
$(\R^n , h_\ssl)$. Hence, the spaces  $(\R^n,h_\ssl)$ have center $0\in \R^n$
and have the same ray structure as $\R^n$. Moreover
the function $r$ (distance to the origin $0\in\R^n$) is the same on the spaces $(\R^n,h_\ssl)$.
\end{enumerate}}

\noindent {\bf Proof.}
Items (i), (ii), (v) and (vi) follow from the definition of $h\0{\lambda}$ and 
5.2.1. Item (iv) is just the definition of $h$ on $B\0{\lambda}$.
Item (iii) also follows  from the definition of $h$ and the fact that
$\rho\0{\lambda-d,d}$ is equal to 1 on the interval $[\lambda-\frac{d}{2},
\lambda]$ (see 3.4).
Item (vii) follows from the previous items. Item (viii) follows from the fact that,
by construction, the metrics $h\0{\lambda}$ are variable metrics (see also
Proposition 3.4.5, property (5),  and remark 4.4.3 (2). This proves the proposition.\\

Recall we are assuming that the family $\{g\0{\lambda}\}$ has cut limit metric
$\hg\0{\infty}$ at $b=0$. Fix a finite atlas $\cA\0{\bS^n}$ for $\bS^n$ as in 3.3 (see remark 1 at the beginning of 3.3).
Let $c\0{\infty}$ be (a fixed constant) such that
the metric $\hg\0{\infty}$ is $c\0{\infty}$-bounded (with respect to
$\cA\0{\bS^n}$).
Similarly, let $c\0{\lambda}$ be such that
the metric $\widehat{(g\0{\lambda})}\0{\lambda}$ is $c\0{\lambda}$-bounded
(with respect to $\cA\0{\bS^n}$).
And  let $c\0{\bS^n}$ be such that
the canonical metric $\sigma\0{\bS^n}$ is $c\0{\bS^n}$-bounded
(see remark 3.4).
\\

\noindent {\bf Proposition 5.2.3.} {\it There is a constant 
 $\lambda\0{\infty}$ such that the following holds.
If the metric $g\0{\lambda}$ is warped $\epsilon$-hyperbolic outside $\B_{\lambda-(1+\xi)}(\R^{n+1})$
with charts of excess $\xi>1$, then the metric $h_\ssl$ is $\eta$-hyperbolic
outside $\B_r$
with charts of excess $\xi-1$, where } $$\eta\geq
C_1\,\Big(\,e^{-r}+\,\frac{12}{d}\epsilon\0{0}\,\Big)\,+\,
C_6\,\Big(\,e^{-2\lambda}\,+\,\epsilon \Big)$$
\noindent {\it provided $\lambda>\lambda\0{\infty}$ and $d>2+4\,\xi$.
Here $C_1=C_1(c\0{\bS^n}+c\0{\infty},n, \xi)$ (see 3.3.3) and $C_6=C_6(\xi)$.}\\

\noindent {\bf 5.2.4. Remarks.}\\
\noindent {\bf 1.} 
We can take $C_6=e^{18+6\xi}$.\\
\noindent {\bf 2.} The number $\lambda\0{\infty}$ can be defined by the following property:
for every $\lambda>\lambda\0{\infty}$ we have      
$|\widehat{(g\0{\lambda})}\0{\lambda}-\hg\0{\infty}|
<\frac{1}{2^{n+1}\,(n!)^2\, c^{n+1}\0{\infty}}$.  Note that $\lambda\0{\infty}$ exists
because, by assumption, $\{g\0{\lambda}\}$ has cut limit at $b=0$,
i.e $\widehat{(g\0{\lambda})}\0{\lambda}\stackrel{C^2}
{\longrightarrow}\hg\0{\infty}$.  \\
\noindent {\bf 3.} Even though $h$ is hyperbolic on $B\0{\lambda-d}$
it is not $\epsilon$-hyperbolic near $o$, for $\epsilon$ small (see example in 3.1).
Hence, rigorously,  Proposition 5.2.4 should say that $h\0{\lambda}$ is $\eta$-hyperbolic
{\it outside $B\0{1+\xi}$}.\\
\noindent {\bf 4.} It is important not to mix the constants
$C_5$ and $C_6$ in the formula that appears in the statement of 5.2.4.
The constant $C_5$ depends of $c\0{\infty}$, and it is essential that
the constant $C_6$ does not depend on $c\0{\infty}$.\\

\noindent {\bf Proof of Proposition 5.2.3.} 
Let $p\in \R^{n+1}-B\0{1+\xi}$ and write $p=(x\0{0},t\0{0})\in
\bS^n\times \R^+$. We have two cases.\\

\noindent {\bf First case. $t\0{0}>\lambda-(1+\xi)$.} \\
Hence  $p\notin \bar{B}\0{\lambda-(1+\xi)}$.
By proposition 4.4.2 the metric $\cW\0{\lambda}\big( g\0{\lambda}\big)$
is $\eta'$-hyperbolic outside $B\0{\lambda-(1+\xi)}$  with charts of excess $\xi-1$,
provided
\begin{equation*}
\eta'\geq e^{18+6\xi}\,\Big(\,  e^{-2\lambda}\,+\, \epsilon\,\Big)
\tag{1}
\end{equation*}

\noindent In particular we can find an $\eta'$-hyperbolic chart $\phi$ for
$\cW\0{\lambda}\big( g\0{\lambda}\big)$ centered at $p$, with excess $\xi-1$.
But by item (iii) of 5.2.2 and the hypothesis $2+4\,\xi<d$
we get that $h\ssl=\cW\0{\lambda}\big( g\0{\lambda}\big)$
outside $B\0{\lambda-2(1+\xi)}$.
Therefore $\phi$ is also a
an $\eta'$-hyperbolic chart for $h\ssl$ centered at $p$, with excess $\xi-1$.\\

\noindent {\bf Second case. $t\0{0}\leq\lambda-(1+\xi)$.} \\
Hence  $p\in \bar{B}\0{\lambda-(1+\xi)}$.
By corollary 3.4.2 the metric
$\cH_{_{(\lambda-d),d}}\Big( {\overline{(g_\ssl )}}_\ssl\,\Big)$
is $\eta''$-hyperbolic outside $B_r$ with charts of excess $\xi$, provided
\begin{equation*}
\eta''\geq C_1\,\Big(\,  e^{-r}\,+\, \frac{12}{d}\epsilon\0{0}\,\Big)
\tag{2}
\end{equation*}

\noindent 
 In particular we can find an $\eta''$-hyperbolic chart $\phi'$ for
$\cH_{_{(\lambda-d),d}}\Big( {\overline{(g_\ssl )}}_\ssl\,\Big)$ centered at $p$, with excess $\xi$.
But by item (iv) of 5.2.2 
we have that $h\ssl=\cH_{_{(\lambda-d),d}}\Big( {\overline{(g_\ssl )}}_\ssl\,\Big)$
on $B\0{\lambda}$.
This together with the assumption $t\0{0}\leq\lambda-(1+\xi)$ imply that $\phi'$ is also a
an $\eta''$-hyperbolic chart for $h\ssl$ centered at $p$, with excess $\xi$.\\

\noindent Hence $h\0{\lambda}$ is $\eta$-hyperbolic with charts of excess
$\xi-1$ provided $\eta\geq \eta',\, \eta''$. Thus from (1) and (2) we 
see that can take 

$$\eta\geq
C_1\,\Big(\,e^{-r}+\,\frac{12}{d}\epsilon\0{0}\,\Big)\,+\,
C_6\,\Big(\,e^{-2\lambda}\,+\,\epsilon \Big)$$
\noindent were $C_1=C_1(c\0{\bS^n}+c\0{\infty},n, \xi)$
 and $C_6=C_6(\xi)=e^{18+8\xi}$.  Also $\epsilon\0{0}$ is as in 3.4.2.
By 3.4.3 we have $\epsilon\0{0}=\epsilon\0{0}(c\0{\lambda})$,
therefore we can write $C'=C'(c\0{\lambda},n,\xi)$.
To complete the proof of the proposition we need the following claim.\\

\noindent {\bf Claim.} {\it Let  $\lambda\0{\infty}$ be such that
$\lambda>\lambda\0{\infty}$ implies
$|\widehat{(g\0{\lambda})}\0{\lambda}-\hg\0{\infty}|
<\frac{1}{2^{n+1}\,(n!)^2\, c^{n+1}\0{\infty}}$. Then we can take $c\0{\lambda}<2c\0{\infty}$, that is 
$\widehat{(g\0{\lambda})}\0{\lambda}$ is $(2c\0{\lambda})$-bounded.}\\

\noindent {\bf Proof of claim.} Write $a=\widehat{(g\0{\lambda})}\0{\lambda}$,
$b=\hg\0{\infty}$ and $c=c\0{\infty}$. We have $|a-b|<\frac{1}{2^{n+1}(n!)^2c^{n+1}}$. Therefore $|a|\leq |a-b|+|b|<\frac{1}{2^{n+1}(n!)^2c^{n+1}}+c<2c$ (because $c>1$).\\

Now, we also have $|b|\,|b^{-1}a-I|<\frac{1}{2^{n+1}(n!)^2c^{n+1}}$.
This together with $|b^{-1}|<(n-1)!c^{n+1}$ (because $b$ is $c$-bounded)
imply $|b^{-1}a-I|<n\, (n-1)! c^{n+1}\frac{1}{2^{n+1}(n!)^2c^{n+1}}=
\frac{1}{2^{n+1}n!}$. But this implies that $det\,(b^{-1}a)\,>\, 1-n!\,2^n\frac{1}{2^{n+1}n!}=\frac{1}{2}$. Therefore $det\, a>\frac{1}{2c}$.
Consequently $a$ is $(2c)$-bounded. This proves the claim.\\

It follows from the claim that, if  $\lambda>\lambda\0{\infty}$, we have
 $C'(c\0{\lambda},n,\xi)<C_5(c\0{\infty},n,\xi)$,
where $C_5$ is obtained from $C'$ by replacing $\epsilon\0{0}(c\0{\lambda})$
by $\epsilon\0{0}(2c\0{\infty})$ (note that $\epsilon\0{0}(c)$
is increasing on $c$, see 3.4.3) and $c\0{\lambda}$ by
$2c\0{\infty}$  in $C_1'$.
Therefore we can take $\eta\geq C_5\, \Big(\,e^{-(\lambda-d)} \,+\,\epsilon\,+\,\frac{1}{d}\, \Big)$. This completes the proof of the proposition.\\

In the next corollary we use the constants $a\0{n}(\epsilon,\xi)$ defined in 3.3.10,
$a'\0{n}(c\0{*},\epsilon,\xi)$ and $d\0{n}(c\0{*},\epsilon,\xi)$ defined in (3.4.4).\\

\noindent {\bf Theorem 5.2.5.} {\it Let $\{g\ssl\}$ be a $\odot$-family and assume it
has a cut limit at $b=0$. Then there is a constant 
 $\lambda\0{\infty}$ such that the following holds.
If the metric $g\0{\lambda}$ is warped $\epsilon$-hyperbolic outside $\B_{\lambda-(1+\xi)}(\R^{n+1})$
with charts of excess $\xi>1$, then the metric $h_\ssl$ is $(3\epsilon'+C_6\epsilon)$-hyperbolic outside $\B_{a\0{n+1}(\epsilon',\xi)}(\R^{n+1})$
with charts of excess $\xi-1$, provided
\begin{enumerate}
\item[(i)] $\lambda-d>\lambda\0{\infty}+\lambda\0{n+1}(c\0{\infty},\epsilon,\epsilon',\xi)$  
\item[(ii)] $d\geq d\0{n+1}(c\0{\infty},\epsilon',\xi)$.
\end{enumerate}

\noindent where $\lambda\0{n+1}(c\0{*},\epsilon,\epsilon',\xi)=
a\0{n+1}(c\0{n+1},\epsilon',\xi)+\frac{1}{2}ln\big(\frac{C_6}{\epsilon}\big)$.}\\

\noindent {\bf Proof.} Take $r=\lambda-d-(1+\xi)$. Hypothesis (i) implies
(1) $e^{-r}<\frac{C_1}{\epsilon'}$, \, (2)   $e^{-2\lambda}<\frac{C_6}{\epsilon}$.
This together with hypothesis (ii) and 5.2.3 imply that outside $\B_r$ the metric
$h\ssl$ is $(3\epsilon'+C_6\epsilon)$-hyperbolic. Note that from the definitions
we get that $r>a\0{n+1}(\epsilon',\xi)+(1+\xi)$. Finally note that on the
ball of radius $r+(1+\xi)=\lambda-d$ the metric is hyperbolic and use 3.3.10.
This proves the corollary.\\

\vspace{.8in}

\noindent {\bf 5.3. Hyperbolic extensions of continued $\odot$-families.}\\
Fix $\xi>0$. Let $M^n$ have metric $h$ and center $o$. Hence we can identify $M-\{o\}$ with $\bS^{n-1}\times\R^+$ and 
$M$ with $\R^n$. Let $\{ h_\ssl\}_{\lambda>(2+\xi)}$ be a family of variable metrics on $M$, that is, $o$ is a center for all $h_\ssl $
and the ray structure (hence also the geodesic spheres at $o$)
of all $h_\ssl$ coincide with that of $h$.
Assume that the family  $\{ h_\ssl\}$
has cut limits on $J_c=(-\infty, c]$, $c>0$. Write as before $h_\ssl=\big(h_\ssl\big)_r+dr^2$.
We also assume that the family has {\it eventually constant (unwarped) cuts} that is there is $d>0$ such that
\begin{equation*}
\big( \hh\0\ssl \big)_{_{\lambda-b}}\,=\, \sigma\0{\bS^{n-1}}\,\,\,\,\,\,\,\,\,\,\,{\mbox{for all}}\,\,\,\,\,\,\,\,\,\,b>d
\tag{1}
\end{equation*}

Hence the metrics $h_\ssl$ are canonically hyperbolic near the origin. 
An example of such family is a continued family corresponding to a $\odot$-family that has cut limits on $[0,c)$ (see 5.2).\\

We consider now the hyperbolic extensions $f=\cE_k(h)$ and $f_\ssl=\cE_k(h_\ssl)$ on $\HH^k\times M$. 
Since the function $r$ (distance to $o\in M$) is the same on the spaces $(M,h_\ssl)$ we can write 
\begin{equation*}
f_\ssl\,=\, \big( cosh\,^2(r) \big)\, \sigma\0{\HH^k}\, +\, h_\ssl
\tag{2}
\end{equation*}

\noindent We write $\cE_k^\ssl(M)=\big( \HH^k\times M, f_\ssl \big)$
and  $\cE_k(M)=\big( \HH^k\times M, f \big)$.
Let $o\in\HH^k\sbs \HH^k\times M$. By remark 2.2.3 (1) we have that $o$
is a center for  $\cE_k(M)$ and all  $\cE_k^\ssl(M)$.
Moreover we have that all geodesics of $\cE_k^\ssl(M)$ emanating from $o$ coincide with the corresponding
geodesics of $\cE_k(M)$, that is, the ray structures
of all  $\cE_k^\ssl(M)$ coincide with the ray structure of $\cE_k(M)$ (see remark 2.2.3 (2)). Hence we can identify $\cE_k^\ssl(M)-\{ o\}$ with $\bS^{k+n-1}\times\R^+$ canonically (i.e independently of $\lambda$). 
Therefore we can write 
$$
f_\ssl\,=\,\big(f_\ssl\big)_s +ds^2
$$
Moreover, the $\Xi$-coordinates on (a subset of) the spheres $\bS_s\big(\cE_K(M)\big)=\bS_s^\ssl\big(\cE_K(M)\big)$ defined on section 2.4
are independent of $\lambda$. We will use these coordinates and all notation from section 2.4. (Recall they are independent of $\lambda$.)\\

Fix $\beta_0>0$ and $a>0$ such that $a+\beta_0\in (0,\pi/2)$ and denote by $S=S_{_{a+\beta_{_0}}}\big( \bS^{k-1} \big)$
the spherical neighborhood of $\bS^{k-1}=\bS^{k+n-1}\cap\HH^k$ \, in \,$\bS^{k+n-1}$ of width $a+\beta_0$.
Using the $\Xi$-coordinates on $\bS^{k+n-1}$ we can write $S-\bS^{k-1}=\bS^{k-1}\times\bS^{n-1}\times (0,\beta_0)$.
We want to determine whether $\{ f_\ssl\}$ has cut limits or not. For this we need to change the indexation of the family $\{ f_\ssl\}$.\\

Let $\lambda(s)=sinh^{-1}\big( sinh(s)\, sin (\beta_0)  \big)$ (see remark 2.4.1). That is, $\lambda(s)$ is the length of
the side opposite to the angle $\beta_0$ of the right hyperbolic triangle whose hypotenuse is $s$. Define
$$j\0{s}\,=\, f\0{\lambda(s)}$$

\noindent {\bf Proposition 5.3.1.} {\it The family $\{j_s \}$ has cut limits over $S$, on $J_{c'}$, where  $c'<c-ln\big(\frac{sin(\beta_0+a) }{sin(\beta_0)}\big)$.}\\

\noindent {\bf Proof.} From corollary 4.2.2. we can express $\big(\hf_\ssl\big)_s$ in $\Xi$-coordinates:
$$\Xi_*\big(\,\, (\hf_\ssl)_{s+b}\,\,\big)\,=\, cos^2\, (\beta)\, \sigma\0{\bS^{k-1}}\, \, +\,\,
sin^2\, (\beta)\,(\hh\0\ssl)_{r(s+b,\beta)}\,\,+\,\,\,d\beta^2$$
\noindent  where $r=r(s,\beta)$ is as in remark 2.4.1. Therefore 
$$\Xi_*\big(\,\, (\hj_s)_{s+b}\,\,\big)\,=\, cos^2\, (\beta)\, \sigma\0{\bS^{k-1}}\, \, +\,\,
sin^2\, (\beta)\,\big(\hh\0{\lambda(s)}\big)_{r(s+b,\beta)}\,\,+\,\,\,d\beta^2$$

\noindent We want to find the limit $\lim_{s\ra\infty}\big(\hh\0{\lambda(s)}\big)_{r(s+b,\beta)}$. Take the inverse
of $\lambda=\lambda(s)$, and we get $s=s(\lambda)=sinh^{-1}\big( \frac{sinh(\lambda)}{ sin (\beta_0)}  \big)$.
Hence

$$\lim_{s\ra\infty}\big(\hh\0{\lambda(s)}\big)_{r(s+b,\beta)}\,\,=\,\, 
\lim_{\lambda\ra\infty}\big(\hh\0{\lambda}\big)_{\vartheta(\lambda,\beta, b)}$$
\noindent where 
$$\vartheta(\lambda,\beta, b)=r\big(s(\lambda)+b,\beta\big)=
sinh^{-1}\bigg(\,sinh\,\bigg\{ \,\,b+ sinh^{-1}\bigg( \frac{sinh(\lambda)}{sin(\beta_0)} \bigg) \,\,\,\bigg\}\,\,sin\, (\beta)\,    \bigg)$$
\noindent and a straightforward calculation shows 
\begin{equation*}
\lim_{\lambda\ra\infty}\Big(\,\vartheta (\lambda,\beta, b)\,-\, \lambda\,\Big)\,=\, b\,+\, ln\big( \,\frac{sin(\beta)}{sin(\beta_0)}\,\big)
\tag{3}
\end{equation*}
\noindent but this limit is not uniform on $\beta$ (the problem is when $\beta \ra 0)$ . 
To solve this problem we use assumption (1), i.e the fact that we are assuming
$\{ h_\ssl\}$ has eventually constant cuts. That is, for $\beta$ less than certain $\beta_1=\beta_1(s)$ we will have that
$\big(\hh\0{\lambda(s)}\big)_{r(s+b,\beta)}=\sigma\0{\bS^{n-1}}$. To compute $\beta_1$ we use the hyperbolic law of sines and get
$$
sin(\beta_1)\,=\, sin(\beta_0)\,\frac{sinh\,(\lambda - d)}{sinh(\lambda)}\,\,\,\,\stackrel{\lambda\ra\infty}{\longrightarrow} \,\,\,e^{-d}\, sin(\beta_0)
$$

\noindent which is greater than zero. Therefore the limit is uniform on $\beta$. Finally, to obtain the condition on $c'$ note that we want
$$
\lambda\,+\,\bigg(\, b\,+\, ln\big( \,\frac{sin(\beta)}{sin(\beta_0)}\,\big)\,\bigg)\,\, <\,\,\lambda\,+\, c
$$
\noindent where we are using equation (3). Since $\beta<\beta_0+a$, this inequality implies that we want
$b<c-ln\big(\frac{sin(\beta_0+a) }{sin(\beta_0)}\big)$. This proves the proposition.\\

It is important to remark that the index $s$ in the statement of proposition 5.3.1 depends on $\lambda$ and also on $\beta_0$.\\

Let us recapitulate. Choosing an index $\lambda$ we get the  metric $h_\ssl$ on $M$. This $h_\ssl$ 
in turn produces uniquely the metric $f_\ssl$ which is just $\cE_k(h_\ssl)$. Note that $f_\ssl$ does not depend on $\beta_0$. 
With a new choice, the angle (or width) $\beta_0$, we get a new indexation of the family $\{ f_\ssl\}$: we write $j_s=f_\ssl$,  
where $s=s(\lambda)$ (or $\lambda=\lambda (s))$.
The relevance of the indices is that they indicate which spherical cut will be taken for the given metric: cut $h_\ssl$ at $\lambda$,
and cut $j_s=f_\ssl=\cE_k(h_\ssl)$ at $s=s(\lambda)$.

\vspace{1in}

\noindent {\bf \large  Section 6. Normal Neighborhoods and Links.}\\

In this section we define  and and give some properties of several types of neighborhoods on
piecewise-spherical and piecewise-hyperbolic complexes, that will be needed in the proof of the main Theorem.\\

We will write
$\R^l_+=(0,\infty)^{l}=\{ (x_1,...,x_{l}),\,\, x_i> 0\}$,
$\bar{\R}^l_+=[0,\infty)^{l}=\{ (x_1,...,x_{l}),\,\, x_i\geq 0\}$ and $\bar{\HH}^l_+=\B_{_{\HH}}^l\cap\bar{\R}^l_+$,
where $\B_{_{\HH}}^l$ is the disc model of $\HH^l$.\\

Decompose $\R^{m+1}=\R^{k+1}\times\R^{l+1}$ (we allow permuting the coordinates) and write $\R^{k+1}=\R^{k+1}\times\{0\}$ and
$\R^{l+1}=\{ 0\}\times\R^{l+1}$.
The spheres $\bS^k=\R^{k+1}\cap\bS^m$ and $\bS^l=\R^{l+1}\cap \bS^m$ are linked (and disjoint)  and note that $\bS^l=\{ u\in\bS^m,
d_{\bS^m}(u, \bS^k)=\pi/2\} $ and analogously for $\bS^k$. Moreover, as mentioned in item (d) of 2.1, $\bS^m-\bS^l$ is
the open normal neighborhood of width $\pi/2$ of $\bS^k$ in $\bS^m$, and we write $\bS^m-\bS^l=\stackrel{\circ}{\s{N}}_{_{\pi/2}}(\bS^k;\bS^m)$.
Since the normal bundle $\stackrel{\circ}{\s{N}}(\bS^k;\bS^m)\ra\bS^k$ is trivial
we can give $\bS^m-\bS^l$ the {\it normal 
coordinates}: $(u,b)\in\bS^k\times \B_{_{\pi/2}}^{l+1}$, where the fibers $\{ u\}\times\B^{l+1}_{_{\pi/2}}$ intersect orthogonally
$\bS^k$ at $(u,0)$, the lines $t\mapsto (u,tw)$, $|w|=1$, are speed-one geodesics emanating orthogonally from $\bS^k$
and the canonical basis $\{ e_i\}$ of tangent vectors $e_i\in T_{_0}\B^{l+1}=\R^{l+1}$ is orthonormal. 
More generally we will write 
${\sf{N}}_{_{\beta}}(\bS^k;\bS^m)=\bS^k\times {\bar{\B}}_{_{\beta}}^{l+1}$, for $\beta<\pi/2$, which is the {\sf closed} normal neighborhood of 
$\bS^k$ in $\bS^m$ of width $\beta$. \\

Since we can identify $\B_{_{\pi/2}}^{l+1}-\{0\}$ with $\bS^l\times (0,\pi/2)$, we can also write
$$\bS^m-\big(\bS^k\cup\bS^l\big)\,\,=\,\,\stackrel{\circ}{\s{N}}_{_{\pi/2}}(\bS^k)-\{\bS^k\}\,\,=
\,\,\bS^k\times\bS^l\times (0,\pi/2)$$
\noindent hence we can assign to a point in
$\bS^m-\big(\bS^k\cup\bS^l\big)$ its {\it polar normal coordinates} $(u,w,r)\in \bS^k\times\bS^l\times (0,\pi/2)$. 
Note the similarity of these coordinates with the
$\Xi$-coordinates of section 2.4.\\

As mentioned in item (d) of section 2.1, we  can write the metric of $\sigma\0{\bS^m}$ on $\bS^m-\big(\bS^k\cup\bS^l\big)$
using the polar normal coordinates as the following doubly warped metric:\\

\noindent {\bf (6.0.1.)}\hspace{1.4in}
$\sigma\0{_{\bS^m}}\,=\, cos^2(r)\, \sigma\0{\bS^k}\,+\,sin^2(r)\,\sigma\0{\bS^l}\,+\, dr^2$\\

\noindent where $r$ is the distance to $\bS^k$.\\

All concepts above can be extended to the sphere $\bS^k_\mu$ of radius $\mu$, and we will write the sub-index $\mu$
in these cases; for instance ${\sf{N}}_{_{\beta}}(\bS^k_\mu;\bS^m_\mu)$ is the closed normal neighborhood of 
$\bS^k_\mu$ in $\bS^m_\mu$ of width $\beta$. 

\vspace{.8in}

\noindent {\bf 6.1. Normal neighborhoods on the sphere $\bS^m$.}\\
As in \cite{ChD} the {\it canonical all-right spherical $k$-simplex} is the spherical simplex with all edges of length
$\pi/2$. We will denote it by $\sd{m}$. Alternatively $\sd{m}=\bS^k\cap\bar{\R}^{m+1}_+$, $k\geq m$. 
The {\it open canonical all-right spherical $k$-simplex} is
$\dDelta^k=\bS^k\cap\R^{m+1}_+$, $k\geq m$. 
The {\it center of mass} of $\sd{m}$ is \,\,$o_{_{m}}=o_{_{\bS^m}}=(\frac{1}{\sqrt{m+1}}, \frac{1}{\sqrt{m+1}},...,\frac{1}{\sqrt{m+1}})$.\\

For $l\leq m$, we have
that  the $l$-face $\sd{m}\cap\R^{l+1}$ of $\sd{m}$ is
(isometric to) the canonical spherical simplex $\sd{l}$. Since there are $C^{m+1}_{l+1}$ ways of canonically embed
$\R^{l+1}$ in $\R^{m+1}$ (by completing the $l+1$ coordinates with zeros) we have that $\sd{m}$ has $C^{m+1}_{l+1}$ $l$-faces.
In fact the poset of faces of $\sd{m}$ is naturally the same
as the poset of faces of the canonical euclidean simplex $\Delta_{_{\R^{m}}}$ and there is a homeomorphism between $\Delta_{_{\R^{m}}}$ and
$\sd{m}$ sending $l$-faces to $l$-faces (use radial projection from the origin
$0\in\R^{m+1}$). 
We will also denote by $\sd{m}$ the poset of all $l$-faces of
$\sd{m}$, $0\leq l\leq m$. Note that every simplex in $\sd{m}$ has a well defined center of mass. As in the euclidean case, for a face $\Delta^j$ of $\Delta^m=\sd{m}$ we can define its {\it opposite face} as the unique
face $\Delta^l$ of dimension $l=m-j-1$ of $\sd{m}$ disjoint from $\Delta^j$.
Hence, also as in the euclidean case, we can identify $\Delta^m$ with the
join $\Delta^j\ast\Delta^l$. This identification can be done in the following way.
Given $x\in\Delta^m$ let $w\in\Delta^j$ the closest point to $x$ in $\Delta^j$, 
and let $u\in\Delta^l$ the closest point to $x$ in $\Delta^l$. Let $[w,u](t)$,
$t\in [0,\pi/2]$ be the (spherical) geodesic segment in $\Delta^m$ from $w$ to $u$.
Then $x=[w,u](\beta)$, where $\beta=\beta(x)$ is the angle between
$w$ and $x$. Note that the geodesic $[w,u]$ is perpendicular to both $\Delta^j$
(at $w$)
and $\Delta^l$ (at $u$).\\

The intersection of $\bS^m\sbs\R^{m+1}$ with the $2^{m+1}$   sets $\R^{m+1}_{j_1,...,j_{m+1}}=\{
(x_1,...,x_{m+1}),\, j_ix_i\geq 0\}$, where $j_i\in\{ -1,1\}$, gives a decomposition of $\bS^m$ in $2^{m+1}$
spherical $m$-simplices; and taking these $m$-simplices together with all their faces we obtain the
{\it canonical all-right  triangulation } of $\bS^m$. We denote the all-right  triangulation also by $\bS^m$.
We will use the symbol $\Delta$ to denote an element in $\bS^{m}$ and we will write $\Delta^k$ if we want to specify the dimension $k$ of $\Delta$.
The notation $\Delta^i<\Delta^k$ means that $\Delta^i$ is a face of $\Delta^k$, but sometimes we will also write $\Delta^i\sbs\Delta^k$.
Let $\Delta=\Delta^k\in\bS^m$, $\Delta\sbs\bS^k$, and let $\beta\in (0,\pi/2)$. The {\it closed normal neighborhood of} $\Delta$ 
of width $\beta$ is 
${\sf{N}}_{_{\beta}}(\Delta;\bS^m)=\Delta\times {\bar{\B}}^{m-k}_{_{\beta}}$, 
where we are using normal coordinates on 
$\stackrel{\circ}{\s{N}}_{_{\pi/2}}(\bS^k;\bS^m)$.
For $x\in \dDelta^k=\dDelta$ the {\it $\beta$-link of $\Delta$ in $\bS^m$ at $x$}
is given by $\sL\0{\beta}(\Delta,\bS^m)=\{x\}\times\bS^{m-k-1}$, where we are using
normal coordinates. Most of the time we will just write  $\sL(\Delta,\bS^m)$
and will not mention $x$ and $\beta$, unless it is necessary. In the special case
when $x=o\0{\Delta}$ we shall call the link $\sL(\Delta,\bS^m)$
the {\it canonical linked sphere of $\Delta^k$ of radius} $\beta$ and write
$\s{S}_{_{\Delta}}=\s{S}_{_{\Delta}}(\beta)=\{ o_{_{\Delta}}\} \times\bS^{m-k-1}_{_{\beta}}$, and sometimes we will write 
$\s{S}_{_{\Delta}}=\bS^{m-k-1}_{_{\beta}}$.
For $x\in\dDelta^k<\Delta^j$ the {\it $\beta$-link of $\Delta^k$ in $\Delta^j$ at $x$},
$\sL_\beta(\Delta^k,\Delta^j)$,
is the intersection of $\sL_\beta(\Delta^k,\bS^m)$, at $x$, with $\Delta^j$.
Hence $\sL_\beta(\Delta^k,\bS^m)=\bigcup\0{\Delta^k<\Delta^j} \sL_\beta(\Delta^k,\Delta^j) $ (as sets and complexes).\\

For a given $\Delta^k\in\bS^m$, the intersections $\s{S}_\Delta\cap \Delta'$, with $\Delta<\Delta'$ 
provide a triangulation for  $\s{S}_{_{\Delta}}=\bS^{m-k-1}_{_{\beta}}$, which, after rescaling using the dilation
$\bS_{_{\beta}}^{m-k-1}\ra\bS^{m-k-1}$, is just the canonical all-right  triangulation of $\bS^{m-k-1}$.\\

Note that from (6.0.1) at the beginning of section 6, we  can identify $\s{N}_{_\beta}(\Delta^k,\bS^m)-\Delta^k$
with $\Delta^k\times\bS^{l}\times (0,\beta)$, $l=m-k-1$,  with the doubly warped metric in
(6.0.1).\\

\noindent {\bf Remark 6.1.1.} The definition of link given above is that of the ``geometric" link. Other two used definitions are the ``simplicial" one and the one  given using directions. For $\Delta^k<\Delta^j$ the {\it simplicial link of
$\Delta^k$ in $\bS^m$} is the subcomplex of $\bS^m$ (with its canonical
all-right triangulation) formed by all simplices $\Delta^j$ such that (1) 
$\Delta^j$ is disjoint from $\Delta^k$, (2) $\Delta^k$ and $\Delta^j$
span a simplex. And, for $\Delta^k<\Delta^j$, the simplicial link of $\Delta^k$ in $\Delta^j$ is just the opposite face of $\Delta^k$ in $\Delta^j$.
Note that if we continue a geodesic $[x,u]$, with $u$ in the (geometric) $\beta$-link of $\Delta^k$ in $\Delta^j$ at $x$, we will hit a unique point in the opposite face of $\Delta^k$ in $\Delta^j$. This radial geodesic projection gives a relationship between 
geometric links and simplicial links. 
For $x\in\dDelta^k$ the {\it direction link of $\Delta^k$ in $\bS^m$ at $x$}
is the set of all vectors at $x$ perpendicular to $\Delta^k$. Using geodesics
emanating from $x$, perpendicularly to $\Delta^k$ we also get a relationship
between geometric links and the direction links.
These different definitions of link all come with natural all-right spherical metrics:
the geometric link with the rescaled induced metric (see paragraph above, before
this remark), the simplicial link with the induced metric and the direction link with the
angle metric. The relationships between the different definitions of link mentioned
above all respect the metrics.\\

All concepts above can be extended to the sphere $\bS^k_\mu$ of radius $\mu$, and we will write the sub-index $\mu$
in these cases; 
for instance we also denote by  $\bS^m_\mu$ is the all-right  triangulation of the sphere $\bS^m_\mu$ of radius $\mu$ and
$\Delta_\mu$ denotes an element of  $\bS^m_\mu$.\\

\vspace{.8in}

\noindent {\bf 6.2.  Sets of widths of spherical normal neighborhoods.}\\
Let ${\sf{B}}=\{\beta_k\}_{k=0,1,2...}$ be an indexed set of real numbers with $\beta_k\in(0,\pi/4)$
and $\beta_{k+1}<\beta_k$. We write $\s{B}(m)=\{\beta_0,...\beta_{m-1}\}$.
The set $\s{B}$ determines {\it the set of spherical} \,$\s{B}$-{\it neighborhoods}\,
${\sf{N}}_{\sf{B}}(\bS^m)={\sf{N}}_{\s{B}(m)}(\bS^m)=\{ {\sf{N}} _{_{\beta_k}}(\Delta^k,\bS^m) \}_{\Delta^k\in\bS^m,\, k<m}$,
for any sphere $\bS^m$ (of any dimension). Note that the normal neighborhoods of all
$k$-simplices $\Delta^k$ have the same width $\beta_k$.
The set $\s{B}$ is called a {\it set of widths of spherical normal neighborhoods}   or simply {\it a set of widths}. The set $\s{B}(m)$ is a {\it finite set of widths
of length $m$}.\\

The definitions above still make sense if we replace $\bS^m$ by $\bS^m_\mu$,
the $m$-sphere of radius $\mu$.\\

We are interested in pairs of sets of widths  $\big(\s{B},\s{A}\big)$,
$\s{B}=\{\beta_k\}$ and $\s{A}=\{\alpha_j\}$, having the following 
{\it disjoint neighborhood property}:\\

\noindent {\bf (6.2.1.)\, DNP:}\,\,\,\,\,{\it For every $k$ and $m$, the following sets 
are disjoint}
$$\bigg\{ \, \s{N}_{\beta_k}(\Delta^k,\bS^m)\, \,-\, \,\bigcup_{j<k} \s{N}_{\alpha_j}(\Delta^j,\bS^m) \,\bigg\}_{\Delta^k\in\bS^m}$$
\vspace{.2in}

\noindent The disjoint neighborhood property obtained by fixing $k$ and $m$
above will be denoted by {\bf DNP($k,m$)}. 
In this case we allow the sets of widths to be finite of length at least $k+1$.
Note that the ordering of the pair $\big(\s{B},\s{A}\big)$
is important.
Alternatively, for every $k$ and $m$ we have the following property. For any two different $k$-simplices $\Delta^k_1$ and $\Delta^k_2$
we have $$\s{N}_{\beta_k}(\Delta_1^k,\bS^m)\,\bigcap\, \s{N}_{\beta_k}(\Delta_2^k,\bS^m)\,\,\,\,
{\mbox{{\LARGE $\sbs$}}}\,\,\, \,\,\bigcup_{j<k} \s{N}_{\alpha_j}(\Delta^j,\bS^m) $$
\noindent That is, the $\s{B}$-neighborhoods of different $k$-simplices intersect only inside the $\s{A}$-neighborhood of the $(k-1)$-skeleton (which is equal to \,$\bigcup_{j<k} \s{N}_{\alpha_j}(\Delta^j,\bS^m)$).\\

\noindent {\bf Proposition 6.2.2.} {\it The pair of
(finite or infinite) sets of widths $\big(\s{B},\s{A}\big)$
satisfy} \,{\bf NDP($k,m$) } {\it if and only if} $$\frac{sin\,\beta\0{k}}{sin\,\alpha\0{k-1}}\leq\frac{\sqrt{2}}{2}$$\\

\noindent Note that the inequality condition is independent of $m$.
The proposition follows directly from lemmas 6.2.3 (taking $k=l$ and $\beta=\gamma$) and 6.2.4 given below, and the fact that $\{\alpha\0{k}\}$ is decreasing.\\

\noindent{\bf Lemma 6.2.3.} {\it  Let
 $\Delta^k$, $\Delta^l\in\bS^m$
and $\Delta^j=\Delta^k\cap\Delta^l$. Let $\alpha, \beta, \gamma\in (0,\pi/2)$
such that \,$\frac{sin\,\beta}{sin\,\alpha}, \,\frac{sin\,\gamma}{sin\,\alpha}\,\leq\,\frac{\sqrt{2}}{2}$. Then }

$$\s{N}_{\beta}(\Delta^k,\bS^m)\,\cap\, \s{N}_{\gamma}(\Delta^l,\bS^m)\,\,\,\,
{\mbox{{\Large $\sbs$}}}\,\,\, \s{N}_{\alpha}(\Delta^j,\bS^m) $$

\noindent {\bf Proof.} 
In this proof $\sL(\Delta,\bS^m)$ shall denote the simplicial link (see 6.1.1)
and $\s{Star}(\Delta,\bS^m)$ the simplicial star, which is the union of all
simplices containing $\Delta$. Note that $\s{N}_\beta(\Delta,\bS^m)\sbs
\s{Star}(\Delta,\bS^m)$,  for every $\Delta\in\bS^m$.\\

Write $\s{S}=\sL(\Delta^j,\bS^m)$, $\Delta'_1=\s{S}\cap\Delta^k$ and $\Delta'_2=\s{S}\cap\Delta^l$.
Then $\Delta'_i$ is a simplex in the all-right triangulation of $\s{S}$.
Also $\Delta'_1$ and  $\Delta'_2$ are disjoint.
Hence their distance in $\s{S}$ is at least $\frac{\pi}{2}$.\\

Suppose there is $q\in\s{N}_{\beta}(\Delta^k,\bS^m)\cap
\s{N}_{\gamma}(\Delta^l,\bS^m)$. Since both of these neighborhoods lie in
$\s{Star}(\Delta^j,\bS^m)$ there is geodesic segment $[p,q]$ in $\s{Star}(\Delta^j,\bS^m)$
with $p\in \Delta^j$ and $[p,q]$ perpendicular to $\Delta^j$ (note that
$p$ may lie in $\p\Delta^j$ and that the geodesic may not be unique
if $q\in\s{S}$). Write $\alpha'=d_{\bS^m}(p,q)$.
We have to prove $\alpha'\leq\alpha$. We assume $\alpha'>\alpha$ and get
a contradiction.
Let $q_1$ be the closest point in $\Delta^k$ to $q$
and $q_2$ be the closest point in $\Delta^l$ to $q$. We have
$a_1=d_{\bS^m}(q_1,q)\leq \beta$ and $a_2=d_{\bS^m}(q_2,q)\leq \gamma$.
We get a right (at $q_i$) spherical triangle with one side equal to $a_i$ and hypotenuse equal to $\alpha'$. Let $\theta_i$ be
the angle at $p$, that is, the angle opposite to the side of length $a_i$. Then by the spherical law of sines we get
$$
sin\,\theta_1\,=\,\frac{sin\, a_1}{sin\,\alpha'} < \frac{sin\,\beta}{sin\,\alpha}
\leq\frac{\sqrt{2}}{2}
$$ 
Consequently $\theta_1<\frac{\pi}{4}$. Similarly we get $\theta_2<\frac{\pi}{4}$.
Let $z_i$ be the intersection of $\s{S}$ with the ray at $p$ with direction $q_i$.
Analogously let $q'$ be the intersection of $\s{S}$ with the ray at $p$ with direction $q$. Note that $z_i\in\Delta'_i$. Also note that
$d_\s{S}(q,z_i)$ is equal to the length of the arc
$qz_i$ in $\s{S}$. This together with the fact that $\s{S}$ is a 
$(m-j-1)$-sphere of radius one imply that $d_\s{S}(q,z_i)=\theta_i$.
Therefore 
$$\frac{\pi}{2}\,\leq\, d_\s{S}(\Delta'_1,\Delta'_2)\,\leq\,
d_\s{S}(z_1,z_2)\,\leq
d_\s{S}(z_1,q')+d_\s{S}(q',z_2)\,\leq \,\big(\theta_1+\theta_2\big)\,
$$
Hence $\frac{\pi}{2}\leq \theta_1+\theta_2<\frac{\pi}{4}+\frac{\pi}{4}=\frac{\pi}{2}$ which is a contradiction. This proves the lemma. \\

\noindent {\bf Lemma 6.2.4.} {\it  Let
$\Delta^k_1$, $\Delta^k_2\in\bS^m$ be two different $k$-simplices,
and $\Delta^{k-1}=\Delta^k_1\cap\Delta^k_2$. Moreover the $\Delta^k_i$
span a simplex.
Let $\alpha, \beta\in (0,\pi/2)$ Suppose that }

$$\s{N}_{\beta}(\Delta_1^k,\bS^m)\,\cap\, \s{N}_{\beta}(\Delta_2^k,\bS^m)\,\,\,\,
{\mbox{{\Large $\sbs$}}}\,\,\, \s{N}_{\alpha}(\Delta^j,\bS^m) $$
\noindent {\it Then $\frac{sin\,\beta}{sin\,\alpha}\leq\frac{\sqrt{2}}{2}$. }\\

\noindent {\bf Proof.} The lemma is certainly true for $\bS^1$.
Using the spherical law of sines it is straightforward to verify the lemma for
$\bS^2$. The case $\bS^m$, $m>2$, can be reduced to the case $m=2$ using
the orthogonal sphere $S$ to $\Delta^{k-1}$ in $\bS^m$ at the center
$o$ of $\Delta^{k-1}$, as in the proof of lemma 6.2.3. This proves the lemma. \\

The next result says that {\bf DNP} implies a seemingly stronger version of itself.\\

\noindent {\bf Lemma 6.2.5.} {\it Suppose the pair of sets of widths
$(\sB, \sA)$ satisfies} {\bf DNP}. {\it Let $\Delta^j=\Delta^k\cap\Delta^l$,
$j<k\leq l$. Then}

$$\s{N}_{\beta_k}(\Delta^k,\bS^m)\,\bigcap\, \s{N}_{\beta_l}(\Delta^l,\bS^m)\,\,\,\,
{\mbox{{\LARGE $\sbs$}}}\,\,\, \,\,\bigcup_{i<k} \s{N}_{\alpha\0{i}}(\Delta^i,\bS^m) $$

\noindent {\bf Remark.} Note that the condition
$\Delta^j=\Delta^k\cap\Delta^l$, $j<k,l,$\, is equivalent to\,
$\Delta^k\not\sbs\Delta^l$ and $\Delta^l\not\sbs\Delta^k$,
where the empty set is considered a simplex of dimension -1.\\

\noindent {\bf Proof of lemma 6.2.5.}
From proposition 6.2.2 we have $\frac{sin\, \beta_l}{sin\,\alpha_j}\leq\frac{sin\, \beta_k}{sin\,\alpha_j}<\frac{\sqrt{2}}{2}$. And from lemma 6.2.3 we now get
$$\s{N}_{\beta_k}(\Delta^k,\bS^m)\,\bigcap\, \s{N}_{\beta_l}(\Delta^l,\bS^m)\,\,\,\,
{\mbox{{\LARGE $\sbs$}}}\,\,\,  \s{N}_{\alpha_j}(\Delta^j,\bS^m)
\,\,\,\,{\mbox{{\LARGE $\sbs$}}}\,\,\, \,\,
\bigcup_{i<k} \s{N}_{\alpha\0{i}}(\Delta^i,\bS^m) $$
\noindent This proves the lemma.

\vspace{.8in}

\noindent {\bf 6.3. Natural neighborhoods on the sphere $\bS^m$.}\\
Let ${\sf{B}}=\{\beta_k\}$ be a set of widths.  
Let $\Delta=\Delta^k\in\bS^m$. As in 6.1, let  
$\s{S}_\Delta=\bS^{m-k-1}_{\beta\0{k}}$ be the linked sphere of $\Delta$ of radius $\beta_k$ (at the center $o\0{\Delta}$ of $\Delta$). By intersecting with $\s{S}_\Delta$, the set of neighborhoods $\s{N}_\s{B}(\bS^m)$ determines the set  
$\s{N}(\s{S}_\Delta,\s{B})\,=\,\{ {\s{S}_\Delta\,\cap\,\sf{N}} _{_{\beta_j}}(\Delta^j,\bS^m) \}_{\Delta^j\in {\sf{D}}(\bS^m)}$.
It is straightforward to verify that for each $\Delta^j$, with $\Delta<\Delta^j$, there is 
some $\beta'_{j-k-1}\in (0,\pi/2]$ such that
$$ \s{S}_\Delta\,\cap\,\s{N} _{_{\beta_j}}(\Delta^j,\bS^m)\,=\, \s{N} _{_{\beta'_{j-k-1}}}(\s{S}_\Delta\cap\Delta^j,\s{S}_\Delta)$$
\noindent where the last term is the  $\beta'_{j-k-1}$-normal neighborhood of
the simplex $\s{S}_\Delta\cap\Delta^j$ in $\s{S}_\Delta$.
Hence we can write $\s{N}(\s{S}_\Delta,\s{B})\,=\,\s{N}_{\s{B'}(m-k-1)}(\s{S}_\Delta)$
where $\s{B'}(m-k-1)=\{\beta'_0,...,\beta_{m-k-2}\}$ and we also say that $\s{N}_{\s{B}'(m-k-1)}(\s{S}_\Delta)$ is the set of $\s{B'}(m-k-1)$-neighborhoods
of $\s{S}_\Delta$. [The set $\s{B'}$ is not strictly speaking a set of widths because
some $\beta'_i$ may be $\geq \pi/4$.]
Note that $\s{B}'(m-k-1)$ depends only on $\s{B}$ and the dimension $k$ of $\Delta^k$.
The next lemma gives this relation explicitly.\\

\noindent {\bf Lemma 6.3.1.} {\it For $i=0,...,m-k-2$
we have the following identity}\\  $$sin\,\bigg(\,\frac{\beta'_i}{sin\,\beta_k}\,\bigg)\,\,=\,\, \frac{sin(\beta_{k+1+i})}{sin(\beta_k)}$$

\noindent {\bf Proof.} Let $p\in\s{S}_{\Delta}\cap\big(\s{S}_{_{\beta_j}}(\Delta^j,\bS^m)\big)$, where $\Delta=\Delta^k<\Delta^j$. Then there is a
$q\in\Delta^j$ such that $d_{_{\bS^m}}(p,q)=d_{_{\bS^m}}(p,\Delta^j)=\beta_j$. Since  $d_{_{\bS^m}}(o_{_{\Delta}},p)=\beta_k$
we get a right (at $q$) spherical triangle with one side equal to $\beta_j$ and hypotenuse equal to $\beta_k$. Let $\theta$ be
the angle opposite to the side of length $\beta_j$. Then, by the spherical law of sines we get
$$
sin(\theta)=\frac{sin(\beta_j)}{sin(\beta_k)}
$$

\noindent But $d_{_{\s{S}_\Delta}}(p,\s{S}_\Delta\cap\Delta^j)=\theta\, sin\,\beta_k$. This proves the lemma.\\

By rescaling $\s{S}_\Delta$ using the dilation $\s{S}_\Delta=\bS^{m-k-1}_{_{\beta_k}}\ra\bS^{m-k-1}$ we get the set of neighborhoods
$\s{N}_{\s{B''}(m-k-1)}(\bS^{m-k-1})$ on $\bS^{m-k-1}$, where $\s{B''}(m-k-1)=\{ \beta''_0,...\beta''_{m-k-2}  \}$ and $\beta''_j=\frac{\beta'_j}{sin\,\beta_k}$.\\

Therefore the set of neighborhood widths $\s{B}$ on $\bS^m$ induces the set  $\s{B''}(m-k-1)$ on $\bS^{m-k-1}$,
by considering $\bS^{m-k-1}$ as the (rescaled) link of the $k$ simplices in the all-right  triangulation of $\bS^m$. And lemma 6.3.1 implies that the relationship
between $\s{B}(m)=\{ \beta_0,...,\beta_{m-1}\}$ and $\s{B''}(m-k-1)=\{\beta''_0,...,\beta''_{m-k-2}\}$ is\\

\noindent {\bf (6.3.2.)}\hspace{1.5in}
$sin(\beta''_j)=\frac{sin(\beta_{k+j+1})}{sin(\beta_k)}$\\

Let $\s{B}=\{ \beta_{_i}\}_{i=0,1,...}$ be a set of widths. 
We say that $\s{B}$ is a {\it natural set of neighborhood widths for all spheres} if $\s{B}(m-k-1)=\s{B''}(m-k-1)$ for all $m $ and $k$ with $m>k$.\\

\noindent {\bf Corollary 6.3.3.} {\it The set of widths $\s{B}=\{\beta_{_i}\}$ is natural if and only if
$sin(\beta_i)=sin^{i+1}(\beta_{_0})$ and $\beta\0{0}<\pi/4$.}\\

\noindent {\bf Proof.} It follows directly from equation 6.3.2 with $j=0$ that 
$sin(\beta_{_{k+1}})=sin(\beta_{_{k}})\,sin(\beta_{_0})$. This proves the corollary.\\

Given $\varsigma\in (0,\frac{\sqrt{2}}{2})$ we define $\s{B}(\varsigma)=\{ \beta_{_i}\}$ by $\beta_{_i}=sin^{-1}\big( \varsigma^{i+1} \big)$.
Hence the corollary says that $\s{B}$ is natural if and only if $\s{B}=\s{B}(\varsigma)$, for some $\varsigma\in (0,\frac{\sqrt{2}}{2})$. 
In fact, in this case we have $\varsigma=sin(\beta_{_0})$.\\

Let $\varsigma\in (0,\frac{\sqrt{2}}{2})$ and $c>1$. We denote by $\sB(\varsigma;c)=\{\gamma_i\}$ the set defined by  $\gamma_i=sin^{-1}\big( c\,\varsigma^{i+1}\big)$.
Note that $\s{B}(\varsigma;c)$ is a set of widths provided $c\varsigma<\frac{\sqrt{2}}{2}$. Proposition 6.2.2 implies the next corollary.\\

\noindent {\bf Corollary 6.3.4.} {\it The pair of sets of 
widths $\big(\s{B}(\varsigma;c),\s{B}(\varsigma;c')\big)$ satisfy}  {\bf DNP} {\it
provided $\frac{c}{c'}\,\varsigma<\frac{\sqrt{2}}{2}$.}\\

Note that we can take $c=c'=1$ which implies that a natural set of widths $\s{B}(\varsigma)$
satisfy {\bf DNP} with $\s{A}=\s{B}=\s{B}(\varsigma)$.
\vspace{.8in}

\noindent {\bf 6.4. Neighborhoods in piecewise spherical complexes.}\\
We use the definitions of spherical complex and all-right  spherical complex 
given in section 1 of \cite{ChD}. 
The set of all-right  spherical simplices of an all-right  spherical complex $P$ will also be denoted by
$P$. For simplicity we shall assume in this paper that all spherical complexes satisfy the
``intersection condition" of simplicial complexes: {\sf every two simplices intersect in a (possibly empty) common face}.
Next we give another (not too different) way of defining
links, using the notions in 6.3.\\

Let $P$ be an all-right spherical complex and
let $\Delta^j$ be a simplex in the all-right  spherical complex $P$. Then the set of all 
$\s{Link}(\Delta^j,\Delta^k)$, with $\Delta^j<\Delta^k$ is also an all-right spherical complex which we denote
by $\s{Link}(\Delta^j, P)$. Thus $\sL(\Delta^j,P)=\bigcup_{\Delta^j<\Delta^k\in P}\sL(\Delta^j,\Delta^k)$
as sets and complexes.
Hence the set $\{\Delta^k\}_{\Delta^j\,<\,\Delta^k\in P}$ is in one-to-one correspondence
with the set of spherical simplices of $\sL(\Delta,P)$, that is $\Delta^k$ corresponds to $\sL(\Delta^j,\Delta^k)$,
which is an all-right  spherical simplex of dimension $k-j-1$ in $\sL(\Delta^j,P)$. 
The all-right  spherical metric on $\s{Link}(\Delta^j,P)$
will be denoted by $\sigma_{_{\s{Link}(\Delta^j,P)}}$.\\


\noindent {\bf Remark 6.4.1.} In the paragraph above we did not specify
the type of link we were using, i.e. geometric, simplicial or direction (see 6.1.1).
Any of these will lead to corresponding definition of $\sL(\Delta,P)$, but
they are all equivalent as metric complexes (note that we are assuming $P$
has the ``intersection condition".  
We will use any of the definitions depending on the situation.
For the geometric link sometimes we will need to specify the base point $x$
and the radius $\beta$ (see 6.1.1).\\

\noindent {\bf Lemma 6.4.2.} {\it We have that}
$$
\sL\bigg( \sL (\Delta^j,\Delta^k), \sL(\Delta^j, P) \bigg)\,\,\,=\,\,\,\sL (\Delta^k, P)
$$
\noindent {\it provided $\Delta^j\,<\,\Delta^k$. (This is an equality of 
all-right spherical metric complexes.)}\\

\noindent {\bf Remark.} If we use the simplicial definition of links this equality
is actually an equality of sets (see 6.1.1).\\

\noindent {\bf Proof.} Let $\Delta^j\,<\,\Delta^k$. Let $\Delta^l$ be the opposite face of $\Delta^j$ in $\Delta^k$.
The statement in the lemma written using the simplicial definition of link is:
$$
\sL\big( \Delta^l, \sL(\Delta^j, P) \big)\,\,\,=\,\,\,\sL (\Delta^k, P)
$$
We prove this. We have that $\Delta^i\in \sL\big( \Delta^l, \sL(\Delta^j, P) \big)$ is equivalent to (recall $\Delta^j\cap\Delta^l=\emptyset$)
(1) $\Delta^i\cap\Delta^j=\emptyset$, (2) $\Delta^i\cup\Delta^j$ is contained in a simplex, 
(3) $\Delta^i\cap\Delta^l=\emptyset$, (4) $\Delta^i\cup\Delta^l\cup\Delta^j$ is contained in a simplex.
On the other hand  $\Delta^i\in \sL (\Delta^k, P)$ is equivalent to
(a) $\Delta^i\cap\Delta^k=\emptyset$, (b) $\Delta^i\cup\Delta^k$ is contained in a simplex. 
Since  $\Delta^l$ is opposite to $\Delta^j$ in $\Delta^k$ we have that (a)+(b) if
and only if (1)-(4). This proves the Lemma.\\

Define
the {\it closed normal neighborhood of $\Delta^j$ in $\Delta^k$ of width $\beta$ } as
\,$\s{N}_{_\beta}(\Delta^j,\Delta^k)=\s{N}_{_\beta}(\Delta^j,\bS^m)\cap\Delta^k$. As mentioned
at the end of section 6.1, equation (6.0.1) implies that
$\s{N}_{_\beta}(\Delta^j,\Delta^k)-\Delta^j$ is isometric to $\Delta^j\times\sd{k-j-1}\times (0,\beta)$
with the doubly warped metric of (6.0.1).
If $\Delta^j$ is a simplex in the all-right  spherical complex $P$, we define
the {\it closed normal neighborhood of $\Delta^j$ in $P$ of width $\beta$ } as
$$\s{N}_{_\beta}(\Delta^j,P)=\bigcup_{\Delta^j\,<\,\Delta^k\in P}\s{N}_{_\beta}(\Delta^j,\Delta^k)$$
\noindent and we can identify $\s{N}_{_\beta}(\Delta^j,P)-\Delta^j$ with
$\Delta^j\times\s{Link}(\Delta^j,P)\times (0,\beta)$ with the doubly warped metric (see item F of section 1)\\

\noindent {\bf (6.4.3.)}\hspace{1in}
$ cos^2(r)\, \sigma_{_{\bS^j}}\,+\,sin^2(r)\,\sigma_{_{\s{Link}(\Delta^j,P)}}\,+\, dr^2$\\

\noindent where $r$ is the distance to $\Delta^j$, and $\sigma_{_{\s{Link}(\Delta^j,P)}}$
is the spherical all-right  metric on $\s{Link}(\Delta^j,P)$.\\

Finally, let $\s{B}=\{\beta_k\}$ be a set of widths (see 6.2). Then, for any all-right  spherical complex $P$
the set $\s{B}$ induces the set of
neighborhoods $\s{N}_{_\sB}(P)=\big\{ \s{N}_{_{\beta_k}}(\Delta^k,P) \big\}_{\Delta\in P}$. The next lemma is the spherical complex version
of lemma 6.2.3.\\

\noindent{\bf Corollary 6.4.4.} {\it  Let
 $\Delta^k$, $\Delta^l\in P$
and $\Delta^j=\Delta^k\cap\Delta^l$. Let $\alpha, \beta, \gamma\in (0,\pi/2)$
such that \,$\frac{sin\,\beta}{sin\,\alpha}, \,\frac{sin\,\gamma}{sin\,\alpha}\,\leq\,\frac{\sqrt{2}}{2}$. Then }

$$\s{N}_{\beta}(\Delta^k,P)\,\cap\, \s{N}_{\gamma}(\Delta^l,P)\,\,\,\,
{\mbox{{\Large $\sbs$}}}\,\,\, \s{N}_{\alpha}(\Delta^j,P) $$

\noindent {\bf Proof.} The proof is the same as the proof of lemma 6.2.5
with minor obvious changes. Just recall that we are assuming $P$ to have the
intersection condition. Also note that if $p$, $q$, $q_i$ are as in the proof of
6.2.5 the they all three lie in an all-right simplex in $P$. This proves the
corollary.\\

As in section 6.2, the next two results follow directly from corollary 6.4.4. 
The first is a version of {\bf DNP} (see 6.2.1) for $P$, obtained by replacing $\bS^m$ by $P$.\\

\noindent {\bf Corollary 6.4.5.} {\it Let the pair of sets of widths $\big(\s{B},\s{A}\big)$
satisfy } {\bf DNP}. {\it Then for any all-right spherical complex $P$ and
dimension $k$ we have that the following sets are disjoint}

$$\bigg\{ \, \s{N}_{\beta_k}(\Delta^k,P)\, \,-\, \,\bigcup_{j<k} \s{N}_{\beta_j}(\Delta^j,P) \,\bigg\}_{\Delta^k\in P}$$\\

\noindent The next is a version of lemma 6.2.5 for general $P$.\\

\noindent {\bf Corollary 6.4.6.} {\it Let the pair of sets of widths $\big(\s{B},\s{A}\big)$ satisfy } {\bf DNP}. {\it Then for any all-right spherical complex $P$ and
$\Delta^j=\Delta^k\cap\Delta^l$, $j<k\leq l$, simplices in $P$, we have that}

$$\s{N}_{\beta_k}(\Delta^k,P)\,\bigcap\, \s{N}_{\beta_l}(\Delta^l,P)\,\,\,\,
{\mbox{{\LARGE $\sbs$}}}\,\,\, \,\,\bigcup_{i<k} \s{N}_{\alpha\0{i}}(\Delta^i,P) $$

\vspace{.8in}

\noindent {\bf 6.5. Neighborhoods in piecewise hyperbolic cones.}\\
Recall that we are denoting by  $\bar{\HH}^{k+1}_+$  the subspace $\B_{_{\HH}}^{k+1}\cap\bar{\R}^{k+1}_+$ of $\HH^{k+1}$,
where $\B_{_{\HH}}^{k+1}$ is the disc model of $\HH^{k+1}$. We can identify 
$\bar{\HH}^{k+1}_+-\{ o\}$ 
with $\sd{k}\times\R^+$ with metric $sinh^2(s)\,\sigma\0{\bS^k}+ds^2$, where $s$ is the distance to the ``vertex" $o$.
We say that $\bar{\HH}^{k+1}_+$ is the {\it infinite hyperbolic cone of} $\sd{k}$ and write $\rC\sd{k}=\bar{\HH}^{k+1}_+$.\\

Let $\Delta^j<\Delta^k$ and let $\Delta^l$ be the opposite face of $\Delta^j$  in $\Delta^k$. We have $l=k-j-1$. Since $\rC\Delta^j=\bar{\HH}_+^{j+1}\sbs\HH^{j+1}$, $\rC\Delta^l=\bar{\HH}_+^{l+1}\sbs\HH^{l+1}$
and $\rC\Delta^k=\bar{\HH}_+^{k+1}=\bar{\HH}_+^{j+1}\times \bar{\HH}_+^{l+1}
\sbs\HH^{k+1}$ we can write \\

\noindent {\bf (6.5.1.)}\hspace{1.6in}$\rC\Delta^k=\rC\Delta^j\times\rC\Delta^l$\\

\noindent  with warped metric $cosh^2(r)\,\sigma_{_{\HH^{j+1}}}+\sigma_{_{\HH^{l+1}}}$,
where $r$ is the distance in $\HH^{l+1}$ to $o$.
Note that the order of the decomposition here is important (see section 1 and proposition 2.3.3). The identification above can be done explicitly in the 
following way. Let $p\in\rC\Delta^k\sbs\HH^{k+1}=\cE_{\HH^{j+1}}(\HH^{l+1})$. 
We use the functions (or coordinates) given in section 2.3: $s$, $r$, $t$, $y$, $v$, $x$,
$u$, $w$. Then $p=sx\in\rC\Delta^k$, $(s,x)\in\R^+\times\Delta^k$, corresponds to
$(y,v)=(tw,ru)\in\rC\Delta^j\times\rC\Delta^l$, $(t,w)\in\R^+\times\Delta^j$,
$(r,u)\in\R^+\times\Delta^l$. Note that $x=[w,u](\beta)$, where
$\beta$ is as in section 2.3, i.e. it is the angle between $w$ and $x$ (see 
also section 6.1 and recall that $\Delta^k=\Delta^j\ast\Delta^l$).\\

Let $P$ be an all-right  spherical complex. Recall that $P$ is constructed by
gluing the all-right spherical simplices $\Delta\in P$. The {\it infinite piecewise hyperbolic cone of $P$} is the space  $\rC P$ 
obtained by gluing the hyperbolic cones $\rC\Delta$, $\Delta\in P$ using the same rules used for obtaining $P$.
Note that all vertex points of the $\rC \Delta$ get glued to a unique {\it vertex} \, $o=o_{_{\rC P}}$.
The piecewise hyperbolic metric on $\rC P$ shall be denoted by $\sigma_{_{\rC P}}$ and its corresponding geodesic metric by
$d_{_{\rC P}}$. Note that all (constant speed) rays emanating from $o$ are length minimizing geodesics defined on $[0,\infty)$.
Then we can identify $\rC P-\{ o\}$ with $P\times\R^+$ with warped metric $sinh^2(r)\,\sigma_{_P}+dr^2$, where $r$ is the distance
to the vertex $o$. (See item F of section 1.)\\

Even though $\sigma\0{\rC P}$ is not (generally) smooth, the set 
of rays emanating from the vertex $o\0{\rC P}$ gives  a well defined
ray structure as in section 1D (see also 2.2).\\

For $s\geq 0$ we denote the {\it open ball of radius s of} $\rC P$ {\it centered at o } by $\B_{_s}(\rC P)$. Note that
this ball is the ``finite open cone" $P\times (0, s)\cup\{o\}$, where we are using the identification above.
The {\it closed ball} will be denoted by ${\bar{\B}}_{_s}(\rC P)$ and the {\it sphere of radius s}, $s>0$, will be denoted
by $\bS_s(\rC P)$, which we shall sometimes identify with $P\times\{ s\}$ or simply with $P$.\\

We say that the cones $\rC \Delta$, $\Delta\in P$, are the {\it cone simplices of}\, $\rC P$ and the
{\it faces } of the cone simplex $\rC \Delta$ are the $\rC \Delta'$, with $\Delta'<\Delta$. The set of all cone simplices 
will also be denoted by $\rC P$. 
The complex $\rC P$ (i.e. $\rC P$ together with its cone faces) is an {\it all-right hyperbolic cone
complex}. Note that the metric $\sigma\0{\rC P}$ can be deduced from the
hyperbolic cone structure. 
For $\Delta\in P$, $\beta>0$ and
$x\in \rC\dDelta-o$, the {\it $\beta$-link $\sL_\beta(\rC \Delta, \rC P)$ of $\rC\Delta$
in $\rC P$, at $x$,} is the set of all end points of geodesics of length 
$\beta$ emanating
perpendicularly to $\rC\Delta$ at $x$. As in section 6.4, by rescaling 
$\sL_\beta(\rC \Delta, \rC P)$ has a natural all-right spherical structure.\\

\noindent {\bf Remarks 6.5.2.}\\
\noindent {\bf 1.} The definition of link above is that of the {\it geometric link}
(see remark 6.1.1). 
We can also define the {\it direction link} of $\rC \Delta$
in $\rC P$ at $x$, which is an all-right spherical complex with the angle metric.
Of course there is a natural identification between these,  as all-right spherical complexes.\\
\noindent {\bf 2.} We also have the following natural identification
between all-right spherical complexes
$$
\sL(\rC\Delta,\rC P)\,=\,\sL(\Delta, P)
$$
\noindent {\bf 3.} In the following chapters we will need an identification as the one above but with with the geometric links and as actual subsets of $\rC P$. For this we will use
the exponential map. We do this in section 7.\\ 
\noindent {\bf 4.} If $\rC P$ is an all-right hyperbolic cone complex
and $\Delta$ is an spherical simplex then $\rC\Delta\times \rC P$ is naturally
an all-right hyperbolic cone complex. In fact the hyperbolic cone structure
is given by the identification

$$
\rC\Delta\times\rC P\,=\, \rC\big(\,    \Delta\,\ast\, P      \,\big)
$$

\noindent where $\Delta\,\ast\, P$ is the all-right spherical complex whose
simplices are the joins $\Delta\ast\Delta'$, $\Delta'\in P$. The identification
above is done using identification 6.5.1 simplexwise. (Note that $\Delta'$ is the opposite face of $\Delta$ in $\Delta\ast\Delta'$.)\\ 
\noindent{\bf 5.}
For $\Delta=\Delta^k\in P^m$, $k<m$,  we shall use the notation
$\s{CS}(\Delta,P)=\rC \big( \s{Star}(\Delta,P) \big)$,
where  $\s{Star}(\Delta,P)$ here is the {\it simplicial star}, formed by all the simplices
in $P$ that contain $\Delta$ as a face.
We have that $\s{CS}(\Delta,P) $ can be identified with
$\rC\Delta\times \rC\big(\sL(\Delta, P)\big)$ with metric $$cosh^2(r)\,\sigma_{_{\HH^{k+1}}}\,+\, \sigma_{_{\rC\sL(\Delta, P)}}$$
\noindent where $r$ is the distance in $\rC\big(\sL(\Delta, P)\big)$ to the vertex $o\in \rC\big(\sL(\Delta, P)\big)$. The identification here is an identification of
all-right hyperbolic cone complexes and it is obtained using remark 4,
because $\Delta\ast\sL(\Delta,P)=\s{Star}(\Delta,P)$.
We can do this identification explicitly in the following way. We use (6.5.1)
and the variables $s$, $r$, $t$, $y$, $v$, $x$,
$u$, $w$ given in 2.3. Using the explanation given in the paragraph following (6.5.1)
in a simplexwise fashion we see that an element $p=sx\in\s{CS}(\Delta^k,P)$ can be written as $(tw,ru)\in\rC\Delta^k\times\rC\Delta^l\sbs\rC \Delta^k\times\rC\sL(\Delta^k,P)$, where $\Delta^k\ast \Delta^l$ is a simplex
in the (simplicial) star $\s{Star}(\Delta^k,P)$; that is, $\Delta^l$ is a simplex
in $\sL(\Delta^k,P)$. Since we can write $x=[w,u](\beta)$, $\beta$ is the angle
between $w$ and $x$ (see paragraph
after (6.5.1) and 6.1), the identification is given by
$$s\,\big[w,u\big](\beta)\,=\, \Big( t\,w\,,\,r\, u \Big)
$$
\\
\noindent {\bf 6.} 
 As mentioned above, even though $\sigma\0{\rC P}$
is not in general smooth it has a well defined ray structure. Furthermore,
we can consider $o\0{\rC P}$ as the {\it center} of $\rC P$ (see section 2.2).
Hence it makes sense to consider $\cE_k(\rC P)$ (with metric $\cE_k(\sigma\0{\rC P})$),
where $P$ is a spherical all-right complex.\\
\noindent {\bf 7.} For $\Delta=\Delta^k\in P$, recall that $\rC \Delta^k=
\bar{\HH}^{k}_+
\sbs\HH^k$, therefore we can write 
$$\s{CS}(\Delta,P)=\cE\0{\rC \Delta}\bigg(\rC\big( \sL(\Delta,P)\big)\bigg)
\,\,\sbs\,\, \cE\0{k}\bigg(\rC\big( \sL(\Delta,P)\big)\bigg)$$
\noindent where we consider $\s{CS}(\Delta,P)$ with metric $\sigma\0{\rC P}$
and $\rC\big(\sL(\Delta,P)\big)$ with metric $\sigma\0{\rC\,\sL(\Delta,P)}$.\\
\noindent {\bf 8.} Also, $\s{CS}(\Delta,P)-\rC \Delta$ can be identified with
$\rC\Delta\times \sL(\Delta, P)\times\R^+$ with metric $$cosh^2(r)\,\sigma_{_{\HH^{k+1}}}\,+\, sinh^2(2)\,\sigma_{_{\sL(\Delta, P)}}\,+\, dr^2$$
\noindent  where $r$ is as above.\\\\

For a cone simplex $\rC\Delta\in \rC P$, we define its {\it closed normal neighborhood of width $s$} by (see remark 4 above)\\

\noindent ${\bf (6.5.3.)} \hspace{1in}
\sN_{_s}(\rC\Delta,\rC P)\,=\, \rC\Delta\times{\bar{\B}}_{_s}(\rC \sL(\Delta,P))\,\,\sbs\,\, 
\s{CS}(\Delta,P)
$\\

\noindent where we are using the identification given in remark 5.
Hence $\sN_{_s}(\rC\Delta,\rC P)$ is the union of (the images of)
all geodesics of length $s$ emanating perpendicularly to and from $\rC \Delta$.\\

\noindent {\bf Lemma 6.5.4.} {\it \, Let  $\Delta^j<\Delta^k\in P$. Then}\
$$
\s{CS}(\Delta^k, P)\,=\,\rC\Delta^k\times\rC\sL(\Delta^k,P)=
\rC\Delta^j\times\s{CS}(\Delta^l,\sL(\Delta^j,P))\sbs
\rC\Delta^j\times\rC\sL(\Delta^j,P)
$$

\noindent  {\it where $\Delta^l=\sL(\Delta^j,\Delta^k)$.
A similar statement holds if we replace $\sN$ by $\stackrel{\circ}{\sN}$.}\\

\noindent {\bf Remark.} The first equality is given by remark 5. The last 
inclusion follows from the inclusion $\s{Star}(\Delta^l,\sL(\Delta^j,P))\sbs\sL(\Delta^j,P)$.
The middle equality in the statement of the lemma is an equality of hyperbolic 
cone complexes.\\

\noindent {\bf Proof.}
We have
$$\rC\Delta^k\times\rC\sL(\Delta^k,P)\,=\,
\rC\Delta^j\times\bigg(\rC\Delta^l\times\rC \sL(\Delta^k,P)\bigg)\,=
\rC\Delta^j\times\s{CS}(\Delta^l,\sL(\Delta^j,P))
$$
\noindent where the first equality follows from 6.5.1 and the second equality
from 6.4.2 and remark 5 above. This proves the lemma.\\

Here is a metric version of lemma 6.5.4.
Let $\Delta^j$, $\Delta^k$, and $\Delta^l$ be as in lemma 6.5.4. Let $h:\bS^{m-k-1}\ra\sL(\Delta^k,P)$
be a homeomorphism
and consider the cone of $h$, $\rC h:\R^{m-k}\ra\rC\sL(\Delta^k,P)$.
Let $f'$ be a variable metric on $\R^{m-k}$, then $f'$ and $\sigma\0{\R^{m-k}}$
have the same ray structure. The metric $f=h_*f'$ is a metric on
$\rC\sL(\Delta^k,P)$, and it has the same ray structure as 
$\sigma\0{\rC\sL(\Delta^k,P)}$.
We can consider the (restriction of the) metric $\cE_k(f)$ on 
$\cE_{\rC\Delta^k}(\rC\sL(\Delta^k,P))=\rC\Delta^k\times\rC\sL(\Delta^k,P)$.
And, since we have $\sL(\Delta^k,P)=\sL(\Delta^l,\rC\sL(\Delta^j,P))$ (see 6.4.2) 
the metric $f$ is also a metric on $\rC\sL(\Delta^l,\sL(\Delta^j,P))$, and we
can consider the metric $\cE_j(\cE_l(f))$ on $\cE_{\rC\Delta^j}\big(\cE_{\rC\Delta^l}(\rC\sL(\Delta^l,\sL(\Delta^j,P))\big)=\rC\Delta^j\times\rC\Delta^l\times
\rC\sL(\Delta^l,\sL(\Delta^j,P))$.\\

\noindent {\bf Corollary 6.5.5.} {\it Using the identification in 6.5.4 we get}
$$
\cE_k(f)\,=\,\cE_j\big(\cE_l(f))
$$

\noindent {\bf Proof.} The proof follows from proposition 2.3.3 and
the proof of lemma 6.5.4. This proves the corollary.\\

\noindent {\bf Lemma 6.5.6.} {\it \, Let  $\Delta^j<\Delta^k\in P$. Then}\
$$
\sN_s(\rC\Delta^k,\rC P)\,=\,\rC\Delta^j\times \sN_s\bigg( \rC\Delta^l,\rC\sL(\Delta^j,P) \bigg)
$$

\noindent {\it where $\Delta^l=\sL(\Delta^j,\Delta^k)$.
A similar statement holds if we replace $\sN$ by $\stackrel{\circ}{\sN}$.}\\

\noindent {\bf Remark.} Note that the left-hand side of the equality, $\sN_s(\rC\Delta^k,\rC P)$, is a subset of $\s{CS}(\Delta^k,P)$. The right-hand side
is a subset of $\rC\Delta^j\times\rC\sL(\Delta^j,P)$. By lemma 6.5.4 we can
write  $\s{CS}(\Delta^k,P)\sbs\rC\Delta^j\times\rC\sL(\Delta^j,P)$. Lemma 6.5.6
says that under this inclusion $\sN_s(\rC\Delta^k,\rC P)$ corresponds to
$\rC\Delta^j\times \sN_s\big( \rC\Delta^l,\rC\sL(\Delta^j,P) \big)$.\\

\noindent {\bf Proof.} We have
$$
\sN_s(\rC\Delta^k,\rC P)\,=\,\rC\Delta^j\times\bigg(\rC\Delta^l\times{\bar{\B}}_{_s}(\rC \sL(\Delta^k,P))\bigg)\,=
\,\rC\Delta^j\times \sN_s\bigg( \rC\Delta^l,\rC\sL(\Delta^j,P) \bigg)
$$
\noindent where the first equality follows from (6.5.1) and (6.5.3) and the
last from 6.4.2 and (6.5.3). This proves the lemma.\\\\

Recall that $\bs(s,\beta)=sinh^{-1}\big(sinh (s)\,sin(\beta)  \big)$ (see remark 2.4.1).\\

\noindent {\bf Lemma 6.5.7.} {\it \, Let $s>0$, $\beta\in(0,\pi/2)$ and $\Delta\in P$. 
Then}
$$
\sN_{\,s_{_\beta}}\big(   \rC\Delta, \rC P \big)\,\,\,\,{\mbox{\LARGE $\cap$}}\,\,\,\,\bS_s(\rC P)\,\,\,{\mbox{{\LARGE =}}}
\,\,\,\sN_{_\beta}\big( \Delta, P \big)\times\{ s\}
$$
\noindent {\it where $s_{_\beta}=\bs(s,\beta)$ and we are identifying $\bS_s(\rC P)$ with $P\times \{ s\}$
(thus $\sN_{_\beta}\big( \Delta, P \big)\times\{ s\}\,\,\sbs\,\, P\times\{ s\}\,\,=\,\,\bS_s(\rC P)$).}\\

\noindent {\bf Proof.} Denote the vertex of $\rC\big( \sL(\Delta, P) \big)$ by $o'$. Note that both sides of the equality above
are contained in  $\bS_s(\rC P)$. From (6.5.3) and $\beta<\pi/2$ we also get that both sides are contained in $\s{CS}(\Delta,P)$.
 Let $p\in \bS_s(\rC P)$, then $d_{_{\rC P}}(o,p)=s$. From 6.5.3 we can write
$p=(x,y)\in \rC\Delta\times{\bar{\B}}_{_{s'}}(\rC \sL(\Delta,P))$, for some $s'$. Consider the geodesic segments $a=[o,p]$, 
$b=[(x,o'),p]\sbs\{ x\}\times \rC \sL(\Delta,P)$ 
and $c=[o, (x,o')]\sbs\rC\Delta\times\{o'\}$. The length of $a$ is $s$. Since
each $\{x\}\times\rC\sL(\Delta,P)$ is totally geodesic (see remark 5) we get that
the length of $b$ is $s'$. Also since $p\in\s{CS}(\Delta,P)$ we have that
$p\in\rC\Delta^j$ for some $\Delta^j\in P$ containing $\Delta$.
But $\rC\Delta^j=\bar{\HH}_+^{j+1}\sbs \HH^{j+1}$ is totally geodesic in
$\rC P$ hence all three segments $a$, $b$, $c$ are contained in $\rC\Delta^j$.
Therefore we get a hyperbolic geodesic triangle with sides
$a$, $b$, $c$, whose angle at $(x,o')$ is $\pi/2$ (because $\rC\Delta\times\{ o'\}$ and $\{ x\}\times \rC\sL(\Delta,P)$ are perpendicular, see remark 5).
Let $\beta'$ be the angle at $o$. Then $p\in \sN_{\,s_{_\beta}}\big(   \rC\Delta, \rC P \big)$ if and only if
$s'\leq s_{_\beta}$. Also $p\in \sN_{_\beta}\big( \Delta, P \big)$ if and only if $\beta'\leq \beta$. But
$s_{_\beta}=sinh^{-1}\big(sinh (s)\,sin(\beta)  \big)$ and by the hyperbolic law of sines we also get that 
$s'=sinh^{-1}\big(sinh (s)\,sin(\beta')  \big)$. Consequently
$s'\leq s_{_\beta}$ is
equivalent to $\beta'\leq \beta$. This proves the lemma.
\vspace{.8in}

\noindent {\bf 6.6. More neighborhoods in hyperbolic cones.}\\
Let $\xi>0$, $\varsigma>0$ and $c>1$  with $c\,\varsigma<e^{-(4+\xi)}$.
Write  $\sB=\sB(\varsigma; c)=\{\beta_i\}$
and $\sA=\sB(\varsigma)=\{\alpha_i\}$
be set of widths as in 6.3. 
We have $sin\,\beta_i=c\,\varsigma^{i+1}$, $sin\,\alpha_i=\varsigma^{i+1}$. 
Since $e^{-(4+\xi)}<\frac{\sqrt{2}}{2}$, the condition $c\,\varsigma<e^{-(4+\xi)}$ together with corollary 6.3.4 imply that
$\big(\s{B},\s{A}\big)$ and $\big(\s{B},\s{B}\big)$
satisfy {\bf DNP} in section 6.2.\\

Given a number $r>0$ and an integer $k\geq 0$ we define  $r\0{k}=r\0{k}(r)=
sinh^{-1}\big(\frac{sinh(r)}{sin(\alpha\0{k})}\big)$. By convention we also
set $r\0{-1}=r$. (Alternatively we could declare that every set of widths $\{\alpha_k\}$
has a (-1) term $\alpha_{-1}$ always equal to $\pi/2$.)
Let $k$ and $m$ be integers with $m\geq 2$ and $0\leq k\leq m-2$. Define
$s\0{m,k}
=sinh^{-1}\big(\frac{sinh(r)\,sin(\beta\0{k})}{sin(\alpha\0{m-2})}\big)=
sinh^{-1}\big(sinh(r\0{m-2})\,sin(\beta\0{k})\big)$.
We write $r\0{m,k}=r\0{m-k-3}$. Note that $r\0{m,k}<s\0{m,k}$.\\

Let $P=P^m$ be an all-right  spherical complex with $m\leq \xi$, and let $r>(4+\xi)$. For every $\Delta^k\in P$, $0\leq k\leq m-2$, define the following subsets of $\rC P$:\\

{\small $\begin{array}{lll}
\cY (P, \Delta^k, r,\xi,(c,\varsigma))&=&
\stackrel{\circ}{\sN}\0{s\0{m,k}}(\rC\Delta^k, \rC P)\,-\, \Bigg(\bigcup_{j<k}
\s{N}\0{r\0{m,j}}(\rC\Delta^j,\rC P)\Bigg)
\,-\, \B_{r\0{m-2}-(2+\xi)}(\rC P)\\\\
\cY (P, r,\xi,(c,\varsigma))&=&
\rC P\,-\, \Bigg(\bigcup_{j<m-1}
\s{N}\0{r\0{m,j}}(\rC\Delta^j,\rC P)\Bigg)
\,-\, \B_{r\0{m-2}-(2+\xi)}(\rC P)\end{array}$}\\\\


Since $\xi$, $c$ and $\varsigma$ will remain constant, in the rest of this section we will write
$\cY(P,\Delta^k,r)$ instead of $\cY(P,\Delta^m,r,\xi,(c,\varsigma))$.\\

\noindent {\bf Lemma 6.6.1.} {\it For $r>(4+\xi)$ we have the following properties}
\begin{enumerate}
\item[{\bf (i)}] $\cY(P,\Delta^k,r)
\,\,\,{\mbox{{\Large $\sbs$}}}\,\,\, \stackrel{\circ}{\sN}\0{s\0{m,k}}(\rC\dDelta^k, \rC P)\,\,\,{\mbox{{\Large $\sbs$}}}\,\,\,int\,\s{CS}(\Delta^k,P)$
\item[{\bf (ii)}] $\cY(P,\Delta^k,r)
\,\,\,{\mbox{{\Large $\cap$}}}\,\,\,
\sN\0{r\0{m,j}}(\rC\Delta^j, \rC P)=\emptyset$, for $j< k$.
\item[{\bf (iii)}] $\cY(P,\Delta^j,r)
\,\,\,{\mbox{{\Large $\cap$}}}\,\,\,
\,\B_{r\0{m-2}-(2+\xi)}(\rC P)=\emptyset$, 
\item[{\bf (iv)}] $\rC P\,-\,\B_{r\0{m-2}-(2+\xi)}(\rC P)\,\,=\,\,\cY(P,r)\,\,\cup\,\,
{\mbox{\Large $\bigcup$}}\0{\Delta^k\in P,\,k\leq m-2}\cY(P,\Delta^k,r)$
\item[{\bf (v)}] {\it $\Delta^j\cap\Delta^k=\emptyset$ \,\,\,implies\,\,\,
$\sN\0{s\0{m,j}}(\rC \Delta^j,\rC P)\cap\sN\0{s\0{m,k}}(\rC \Delta^k,\rC P)=\emptyset$}
\item[{\bf (vi)}] {\it $\Delta^j\cap\Delta^k=\emptyset$ \,\,\,implies\,\,\,
$\cY(P,\Delta^j,r)\cap\cY(P,\Delta^k,r)=\emptyset$}
\item[{\bf (vii)}] {\it  $\Delta^k=\Delta^i\cap\Delta^j$, \,with $k<i,\,j$,\,\,
implies \,\,, $\cY(P,\Delta^i,r)\cap\cY(P,\Delta^j,r)=\emptyset$}
\item[{\bf (viii)}]  {\it for any two different $k$-simplices\,\, $\Delta_1^k,\,\Delta_2^k$ \,\,\,we have\,\,\, $\cY(P,\Delta_1^k,r)\cap\cY(P,\Delta_2^k,r)=\emptyset$}
\item[{\bf (ix)}] 
{\it $\Delta^j<\Delta^k$ \,\,\,implies\,\,\,
$\cY(P,\Delta^j,r)\cap\cY(P,\Delta^k,r)\,\,\,{\mbox{{\Large $\sbs$}}}\,\,\,
\s{CS}(\Delta^j,P)\,\,\,{\mbox{{\Large $\sbs$}}}\,\,\,\s{CS}(\Delta^k,P)$}
\end{enumerate}\vspace{.2in}

\noindent {\bf Proof.}
The statements (ii) and (iii) follow from the definition of $\cY(P,\Delta,r)$
(recall $\s{CS}(\Delta,P)=\rC \big(\s{Star}(\Delta,P)\big)$, see section 6.2).
To prove (i) note that from the definition of $\cY(P,\Delta,r)$
we have
$\cY(P,\Delta^k,r)\sbs\, \sN\0{s\0{m,k}}(\rC\dDelta^k, \rC P)\sbs\s{CS}(\Delta^k,P)$. If  a point $p\in\p \s{CS}(\Delta^k,P)\cap\stackrel{\circ}{\sN}\0{s\0{m,k}}(\rC\Delta^k, \rC P)$ then its distance to $\p \Delta^k$ is
$<s\0{m,k}$. But it can be checked that $r\0{m,j}>s\0{m,k}$, $j<k$
(this follows from $c\varsigma<e{-(4+\xi)}<1$). Therefore
$p\notin \cY(P,\Delta^k,r)$. This proves (i).
Statement (ix) follows from (i).\\

Since $ P\,\,=\,\, {\mbox{\large $\bigcup$}}\0{\Delta^k\in P}
int\,\Delta^k$ we have 
that $$\rC P\,\,=\,\, {\mbox{\large $\bigcup$}}\0{\Delta^k\in P}
\rC (int\,\Delta^k)\,\,=\,\,
{\mbox{\large $\bigcup$}}\0{\Delta^k\in P}
\stackrel{\circ}{\sN}\0{s\0{m,k}}(\rC\Delta^k, \rC P)$$
\noindent This together with (iii)
and the definition of $\cY(P,r)$  imply that we can prove (iv) by showing, by induction on $k$, that $U=
{\mbox{\large $\bigcup$}}\0{\Delta^l\in P,\, l\leq m-2}\cY(P,\Delta^l,r)$ contains
$\stackrel{\circ}{\sN}\0{s\0{m,k}}(\rC\Delta^k, \rC P)\,\,
-\,\,\B_{r\0{m-2}-(2+\xi)}(\rC P)$
for every $k$-simplex of $P$, $k\leq m-2$. For $k=0$ this statement holds because $\cY(\Delta^0,P)=
\,\,\stackrel{\circ}{\sN}\0{s\0{m,0}}(\rC\Delta^0, \rC P)
-\B_{r\0{m-2}-(2+\xi)}(\rC P)$. Assume $U$ contains
every   $\stackrel{\circ}{\sN}\0{s\0{m,j}}(\rC\Delta^j, \rC P)\,\,
-\B_{r\0{m-2}-(2+\xi)}(\rC P)$, for all $j<k$. By the definition of
$\cY(\Delta^k,P)$ we have that
$\stackrel{\circ}{\sN}\0{s\0{m,k}}(\rC\Delta^k, \rC P)
\,-\,\B_{r\0{m-2}-(2+\xi)}(\rC P)$ is contained in
 
$$\Bigg[\cY(\Delta^k,P)\,\,\,\,\,\,\cup\,\,\,\,\,\,
{\mbox{\large $\bigcup$}}\0{\Delta^j\in P,\, j<k}
\sN\0{r\0{m,j}}(\rC\Delta^j, \rC P)
\Bigg]\,\,\,-\,\,\, \B_{r\0{m-2}-(2+\xi)}(\rC P)$$\\

\noindent This together with the fact that
$s\0{m,k}>r\0{m,k}$ and the inductive hypothesis imply that
$\stackrel{\circ}{\sN}\0{s\0{m,k}}(\rC\Delta^k, \rC P)
\,-\,\B_{r\0{m-2}-(2+\xi)}(\rC P)
\,\,\,{\mbox{{\large $\sbs$}}}\,\,\,U$. This proves (iv).\\

To prove the other two statements we need the following lemma.\\

\noindent {\bf Lemma 6.6.2.}
{\it For $t\geq r\0{m-2}-(2+\xi)$ and $r>(4+\xi)$ we have (see lemma 6.5.7)}
$$\sN\0{r\0{m,k}}\big(   \rC\Delta, \rC P \big)\,\,\,\,{\mbox{\LARGE $\cap$}}\,\,\,\,\bS_t(\rC P)\,\,\,{\mbox{{\LARGE =}}}
\,\,\,\sN\0{\theta\0{m,k}(t)}\big( \Delta, P \big)\times\{ t\}
$$

$$\sN\0{s\0{m,k}}\big(   \rC\Delta, \rC P \big)\,\,\,\,{\mbox{\LARGE $\cap$}}\,\,\,\,\bS_t(\rC P)\,\,\,{\mbox{{\LARGE =}}}
\,\,\,\sN\0{\phi\0{m,k}(t)}\big( \Delta, P \big)\times\{ t\}
$$
\noindent {\it where $\theta\0{m,k}(t)$ and $\phi\0{m,k}(t)$ are defined by 
the equations $sin(\theta\0{m,k}(t))=c''sin(\alpha\0{k})$, $sin(\phi\0{m,k}(t))=c''sin(\beta\0{k})$, with $c''=\frac{sinh(r\0{m-2})}{sinh(t)}<4e^{2+\xi}$.
Moreover $\theta\0{m,k}(t)$ and $\phi\0{m,k}(t)$ are well defined
and less that $\pi/4$.}\\

\noindent {\bf Proof.} From lemma 6.5.7
we have
$$sin (\theta\0{m,k}(t))\,\,=\,\, \frac{sinh(r\0{m,k})}{sinh(t)}\,\,=\,\,
\frac{sinh(r\0{m-2})}{sinh(t)}\,\,\frac{sinh(r\0{m,k})}{sinh(r\0{m-2})}\,\,=\,\,
c''sin(\alpha\0{k})$$
\noindent  and 
$$sin (\phi\0{m,k}(t))\,\,=\,\, \frac{sinh(s\0{m,k})}{sinh(t)}\,\,=\,\,
\frac{sinh(r\0{m-2})}{sinh(t)}\,\,\frac{sinh(s\0{m,k})}{sinh(r\0{m-2})}\,\,=\,\,
c''sin(\beta\0{k})$$
\noindent A simple calculation shows
that $c''<2e^{2+\xi}$, provided $r>4+\xi$ (thus $r\0{m-2}>4+\xi$). Hence
the definitions of $\alpha_k$ and $\beta_k$ at the beginning
of this section imply
$c''sin(\alpha\0{k})\,=\,c''\varsigma^{k+1}<\frac{\sqrt{2}}{2}$
and $c''sin(\beta\0{k})\,=\,c''c\,\varsigma^{k+1}<\frac{\sqrt{2}}{2}$.
This proves the lemma.\\

We now finish the proof of lemma 6.6.1. Statement (v) follows from lemma
6.6.2 and the fact that  $\beta$-neighborhoods, $\beta<\pi/4$, of disjoint
simplices in an all-right spherical complex are disjoint. Statement (vi) follows
from (v). We prove (vii). Note that $c''=c''(m,t)$.
Using items (i), (ii), and  lemmas 6.6.2 and 6.2.5 it is enough to prove that, for fixed $t$ and $m$, 
the pair of sets of widths $\big(\{ \phi\0{m,k}(t) \},\{\theta\0{m,k}(t)  \}\big)$
satisfies {\bf DNP}. But from the definitions we have 
$\{ \phi\0{m,k}(t) \}=\s{B}(\varsigma,c c'')$ and $\{\theta\0{m,k}(t)  \}=\s{B}(\varsigma, c'')$. Therefore lemma 6.3.4  and the condition
$c\,\varsigma<e^{-(4+\xi)}$ imply
$\big(\{ \phi\0{m,k}(t) \},\{\theta\0{m,k}(t)  \}\big)$
satisfies {\bf DNP}. This proves (vii).
Statement (viii) follows from (vii) by taking $i=j$.
This proves Lemma 6.6.1.\\

Define the sets\\ 

\hspace{.5in}$\begin{array}{lll}
\cX(P^m,\Delta^k,r)&=&\cY(P^m,\Delta^k,r)-\B_{r\0{m-2}}(\rC P^m)
\\\\\cX(P^m,r)&=&\cY(P^m,r)-\B_{r\0{m-2}}(\rC P^m)
\end{array}$\\\\

Alternatively, we can define  $\cX(P^m,\Delta^k,r)$ by the same formula 
that defines  $\cY(P^m,\Delta^k,r)$ with just one small change: in the last term
replace the radius $r\0{m-2}-2\xi$ by $r\0{m-2}$. Similarly for
 $\cX(P^m,r)$.\\

\noindent {\bf Lemma 6.6.3.} {\it For $\Delta^j<\Delta^k\in P$ we have}
$$\cY(P,\Delta^k,r)\,\,\,\,{\mbox{\Large $\sbs$}}\,\,\,\,\,\rC\Delta^j\,\,\,\,\times\,\,\,\,
\cX\bigg(\sL\Big(\Delta^j,P\Big), \Delta^l , r \bigg)
$$
\noindent {\it where $\Delta^l=\Delta^k\cap\sL(\Delta^j,P)$.}\\

\noindent {\bf Remark.} The left term in the lemma is a subset of
$\s{CS}(\Delta^k,P)$, thus also a subset of $\s{CS}(\Delta^j,P)$. The right term is a subset $\rC\Delta^j\times
\rC\sL(\Delta^j, P)$. Item (5) in remark 6.5.2 says that we can write
$\s{CS}(\Delta^j,P)=\rC\Delta^j\times
\rC\sL(\Delta^j, P)$. Lemma 6.6.3 says that 
$\cY(P,\Delta^k,r)$ corresponds to a subset of
$\rC\Delta^j\times
\cX\big(\sL\big(\Delta^j,P\big), \Delta^l , r \big)$ under this identification.\\

\noindent {\bf Proof.} Using lemma 6.5.6, the definitions of $\cY(P,\Delta^k,r)$
and $\cX\big(\sL\big(\Delta^j,P\big), \Delta^l , r \big)$,
and (v) of lemma 6.6.1 we see
that it is enough to prove the following three statements
\begin{enumerate}
\item[{\it (1)}]\,\,\,\,$\cY(P,\Delta^k,r)\,\,\,{\mbox{\Large $\sbs$}}\,\,\,
\rC\Delta^j\,\,\times\,\,
\stackrel{\circ}{\s{N}}\0{s\0{m-j-1,l}}\bigg( \rC\Delta^l,\rC\sL\big( \Delta^j,P\big)\bigg)$

\item[{\it (2)}] For $\Delta^i\in\sL(\Delta^j,P)$,
$i<l=m-j-1$, we have
$$\cY(P,\Delta^k,r)\,\,\,{\mbox{\Large $\cap$}}\,\,\,
\Bigg[\rC\Delta^j\,\,\times\,\,\s{N}\0{r\0{m-j-1,i}}\Big( \rC\Delta^i,\rC\sL\big( \Delta^j,P\big)\Big)
\Bigg]\,\,=
\,\,\emptyset$$
\item[{\it (3)}]
$\cY(P,\Delta^k,r)\,\,\,{\mbox{\Large $\cap$}}\,\,\,
\bigg[\rC\Delta^j\,\,\times\,\,\B\0{r\0{m-j-3}}\Big(\rC\sL\big( \Delta^j,P\big)\Big)
\bigg]\,\,=\,\,\emptyset$
\end{enumerate}

Statement (1) follows from (i) of lemma 6.6.1, lemma 6.5.6 and the
equalities $s\0{m,k}=s\0{m-j-1,k-j-1}$ and $l=k-j-1$. 
Statement (2) follows from (ii) of lemma 6.6.1, lemma 6.5.6 and the
statements $r\0{m,i+j+1}=r\0{m-j-1,i}$,  $i+j+1<k$. For (3) note that
(6.5.3) and the definition of $r\0{m,j}$ imply

$$\rC\Delta^j\times\B\0{r\0{m-j-3}}\Big(\rC\sL\big( \Delta^j,P\big)\Big)
\,\,\,\,\,\,\,{\mbox{\Large $=$}}\,\,\,\,\,\,\,
\sN\0{r\0{m,j}}\Big(\rC\Delta^j,\rC P\Big)
$$
This together with (ii) of lemma 6.6.1 imply (3). This proves the lemma.\\

\noindent {\bf Lemma 6.6.4.} {\it For $\Delta^k\in P$, $k\leq m-2$, we have}
$$\cY(P,r)\,\,\,\,{\mbox{\Large $\cap$}}\,\,\,\,
\cY(P,\Delta^k,r)\,\,\,\,{\mbox{\Large $\sbs$}}\,\,\,\,\,\rC\Delta^k\,\,\,\,\times\,\,\,\,
\cX\bigg(\sL\Big(\Delta^k,P\Big), r \bigg)
$$


\noindent {\bf Proof.} Using the definition of 
$\cX\big(\sL\big(\Delta^k,P\big), r \big)$ we see
that it is enough to prove the following three statements
\begin{enumerate}
\item[{\it (1)}]\,\,\,\,$\cY(P,\Delta^k,r)\,\,\,{\mbox{\Large $\sbs$}}\,\,\,
\rC\Delta^k\,\,\times\,\,
\rC\sL\big( \Delta^k,P\big)$

\item[{\it (2)}] For $\Delta^j\in P$, $\Delta^k<\Delta^j$,
$l\leq m-k-3$, 
and $\Delta^l$ opposite to $\Delta^k$ in $\Delta^j$, we have
$$\cY(P,r)\,\,\,{\mbox{\Large $\cap$}}\,\,\,
\Bigg[\rC\Delta^k\,\,\times\,\,\s{N}\0{r\0{m-k-1,l}}\Big( \rC\Delta^l,\rC\sL\big( \Delta^k,P\big)\Big)
\Bigg]\,\,=
\,\,\emptyset$$
\item[{\it (3)}]
$\cY(P,r)\,\,\,{\mbox{\Large $\cap$}}\,\,\,
\bigg[\rC\Delta^k\,\,\times\,\,\B\0{r\0{m-k-3}}\Big(\rC\sL\big( \Delta^k,P\big)\Big)
\bigg]\,\,=\,\,\emptyset$
\end{enumerate}

Statement (1) follows from (i) of lemma 6.6.1, 6.5.2 (5).
Statement (2) follows lemma 6.5.6, the identities
$r\0{m-k-1,j-k-1}=r\0{m,j}$, $k+l+1=j$, the fact that $l\leq m-k-3$ if and only if
$j\leq m-2$,
and the definition of $\cY(P,r)$. Finally (3) 
follows from (6.5.1), the definition of $r\0{m,k}$
and the definition of  $\cY(P,r)$. This proves the lemma.
\vspace{.8in}

\noindent {\bf 6.7. Radial stability of the sets $\cY(P,\Delta^k,r)$.}\\
In section 8 we will need a sort of a stable property for the sets
$\cY(P,\Delta^k,r)$, $\cY(P,r)$. We use the objects and notation used
in section 6.6. Recall that $\s{Star}(\Delta,P)$ is the simplicial star of
$\Delta$ in $P$. Also recall that an element in $\rC P$ can be written as
$sx$, $s\in [0,\infty)$, $x\in P$.\\

\noindent {\bf Lemma 6.7.1.} {\it Fix $b\in \R$. Let $\Delta=\Delta^k\in P$,
$x\in\s{Star}(\Delta,P)$ and $a>0$. Then $(s+b)x\in\s{N}(\rC \Delta,\rC P)$
if and only if}
$$
sin(\gamma)\,\frac{sinh(s+b)}{sinh (s)}\,\,\leq\,\,sin(\alpha) 
$$
\noindent {\it where $s>0$, $\gamma=\gamma(x)=d\0{P}(x,\Delta)$ and 
$\alpha=sin^{-1}(\frac{sinh(a)}{sinh(s)})$.}\\

\noindent {\bf Proof.} Note that $\gamma$ is the angle opposite to
the cathetus of length $d=d\0{\rC P}((s+b)x,\rC \Delta)$ of the right
hyperbolic triangle with hypotenuse $(s+b)$. Also 
$\alpha$ is the angle opposite to
the cathetus of length $a$ of the right
hyperbolic triangle with hypotenuse $s$. Using the hyperbolic law of
sines we see that the inequality in the statement of the lemma
is equivalent to $d\leq a$, which, in turn, is equivalent to
$(s+b)x\in\s{N}(\rC \Delta,\rC P)$. This proves the lemma.\\

Write $R(s)=R_{x,b}(s)=(s+b)x$.\\

\noindent {\bf Lemma 6.7.2.} {\it Let $\Delta$, $P$ and $x$ as in lemma 6.7.1.
We have the following three mutually exclusive possibilities}
\begin{enumerate}
\item[{\it (i)}] {\it $e^b sin(\gamma)<sin(\alpha)$, which implies that
$R(s)\in\s{N}(\rC \Delta,\rC P)$, $s\geq s_0$, for some $s_0$.}
\item[{\it (ii)}] {\it $e^b sin(\gamma)>sin(\alpha)$, which implies that
$R(s)\notin\s{N}(\rC \Delta,\rC P)$, for all $s>0$.}
\item[{\it (iii)}] {\it $e^b sin(\gamma)=sin(\alpha)$, which implies that
$R(s)\notin\s{N}(\rC \Delta,\rC P)$, for all $s>0$.}
\end{enumerate}

\noindent {\bf Proof.} The lemma follows from lemma 6.7.1 and the following
two facts: (1) the function $s\mapsto\frac{sinh(s+b)}{sinh(s)}$ is strictly
decreasing for $s>0$, and (2) $\lim_{s\ra \infty}\frac{sinh(s+b)}{sinh(s)}=e^b$.
This proves the lemma.\\

If $R=R_{x,b}$ satisfies (i) we say that $R$ is {\it eventually in}
$\s{N}(\rC \Delta,\rC P)$. If $R=R_{x,b}$ satisfies (ii) we say that $R$ is {\it 
stably disjoint from}
$\s{N}(\rC \Delta,\rC P)$. If $R=R_{x,b}$ satisfies (iii) we say that $R$ is {\it 
unstably disjoint from}
$\s{N}(\rC \Delta,\rC P)$. Note that the first two cases are {\it stable} in the
following sense: it $x'$ and $b'$ are sufficiently close to $x$ and $b$, respectively,
and $R_{x,b}$ is eventually in $\s{N}(\rC \Delta,\rC P)$, then $R_{x',b'}$ is
also eventually in $\s{N}(\rC \Delta,\rC P)$. Similarly for $R$ being stably
disjoint from $\s{N}(\rC \Delta,\rC P)$. The next result is the main purpose of this section.\\

\noindent {\bf Proposition 6.7.3.}  {\it Fix $b\in \R$ and let $x\in P$.
Then at least one of the following conditions hold.}
\begin{enumerate}
\item[{\it (1)}] {\it There is $\Delta^k\in P$, $k\leq m-2$, such that
$R_{x,b}(r_{m-2})\in\cY(P,\Delta^k,r(r_{m-2}))$, for all $r_{m-2}>r'$, for some 
$r'$.}
\item[{\it (2)}]  {\it We have that
$R_{x,b}(r_{m-2})\in\cY(P,r(r_{m-2}))$, for all $r_{m-2}>r'$, for some 
$r'$.}
\end{enumerate}

\noindent {\it Moreover, these two conditions are stable in the following sense.
If $x'$ and $b'$ are sufficiently close to $x$ and $b$, respectively,
and $R_{x,b}$ satisfies (i) then $R_{x',b'}$ also satisfies (i) (with the same
$r'$). Similarly for condition (ii).}\\

\noindent {\bf Proof.} By induction. Suppose $R=R_{x,b}$ is eventually in
$\s{N}\0{r\0{m,0}}(\Delta^0,P)$, for some $\Delta^0$.
Then, since $\s{N}\0{r\0{m,0}}(\Delta^0,P)\sbs\cY(P,\Delta^0,r)$ we see that $R$ satisfies (1) for $\cY(P,\Delta^0,r)$ and we are done.
Suppose $R$ is unstably disjoint from $\s{N}\0{r\0{m,0}}(\Delta^0,P)$, for some $\Delta^0$. Then $x\in\s{Star}(\Delta^0,P)$ and, by 6.7.2, $e^b sin(\gamma)=sin(\alpha_0)$, where $\gamma=\gamma(x)$ is as above.
Then $e^b sin(\gamma)<sin(\beta_0)$, hence, by 6.7.1 and 6.7.2
$R$ is eventually in $\stackrel{\circ}{\s{N}}\0{s\0{m,0}}(\Delta^0,P)$, and also follows 
that $R$ satisfies (1) for $\cY(P,\Delta^0,r)$ and we are done.
Note that these two cases are stable.\\

Now suppose that $R$ is stably disjoint from all $\s{N}\0{r\0{m,0}}(\Delta^0,P)$.
As before we have three possibilities. First if
$R=R_{x,b}$ is eventually in
$\s{N}\0{r\0{m,1}}(\Delta^1,P)$, for some $\Delta^1$.
Then, using the assumption that $R$ is stably disjoint from all
$\s{N}\0{r\0{m,0}}(\Delta^0,P)$ together with the definition of
$\cY(P,\Delta^1,r)$ we see that $R$ satisfies (1) for $\cY(P,\Delta^1,r)$ and we are done.
Suppose $R$ is unstably disjoint from $\s{N}\0{r\0{m,1}}(\Delta^1,P)$, for some $\Delta^1$. Using the same argument as in the $\Delta^0$ case 
(when we assumed $R$ unstably disjoint from some $\s{N}\0{r\0{m,0}}(\Delta^0,P)$)
we get that
$R$ is eventually in $\stackrel{\circ}{\s{N}}\0{s\0{m,1}}(\Delta^1,P)$, and it also follows 
that $R$ satisfies (1) for $\cY(P,\Delta^1,r)$ and we are done.
The third case is that $R$ is stably
disjoint from all  $\s{N}\0{r\0{m,k}}(\Delta^k,P)$, $k\leq 1$.
Again, note that these three cases are stable.
Proceeding in this way we obtain that
either $R$ satisfies (1), for some $\Delta^k$, $k\leq m-2$ or $R$ is stably
disjoint from all  $\s{N}\0{r\0{m,0}}(\Delta^k,P)$, $k\leq m-2$. Hence (2)
holds for $R$. Moreover it does so stably. This proves the proposition.

\vspace{1.2in}

\noindent {\bf \large  Section 7. Link Smoothings and Smooth Cubifications.}\\

For the basic definitions and results about cubic and all-right spherical complexes we refer to \cite{ChD} (see also section 6).
Given a (cube or all-right spherical) complex $K$ we use the same notation $K$ for the complex itself (the collection of all
closed cubes or simplices) and its realization (the union of all cubes or simplices).
We assume that all complexes here satisfy the following condition: any two closed cubes or
simplices intersect
on at most one (possibly empty) common subcube or subsimplex.\\

Throughout this section $M^n$ will denote a smooth manifold of dimension $n$, unless otherwise stated. Also
$K$ will denote a cubical or all-right complex and $f:K\ra M$ a non-degenerate
$PD$ homeomorphism \cite{MunkresLectures}, that is, for all $\sigma\in K$  we have $f|_{\sigma}$ is a smooth embedding.
Recall that in the cubical case the pair $(K,f)$ is called
a {\it smooth cubification} of  $M$. Sometimes we will write $K$ instead
of $(K,f)$. For $\sigma\in K$ we denote its interior by $\dsigma$.\\

Every cube complex has a natural $PL$ structure
given by the simplicial structure obtained by subdividing the cubes properly. It follows then that for a cubification
$(K,f)$ of $M$, this natural $PL$ structure on $K$ induces a $PL$ structure on $M$ that is Whitehead compatible with $M$.\\

In this section 
we consider $\sL(\sigma^{j},K)$,  the (geometric) link of an open $j$-cube or $j$-all-right simplex $\sigma^j$,
as the union of the end points of straight segments of small length $\epsilon>0$
emanating perpendicularly (to $\dsigma^j$) from some point  $x\in \dsigma^j$
(see 6.1 and remark 6.1.1). 
We say that the link is {\it based at $x$.} And the star $\s{Star}(\sigma,K)$ as the union of such segments. We can identify the star with the cone of the link
$\rC\sL (\sigma,K)$  (or $\epsilon$-cone) defined as
$$\rC \, \sL(\sigma, K)=\sL(\sigma,K)\times [0,\epsilon)\,/\, \sL(\sigma,K)\times\{0\}$$ 
\noindent Thus a point $x$ in $\rC \, \sL(\sigma, K)$, different from the {\it cone point} $o\0{\rC \, \sL(\sigma, K)}$, can be written as
$x=t\,u$, $t\in (0,\epsilon)$, $u\in \sL(\sigma, K)$. For $s>0$ we get a
{\it the cone homothety} $x\mapsto sx=(st)u$ (partially defined if $s>1$). 
If we want to make explicit the dependence of the link or the cone on $\epsilon$ we shall write $\sL_\epsilon(\sigma,K)$ or $\rC _\epsilon\,\sL(\sigma,K)$ respectively. 
Also, we will always take $\epsilon <1/2$ ($<\pi/4$ in the spherical case) and it can be verified
that all results in this section (unless otherwise stated) are independent
of the choice of the $\epsilon$'s. 
As usual we shall identify the $\epsilon$-neighborhood
of $\dsigma$ in $K$ with $\rC _\epsilon\,\sL(\sigma,K)\times
\dsigma$ (or just $\rC \,\sL(\sigma,K)\times
\dsigma$). Hence a cone homothety induces a {\it neighborhood homothety}
obtained by crossing it with the identity $1\0{\dsigma}$.\\

 Recall that the link $\sL(\sigma^i,K)$,  $\sigma^i\in K$, has a natural
all-right piecewise spherical structure, which induces a simplicial structure and thus a $PL$ structure on  $\sL(\sigma^i,K)$.
Since the $PL$ structure on $M$ induced by $K$ is Whitehead compatible with $M$ we have that
the link $\sL(\sigma^i,K)$ is $PL$ homeomorphic to $\bS^{n-i-1}$.
A {\it link smoothing for}  $\dsigma^i$ (or $\sigma^i$)
is just a homeomorphism $h_{\sigma^i}:\bS^{n-i-1}\ra\sL(\sigma^i,K)$.
The {\it cone} of $h_{\sigma^i}$ is the map 
$$\rC \,h_{\sigma^i}:\D^{n-i}\longrightarrow \rC\sL (\sigma^i, K)$$
\noindent given by $t\,x= [x,t]\mapsto t\,h_{q^i}(x)=[h_{q^i}(x),  \,t]$,
 where we are canonically identifying the $\epsilon$-cone  of $\bS^{n-i-1}$ with the disc $\D^{n-i}$.\\

A link smoothing $h_{\sigma^i}$  induces the following {\it neighborhood smoothing for}
$\dsigma^i$:
$$h^\bullet_{\sigma^i}=f\,\, \circ\,\,\Big(\rC \,h_{\sigma^i}\times 1_{\dsigma^i}\Big):\D^{n-i}\times \dsigma^i\longrightarrow M $$

Here is the main result of section 7.\\

\noindent {\bf Theorem 7.1.} {\it Let $M^n$ be a smooth manifold with smooth structure $\cS$. Let $K$ be a cubical or all-right complex and  $f:K\ra M$ a non-degenerate PD homeomorphism.
Then there is a smooth structure $\cS'$
and link smoothings $ h_{\sigma^i}$, for all $\sigma^i\in K$, such that all induced
neighborhood smoothings $h^\bullet_{\sigma^i}:\D^{n-i}\times 
q^i\longrightarrow \big(M,\cS'\big)$ are smooth embeddings. 
The smooth structure $\cS'$ is diffeomorphic to $\cS$.}\\

\noindent {\bf Remarks.}

\noindent {\bf 1.} Note that $\cA=\big\{\,(\,h^\bullet_{\sigma^i}\, ,\,\D^{n-i}\times 
\dsigma^i\,)\,\big\}_{\sigma^i\in K}$ is a differentiable atlas for $(M,\cS')$.
Sometimes will just write $\cA=\big\{\,h^\bullet_{\sigma^i}\,\big\}_{\sigma^i\in K}$.
Also  the maps $f|_{\dsigma^i}:\dsigma^i\ra (M,\cS')$ are embeddings. 
The topological atlas $\cA$ depends uniquely on the
the complex $K$, the map $f$ and the collection of link smoothings 
$\{h_{\sigma}\}_{\sigma\in K}$. 
To express the dependence of the atlas on the set of links smoothings
we shall write $\cA=\cA\big(\{h_{\sigma}\}_{\sigma\in K}\,\big)$
(this is different from $\cA=\big\{\,h^\bullet_{\sigma^i}\,\big\}_{\sigma^i\in K}$,
as given above).
Note that not every collection of link smoothings induce a smooth atlas, but the atlas $\cA$  given by Theorem 7.1 is smooth.
We call $\cA$ a {\it normal atlas for } $K$ and the corresponding smooth structure $\cS'$
a {\it normal smooth structure} on $M$ for $K$. 
In section 7.7 we  show that the
atlas $\cA\big(\{h_{\sigma_i}\}\big)$ is smooth if and only if
the set of smoothings $\{h_{\sigma_i}\}$ is ``smoothly compatible".
The most important feature about these normal
atlases is that they preserve the radial and sphere (link) structure given by $K$.

\noindent {\bf 2.} 
It can be checked from th proof that in fact $\cS'$ is isotopic to $\cS$ for 
$n\neq 4$. That is, there is a diffeomorphism $\phi : (M,\cS)\ra (M,\cS')$ that is (topologically) isotopic to the identity map $1_M$. With some care we can probably include
the case $n=4$.

\noindent {\bf 3.} Note that the image of the chart $h^\bullet_{\sigma}$
is the open normal neighborhood $\stackrel{\circ}{\s{N}}_\epsilon(\dsigma,K)$
of width $\epsilon$ of $\dsigma$ in $K$.
Even though we are assuming, for simplicity, that $\epsilon<1/2$ ($\epsilon<\pi/4$
in the spherical case) it can be checked from the proof of 7.1 that we can
actually take $\epsilon=1$ ($\epsilon=\pi/2$) for the charts in 7.1. This is not
an essential fact but it will simplify an argument given in section 8.

\noindent {\bf 4.} We can not expect the maps $f|_{\sigma^i}$ to be embeddings, nor the map $f:K\ra (M,\cS')$ to be PD.
This does not even happen in the next simple example.\\

\noindent{\bf Example.} Consider $\R^2$ with its canonical cube structure, that is, the vertices of the cubes are points in $\R^2$
with integer coordinates. Consider now $\R^2$ with a cube structure given by pulling back the canonical one using
the map $(r,\theta)\mapsto (r,\ell \,\theta)$ where we are using polar coordinates, and $\ell\neq \pm 1$ is an integer.
Note that this map is not differentiable at the origin, but all partial derivatives exist.
Let $K_\ell$ be this new cube complex and $f:K_\ell\ra \R^2$ be the inclusion.
Also note that the differentiable structure
$\cS'$ mentioned in Theorem 7.1 for $K_\ell$ coincides with the canonical structure, provided we choose the link
smoothing at the origin in the obvious way.
Then for any closed 2-cube $\sigma^2$ containing the origin we have that  $f|_{\sigma^2}:\sigma^2\ra \R^2$ is not differentiable at 
the origin, hence $f$ is not $PD$.\\

The neighborhood smoothings $h^\bullet_{\sigma^i}$ in Theorem 7.1 preserve radial structure by definition, and this fact is
needed for our geometric constructions in other sections. But the example above shows that the $PL$ structure given by $K$ has to be sacrificed.\\

\noindent {\bf Proof of Theorem 7.1.}\\
 
To simplify our notation we will write the proof for cubical complexes $K$; the proof for the all-right complex case 
is similar. We will denote an open $i$-cube by $q^i$ and the corresponding closed cube
by $\bq^i$. We write $\sL(q,K)=\sL(\bq,K)$.
Also we denote by $K_j$ the $j$-skeleton of $K$, that is the union
of all cubes of dimension $\leq j$. Also write $M_j=f(K_j)$. We will prove the Theorem by induction on $k=n-j$. Consider the following
statement:\\

\begin{enumerate}
\item[{\bf S($k$)}]
{\it There is a differentiable structure $\cS_{k}$ on $M-M_j$ and link smoothings for all $i$-cubes, $i> j $, such that
all induced
neighborhood smoothings $h^\bullet_{q^i}:\D^{n-i}\times 
q^i\longrightarrow \big(M-M_j,\cS_{k}\big)$ are smooth embeddings}
\end{enumerate}

\noindent {\bf Remark.} If {\bf S($k$)} holds, then $\big\{\,(\,h^\bullet_{q^i}\, ,\,\D^{n-i}\times 
q^i\,)\,\big\}_{\bq^i\in K,\, i>j}$ is a differentiable atlas for $(M-M_j,\cS_k)$.
Moreover, the smooth structure $\cS_{k}|_{M-M_{j+1}}$ coincides with $\cS_{k-1}$, provided the link smoothings of both structures for all the $q^i$
coincide, $i>j+1$.
Also  the maps $f|_{q^i}:q^i\ra (M-M_j,\cS_{k})$ are embeddings, $i>k$.
\\

\noindent {\bf Lemma 7.2.} {\it Assume {\bf S($k$)} holds. Then  $f\big(\sL(q^{j},K))$ is a smooth submanifold of $(M-M_j,\cS_{k})$,
for every $j$-cube $q^{j}$ of $K$.}\\

\noindent {\bf Addendum to Lemma 7.2.} {\it Let $h_{q^j}:\bS^{k-1}\ra \sL(q^j,K)$ be a link smoothing such that $f\circ h_{q^j}$
is a diffeomorphism. Then the corresponding neighborhood smoothing
$$h^\bullet_{q^j}|_{\D^k\times q^j-q^j}:\D^k\times q^j-q^j\ra \big(M-M_j,\cS_k\big)$$ 
\noindent is a smooth embedding.}\\


\noindent {\bf Proof of Lemma 7.2.} By the remark above it is enough to verify that $\sL(q^{j},K)$ is a smooth submanifold in
every chart  $h^\bullet_{q^i}:\D^{n-i}\times 
q^i\longrightarrow \big(M-M_j,\cS_{k}\big)$. We can assume that $\bq^j$ is a subcube of $\bq^i$.
To be specific consider $\sL(q^{j},K)=\sL\0{\epsilon}(q^{j},K)$,
$\sL(q^{i},K)=\sL\0{\delta}(q^{i},K)$ and
assume that $\sL(q^{j},K)$ is based at $x\in q^j$.  Write $S=(h^\bullet_{q^i})^{-1}
\big(f(\sL(q^{j},K))\big)$ and $\bq^i=\bq^l\times \bq^j$. Then $S\sbs \D^{n-i}\times q^l\times\{x\}$. Identify $\bq^l$ with
$\bq^l\times\{ x\}$ and $x$ with a vertex of $\bq^l$.
Since $h^\bullet_{q^i}$ is, by definition, radial in the first coordinate, we have that $S$ is the set of points  $(p,q)\in\D^{n-1}\times q^l$ 
such that the segment $[x,h^\bullet_{q^i}(p,q)]$ has length $\epsilon$ (and is perpendicular to $q^j$). Hence $S$ is the set of points $(p,q)$
that satisfy the equation 

\begin{equation*}\epsilon^2\,=\,\delta^2\,d_{\D^{n-i}}^2(p,0)\,+\,d_{\bq^l}^2(q,x)\tag{1}\end{equation*}

\noindent where $d_{\D^{n-i}}$ and $d_{\bq^l}$ are the euclidean distances on $\D^{n-i}\sbs\R^{n-i}$ and
$\bq^l\sbs\R^l$, respectively. This proves the lemma.\\

Now we prove the addendum. Note that
$ \D^{n-i}\times \bq^i\sbs \R^{n-i}\times\R^i=\R^n$, hence
using linear (thus smooth) parametrized homotheties  $\R^{n}\times (0,1) \ra\R^n$, $(x, s)\mapsto sx$,  we can show that the 
parametrized neighborhood homothety $$f\Big(\rC \sL(q^{j},K)\times\dsigma\Big)\times (0,1)\ra f\Big(\rC \sL(q^{j},K)\times \dsigma\Big),\,\,\,\,\,\,\,f(x,y,s)\mapsto f(sx,y)$$
\noindent is smooth. Moreover, for fixed $s\in (0,1)$ the map above is
an embedding.
To prove the addendum just note that for $x\neq 0\in \D^k$ we have 
$$
h^\bullet_{q^j}(x,y)\,=\,f\,\Big(\,|x|\,h_{q^j}\big(\frac{1}{|x|}x    \big)\,,\,
y   \,\Big)
$$
\noindent Hence $h^\bullet_{q^j}$ is smooth. It is straightforward
to prove that it is a smooth embedding. This proves the addendum.\\

\noindent {\bf Remark.} Equation (1) in the all-right spherical
case can be obtained using the spherical law of cosines. The equation is
$$
\epsilon\,=\, cos\big( \delta\,d_{\D^{n-i}}(p,0)\big)
\,cos\big(\,d_{\bq^l}(q,x) \big)
$$\\

This lemma implies that the smooth structure $\cS_k$ induces a smooth structure on the $PL$ sphere $f\big(\sL(q^{j},K)\big)$, but the
smooth structure and the $PL$ structure are not necessarily compatible (see remark 4 above). If $k\leq 7$, $k\neq 5$
the $(k-1)$-sphere $f\big(\sL(q^{j},K)\big)$,   with the induced smooth structure, is a canonical sphere. The next lemma says
that the same is true for $k=5$.\\

\noindent {\bf Lemma 7.3.} {\it Let $j=n-5$ and consider the 4-sphere $f\big(\sL(q^{j},K)\big)$ with the induced smooth structure
from lemma 7.2. Then  $f\big(\sL(q^{j},K)\big)$ is diffeomorphic to $\bS^4$.}\\

The proof of this lemma is given in Appendix C.\\

We now prove Theorem 7.1.  Recall we are writing $j+k=n$.
We begin with $k=1$ and define $\cS_1=\cS|_{M-M_{n-1}}$.
Then  $\cA_1=\big\{\,(\,h^\bullet_{q^n}\, ,\,\D^{0}\times 
q^n\,)\,\big\}_{\bq^n\in K}$ is a differentiable atlas for $\cS_1$.
Using lemma 7.2 and 7.3 (for the case k=5) we can construct $\cS_k$, $k\leq 7$ inductively:
$\cS_{k+1}$ has an atlas $\cA_{k+1}$ which  is obtained from $\cA_k$ by adding the
charts  $\big\{\,(\,h^\bullet_{q^{j}}\, ,\,\D^{n-j}\times 
q^j\,)\,\big\}_{\bq^j\in K}$, for some smoothings $h_{q^j}$ of the links of the $j$-cubes,
such that $h_{q^j}:\bS^{k-1}\ra f(\sL(q^j,K))$ is a diffeomorphism.
Hence, by construction, we have that $\cA_k= \big\{\,(\,h^\bullet_{q^i}\, ,\,\D^{n-i}\times 
q^i\,)\,\big\}_{\bq^i\in K,\, i>j}$, for $k\leq 7$.\\

This proves the Theorem for $n\leq 3$, since (diffeomorphism classes of) smooth structures are unique for $n\leq 3$. For $n=4$ it only remains to prove that the smooth
structures $\cS'=\cS_5$ and $\cS$ are diffeomorphic, and this is done in appendix C
(see C.2.2).
Hence from now on we assume $n\geq 5$. We will prove by induction the following stronger statement. Let
{\bf S'($k$)}, $k\geq 4$, be the statement obtained from {\bf S($k$)} by adding the extra condition: \\

\begin{enumerate}
\item[{\bf C($4$)}:]
\,\,\, {\it $\cS_{4}$ extends to a smooth structure $\cS_{4}'$ on the whole $M$, 
and $\cS_{4}'$ is diffeomorphic to $\cS$. }
\item[{\bf C($k$)}:]
\,\,\, {\it For $k>4$ we have that $\cS_{k}$ extends to a smooth structure $\cS_{k}'$ on the whole $M$, 
and $\cS_{k}'$ is diffeomorphic to $\cS_{k-1}'$.}
\end{enumerate}

Recall that we have proved {\bf S($k$)} for $k\leq 7$. The next propsoition
shows our strategy to prove Theorem 7.1.\\

\noindent {\bf Proposition 7.4.} {\it We have that}
\begin{enumerate}
\item[{\bf (1)}] {\it statement} {\bf C}(4) {\it holds}. 
\item[{\bf (2)}] {\it Statements} {\bf C}($k$) {\it hold for $4<k\leq 6$}.
\item[{\bf (3)}] {\it Statement} {\bf S'}($k-1$) {\it implies} {\bf S}($k$).
\item[{\bf (4)}] {\it Statements} {\bf C}($k$) and {\bf S}($k+1$)
{\it imply} {\bf C}($k+1$) {\it for $k\geq 6$.}
\end{enumerate}
 
\noindent {\bf Proof of Theorem 7.1 assuming proposition 7.4.}\\
Just note that it follows from (2), (3), (4) and the fact that we have already proved     {\bf S($7$)} that 
statement {\bf S'}($k$) implies {\bf S'}($k+1$), for $k\geq 6$. Thus induction
shows that statement {\bf S'}($n+1$) is true. This proves Theorem 7.1
assuming proposition 7.4.\\

\noindent {\bf Proof of Proposition 7.4.}\\
First we need a definition. A {\it spindle neighborhood map} of
an open cube $q^i$ in $K$ is a topological embedding $\alpha:q^i\times \D^{n-i}\ra K$ such that:  (1) $\alpha(x,0)=x$ and \, (2) the diameters of
the fibers $\alpha (\{x\}\times\D^{n-i})$ (with respect to any compatible metric on $K$) tend uniformly to zero as $x$
tends to $\p q^i$. (This terminology is essentially due to J. Munkres \cite{Munkres}). 
The image of $\alpha$ is a {\it spindle neighborhood} of $q^i$. \\

The proof of  (1) uses some ingredients of the proof of lemma 7.2, and
it is given in the second part of appendix C.\\

We prove (2) next. Assume {\bf C}($k-1$) holds, where $5\leq k\leq 6$.
Choose a spindle neighborhood for each $(j+1)$-cube $q^{j+1}$. We can assume all spindle neighborhoods
to be disjoint.
Since we are constructing $\cS_k$ and $\cA_k$ by induction, the smooth structure $\cS_{k}|_{M-M_{j+1}}$ coincides with $\cS_{k-1}$.
Hence the smooth structures $\cS_k$ and $\cS_{k-1}'$ coincide on $M-M_{j+1}$.
Let $q^{j+1}$ be a $(j+1)$-cube and pull back these two smooth structures to $ \D^{k-1}\times
q^{j+1}$ via the corresponding
spindle map $\alpha_{q^{j+1}}$. Call these structures $\cB$ and $\cC$. They coincide outside $q^{j+1}=\{0\}\times
q^{j+1}$. Since $\pi_{k-1}(TOP/O)=0$, 
for $4\leq k-1 <7$, we can apply the Classification Theorem 10.1 of \cite{KiSi}, p.194
(see also Theorem 4.1 in p.25 of \cite{KiSi}) and obtain a diffeomorphism 
$\beta:( \D^{k-1}\times q^{j+1}, \cB)\ra ( \D^{k-1}\times q^{j+1},\cC)$, which is the identity outside a pre-chosen small spindle
neighborhood of $q^{j+1}$ in $ \D^{k-1}\times q^{j+1}$. Moreover $\beta$ is isotopic to the identity modulo this neighborhood.
Using these maps $\beta$ we can define a diffeomorphism
$\phi: (M-M_j,\cS_k)\ra (M-M_j,\cS_{k-1}')$ and $\phi$ extends to a homeomorphism $\phi : M\ra M$
by declaring $\phi(x)=x$, for $x\in M_j$. Take now $\cS_k'$ to be the pull back of $\cS_{k-1}'$ by $\phi$.
This proves (2).\\

To prove (3) we need the following result. Recall we are assuming
$n\geq 5$.\\

\noindent {\bf Lemma 7.5.} {\it Assume {\bf S'($k$)} holds.
Consider the $(k-1)$-sphere $f(\sL(q^{j},K))$ with the induced smooth structure
from lemma 7.2. Then  $f(\sL(q^{j},K))$ is diffeomorphic to $\bS^{k-1}$.}\\

\noindent {\bf Proof.} 
 Let $h_{q^j}:\bS^{k-1}\ra \sL(q^j,K)$
be a homeomorphism. 
Pull back the differentiable structure  $\cS'_k$ by the map $h^\bullet_{q^j}$ to
obtain a differentiable structure $\cB$ on $\D^{n-j}\times q^j$. By the addendum to
Lemma 7.2 the differential structure $\cB$ is a product structure outside 
$q^j=\{ 0\}\times q^j\sbs \D^k\times q^j$. By lemma 7.3 we can assume 
$k-1> 6$ and we can apply the Product Structure Theorem (see Theorem 5.1 in p.31 of \cite{KiSi}) to obtain a diffeomorphism
$\big( \D^k\times q^j, \cC\times \cS_{\R^j} \big)\ra\big(\D^k\times q^j,\cB \big)$
for some differentiable structure $\cC$ on the disc $\D^k$. Furthermore, we can assume this
diffeomorphism to be the identity outside a small neighborhood of $q^j$.
Hence the differentiable structure $\cB|_{\bS^{k-1}}$ extends to a differentiable structure
on the whole disc. We can apply now the smooth $h$-cobordism Theorem to conclude that
$\cB|_{\bS^{k-1}}$ is diffeomorphic to the canonical structure on the sphere. This
proves the lemma.\\

We can prove (3) now:
just use lemma 7.5 and the addendum to lemma 7.2 to construct $\cA_k$ and $\cS_k$ as done before. 
Therefore {\bf S'($k-1$)} implies {\bf S($k$)}. \\

It remains to prove (4), that is we have to prove that  {\bf C($k-1$)} and {\bf S($k$)}
imply  {\bf C($k$)}, i.e. that $\cS_k$ extends to the whole $M$. For this we
may have to change the smoothings of the links of $j$-cubes given in the proof of lemma 7.5.
The following lemma will indicate how to change the smoothing $h_{q^j}$.
Recall that for $f:\bS^k\ra\bS^k$ the map $\rC f:\D^{k+1}\ra\D^{k+1}$
is the cone of $f$ defined by $\rC f (x)=|x|f(\frac{x}{|x|})$, $x\neq 0$ and
$f(0)=0$.\\

\noindent {\bf Lemma 7.6.} {\it Let $\cD$ be a smooth structure on the disc $\D^k$, $k \geq 6$.
Then there is a diffeomorphism $g:\bS^{k-1}\ra\bS^{k-1}$ such that the smooth structure
$(\rC \,g)^*\cD$ is diffeomorphic to the canonical one, by a diffeomorphism that is the identity outside a small neighborhood
of the origin.}\\

\noindent {\bf Proof.} Let $G:\D^k\ra(\D^k,\cD)$ be a diffeomorphism that is radial outside
a small neighborhood of the origin. Take $g=G|_{\bS^{k-1}}$
and note that $((\rC g)^{-1}\circ G)^*(\rC g)^*\cD$ is the canonical smooth
structure on $\D^k$.  Moreover $(\rC g)^{-1}\circ G$ is the identity outside
a small neighborhood of the origin.
This proves lemma 7.6.\\

We now prove (4). Assume {\bf C($k$)} and {\bf S($k+1$)} hold.
Let $\cB=\cB(h_{q^j})$ be the pull back of the smooth structure $\cS_{k}'$ by $h^\bullet_{q^j}$.
Then $\cB$ is a smooth structure on $\D^k\times q^j$ that coincides with the canonical one
outside $q^j$. For a self-homeomorphism $g$ on $\bS^{k-1}$ 
write $g^\bullet=(\rC \, g)\times 1_{q^j}$, which is a self-homeomorphism 
on $\D^{k-1}\times q^j$.\\

\noindent {\bf Claim.} {\it There is a self-homeomorphism $g$ on $\bS^{k-1}$
such that $\cB( h_{q^j}\circ g)=(g^\bullet)^*\cB$ is diffeomorphic to the canonical
smooth structure on $\D^k\times q^j$, via a diffeomorphism that is the
identity outside a small spindle neighborhood of $q^j$ in $\D^k\times q^j$.}\\

\noindent {\bf Proof of claim.}
Using the Product
Structure Theorem (see Theorem 5.1 in \cite{KiSi}, p.31), we get a diffeomorphism
$$\gamma: \big( \D^k\times q^j ,\cD\times\cS_{\R^j}\big)\ra \big( \D^k\times q^j,\cB \big)$$
\noindent for some smooth structure $\cD$ on the disc $\D^k$. Furthermore,
$\gamma$ is the identity outside a small spindle neighborhood of $q^j$ and isotopic to the identity (by an isotopy that is constant outside a small spindle neighborhood of $q^j$). Apply lemma 7.6 to $\cD$ to obtain a diffeomorphism $g$. The smooth structures $(\gamma\circ g^\bullet)^*\cB$ and $(g^\bullet)^*\cB$  
are diffeomorphic via a diffeomorphism that is the identity on the complement of a small neighborhood of $q^j$ (the diffeomorphism is $(g^\bullet)^{-1}\circ\gamma\circ g^\bullet$).
But by lemma 7.6 we have that $(\gamma\circ g^\bullet)^*\cB=(g^\bullet)^*(\cD\times \cS_{\R^j})=(g^*\cD)\times\cS_{\R^j}$
is diffeomorphic to the canonical smooth structure on $\D^k\times q^j$,
via a diffeomorphism that is the identity on the complement of a small neighborhood of $q^j$. Moreover, by contracting
the fibers $\D^k\times\{ x\}$, as $x$ approaches $\p q^j$ we can  assume that 
that the aforementioned diffeomorphism is the identity outside a small spindle neighborhood of $q^j$.
Hence the same is true for $(g^\bullet)^*\cB$. This proves the claim.\\

From the claim it follows that
if we replace $h_{q^j}$ by $ h_{q^j}\circ g$ and take $\cB=\cB(h_{q^j}\circ g)$
instead of $\cB=\cB(h_{q^j})$
we get that the new structure $\cB$ is now diffeomorphic to the canonical one
 modulo the complement of a small spindle neighborhood of $q^j$. That is
there is a diffeomorphism $\beta=\beta_{q^j}:\D^k\times q^j\ra (\D^k\times q^j,\cB)$ that
is the identity outside a small spindle neighborhood of $q^j$.
We can now 
proceed as in the case $k\leq 7$ and use the maps $\beta_{q^j}$ to obtain
a diffeomorphism  
$\phi: (M-M_j,\cS_k)\ra (M-M_j,\cS_{k-1}')$ and $\phi$ extends to a homeomorphism $\phi : M\ra M$
by declaring $\phi(x)=x$, for $x\in M_j$. Take now $\cS_k'$ to be the pull back of $\cS_{k-1}'$ by $\phi$. This proves (4) and completes the proof of proposition 7.4
\vspace{.8in}

\noindent {\bf 7.7. Induced Link Smoothings.}\\
Let $K$ be a cubical or all-right spherical complex. Then the links of $\sigma\in K$ are
all-right spherical complexes. Furthermore the all-right spherical structure on 
$\sL(\sigma,K)$ induced by $K$ has all-right spherical simplices
$$\big\{ \,\tau\,\,\cap\,\, \sL(\sigma,K)\,\,  ,\,\, \tau\in K\big\} $$ 
\noindent Note that $\tau\,\,\cap\,\, \sL(\sigma,K)$ is non-empty only when $\sigma\sbs \tau$. We have the following identifications

$$\sL\,\Big(\, \tau\,\cap\,\sL\,(\sigma,K)  , \sL\,(\sigma,K)  \,\Big)\,=\,\sL\,\big( \, \tau   , K  \,\big)
$$

\noindent and

$$\rC\sL \,\Big(\, \tau\,\cap\,\sL\,(\sigma,K)  , \sL\,(\sigma,K)  \,\Big)\,=\,\rC\sL\,\big( \, \tau   , K  \,\big)
$$\\

\noindent In particular we can write $\rC\sL\,\big( \, \tau   , K  \,\big)\sbs
 \sL\,(\sigma,K)$. \\

\noindent {\bf Remark 7.7.1.} The identifications above can be done in many ways.
In 6.4.2 we proved the first equation for all-right spherical complexes.
The proof for cubical complexes is the same. But in section 7.7 and in section 9
we will need a different identification, one using radial projections. 
For example the inclusion  $\rC\sL\,\big( \, \tau   , K  \,\big)\sbs
 \sL\,(\sigma,K)$ in the cubical case can be done in the following way.
Let $v\in\rC\sL\,\big( \, \tau   , K  \,\big)=\rC\sL_r\,\big( \, \tau   , K  \,\big)$
at some $p\in \dot{\tau}$ and consider
$ \sL\,(\sigma,K)= \sL_s\,(\sigma,K)$ at some point $q\in \dsigma$, with
$d_K(p,q)=s$ with the segment $[p,q]$ perpendicular to $\sigma$.
Then $v$ corresponds to the point $v'$ in the segment $[q,v]$ at a distance
$s$ from $q$. The (angular) distance in $\sL(\sigma,K)$ from $p$ to $v'$ is
$tan^{-1}(\frac{d_K(v,p)}{s})$.
If we need to specify the type of identification we are using we shall
write $\Re:\rC\sL(\tau,K)\hookrightarrow \sL(\sigma,K)$ for the radial
projection described above. We will write $\Re=\Re\0{p,q,r,s}$ if we want to make explicit the
dependence of $\Re$ on the choices above.\\

\noindent {\bf Remark 7.7.2.} Note that the representation of the radial projection $\Re$ in the chart
$\big(\, h_{\tau}^\bullet, \D^{n-j}\times\dot{\tau}\,\big)$, $\tau$ a $j$-simplex, is smooth.\\

Since we can write $\sL\,\Big(\, \tau\,\cap\,\sL\,(\sigma,K)  , \sL\,(\sigma,K)  \,\Big)\,=\,\sL\,\big( \, \tau   , K  \,\big)$
we can say that a set of link smoothings $\{ h_\sigma\}_{\sigma\in K}$ for
(the links of) $K$
induces, just by restriction, a collection of link smoothings for $\sL(\sigma, K)$,
$\sigma\in K$, given by
$\{ h_\tau\}_{\sigma\lneqq\tau}$. Note that the atlas  $\cA_\sigma=\cA_{\sL(\sigma,K)}=\big\{h^\bullet_\tau\big\}_{\sigma\lneqq \tau}$ 
is a (a priori just topological) normal atlas on $\sL(\sigma,K)$.\\

 Let $\cA=\cA\big(\{h_\sigma\}_{\sigma\in K}\big)$ be as in Theorem 7.1
(and remark 1 after 7.1). This atlas induces the smooth structure $\cS'$. 
Also let $\cS'_\sigma=\cS'_{\sL(\sigma,K)}$
be the smooth structure on $f(\sL(\sigma,K))$ constructed in the proof of lemma 7.2.
Hence the inclusion $\big(\,f(\sL(\sigma,K)),\cS'_\sigma\,\big)\ra\big(\, M,\cS'\,\big)$ is an embedding.
Then it it straightforward to verify that $\cA_\sigma$ is the atlas constructed in the
proof of lemma 7.2. In particular we obtain   the following result.\\

\noindent {\bf Corollary 7.7.3.} {\it The atlas $\cA_\sigma$ is smooth and it induces 
the structure $\cS'_\sigma$.}\\

And, since a link smoothing $h_\sigma$ is the restriction of the embedding $h^\bullet_\sigma$ we get the following corollary.\\ 

\noindent {\bf Corollary 7.7.4.} {\it For every $\sigma^k\in K^n$, the link smoothing $h_{\sigma}:\bS^{n-k-1}\ra \big(\sL(\sigma,K),\cS'_\sigma\big)$ is a diffeomorphism.}

\vspace{.8in}

\noindent {\bf 7.8. Change of charts and smooth compatibility.}\\
Let $K$ be a cube complex or an all-right spherical complex, $M$ a manifold and $f:K\ra (M,\cS)$
a non-degenerate $PD$ homeomorphism. Theorem 7.1 states that there is a (diffeomorphic to $\cS$) differentiable structure $\cS'$
on $M$ so that we get nice normal charts that respect the radial and link structure of $K$ that we called a normal differentiable structure on $M$ for $K$.
But with this new structure $\cS'$ the maps $f|_\sigma$, $\sigma\in K$, are not necessarily embeddings (here $\sigma$ is a {\sf closed} cube or simplex).
Furthermore the natural $PL$ triangulation of the links $\sL(\sigma,K)$ is not necessarily smooth (with respect to $\cS'_\sigma$, see 7.6). We need some technical results that will correct, at least in part, the lack
of strong smoothability of the charts on closed cubes (or simplices). 
We will use the results of this section in section 8.2.\\

Let $U$ be a bounded open set of $\R^n$. A smooth map $f:U\ra\R^N$
is {\it polynomially bounded with respect to a subset } $B\sbs \R^n$ if for every
$p\in B$ and every partial derivative $\p^\alpha f$ of $f$
we have constants $C,\, m$ such that $|\p^\alpha f(p)|\leq C d\0{\R^n}(p, B)^m$.
The map is {\it $C^1$ well-bounded} if the norm of the first derivative $|Df|$ is bounded
and bounded away from zero.
For $U\sbs\bS^n$ we write $\rC^+U=\rC U-\{0\}$.\\

\noindent {\bf Lemma 7.8.1.} {\it Let $U\sbs\bS^n$ be open and $f:U\ra f(U)\sbs\R^N$ be smooth. Let $V$ open with $\bar{V}\sbs U$. Then 
$\rC \big(f|_{V}\big):\rC^+V\ra \rC f(V)\sbs \R^N$ is polynomially bounded at 0.
Furthermore $\rC \big(f|_{V}\big)$ is also $C^1$ well-bounded.}\\

\noindent {\bf Proof.} We have $(\rC f) (x)=|x|f(\frac{x}{|x|})$, and the lemma follows from differentiating this equation.
This proves the lemma.\\

\noindent {\bf Remark.} Actually since $\rC f$ is a cone map then $D(\rC f)_x=D(\rC f)\0{\frac{x}{|X|}}$, for $x\neq 0$.
In particular
if $u\in V$, and $t> 0$ then
$D(\rC f)_{tu}u=f(u)$.\\

Let  $\cA=\Big\{ \big(  h_{\sigma^i}^\bullet,\D^{n-i}\times\dsigma^i    \big)\Big\}$
be a normal atlas as in 7.1. Let $S=f(\sL(\sigma^k,K))$ and $U'\sbs S\cap W_{\sigma^i}$,
where $\sigma^i>\sigma^k$ and $W_{\sigma^i}$ is the image of 
$h_{\sigma^i}^\bullet$. Write $U=\big(h^\bullet_{\sigma^k}\big)^{-1}(U')$
and let $V$ open with $\bar{V}\sbs U$. We have the following corollary about
the change of charts $\big(h^\bullet_{\sigma^i}\big)^{-1}\circ h_{\sigma^k}^\bullet:
\rC^+V\ra \D^{n-i}\times\sigma^i\sbs\R^n$, restricted to the cone of $V$.\\

\noindent {\bf Corollary 7.8.2.} {\it The change of charts
$\big(h^\bullet_{\sigma^i}\big)^{-1}\circ h_{\sigma^k}^\bullet$,
restricted to $\rC^+ V$, is polynomially bounded at 0, and
$C^1$ well-bounded.}\\\\

Let $K$ be a cube complex and $f:K\ra M$ a homeomorphism.
Since in what follows of this subsection the function $f$ is not essential,
to simplify our notation we identify $K$ and $M$ via $f$. Let 
$\big\{h_{\sigma}\big\}$ be a set of link smoothings on the topological 
$n$-manifold $M$, for $K$. Recall that this set determines a 
(not necessarily smooth) atlas $\cA=\Big\{h^\bullet_{\sigma}\Big\}$ on  $M$.
As above we denote by $W_{\sigma}$ the image of the chart $h^\bullet_\sigma$.\\

Let $\sigma^k\in K$ and write $S=\sL(\sigma^k,K)$.
Recall that $K$ induces an all-right triangulation on $S$ with simplices
$\sigma^i\0{S}=\sigma^i\cap S$, and open simplices
$\dsigma^i\0{S}=\dsigma^i\cap S$.
Notice that if $\sL(\sigma^k,K)=\sL_s(\sigma^k,K)$ is based at $p\in \sigma^k$, then
$\sigma\0{S}^i$ is the set of points $x$ in $\sigma^i$ with $d\0{\sigma^i}(p,x)=s$ and $[p,x]$
perpendicular to $\sigma^k$ at $p$.\\

In section 7.7 we mentioned that for $\sigma^k\sbs\sigma^j$, using
radial projection, we  can consider $\rC\sL(\sigma^j,K)$ as a subset of $\sL(\sigma^k,K)$
and we write $\Re:\rC\sL (\sigma^j,K)\hookrightarrow \sL(\sigma^k,K)$ (see remark 7.7.1).
We say that the set of smoothings is {\it smoothly compatible} if
for every $\sigma^k<\sigma^j$ the composition maps\\

\noindent {\bf (7.8.3.)}\,\,\,\,
{\small $\D^{n-j} \times \dsigma\0{S}^j\,\,\stackrel{{\mbox{{\tiny inclusion}}}}{\longrightarrow}\,\D^{n-j} \times \dsigma^j
\,\stackrel{{\mbox{{\tiny $ h^\bullet_{\sigma^j}$\,}}}}{\longrightarrow}\,\rC\sL\big(\sigma^j,K\big)\times\dsigma^j\,
\stackrel{{\mbox{{\tiny $\Re$}}}}{\longrightarrow}\, \sL(\sigma^k,K)\,\stackrel{{\mbox{\tiny $(h_{\sigma^k}^\bullet )^{-1}$}}}{\longrightarrow}\,\bS^{n-k-1}$}\\

\noindent are smooth embeddings. Here we are denoting (with a bit of abuse of notation) by $\Re$ the map given by
$\Re(v,p)=\Re\0{p,q,r,s}(v)$, where $q,r,s$ are fixed. Note that the map given in (7.8.3) can also be written as\\

\noindent {\bf (7.8.4.)}\hspace{.8in}
{\small $\D^{n-j} \times \dsigma\0{S}^j\,\,\stackrel{{\mbox{{\tiny $\Re'$}}}}{\longrightarrow}\, S'
\,\stackrel{{\mbox{{\tiny $ h^\bullet_{\sigma^j}$\,}}}}{\longrightarrow}\,S= \sL(\sigma^k,K)\,\stackrel{{\mbox{\tiny $(h_{\sigma^k}^\bullet )^{-1}$}}}{\longrightarrow}\,\bS^{n-k-1}$}\\

\noindent where $S'=(h^\bullet_{\sigma^j})^{-1}(S)$ and $\Re'$ is the representation of $\Re$ in the chart
corresponding to $\sigma^j$ (which is always smooth).
Note that the inverse of the map given in 7.8.3 (or 7.8.4) is just the change of charts
$(h_{\sigma^j}^\bullet )^{-1}\circ h_{\sigma^k}^\bullet $ restricted to an open subset of $\bS^{n-k-1}$,
plus the ``straightening map "  $(\Re')^{-1}$.\\

\noindent {\bf Proposition 7.8.5.} {\it The set of smoothings $\{h_\sigma\}$ is smoothly
compatible if and only if the atlas $\cA\big(\{h_\sigma\}\big)$ is smooth.}\\

\noindent {\bf Proof.} If the atlas $\cA$ is smooth then all chart maps $h^\bullet_{\sigma}$ are embeddings
(with respect to the smooth structure generated by $\cA$).
Therefore  the composition given in 7.8.3  above
is smooth. Note the ``identification" $\Re$
is also smooth, for it is smooth in the chart corresponding to $\sigma^j$ (see remark 7.7.2).
Hence $\{h_\sigma\}$ is smoothly compatible.
We prove the converse by induction on the codimension of the skeleta.\\

Write $k+j=n$.
Suppose the set of smoothings $\{h_\sigma\}$ is smoothly compatible.
Denote by $W_{\sigma^i}$ the image of $h_{\sigma^i}^\bullet$.
Assume that we have proved that the atlas  $\cA_{j}=\Big\{ \big(h^\bullet_{\sigma^i},
\D^{n-i}\times\dsigma^i\big)    \Big\}_{k<i}$ is smooth, and we want to prove that
the atlas  $\cA_{j+1}=\Big\{ \big(h^\bullet_{\sigma^i},
\D^{n-i}\times\dsigma^i\big)    \Big\}_{k\leq i}$ is smooth.
Recall that $\cA_{j}$ is a smooth atlas on the complement $M-M_k$ of the $k$-skeleton. The difference between the atlas $\cA_j$ and the atlas $\cA_{j+1}$
are the charts with maps $h^\bullet_{\sigma^k}$, for all link smoothings
$h_{\sigma^k}$ of $k$-cubes $\sigma^k$. We prove the proposition by
proving that the following maps are embeddings
$$h^\bullet_{\sigma^k}|_{(\D^k\times\dsigma^k)-(\{0\}\times\dsigma^k)}\ra
(M-M_j,\cA_j)$$
\noindent Fix $\sigma^k$. The open sets 
$U_{\sigma^i}=S\cap W_{\sigma^i}$, $\sigma^i>\sigma^k$, form an open cover of $S$. 
Note that $U_{\sigma^i}$ is a normal neighborhood of $\sigma\0{S}^i=\sigma^k\cap\sigma^i$
in $S$.
Write $V_{\sigma^i}=(h_{\sigma^k}^\bullet)^{-1}(\rC^+ U_{\sigma^i})$
(here the vertex of the cone $o$ is the center of $S$). Let $u\in S$. 
Take $i$ so that
$u\in U_{\sigma^i}$. \\

\noindent {\it Claim.  The map \,\,\,
$h_{\sigma^k}^\bullet|\0{V_{\sigma^i}\times\dsigma^k}: V_{\sigma^i}\times
\dsigma^k\ra W_{\sigma^i}\sbs (M-M_k,\cA_{j})$ is an embedding. }\\

\noindent {\it Proof of the claim.} Since  $h^\bullet_{\sigma^i}$ is already an embedding
it is enough to prove that the map $h=\big(h^\bullet_{\sigma^i}\big)^{-1}\circ
h_{\sigma^k}^\bullet|\0{V_{\sigma^i}} $ is an embedding. But 
we can consider $V_{\sigma^i}\sbs \D^{n-k}-\{0\}=\bS^{n-k-1}\times(0,1)$. Write $\sigma^i=\sigma^l\times\sigma^k$,
and let $o$ be the vertex in $\sigma^l$ corresponding to $\sigma^k$ by the projection $\sigma^i\ra\sigma^l$.
For $v=
(t\,u,w)\in V_{\sigma^i}\times \dsigma^k$, $u\in\bS^{n-k-1}$, we can write 
$$h(v)=( \alpha_u(t) ,w)\in (\D^{n-i}\times \dsigma^l)\times\dsigma^k=\D^{n-i}\times\dsigma^i$$

\noindent where $\alpha_u$ is the segment $[o,g(u)]$
and $g$ is the inverse of the map in 7.8.4. This proves the claim.\\

Since the open sets $V_{\dsigma^i}\times\dsigma^k$ cover 
$(\D^{n-k}\times\dsigma^k)-(\{0\}\times\dsigma^k)$ we can conclude that $h_{\sigma^k}^\bullet$ is an embedding
away from $\{0\}\times\dsigma^k$, into the smooth manifold $(M-M_k,\cA_{j+1})$.
This proves the proposition.\\

Let $K$ be a smooth cubification of $M$ and $\cA$ be a normal atlas for $K$, inducing the smooth structure $\cS'$ on $M$. We mentioned several times that for a (closed)
cube $\sigma$ the inclusion $\sigma\hookrightarrow (M,\cS')$ is (almost always) not
a smooth embedding. But the following regularity condition can be easily verified.
Consider $\sigma^j=\sigma^l\times\sigma^k\in K$. Identify a normal
neighborhood of $\sigma^k$ in $\sigma^j$ with 
$\rC\sL(\sigma^k,\sigma^j)\times \sigma^k$ and write an element
in $\rC\sL(\sigma^k,\sigma^j)$ in the form $tu$, $u\in\sL(\sigma^k,\sigma^j)$.
We have an inclusion $\sL(\sigma^k,\sigma^j)\times \sigma^k\sbs
\rC\sL(\sigma^k,\sigma^j)\times \sigma^k$.
Also denote the inclusion $\sigma^j\hookrightarrow (M,\cS')$ by $\iota$.\\

\noindent {\bf Lemma 7.8.6.} {\it 
Let $(u\0{0},p\0{0}), (u_n, p_n)\in \sL(\sigma^k,\sigma^j)\times \dsigma^k\sbs
\sigma^j$ with $(u_n,p_n)\ra (u\0{0},p\0{0})$.
Also let $(0, v\0{0}), (0,v_n)\in T_{(u\0{n},p\0{n})}\Big(\rC\sL(\sigma^k,\sigma^j)\times \sigma^k\Big)$ (hence
they are ``parallel to $\sigma^k$), with $v_n\ra v\0{0}$ and $t_n\ra 0\in [0,\infty)$. Then}
\begin{enumerate}
\item[{\it (i)}] {\it  we have that $D\iota\0{(t_nu_n,p_n)}(u_n,0)\ra D\iota\0{(0,p\0{0})}(u\0{0},0)$,}
\item[{\it (ii)}] {\it  we have that $D\iota\0{(t_nu_n,p_n)}(0,v_n)\ra D\iota\0{(0,p\0{0})}(0,v\0{0})$.}
\end{enumerate}

\noindent {\bf Proof.} Just take that chart $h^\bullet_{\sigma^j}$ and recall that
$h^\bullet_{\sigma^j}$ is a product map and respects the radial structure.
See also remark after lemma 7.1. This proves the lemma.

\vspace{.8in}

\noindent {\bf 7.9. Manifolds with codimension zero singularities.}\\
In this section we treat the case of manifolds with a one point singularity.
The case of manifolds with many (isolated) point singularities
is similar. The results in this section will be used in section 11.\\

Let $Q$ be a smooth manifold with a one point singularity $q$, that is
$Q-\{q\}$ is a smooth manifold and there is a topological embedding
$\rC_1 N\ra Q$, with $o\0{\rC N}\mapsto q$, that is a smooth embedding outside 
the vertex $o\0{\rC N}$. Here $N=(N,\cS_N)$ is a closed smooth manifold (with smooth
structure $\cS_N$).
Also $\rC_1 N $ is the (closed ) cone of width 1 and we identify
$\rC_1 N-\{o\0{\rC N}\}$ with $N\times (0,1]$. We write $\rC_1N\sbs Q$.
We say that the {\it singularity $q$ of $Q$ is  modeled on  $\rC N$.}\\

Assume $(K,f )$ is a smooth cubification of $Q$, that is
\begin{enumerate}
\item[(i)] 
$K$ is a cubical complex.
\item[(ii)] $f:K\ra Q$ is a homeomorphism. Write $f(p)=q$ and $L=\sL(p,K)$.
\item[(iii)] $f|_{\sigma}$ is a smooth embedding for every
cube $\sigma$ not containing $p$.
\item[(iv)]  $f|_{\sigma-\{p\}}$ is a smooth embedding for every
cube $\sigma$ containing $p$.
\item[(v)]  There is $g:L\ra N$ such that $(L,g)$ is a smooth triangulation of $N$.
\end{enumerate}

Many of the definitions and results given in section 7 still hold (with minor changes)
in the case of manifolds with a one point singularity:\\

\begin{enumerate}
\item[{\bf (1)}] A {\it link smoothing} for $L=\sL(p,K)$ (or $p$) is just a
homeomorphism $h_p:N\ra L$.
Since all but one of the links of $K$ are spheres, sets of  link smoothings for $K$
are defined, that is they are sets of link smoothings for the sphere links plus
a link smoothing for $L$.
\item[{\bf (2)}] Given a set of link smoothings for $K$ we get a set of charts
as before. For the vertex $p$ we mean the cone map 
$h_p^\bullet=f\circ\rC h_p:\rC N\ra Q$. We will also denote the restriction
of $h_p^\bullet$ to $\rC N-\{o\0{\rC N}\}$ by the same notation $h_p^\bullet$.
As before $\{h^\bullet_\sigma\}_{\sigma\in K}$ is a {\it normal atlas} for $K$
(or $K-\{p\}$). The atlas is smooth if all change of charts are smooth. A smooth normal atlas on $Q$ for $K$ induces, by restriction,  a smooth normal 
structure on $Q-\{q\}$ for $K-\{p\}$ (this makes sense even though
$K-\{p\}$ is not, strictly speaking, a cube complex).
\item[{\bf (3)}] We say that the set $\{h_\sigma\}$ is smoothly compatible if
condition 7.8.3 holds. For the case $\sigma=p$ we are considering
$h_p^\bullet:\rC N-\{o\0{\rC N}\}\ra Q-\{q\}$ and identifying $\rC N-\{
o\0{\rC N}\}$ with $N\times (0,1]$ with the product smooth structure obtained from
\s{some} smooth structure ${\tilde{\cS}}_N$ on $N$. It can be verified that proposition 7.8.5 holds:  $\{h_\sigma\}$ is smoothly compatible if and only if
$\{h^\bullet_\sigma\}$ is a smooth atlas on $K$. In this case we say that
the smooth atlas  $\{h^\bullet_\sigma\}$ (or the induced smooth structure, or the set $\{h_\sigma\}$) 
is {\it correct with respect to $N$} if $\cS_N$ and ${\tilde{\cS}}_N$ are diffeomorphic.
\item[{\bf (4)}] Theorem 7.1 also holds in this context:
\end{enumerate}

\noindent {\bf Theorem 7.9.1.} {\it  Let $Q$ be a smooth manifold
with one point singularity $q$ modeled on $\rC N$, where $N$ is a closed
smooth manifold. Let $(K,f)$ be a smooth cubification of $Q$. Then
$Q$ admits a normal smooth structure for $K$ which restricted to $Q-\{q\}$ is diffeomorphic to $Q-\{q\}$.
Moreover this normal smooth structure is correct with respect to $N$ if}
\begin{enumerate}
\item[(a)] {\it $dim\, N\leq 4$.}
\item[(b)] {\it $dim\, N>5$ and the Whitehead group $Wh(N)$ of $N$ vanishes. }
\end{enumerate}

\noindent {\bf Proof.} The proof of the existence part is the same as the proof of
Theorem 7.1. Moreover it also follows from the proof of 7.1 that the smooth
structures $\cS$, $\cS'$ on $Q-\{q\}$ are diffeomorphic. Here
$\cS$ is the given smooth structure on $Q-\{q\}$ and  $\cS'$
is the (restriction of) normal smooth structure on
$Q-\{q\}$ for $K$.
From the proof of 7.1 we see that $S=f(\sL(p,K))$ is a smooth submanifold of 
$(Q-\{q\},\cS')$. Let  $\cS'_N$ be 
the smooth structure on $S$ induced by $\cS'$. Note that  $(S,\cS'_N)$ 
is diffeomorphic to  $(N,{\tilde{\cS}}_N)$ (see (3) above). Note also that
$\cS'_N$
is a normal smooth structure on $N$ for $L$, as can be verified from the
proof of 7.1 applied to $N$ and $L$ (actually the link smoothings are the same).
Applying C.2.2 we get that $(S,\cS'_N)$ is $PL$ homeomorphic to $L$.
On the other hand condition (v) above (in the definition of cubification for
a one point singularity manifold) imply that $(N,\cS_N)$ is also $PL$ homeomorphic
to $L$.
It follows that all three $(S,\cS'_N)$, $(N,{\tilde{\cS}}_N)$, $(N,\cS_N)$
are $PL$ equivalent. If we are in dimension four all three spaces above
are diffeomorphic and (a) follows.\\

We prove (b). By hypothesis there is an embedding $(N,\cS_N)\ra (Q-\{q\},\cS)$
whose image lies near $q$. Since $\cS$ is diffeomorphic to $\cS'$, we get that
there is an embedding $(N,\cS_N)\ra (Q-\{q\},\cS')$ whose image lies near $q$.
But near $q$ the smooth structure $\cS'$ is diffeomorphic to the product $I_\epsilon\times S$ with smooth structure $\cS_{I_\epsilon}\times\cS'_N$, where $I_\epsilon=(0,\epsilon)$. And this product is diffeomorphic to 
$I_\epsilon\times N$ with smooth structure $\cS_{I_\epsilon}\times {\tilde{\cS}}_N$. 
Therefore there is an embedding $(N,\cS_N)\ra (I_\epsilon\times N,\cS_{I_\epsilon}\times {\tilde{\cS}}_N)$. Hence, by a standard topological
argument, 
$(N,{\tilde{\cS}}_N)$ and $(N,\cS_N)$ are smoothly h-cobordant. 
Consequently the proof of (b) is obtained by applying the h-cobordism Theorem. This
proves the Theorem.

\vspace{1in}

\noindent {\bf \large  Section 8. Smoothing Cones.}\\

Given an all-right spherical complex $P^m$ of dimension $m$ and a compatible smooth structure $\cS_P$ on $P$,
by Theorem 7.1 (see also remark 1 after the statement of 7.1) we can assume that $\cS_P$ is a normal smooth structure,
and $\cS_P$ has a normal atlas $\cA_P$.  The atlas $\cA_P$ and its induced differentiable structure $\cS_P$ are constructed (uniquely) from a set of link smoothings
$\cL=\{ h\0{\Delta}\}\0{\Delta\in P}$. To express this dependence we will sometimes write $\cA_P=\cA_P\big(\cL\big)$ and $\cS_P=\cS_P\big(\cL\big)$.\\

Recall that the cone $\rC P$ has a piecewise hyperbolic metric induced by the piecewise spherical metric on $P$ (see 6.5).
We denote these metrics by $\sigma\0{\rC P}$ and $\sigma\0{P}$ respectively.
As mentioned in section 6.5, the piecewise hyperbolic metric $\sigma\0{\rC P}$
is piecewise variable, hence it has a well defined ray structure as in
section 1D (see also 2.2). \\

\noindent {\bf (8.0.1.)}\,\,\, Consider the following data.
\begin{enumerate}
\item[{\bf 1.}] A positive number $\xi$.
\item[{\bf 2.}] A collection $\s{d}=\{ d_2, d_3,....\}$ of real numbers, with $d_i>(4+\xi)$.
We write $\s{d}(k)=\{ d_2, d_3,...,d_k\}$
\item[{\bf 3.}] 
A positive number  $r$, with $r>2d_i$, $i=2,...,m+1$, and $m$ as in item 4.
\item[{\bf 4.}] 
Real numbers $\varsigma>0$, $c>1$, 
with $c\,\varsigma<e^{-(4+\xi)}$.
Hence we get sets of widths  (see 6.3 and 6.4)
$\sA=\sB(\varsigma)=\{\alpha_i\}$
and $\sB=\sB(\varsigma;c)=\{\beta_i\}$,
where $sin\,\alpha_i=\varsigma^{i+1}$, $sin\,\beta_i=c\,\varsigma^{i+1}$. 
\item[{\bf 5.}] An all-right spherical complex $P^m$, $dim\, P=m$, with compatible normal atlas $\cA\0{P}\big(\cL_P\big)$, where $\cL_P$ is a (smoothly compatible, see 7.8.5) set
of link smoothings on $P$.
\item[{\bf 6.}] A diffeomorphism $\phi\0{P}=\phi\0{P,\cL\0{P}}: (P,\cS\0{P}(\cL_P))\ra \bS^{m}$ to the standard $m$-sphere. The map $\phi\0{P}$ is called a {\it global
smoothing for $P$, with respect to $\cS_P$ (or $\cA_p$, or $\cL_P$)}.
For $m=1$ we shall take $\phi\0{P}$ in a canonical way
(that is, depending only on $P$).
\end{enumerate}\vspace{.2in}

The smooth atlas $\cA\0{P}(\cL_P)$ on $P$ induces, by coning, a smooth atlas
on $\rC P-\{o\0{\rC P}\}$, and this atlas together with  the coning
$\rC\phi\0{P}:\rC P\ra\R^{m+1}$ of the map
$\phi\0{P} $
induce a smooth atlas $\cA\0{\rC P}=\cA\0{\rC P}\big(\cL\0{P},\phi\0{P}\big)$ on $\rC P$. We denote the corresponding smooth structure by $\cS\0{\rC P}=\cS\0{\rC P}\big(\cL\0{P},\phi\0{P}\big)$. Note that
we get a diffeomorphism
$\rC\phi\0{P}:(\rC P,\cS\0{\rC P})\ra\R^{m+1}$.\\

With the data given in items 1-6 above we will construct for every $P^m$, by induction on the dimension $m$, the {\it smoothed} Riemannian metric
$\cG\big(P,\cL\0{P},\phi\0{P},r,\xi,\s{d},(c,\varsigma)\big)$ on the cone $\rC P$ of $P$, where we consider $\rC P$ with smooth structure $\cS\0{\rC P}$. \\

In sections 8.1 and 8.2 we will assume $\xi$, $\s{d}$,
$c$, $\varsigma$ fixed. In particular we shall assume
$\s{A}$, $\s{B}$ fixed. 
So, to simplify our notation, we shall denote the smoothed metric by $\cG(P,\cL\0{P},\phi\0{P},r)$ or just $\cG(P,r)$ or
$\cG(P)$. In sections 8.3 and 8.4 we need to make explicit the dependence of the smoothed metric on  the other variables, and we will show that, given $\epsilon>0$, we can choose 
$r$ and $d_i$, $i=2,...m$, large so that $\cG\big(P,\cL\0{P},\phi\0{P},r,\xi,\s{d}, (c,\varsigma))$ has curvatures $\epsilon$-pinched to -1, provided the variables
satisfy certain conditions.
Before we begin with dimension 1 we need to discuss induced structures.\\

Let $\Delta=\Delta^k\in P$. The {\it restriction of $\cL_P$ to $\sL(\Delta,P)$}
is the set $\cL_P|_{\sL(\Delta,P)}=\{h_{\Delta'}\}_{\Delta\lneqq\Delta'}$. Sometimes we will
just write $\cL_{\sL(\Delta,P)}$ or, more specifically,
$\cL_{\sL(\Delta,P)}(\cL_P)$. The corresponding induced atlas on $\sL(\Delta,P)$ is $\cA_{\sL(\Delta,P)}(\cL_P)=\{h^\bullet_{\Delta'}\}_{\Delta\lneqq\Delta'}$,
and sometimes we will
simply write $\cA_{\sL(\Delta,P)}$ (see section 7.7). The
smooth structure on $\sL(\Delta,P)$ induced by $\cA_{\sL(\Delta,P)}$ will be denoted
by $\cS_{\sL(\Delta,P)}(\cL_P)$, or simply by $\cS_{\sL(\Delta,P)}$. By corollary 7.7.4
we have that, for $\Delta\in P$, the link smoothing $h_\Delta$ is a global
smoothing for $\sL(\Delta,P)$ with respect to $\cS_{\sL(\Delta,P)}$. Write $\phi_{\sL(\Delta,P)}=\phi_{\sL(\Delta,P)}(\cL_P)=h_\Delta$.
Therefore we obtain
the following {\it restriction rule}:\\

\noindent {\bf (8.0.2.)}\hspace{.9in} $\cL_P\longrightarrow \Big(\,\cL_{\sL(\Delta,P)}
(\cL_P)\,,
\,\phi_{\sL(\Delta,P)}(\cL_P)\, \Big)$ \\

\noindent where $\cL_P$ satisfies 5 in (8.0.1) for $P$, and the objects
$\cL_\Delta$, $\phi_\Delta$ satisfy 5, 6 of (8.0.1) for $\sL(\Delta,P)$. 
The smooth structure on $\rC\sL(\Delta,P)$ constructed from
the data $\big(\cL_\Delta,\phi_\Delta\big)$ will be denoted by
$\cS_{\rC\sL(\Delta,P)}(\cL_P)$, or $\cS_{\rC\sL(\Delta,P)}(\cL_\Delta,\phi_\Delta)$,
or simply by $\cS_{\rC\sL(\Delta,P)}$.\\ 

Let $\Delta^j<\Delta^k\in P$ and let $\Delta^l=\sL(\Delta^j,\Delta^k)$. We have the
identity $\sL(\Delta^l,\sL(\Delta^j,P))$\newline
$=\sL(\Delta^k,P)$ (see 6.4.2). The next lemma says that the restriction rule
(8.0.2) is transitive, that is, it respects the identity above.\\

\noindent {\bf Lemma 8.0.3.} {\it Let $\Delta^j<\Delta^k\in P$ and let $\Delta^l=\sL(\Delta^j,\Delta^k)$. Then we have}
$$
\cL_{\sL\big(\Delta^l,\sL(\Delta^j,P)\big)}\bigg(    \cL_{\sL(\Delta^j,P)} \big(  \cL_P \big)            \bigg)\,=\,\cL_{\sL(\Delta^k,P)}\big(  \cL_P \big)
$$
\noindent and
$$
\phi_{\sL\big(\Delta^l,\sL(\Delta^j,P)\big)}\bigg(    \cL_{\sL(\Delta^j,P)} \big(  \cL_P \big)            \bigg)\,=\,\phi_{\sL(\Delta^k,P)}\big(  \cL_P \big)
$$\\

\noindent {\bf Proof.} If we use the simplicial definition of link (see 6.1.1)
the identity $\sL(\Delta^l,\sL(\Delta^j,P))=\sL(\Delta^k,P)$ is an equality
of sets of simplices of $P$. Therefore the lemma follows from the definition
of $\cL$ and $\phi$. This proves the lemma.\\

Recall that we have an identification $\s{CS}(\Delta,P)=\rC\Delta\times\rC\sL(\Delta,P)$
(see 6.5.2, remark 5). The ``open" version of this identification is
$\stackrel{\circ}{\s{CS}}(\Delta,P)=\rC\dDelta\times\rC\sL(\Delta,P)$, where
$\stackrel{\circ}{\s{CS}}(\Delta,P)=\rC(\stackrel{\circ}{\s{Star}}(\dDelta,P))$.
Here $\stackrel{\circ}{\s{Star}}(\dDelta,P)=\,\stackrel{\circ}{\s{N}}_{\pi/2}(\dDelta,P)$. Note that  $\stackrel{\circ}{\s{CS}}(\Delta,P)$ as an open subset of
$\rC P$ has the induced smooth structure $\cS_{\rC P}|_{\stackrel{\circ}{\s{CS}}(\Delta,P)}$, and, for simplicity, we will just write $\cS_{\rC P}$. On the other hand note
that $\rC\dDelta=\HH^{k+1}_+\sbs\HH^{k+1}$ has the natural smooth structure
$\cS_{\HH^{k+1}}$, and $\rC\sL(\dDelta,P)$ has the smooth structure
$\cS_{\rC\sL(\Delta,P)}$. Therefore we can give
$\rC\dDelta\times\rC\sL(\Delta,P)$ the ``product" smooth structure
$\cS_\times=\cS_{\rC\dDelta\times\rC\sL(\Delta,P)}$.\\

\noindent {\bf Lemma 8.0.4.} {\it The following identification is a diffeomorphism}
$$
\Big(\,\stackrel{\circ}{\s{CS}}(\Delta,P)\,,\, \cS_{\rC P}\,   \Big)\,\,=\,\, \Big(  \,\rC\dDelta\times\rC\sL(\Delta,P)\,,\,   \cS_\times \, \Big)
$$

\noindent {\bf Proof.} We use the variables $s$, $t$, $r$, $y$, $v$, $x$, $w$, $u$, $\beta$ defined
in section 2.3 (also see 6.5.1 and 6.5.2 (5)).\\

We assume that the image of the chart $h^\bullet_\Delta$
is $\stackrel{\circ}{\s{N}}_{\pi/2}(\Delta,P)$ (see remark 3 after 7.1). By rescaling,
and using the notation in 6.5.2 (5) (see also 6.1 and 6.5.1)
we can write 

\begin{equation*}\begin{array}{lccc}
h^\bullet_\Delta\,\,:&\D^{m-k}(\pi/2)\times\dDelta&\longrightarrow&P\\\\
&(\,\beta\,u'\,,\, w\,)&\mapsto&\big[\,w,  h_\Delta(u')  \,\big](\beta)
\end{array}
\tag{1}
\end{equation*}\\

\noindent where $\D^{m-k}(\pi/2)$ is the disc of radius $\pi/2$, and 
we are expressing and element $\D^{m-k}(\pi/2)$ as $\beta u'$, with
$\beta\in [0,\pi/2)$, $u'\in\bS^{m-k-1}$. A chart for
$(\,\stackrel{\circ}{\s{CS}}(\Delta,P)\,,\, \cS_{\rC P}\,  )$
is the cone of $h^\bullet_\Delta$, which we shall denote by $h^\ast_\Delta$. Explicitly,
from (1) we have

\begin{equation*}\begin{array}{lccc}
h^\ast_\Delta\,\,:&\R_+\times\D^{m-k}(\pi/2)\times\dDelta&\longrightarrow&\rC P\\\\
&(\,s\,,\,\beta\,u'\,,\, w\,)&\mapsto&s\,\big[\,w,  h_\Delta(u')  \,\big](\beta)
\end{array}
\tag{2}
\end{equation*}\\

\noindent And for $(  \,\rC\dDelta\times\rC\sL(\Delta,P)\,,\,   \cS_\times \, )$  we can take the following chart

\begin{equation*}\begin{array}{lccc}
h^\dag_\Delta\,\,:&\R_+\times \R^{m-k}\times\dDelta&\longrightarrow&P\\\\
&(\,t\,,\,r\,u'\,,\, w\,)&\mapsto&\big(\,t\,w\,,\,  r\,h_\Delta(u')  \,\big)
\end{array}
\tag{3}
\end{equation*}\\

\noindent where we are expressing an element in $\R^{m-k}$ as
$r u'$, $r\in[0,\infty)$, $u'\in\bS^{m-k-1}$. From (2) and (3)  and
(5) of 6.5.2 we get

\begin{equation*}
\Big(h^\dag_\Delta\Big)^{-1}\circ h^\ast_\Delta\big(\,s\,,\, \beta\,u'   \,,\, w  \,\big)
\,\,=\,\, \Big( t\,,\,   r\,u'\,,\, w  \,\Big)
\tag{4}
\end{equation*}\\

\noindent and recall that (see 2.3) the relationship  between the variables
$s$, $\beta$, $t$, $r$ is the following. There is a right hyperbolic triangle
with catheti of length $t$, $r$, hypotenuse of length $s$ and angle $\beta$
opposite to the cathetus of length $r$. Using hyperbolic trigonometry we can
find an invertible transformation $(s,\beta)\ra (t,r)$. In particular 
$r=sinh^{-1}(sin\beta\, sinh(s))$.  The variables $s$ and $t$ are never zero,
but $\beta$ and $r$ could vanish. Note that $\beta=0$ if and only if $r=0$. 
To get differentiability at $\beta=0$ note that the map $(s,\beta u')\ra ru'$
can be rewritten as $(s,z)\ra(\frac{r(s,\beta)}{\beta}z)$, $\beta=|z|$, which is smooth
because $\frac{r(s,\beta)}{\beta}$ is a smooth even function on $\beta$. Similarly,
the smoothness of the inverse of the map in (4) follows from the fact that
the map $(t,r)\ra\frac{\beta(r,t)}{r}$ is a smooth even function on $r$.
This proves the lemma.

\vspace{.8in}

\noindent {\bf 8.1 Dimension one.} \\
An all-right spherical complex $P^1$ of dimension one satisfying (8.0.1) (6) 
is formed by a finite number
$k'$ of segments of length $\pi/2$ glued successively forming a circle. Hence we shall canonically take (up to rotation)
$\phi=\phi\0{P}:P\ra\bS^1$ so that $\phi$ maps each 1-simplex to an arc of length $2\pi/k'$, and it does so with constant speed. Using $\phi$  we shall identify $P$ with a circle of length $k'\pi/2$.   Therefore $P$ with metric $\sigma\0{P}$, is isometric to $\bS^1$
with metric $k\sigma\0{\bS^1}$, $k=k'/4$. Consequently we identify $\rC P$ with $\R^2$,
and $\rC P-\{ o\0{\rC P}\}$ to $\R^2-\{0 \}$ with hyperbolic metric
$\sigma\0{\rC P}=sinh^2(t)\,k\,\sigma\0{\bS^1}+dt^2$. Notice that this metric is smooth on $\R^2$ away from the cone point $o\0{\rC P}=0\in\R^2$,
and it does have a singularity at 0 unless $k=1$.\\

Recall we are assuming  $r> d_2$, where both $r$ and $d_2$ are given. Next we give two constructions of $\cG(P)$. The first one
is very explicit (see Gromov-Thurston \cite{GT}) and does not use 
(directly) any of the methods introduced
previously.  The second one looks more like
the inductive construction in 8.2, and uses the continuation metric construction
given in section 5.2. These two constructions are slightly different but both
satisfy the two properties {\bf P'1}, {\bf P'2} given below. Here is the first construction.\\

 Let $\rho$ be the function in section 3.4 (and 4.3).
Define $$\mu(t)=\mu\0{d_2,r,k}(t)=k\,\rho({\mbox{\scriptsize $\frac{t}{d_2}-\frac{r-d_2}{d_2}$}})\, + \,\big(1-\rho({\mbox{\scriptsize $\frac{t}{d_2}-\frac{r-d_2}{d_2}$}})\big)$$ 
\noindent hence
$\mu(t)=1$, for $t\leq r-d_2$ and $\mu(t)=k$ for $t\geq r$. Define
$$\cG(P,r)\,\,=\,\,sinh^2(t)\,\mu(t)\,\sigma\0{\bS^1}\,+\,dt^2
$$

Since the metric $\cG(P,r)$ is equal to the canonical hyperbolic warped metric
$sihn^2(t)\sigma\0{\bS^1}+dt^2$ on the ball of radius $r-d_2$, we can
extend $\cG(P,r)$ to the cone point $o\0{\rC P}=0\in\R^2$.
It is straightforward to verify that $\cG(P,r)$ satisfies the following three properties:
\begin{enumerate}
\item[{\bf P'1.}] The metrics $\cG(P,r)$ and $\sigma\0{\rC P}$ have the same ray structure
(see section 1D and lemma 1.1).
\item[{\bf P'2.}] The metric $\cG(P,r)$ coincides with $\sigma\0{\rC P}$ outside
the ball of radius $r$.
\item[{\bf P'3.}] The metric $\cG(P,r)$ coincides with $sinh^2(t)\sigma\0{\bS^1}
+dt^2$ on
the ball of radius $r-d_2$.
\item[{\bf P'4.}] The family of metrics $\{\cG(P,r)\}\0{r>d_2}$ has cut limits
everywhere (that is, it has cut limits on $I=\R$, see 5.1). 
Notice that $d_2$ is fixed and $r$ is the index of the family.\end{enumerate}

\noindent Actually from the definition of cut limit we have that the cut
limit of $\cG(p,r)$ at $b$ is \\

\noindent {\bf (8.1.1.)}\hspace{.5in}$\Big(\,\lim_{r\ra\infty}\mu\0{d_2,r,k}(r+b)\,\Big)\,\sigma\0{\bS^1}\,=\,\Big(\, 1\,+\,\big(k-1\big)\rho\big( 1+\frac{b}{d_2}  \big)\,\Big)\,\sigma\0{\bS^1}$\vspace{.1in}

Here is the second construction.  Recall that we are identifying $\rC P-\{o\0{\rC P}\}$ with
$\R^2-\{0\}$. Consider the constant $\odot$-family of metrics $\{g\0{r}\}_{r-\frac{1}{2}>d_2}$
given by $g\0{r}=k\sigma\0{\bS^1}+dt^2=\sigma\0{\rC P}$. Now just define (see section 5.2)
$$\cG(P,r)\,=\,   \cG(P,r,\s{d})  \,=\, \cC\0{d\0{2}-\frac{1}{2}}(g\0{r-\frac{1}{2}})$$ 
\noindent In this case we also have that $\cG(P)$ satisfies 
{\bf P'1}, {\bf P'2},  {\bf P'3} and {\bf P'4} (see proposition 5.2.2).
\vspace{.8in}

\noindent {\bf 8.2 The inductive step .} \\
Recall that in this section we are assuming $\xi$, $c$, $\varsigma$ (hence $\s{A}$, $\s{B}$),
and $\s{d}$ constant.
With the data $\xi$, $\s{A}$, $\s{B}$, $r>0$ and an all-right spherical complex $P$ we constructed in section 6.6 the numbers $r\0{k}=r\0{k}(r)$ and the sets 
$\cY(P,\Delta^k,r)$, $\cY(P,r)$, $\cX(P,\Delta^k,r)$, $\cX(P,r)$, 
where $\Delta^k\in P$. The inverse of the function
$r\0{k}=r\0{k}(r)$ shall be denoted by $r=r(r\0{k})$.
Recall also that in section 6.5 (see 6.5.2 (7)) we identified $\s{CS}(\Delta^k, P)$, with metric $\sigma\0{\rC P}|\0{\s{CS}(\Delta^k, P)}$, with
$\cE_{\rC\Delta^k}(\rC\sL(\Delta^k))$, with metric $\cE_k(\sigma\0{\rC\sL(\Delta^k,P)})$.
We will use these objects in this section.\\

Let $m\geq 2$ and
suppose that for every triple $(P,\cL\0{P},\phi\0{P})$,  $j=dim\,P\leq m-1$, as in items 
5 and 6 of (8.0.1) above, and $r>d_i$, $i=2,...,m+1$\,
there are couple of Riemannian metrics: the {\it smoothed} metric 
$\cG(P,\cL\0{P},\phi\0{P},r,\xi,\s{d},(c,\varsigma))$ (or simply $\cG(P,\cL\0{P},\phi\0{P},r)$, or even $\cG(P,r)$), 
and the {\it patched} metric $\ccP(P,\cL\0{P},r)$ (or just $\ccP(P,r)$), satisfying the 
following properties
\begin{enumerate}
\item[{\bf  P1.}] The smoothed metric $\cG(P,r)$ is a Riemannian metric on $(\rC P,\cS_{\rC P})$, and it has
the same ray structure as $\sigma\0{\rC P}$.
\item[{\bf  P2.}] The patch metric $\ccP(P,r)$ is a Riemannian metric defined on 
$\rC P-\B\0{r\0{j-2}-(2+\xi)}(\rC P)$ (with smooth structure $\cS_{\rC P}$), and it has
the same ray structure as $\sigma\0{\rC P}$.
\item[{\bf  P3.}] On $\cY(P,\Delta^k,r)$, $k\leq j-2= dim\,P-2$,
 the patched metric $\ccP(P,r)$ coincides with the
metric    (see section 2.2)
$$\cE_{\rC \Delta^k}\bigg(\cG\Big(\sL(\Delta^k,P),r\Big)\bigg)$$
where $\cG\Big(\sL(\Delta^k,P),r\Big)=\cG\Big(\sL(\Delta^k,P),\cL_{\Delta^k}(\cL_P),
\phi_{\Delta^k}(\cL_P), r\Big)$
is defined on the space $(\s{CS}(\Delta,P), \cS_{\rC P})$.
(Recall that $ \cY(P,\Delta^k,r)\sbs\s{CS}(\Delta^k,P)=\rC\Delta^k\times\rC\sL(\Delta^k,P)$, see 6.5.2 (5), (7) and lemma 8.0.4.)

\item[{\bf  P4.}] On $\cY(P,r)$
 the patched metric $\ccP(P,r)$ coincides with the
metric  $\sigma\0{\rC P}$.
\item[{\bf  P5.}] The metrics $\cG(P,r)$ and $\ccP(P,r)$ coincide on $\rC P-\B\0{r\0{j-2}}(\rC P)$.
\end{enumerate}\vspace{.2in}

\noindent Note that the patched metric $\ccP(P,\cL_P,r)$ does not depend on
$\phi\0{P}$.
Properties {\bf P3}, {\bf P4}, {\bf P5} and the definition of the sets
$\cX(P,\Delta^k,r)$, $\cX(P,r)$ imply\\

\begin{enumerate}
\item[{\bf  P6.}] On $\cX(P,\Delta^k,r)$, $k\leq j-2= dim\,P-2$,
 the smoothed metric $\cG(P,r)$ coincides with the
metric 
$$\cE_{\rC \Delta^k}\bigg(\cG\Big(\sL(\Delta^k,P),r\Big)\bigg)$$
where $\cG\Big(\sL(\Delta^k,P),r\Big)=\cG\Big(\sL(\Delta^k,P),\cL_{\Delta^k}(\cL_P),
\phi_{\Delta^k}(\cL_P), r\Big)$
is defined on the space $(\s{CS}(\Delta,P), \cS_{\rC P})$.
\item[{\bf  P7.}] On $\cX(P,r)$
the smoothed metric $\cG(P,r)$ coincides with the
metric  $\sigma\0{\rC P}$.
\end{enumerate}\vspace{.3in}

Note that the metrics $\cG(P^1,r)$ constructed for spherical all-right 
1-complexes in section 8.1, together with the choice $\ccP(P^1,r)=\sigma\0{\rC P^1}$ 
satisfy properties {\bf P1-P5}. Indeed {\bf P1'} implies {\bf P1},
{\bf P2'} implies {\bf P5} (recall $r_{-1}=r$, see 6.6) and 
{\bf P2}, {\bf P3}, {\bf P4} are trivially satisfied.\\

Now, assume we are given the data: $P$,
$dim\, P=m$, $\cL\0{P}$, $\phi\0{P}$, $r$ as  items 5 and 6 in (8.0.1) at the beginning of section 8. Define the patched metric $\ccP(P,r)=\ccP(P,\cL\0{P},r)$
on $\rC P-\B\0{r\0{m-2}-(2+\xi)}(\rC P)$ as in {\bf P3} and {\bf P4} above.
That is, we define $\ccP(P,r)$ by demanding that:

\begin{enumerate}
\item[{\bf  P''3.}] On $\cY(P,\Delta^k,r)$, $k\leq dim\,P-2$,
 the patched metric $\ccP(P,r)$ coincides with the
metric   $\cE_{\rC \Delta^k}\bigg(\cG\Big(\sL(\Delta^k,P),r\Big)\bigg)$.
\item[{\bf  P''4.}] On $\cY(P,r)$, the patched metric $\ccP(P,r)$ coincides with the
metric  $\sigma\0{\rC P}$.
\end{enumerate}\vspace{.2in}

\noindent {\bf Lemma 8.2.1.} {\it The patched metric $\ccP(P,r)$ defined by
properties {\bf P''3} and {\bf P''4} is well defined.}\\

\noindent {\bf Proof.} The metric $\ccP(P,r)$ is defined on the ``patches''
$\cY(P,\Delta,r)$, $\Delta\in P$, and $\cY(P,r)$. We have to prove that these definitions coincide on the
intersections $\cY(P,\Delta^k,r)\cap\cY(P,\Delta^j,r)$,
$\cY(P,r)\cap\cY(P,\Delta^j,r)$. If $\Delta^j\cap\Delta^k=
\emptyset$ then (vi) of lemma 6.6.1 implies $\cY(P,\Delta^j,r)\cap\cY(P,\Delta^k,r)=
\emptyset$. Also if $\Delta^j\not\sbs\Delta^k$ and $\Delta^k\not\sbs\Delta^j$ by (vii) of lemma 6.6.1, we also get
$\cY(P,\Delta^j,r)\cap\cY(P,\Delta^k,r)=\emptyset$. Therefore we assume
$\Delta^j<\Delta^k$, $j<k$.\\

Recall that $\cY(P,\Delta^j,r)\sbs\s{CS}(\Delta^j,r)$
and  $\cY(P,\Delta^k,r)\sbs\s{CS}(\Delta^k,r)$ (see 6.6.1 (i)).
The metrics  

\begin{equation*}
\begin{array}{ccc}
h
&=&\cE_{\rC \Delta^j}\bigg(\cG\Big(\sL(\Delta^j,P),\cL_{\sL(\Delta^j,P)}(\cL_P), \phi_{\sL(\Delta^j,P)}(\cL_P),r\Big)\bigg)
\end{array}
\tag{1}
\end{equation*}

\begin{equation*}
\begin{array}{ccc}g
&=&\cE_{\rC \Delta^k}\bigg(\cG\Big(\sL(\Delta^k,P),\cL_{\sL(\Delta^k,P)}(\cL_P), \phi_{\sL(\Delta^k,P)}(\cL_P),r\Big)\bigg)
\end{array}
\tag{2}
\end{equation*}

\noindent are defined on the whole of $\s{CS}(\Delta^j,P)$ and $\s{CS}(\Delta^k,P)$,
respectively. 
From remark 5 in 6.5.2 we have that $\s{CS}(\Delta^j,P)=\rC\Delta^j\times
\rC\sL(\Delta^j,P)$. And from lemma 6.6.3 we have that
$\cY(P,\Delta^k,r)\sbs\rC\Delta^j\times
\cX\Big(\sL\big(\Delta^j,P\Big), \Delta^l , r \Big)$, where $\Delta^l=\Delta^k\cap
\sL(\Delta^j,P)$ (alternatively $\Delta^l$ is opposite to $\Delta^j$ in $\Delta^k$, or $\Delta^l=\sL(\Delta^j,\Delta^k)$).
Hence it is enough to prove that the metrics $h$ and $g$ coincide on 
$\rC\Delta^j\times \cX\Big(\sL\big(\Delta^j,P\Big), \Delta^l , r \Big)$.
But (2) and corollary 6.5.5 (see also proposition 2.3.3) imply 

\begin{equation*}
g=\cE_{\rC \Delta^j}\bigg[
\cE_{\rC \Delta^l}\bigg(\cG\Big(\sL(\Delta^k,P),  \cL_{\sL(\Delta^k,P)}(\cL_P), \phi_{\sL(\Delta^k,P)}(\cL_P),       r\Big)\bigg)\bigg]
\tag{3}
\end{equation*}

\noindent Note that
the inductive hypothesis (that is, properties {\bf P3}, {\bf P5}, which imply {\bf P6}) applied to the data $\sL(\Delta^j,P)$ and $\Delta^l$
gives us that on the set $\cX\Big(\sL\big(\Delta^j,P\big), \Delta^l , r \Big)$
we have 

\begin{equation*}\cG\Big(\sL(\Delta^j,P),
\cL_{\sL(\Delta^j,P)}(\cL_P),\phi_{\sL(\Delta^j,P)}(\cL_P),
r\Big)\,=\,\cE_{\rC \Delta^l}(f)
\tag{4}
\end{equation*}

\noindent where

{\small\begin{equation*}f=
\cG\bigg(\sL(\Delta^l,\sL(\Delta^j,P)),
\cL_{\sL(\Delta^l,\sL(\Delta^j,P))}(\cL_{\sL(\Delta^j,P)}),
\phi_{\sL(\Delta^l,\sL(\Delta^j,P)}(\cL_{\sL(\Delta^j,P)}),
r\bigg)
\tag{5}
\end{equation*}}

\noindent From (5) and transitivity of the restriction rule (see 8.0.3) we get

\begin{equation*}f=
\cG\Big(\sL(\Delta^k,P),
\cL_{\sL(\Delta^k,P)}(\cL_{P}),
\phi_{\sL(\Delta^k,P)}(\cL_{P}),
r\Big)
\tag{6}
\end{equation*}

\noindent Putting together (1), (4) and (6) we obtain an equation with the same right-hand side as in (3) but with $h$ instead of $g$ on the left-hand side. This proves
that $g=h$ on $\cY(P,\Delta^j,r)\cap\cY(P,\Delta^k,r)$. \\

The proof that the patched metric is well defined on $\cY(P,\Delta^k,r)\cap\cY(P,r)$
uses a similar argument and it follows from lemma 6.6.4, the inductive hypothesis
applied to $\sL(\Delta^k,P)$ (that is properties {\bf P4}, {\bf P5} which imply {\bf P7})
and remarks 6 and 7 in 6.5.2. This proves the lemma.\\\\

Recall that $r\0{m-2}=r\0{m-2}(r)$. Let $r=r(r\0{m-2})$ be the inverse,
where we consider $r\0{m-2}$ as a large real variable.
For $P=P^m$ using $\rC\phi\0{P}$ we get an identification between
$\rC P$ and $\R^{m+1}$. Therefore we can consider the
family of metrics
$\Big\{\ccP\big(P,r(r\0{m-2})\big)\Big\}_{r\0{m-2}}$ as a $\odot$-family of
metrics on $\R^{m+1}$. Because of the unpleasant constant 1/2 in
the warp forcing process (see 4.4 and 5.2) we rescale this family and
now consider the $\odot$-family of metrics\vspace{.1in}

\noindent {\bf (8.2.2.)}\hspace{1.5in} $\Big\{\ccP\big(P,r(r\0{m-2})\big)\Big\}_{r\0{m-2}-
{\mbox{\tiny$\frac{1}{2}$}}}$ 

\noindent (recall that the indexation
of the family tells us where we take the spherical cuts).
We finally define\\

\noindent {\bf (8.2.3.)}\hspace{1in}
$\cG(P,r)\,\,=\,\,\cC\0{d\0{m+1}-\frac{1}{2}}\Big(\ccP\big(P,r(r\0{m-2})\big)   \Big)$\\

\noindent where, as mentioned above, we are applying the continuation process
of section 5.2 to the $\odot$-family in (8.2.2).
Then, by construction, this metric satisfies {\bf P3} and {\bf P4}. 
Properties {\bf P1} and {\bf P2} can be proved by induction on the dimension of
$P$ (using properties {\bf P1-P5}), together with remark 2.2.3.
Property {\bf P5} holds by construction and by (ii) of 5.2.2. 
This concludes the construction of the smoothed
Riemannian metric $\cG(P,r)=\cG(P,\cL_P,\phi_P,r,\xi,\s{d}, (c,\varsigma))$.\\

Hence, by construction, we have the following properties.\\\\

\noindent {\bf P8.} The smoothed metric $\cG(P^m,r)$ is hyperbolic on 
$\B_{r\0{m-2}-d\0{m+1}}(\rC P)$.\\

\noindent {\bf P9.} The patched metric $\ccP(P^m,r)$ does not depend 
$d_i$, $i>m$.\\

\noindent {\bf P10.} The smoothed metric $\cG(P^m,r)$ does not depend 
$d_i$, $i>m+1$.\\

\vspace{.8in}

\noindent {\bf 8.3. On the dependence of $\cG(P,r)$ on the variables $c$ and $\xi$.}\\
In this section we show that the smoothed metric $\cG(P,r)=\cG(P,\cL_P,\phi_P,\xi,
r,(c,\varsigma))$ does not depend on the variables $\xi$ and $c$, provided
$\varsigma$ is fixed and $c$, $\xi$ and $\varsigma$ satisfy certain relation. In the next section we will show that, assuming $\s{d}$ and $r$
large, the metric $\cG(P,r)$ is $\epsilon$-hyperbolic. However the excess of the 
$\epsilon$-hyperbolic charts does depend on the variables $c$ and $\xi$.
In the next result assume $\varsigma$ and $\s{d}$ fixed. We shall write
$\cG(P,r,\xi,c)=\cG(P,\cL_P,\phi_P,r,\xi,(c,\varsigma))$. \\

\noindent {\bf Proposition 8.3.1.} {\it Let $c'>c>1$ and $\xi'>\xi>0$ be such that  $c'\varsigma <e^{-(4+\xi')}$
Then on $\rC P-\B_{r\0{m-2}-(2+\xi)}(\rC P)$,  for $r>(1+\xi)$  we have}
$$
\cG(P,r,\xi',c')\,\,=\,\,\cG(P,r,\xi,c)
$$\\

\noindent {\bf Proof.} Write $\s{A}'=\s{B}(c',\varsigma)$. Denote by
$\cY'(P,\Delta,r)=\cY(P,\Delta,r,\xi',(c',\varsigma))$ the sets obtained by replacing $c$ and $\xi$ in the definition
of $\cY(P,\Delta,r)=\cY(P,r,\xi,(c,\varsigma))$ (see 6.6) by $c'$ and $\xi'$, respectively. 
Similarly we obtain $\cY'(P,r)$. 
We have 
\begin{equation*}
\cY(P,\Delta,r)\sbs\cY'(P,\Delta,r) {\mbox{\,\,\,and\,\,\,}}\cY(P,r)\sbs \cY'(P,r)
\tag{1}
\end{equation*}

We will prove
the proposition by induction on the dimension $m$ of $P^m$.
It can be checked from section 8.1 that the case $m=1$ does not depend
on the variables $c$ and $\xi$.
Assume $\cG(P^k,r,\xi',c'))\,\,=\,\,\cG(P^k,r,\xi,c)$, for every $P^k$, $k<m$.
Consider $P^m$. First we prove that the corresponding patched metrics
$\ccP(P,r,\xi',c')$ and $\ccP(P,r,\xi,c)$ coincide. But it follows from properties {\bf P3} 
and {\bf P4} applied to both metrics, the inductive hypothesis and (1) that 
$\ccP(P^m,r,\xi',c'))\,\,=\,\,\ccP(P^m,r,\xi,c)$ on $\cY(P,\Delta^k,r)$, for all $\Delta^k\in P$,
$k\leq m-2$, and on $\cY(P,r)$. Therefore, by  6.6.1 (iv), 
the metrics $\ccP(P^k,r,\xi',c'))$, $\ccP(P^k,r,\xi,c)$ coincide on 
$\rC P-\B_{r\0{m-2}-(2+\xi)}(\rC P)$.
Finally note that the smoothed metrics $\cG(P,r,\xi,c)$, $\cG(P,r,\xi',c')$ are obtained
from the corresponding patched metrics by using the continuation process of
section 5.2 (see also (8.2.3)). But this process depends only on $\s{d}$ and $r\0{m-2}=sinh^{-1}(\frac{sinh(r)}{sin(\alpha\0{m-2})})$. The former
is fixed and the later, since $sin(\alpha\0{m-2})=\varsigma^{m-1}$ (see 6.6),
is independent of $c$, $c'$, $\xi$ and $\xi'$. This proves the proposition.\\

In the next section we will need the following result.
We use the notation in the proof of the previous proposition.
Let $s'\0{m,k}=s\0{m,k}(r)$ be obtained from $s\0{m,k}=s\0{m,k}(r)$ by replacing $c$ by $c'$
(see 6.6).\\

\noindent {\bf Lemma 8.3.2.} {\it If $c'>c$ we have}
$$
s'\0{m,k}\,-\,s\0{m,k}\,\,>\,\,ln\Big(\frac{c'}{c}\Big)\,-\,1$$

\noindent {\it provided $r>1$}.\\

\noindent {\bf Proof.} From the definition at the beginning of 6.6 we have $s\0{m,k}=sinh^{-1}(c\,\frac{sinh(r)}{\varsigma^{m-k-2}})$ and
$s'\0{m,k}=sinh^{-1}(c'\,\frac{sinh(r)}{\varsigma^{m-k-2}})$.
A simple calculation shows that the function $t\mapsto sinh^{-1}(c't)-sinh^{-1}(c t)$
is increasing. And another calculation shows that the value of this function
at $t=1$ has value at least $ln(c')-ln(c)-1$. This proves the lemma.

\vspace{.8in}

\noindent {\bf 8.4. On the warped $\epsilon$-hyperbolicity of $\cG(P,r)$ .} \\
In this section we prove that the smoothed
metrics on $\rC P^m$ are $\epsilon$-hyperbolic, provided $d_{2},...,  d_{m+1}$ and $r$
are large enough, $m\leq \xi$ and $\varsigma$ is small (how small depending only on $\xi$). First we need the following lemma.
Recall that an element of $\rC P$ can be written as $sx$, $s\geq 0$, $x\in P$.\\

\noindent {\bf Lemma 8.4.1.}
{\it The family of metrics 
$$\Big\{\cG\big(P^m,r(r\0{m-2})\big)\Big\}\0{r\0{m-2}}$$
has cut limits on $\R$ (see 5.2).}\\

\noindent {\bf Proof.}
We prove this by induction on the dimension $m$ of
$P^m$. For $m=1$ the lemma follows from (v), (vi) and (vii) of proposition
5.2.2 (alternatively we can use {\bf P4'} and 8.1.1 in section 8.1)\\

\noindent {\bf Claim.} {\it Suppose the $\odot$-family of metrics
$\Big\{\cG\big(\sL(\Delta^k,P^m),r(r\0{m-k-3})\big)\Big\}\0{r\0{m-k-3}}$
has cut limits on $\R$. Then the $\odot$-family of metrics} $$
\Big\{\,\,\cE_{\rC\Delta^k}\Big(\,\cG\big(\sL(\Delta^k,P),r(r\0{m-2})\,\Big)\,\,\Big\}\0{r\0{m-2}}
$$
\noindent {\it also has cut limits on $\R$ over the set $\sN(\Delta^k,\alpha)$, for any $\alpha\in(0,\pi/2)$ (see 5.1).}\\

\noindent {\bf Proof of claim.} Without loss of generality take $\alpha>\alpha_{k}$
(recall $\alpha_{k}\in \s{A}$ is fixed, see 8.0.1 (4)).
Note that, by construction  (see 5.2.2 (i) or {\bf P8}), the family
$\Big\{\cG\big(\sL(\Delta^k,P),r(r\0{k})\big)\Big\}\0{r\0{k}}$
satisfies the hypothesis of proposition 5.3.1 stated at the beginning of 5.3:
the family has eventually constant cuts.
In proposition 5.3.1 take $\beta_0=\alpha_{k}$, $a=\alpha-\alpha_{k}$ 
and $c$ to be any positive real number.
Since $\alpha_{k}>0$ is fixed
and $r\0{m-k-3}=sinh^{-1}(sinh(r\0{m-2}) sin(\alpha\0{k}))$ (see beginning of 6.1)
we get that the claim follows from proposition 5.3.1. This proves the claim.
\\

We continue with the proof of lemma 8.3.1.
Assume the lemma holds for $P^k$, $k<m$. Let $P=P^m$
and $\Delta^k\in P$. Suppose that the lemma does not hold for 
the family ${\cal{F}}=\Big\{\cG\big(P^m,r(r\0{m-2})\big)\Big\}\0{r\0{m-2}}$.
We will show a contradiction. To simplify our notation write $s=r_{m-2}$
and $g_s=\cG\big(P^m,r(s)\big)$. We have $g_s=sinh^2(r)(g_s)_s+dr^2$,
where $r$ is the distance to $o\0{\rC P^m}$.
Since $\cF$ does not have cut limits on $\R$ the is a bounded closed interval
$I\sbs \R$ such that $\cF$ does not have cut limits on $I$. For $(x,b)\in P^m\times I$ write $g^\ast_s(x,b)=(g_s)_{s+b}(x)$. Note that
\begin{equation*}
sinh^2(s+b)g^\ast_s(x,b)+dt^2\,=\, g_s\big((s+b)x\big)
\tag{1}
\end{equation*}

\noindent Since we are assuming that $\cF$ does not have cut limits on $I$ we have that
the sequence $\{g^\ast_s\}$ defined on $P\times I$
does not converge in the $C^2$ topology.  Hence there is a derivative $\p^J$,
for some multi-index of order $\leq 2$, and sequences $s_n\ra\infty$, $x_n\ra x$,
$b_n\ra b$ such that $|\p^Jg^\ast_{s_n}(x_n,b_n)-\p^Jg^\ast_{s_{n+1}}(x_n,b_n)|\geq a$
for some fixed $a>0$, and $n$ even. By proposition
6.7.3 we have that $R_{x,b}(s)=(s+b)x\in\cY(P,\Delta^k,r(s))$, for some $\Delta^k$,
$k\leq m-2$, and $s>s'$, for some $s'$; or 
$R_{x,b}(s)=(s+b)x\in\cY(P,r(s))$, $s>s'$, for some $s'$.
We assume the first case: $R_{x,b}(s)=(s+b)x\in\cY(P,\Delta^k,r(s))$, for some $\Delta^k$,
$k\leq m-2$. The other case is similar.  Since $s$ is large, we also get that 
$R_{x,b}(s)=(s+b)x\in\cX(P,\Delta^k,r(s))$, $s>s'$.
Moreover, also by 6.7.3, we can assume
$R_{s_n,b_n}(s)=(s+b_n)x_n\in\cX(P,\Delta^k,r(s))$, for $s>s'$. 
But by property {\bf P6}, on $\cX(P,\Delta^k,r(s))$ the metric $g_s$
is equal to $\cE_{\rC\Delta^k}\Big(\,\cG\big(\sL(\Delta^k,P),r(s)\,\Big)$.
Consequently the family of metrics $\big\{ \cE_{\rC\Delta^k}\Big(\,\cG\big(\sL(\Delta^k,P),r(s)\,\Big) \big\}_s$ does not have cut limits on $I$ either.
But the claim, together with the inductive hypothesis, imply that this family
does have cut limits, which leads to a contradiction. This proves the lemma.\\

For a positive real number $\xi$ and a positive integer write $\xi\0{k}=\xi-k+\frac{1}{k}$. Note that $\xi\0{1}=\xi$.\\


\noindent {\bf Proposition 8.4.2.} {\it Let $\xi>0$, $c>0$ and $\varsigma<\frac{\sqrt{2}}{2}$ be such that 
\begin{enumerate}
\item[(i)] $\varsigma<e^{-(15+3\xi)}$
\item[(ii)] $c\geq e^{5+\xi}$
\end{enumerate}
Let $\big(P^m,\,\cL\0{P},\,\phi\0{P}\big)$, and
$\epsilon>0$. Then we have that $\cG\big(P,\cL\0{P},\phi\0{P},\xi,r,\s{d}, (c,\varsigma))$ is warped $\epsilon$-hyperbolic outside the ball of radius $a\0{m+1}(\epsilon,\xi)$
(see 3.3.10), with charts of excess $\xi\0{m}$, provided}
\begin{enumerate}
\item[{\it (a) }] {\it $d_i$ and $r-d_i$, $i=2,...m+1$, are sufficiently large.}
\item[{\it (b) }] {\it $m+1\leq \xi$.}
\end{enumerate}\vspace{.1in}

\noindent {\bf Remarks.} \\
\noindent {\bf 1.} By ``sufficiently large" in (a) we mean that there are
$r_i(P)$ and $d_i(P)$, $i=2,...,m+1$, such that the proposition holds
whenever we choose $r-d_i\geq r_i(P)$ and $d_i\geq d_i(P)$.\\ 
\noindent {\bf 2.} How large $r$ and the $d_i$ have to be depends on $\epsilon$.\\
\noindent {\bf 3.} Note that the choices of $c$, $\xi$ and $\varsigma$ do not depend on
$\epsilon$.\\
\noindent {\bf 4.} 
If we want the smoothed metric on a cone $\rC P^m$
to be $\epsilon$-hyperbolic (outside the ball of radius $a\0{m+1}(\epsilon,\xi)$) we can choose
$\xi=m+1$, $c=e^{5+\xi}$ and $\varsigma=e^{-(16+3\xi)}$. With these choices the method would not work for $P$ of dimension $>m$.\\
\noindent {\bf 5.} The condition $\varsigma=e^{-(15+3\xi)}$ is stronger than the
condition $\varsigma<e^{-(4+\xi)}$. The later is used to construct the
smoothed metric but it is not strong enough to give us $\epsilon$-hyperbolicity.\\
{\bf 6.} The condition ``$\epsilon$-hyperbolic outside the ball of radius $a\0{m+1}(\epsilon,\xi)$" is the best we can hope for because this is what happens in
hyperbolic space (see 3.3.10). As mentioned before this is because, even though
our definition of $\epsilon$-hyperbolicity is well suited for our purposes, it has
some flaws.\\

\noindent {\bf Proof.} We will only mention the relevant objects to our argument
in the notation for the smoothed metrics. That is we will write
$\cG(P,\s{d},r,\xi, (c,\varsigma))$ or just $\cG(P,\s{d},r)$.\\

Our proof is by induction on the dimension $m$ of $P^m$, with $m+1<\xi$. So, we assume $c$, $\xi$ and $\varsigma$ fixed and satisfying (i) and (ii), that is, $\varsigma<e^{-(15+3\xi)}$ and $c\geq e^{15+3\xi}$. We also assume 
$\epsilon >0$ and without loss of generality we can assume
\begin{equation*}
\epsilon\,<\,\frac{1}{(1+\xi)^2}
\tag{1}
\end{equation*}

For $m=1$ we have that the lemma follows from 8.1 and
Theorem 5.2.5 by writing $\lambda=r$, choosing $g_r=\sigma\0{\rC P}$,
replacing $\xi$ by $\xi+1$, and taking $\epsilon'=\epsilon/4$. 
Also, since $g_r=\sigma\0{\rC P}$ is $\epsilon$-hyperbolic, for every $\epsilon$,
we can take $\epsilon=0$ in 5.2.5. With all these choices Theorem 5.2.5 implies that
$\cG(P,d_2,r)$ is $\epsilon$-hyperbolic, with charts of excess $\xi=\xi\0{1}$, provided $r-d_2$ and $d_2$ are large enough.\\

Let $m$ such that $m+1\leq\xi$.
We now assume that the proposition holds for all $k<m$. That is, given $\epsilon>0$
and $P^k$,  the smoothed metric $\cG(P^k,r,\s{d})$ is warped $\epsilon$-hyperbolic
outside the ball of radius $a\0{k+1}(\epsilon,\xi)$,
with charts of excess $\xi\0{k}$,
provided $r-d_i$ and $d_i$, $i=2,...,d_{k+1}$ are large enough. Note that, since we
are assuming $k<m$, we get that $k+1<\xi$. 
We use the following notation

\begin{equation*}A\0{k}=C_3'(m-k,k+1,\xi)\hspace{.5in}
B=C_6(\xi)\hspace{.5in}
\epsilon\0{k}=\frac{\epsilon}{3A\0{k}B}
\tag{2}
\end{equation*}

\noindent where $C_3'$ is as in Theorem 3.5.8 and $C_6$ as in Theorem 5.2.5.\\

\noindent Let $P=P^m$. For $k<m$ write

$$
L_k=\Big\{ \sL(\Delta^k,P) \Big\}\0{\Delta^k\in P}$$

\noindent A generic element in $L_k$ will be denoted by $Q=Q^j$, $j+k=m-1$. By inductive
hypothesis, for each $Q^j$ there are $r(Q^j)$ and $d_i(Q^j)$, $i=2,...,j+1$ such that
$\cG(Q,r(Q^j),\s{d}(Q^j))$ is $\epsilon\0{k}$-hyperbolic
outside the ball of radius $a\0{j+1}(\epsilon,\xi)$, with charts of excess $\xi\0{j}$,
provided $r-d_i\geq r_i(Q^j)$ and $d_i\geq d_i(Q^j)$. For $i\leq m$, let
$d_i(P)$ be defined by\\

\hspace{1.6in}$d_i(P)\,=\, {\mbox{max}}\0{Q^j,\,i\leq j+1}\Big\{ d_i(Q^j)   \Big\}$\\

\noindent We write $\s{d}(P)=\{ d_2(P),...,d_{m}(P),...\}$ where
$d_i(P)$, $i>m+1$, is any positive number. This is just for notational purposes
and the arguments given below will not depend the $d_i(P)$, $i>m+1$.
We do reserve the right to later choose $d_{m+1}(P)$ larger.
Also write\\

\hspace{1.1in}
$ r_i(P)=d_i(P)\,+\, {\mbox{max}}\0{{Q^j,\, i\leq j+1}}\big\{\,4\,ln(m)\,,\, r_i(Q^j)\,  \big\}$\\

\noindent Therefore, from the inductive hypothesis and property {\bf P8} we get that\\

\noindent {\bf (8.4.3.)} {\it 
For every $Q^j\in L_k$, the metric 
$\cG(Q^j,r,\s{d})$ is $\epsilon\0{k}$-hyperbolic outside the ball of 
\newline \hspace*{.61in}radius $a\0{j+1}(\epsilon,\xi)$, with charts of excess 
$\xi\0{j}$, provided $r-d_i\geq r_i(P)$ and $d_i\geq d_i(P)$, 
\newline \hspace*{.61in}$i=2,...,k+1$. Moreover,  $\cG(Q^j,r,\s{d})$
is hyperbolic on $\B_{r\0{j-2}-d_{j+1}}(\rC Q^j)$.}\\

\noindent  By definition we have $r\0{i}(P)\geq 4\,ln(m)$. Hence, if $r>r\0{i}(P)$ and $0\leq j\leq m-1$
we get that $\xi\0{j}-e^{-(r/2)}>\xi-j+\frac{1}{j+1}\geq\xi-(m-1)+\frac{1}{m}$. This together with 
(8.4.3), Theorem 3.5.8 and the definitions given in
(2) imply  that\\

\noindent {\bf (8.4.4.)} {\it 
For every \,$\sL(\Delta^k,P)\in L_k$, the metric 
$\cE_{k+1}\big(\cG(\sL(\Delta^k,P),r,\s{d})\big)$, defined  on
\newline \hspace*{.58in} $\cE_{k+1}\big(\rC\sL(\Delta^k,P)  \big)$, is 
$(\frac{\epsilon}{3\,B})$-hyperbolic
outside the ball of radius $a\0{j+1}(\frac{\epsilon}{3\,B},\xi)$, with
\newline \hspace*{.57in} charts of excess 
$\xi-(m-1)+\frac{1}{m}$, 
provided $r-d_i\geq r_i(P)$ and $d_i\geq d_i(P)$, 
\newline \hspace*{.61in}$i=2,...,k+1$. Moreover 
$\cG(\sL(\Delta^k,P),r\0{P},\s{d})$
is hyperbolic on $\B_{r\0{j-2}-d_{j+1}}(\cE_{k+1}(\rC\sL(\Delta^k,P)))$.}\\

\noindent {\bf Lemma 8.4.5.} {\it The patched metric $\ccP(P,r,\s{d})$ is
warped $(\frac{\epsilon}{3\, B})$-hyperbolic, provided $r-d_i\geq r_i(P)$,
$d_i\geq d_i(P)$, $i=2,...,m$.}\\

\noindent {\bf Proof.} 
Before we prove the proposition we need some preliminaries. Recall that we have functions $r\0{m,k}=r\0{m,k}(r)$
and $s\0{m,k}(r)$. For $\Delta=\Delta^k\in P$ write
$Y_\Delta=\cY(P,\Delta,r,\xi,(c,\varsigma))$ and 
$Y=\cY(P,r,\xi,(c,\varsigma))$ (see 6.6). For $\Delta=\Delta^k$, $k\leq m$ define

\begin{equation*}
N_\Delta\,=\,\s{N}_{s\0{m,k}}(\rC\Delta,\rC P)\,-\, \bigcup_{\Delta^l\in P,\, l<k}\s{N}_{s\0{m,k}}
(\rC\Delta^l,\rC P)\,-\,\B_{r\0{m-2}-(2+\xi)}(\rC P)
\tag{3}
\end{equation*}

\begin{equation*}
N\,=\,\rC P\,-\, \bigcup_{\Delta^l\in P,\, l\leq m-2}\s{N}_{s\0{m,k}}
(\rC\Delta^l,\rC P)\,-\,\B_{r\0{m-2}-(2+\xi)}(\rC P)
\tag{4}
\end{equation*}

\noindent 
For $k=m-1$ in the definition of $N_{\Delta^k}$ above take $s\0{m,m-1}=sinh^{-1}(c'sinh(r)\varsigma)$. Also for $k=m$ take 
$N\0{m,m}(\rC\Delta^m,\rC P)=\rC\Delta^m$.
Write $N_k=\cup_{\Delta^k\in P}N_{\Delta^k}$. Define
$c'=e^{5+\xi}c$ and $\xi'=6+2\xi$. From  hypothesis (i), that is from $c\varsigma<e^{-(15+3\xi)}$, we get that $c'\varsigma<e^{-(4+\xi')}$,
hence we can define the sets $Y'_\Delta=\cY(P,\Delta,r,\xi',(c',\varsigma))$
and $Y'=\cY(P,r,\xi',(c',\varsigma))$ (see 6.6). That is

\begin{equation*}
Y'_\Delta\,=\,\stackrel{\circ}{\s{N}}_{s'\0{m,k}}(\rC\Delta^k,\rC P)\,-\, \bigcup_{\Delta^l\in P,\, l<k}\s{N}_{r\0{m,k}}
(\rC\Delta^l,\rC P)\,-\,\B_{r\0{m-2}-(2+\xi')}(\rC P)
\end{equation*}

\noindent and analogously for $Y'$. Here $s'\0{m,k}$ is defined by replacing
$c$ by $c'$ in the definition of $s\0{m,k}$. Note that if we replace $c$ by 1
in the definition of $s\0{m,k}$ we obtain $r\0{m,k}$. This together with
hypothesis (ii),  the definition of $c'$, and lemma 8.3.2 imply

\begin{equation*}
\begin{array}{lll}
s'\0{m,k}\,-\,s\0{m,k}\,>\,\,5+\xi\\\\
s\0{m,k}\,-\,r\0{m,k}\,>\,\,5+\xi
\end{array}
\tag{5}
\end{equation*}

\noindent It follows from the condition $c'\varsigma<e^{-(4+\xi')}$, lemma
8.3.1 and properties {\bf P3}, {\bf P5} that for each $\Delta^k\in P$, $k\leq m-2$, we have that\\

\noindent {\bf (8.4.6.)} {\it  the metrics $\ccP(P,r,\xi,(c,\varsigma))$,
$\cE_{k+1}\big(\cG(\sl(\Delta^k,P), r,\xi, (c,\varsigma)\big)$ and 
$\cE_{k+1}\big(\cG(\sl(\Delta^k,P), r,\xi', (c',\varsigma)\big)$
\newline\hspace*{.55in} coincide on $Y_{\Delta^k}$}.\\

For $p\in\rC P$ denote the ball of radius $s$ centered at $p$ by
$\B_{s,p}(\rC P)$, with respect to the metric $\sigma\0{\rC P}$.\\

\noindent {\bf Claim 1.}
{\it For $\Delta=\Delta^k$, $k\leq m-2$, we have that}\hspace{.2in}
$d\0{\ccP}\Big(\,  N_\Delta    \,,\, \rC P\,-\, Y'_\Delta      \,\Big)\,\, \geq\, 
\, 5+\xi$\\

\noindent {\bf Remark.} Here $d\0{\ccP}(.,.)$ denotes path distance with respect to the metric  $\ccP(P,r)$.\\

\noindent {\bf Proof of claim.} Write $\Delta=\Delta^k$. From the definitions we have
$N_\Delta\sbs Y_\Delta\sbs Y'_\Delta\sbs\s{CS}(\Delta,P)$. Note that
$\s{CS}(\Delta,P)\sbs \rC P$ but we can also  consider
$\s{CS}(\Delta,P)\sbs E=\cE_{k+1}(\rC\sL(\Delta,P))$.
We claim that we can work on $E$, that is, it is enough to prove that
$d\0{E}\Big(\,  N_\Delta    \,,\, \rC\,-\, Y'_\Delta      \,\Big)>5+\xi$
($E$ with metric $\cE_{k+1}(\cG(\sl(\Delta,P),r))$).
To see this note first that $\overline{N_\Delta}\sbs int\,Y'_\Delta$ (this follows from
$s\0{m,k}<s\0{m,j}$, $j<k$). Hence
if there is $p\in N_\Delta$ and $q\notin Y'_\Delta$
and $\alpha$ a path joining $p$ to $q$ of $\ccP$-length $< 5+\xi$ then
there a restriction $\beta$ of $\alpha$ such that: (1) it begins at $p$, \,(2) it
ends at some point $q'\in\p Y'_\Delta$,\, (3) it is totally contained in $\overline{Y'_\Delta}$,\,
(4) its $\ccP$-length is $< 5+\xi$. By property {\bf P3} the path
$\beta$ (now considered in $Y'\sbs E$, with metric $\cE_{k+1}(\cG(\sL(\Delta,P),r))$) has
the same properties (1)-(4). This shows that we can work on $E$ with metric
$\cE_{k+1}(\cG(\sL(\Delta,P),r))$ instead of $\rC P$ with metric $\ccP(P,r)$.\\

Let $o\in \s{CS}(\Delta,P)\sbs E$ correspond to $o\0{\rC P}$.
Let $p\in N_\Delta$, and $q\notin Y'_\Delta$.
From definition (3) we have three cases.\\

\noindent {\bf Case 1.} $q\notin \s{N}_{s'\0{m,k}}(\rC\Delta,\rC P)$.\\
Since $\s{N}_{s'\0{m,k}}(\rC\Delta,\rC P)=\bar{\HH}^{k+1}_+\times
\B_{s'\0{m,k}}(\rC\sL(\Delta,P))$ and $N_\Delta\sbs\s{N}_{s\0{m,k}}(\rC\Delta,\rC P)=
\bar{\HH}^{k+1}_+\times
\B_{s\0{m,k}}(\rC\sL(\Delta,P))$, this case follows from the first inequality in (5).
Note that the $s$ neighborhoods of $\HH^{k+1}$ in $\cE_{k+1}(\rC\sL(\Delta,P))$
with respect to the metrics $\cE_k(\cG(\sL(\Delta,P),r))$ and $\cE_{k+1}(\sigma\0{\rC\sL(\Delta,P)})$ coincide (see remark 2.2.2 (3,4)).\\

\noindent {\bf Case 2.} $q\in \s{N}_{r\0{m,j}}(\Delta^j,P)$, $j<k$.\\
Since $p\in N_\Delta$ w have that $p\notin \s{N}_{s\0{m,j}}(\Delta^j,P)$,
this case is similar to case 1, but uses the second inequality in (5), instead
of the first.\\

\noindent {\bf Case 3.} {\it $q\in\B_{r\0{m-2}-(2+\xi')}(E)$, where the ball is centered
at $o$.}\\
Since $p\notin \B_{r\0{m-2}-(2+\xi)}(E)$ this case follows from the fact
that $\big(r\0{m-2}-(2+\xi)\big)-\big(r\0{m-2}-(2+\xi')\big)=\xi'-\xi=6+\xi>5+\xi.$
(Also see remark 2.2.2 (3).)
This proves the claim.\\\\

\noindent {\bf Claim 2.}  
{\it We have that}\hspace{.2in}$
d\0{\rC P}\Big(\,  N    \,,\, \rC P\,-\, Y     \,\Big)\,\, >\, 
\, 5+\xi$\\

\noindent {\bf Proof of claim.} First, we can extend the argument given in the
proof of claim 1 one more step, that is for $\Delta=\Delta^{m-1}$. The only problem
here is that we have not defined the sets $\cY(P,\Delta^k,r)$ for
$k=m-1$  (see 6.6). So just define $s'\0{m,m-1}=sinh^{-1}(c'sinh(r)\varsigma)$,
and $Y'_\Delta$ accordingly (as in 6.6). It can be checked, using the results in
section 6, that the argument above goes through and we get
$d\0{\ccP}\Big(\,  N_\Delta    \,,\, \rC P\,-\, Y'_\Delta      \,\Big)\,\, \geq\, 
\, 5+\xi$ (in this case the patched metric is just $\sigma\0{\rC P}$). Since $Y'_\Delta\sbs Y$ it remains to prove that for $\Delta=\Delta^m$
we have $d\0{\ccP}\Big(\,  N_\Delta    \,,\, \rC P\,-\, Y     \,\Big)\,\, \geq\, 
\, 5+\xi$. To deal with this case define $Y^*=\rC P-\bigcup_{k< m}\s{N}_{r\0{m,k}}(\rC\Delta^k,\rC P)$ and let $Y^*_\Delta$ be the component of $Y^*$ contained in
$\rC\Delta$. Here we are taking 
$r\0{m,m-1}=sinh^{-1}(sinh(r)\varsigma)$.
Then $N_\Delta\sbs Y^*_\Delta\sbs Y$. We want to prove that
$d\0{\ccP}\Big(\,  N_\Delta    \,,\, \rC P\,-\, Y^*_\Delta     \,\Big)\,\, \geq\, 
\, 5+\xi$, but this case now
can be reduced to the case $\Delta=\Delta^m\sbs\bS^m\sbs\HH^{m+1}=\rC \bS^m$,
which can easily be dealt with. This proves the claim.\\

We are now ready to prove the lemma.
Assume $p\in N_k$. We prove by induction on $k$ that the patched metric
$\ccP(P,r,\s{d})$ is warped $(\frac{\epsilon}{3\,B})$-hyperbolic, with excess
$\xi''=\xi-(m-1)+\frac{1}{m}$, that is, there is a warped  
$(\frac{\epsilon}{3\,B})$-hyperbolic chart $\phi:\T_{\xi''}\ra \rC P$
centered at $p$. We begin with $k=0$. Assume $p\in N_0$. Then $p\in N_{\Delta^0}$ for some $\Delta^0$. From (8.4.4) there is a warped $(\frac{\epsilon}{3\,B})$-hyperbolic chart $\phi:\T_{\xi''}\ra \cE_k(\rC\sL(\Delta^0,P))$. But it follows from
lemma 3.1.2,  claim 1 and (1) that the image $\phi(\T_{\xi''})\sbs Y'_{\Delta^0}$.
This together with (8.4.6) imply $\phi$ is also  chart for $\ccP(P,r)$. This proves
the case $k=0$. The inductive step  $k\leq m-2$ is similar. It remains to prove
the case $p\in N$. But this case follows from a similar argument as 
above (in this case fitting a chart in $Y$) and using claim 2 and property {\bf P5}. 
This proves lemma 8.4.5.\\\\

We now finish the proof of proposition 8.4.2. Take $\epsilon'=\frac{\epsilon}{6}$
and  apply Theorem 5.2.5 to the family
$\big\{ \ccP(P,r(r\0{m-2}),\s{d}) \big\}_{r\0{m-2}-\frac{1}{2}}$.
Note that we have to use lemma 8.4.1
to satisfy one of the hypothesis of Theorem 5.2.5.  Since $3\epsilon'+B\frac{\epsilon}{3\,B}<\epsilon$
(recall $B=C_6$, see (2)) from Theorem 5.2.5 we
 obtain a number $r\0{m+1}(P)$ and a (possibly larger) number  $d_{m+1}(P)$ such that
$\cG(P,r,\s{d})$ is warped $\epsilon$-hyperbolic, provided $r-d_i\geq r_i(P)$ and
$d_i\geq d_i(P)$, $i=2,...,m+1$. Finally note that the excess of the charts given
by Theorem 5.2.5 is $\xi''-1=(\xi-(m-1)+\frac{1}{m})-1=\xi\0{m}$.
This proves proposition 8.4.2.

\vspace{.8in}

\noindent {\bf 8.5.  Smoothing cones over manifolds.}\\
As in the beginning of section 8, let $P^m$ be an all-right spherical complex
and $\cS_P=\cS(\cL_P)$ a compatible normal smooth structure on $P$.
In the previous sections we have canonically constructed a 
Riemannian metric $\cG(P,\cL_P, \phi_P,r,\xi,\s{d},(c,\varsigma))$ on the cone
$\rC P$. An important assumption was that $(P,\cS_P)$ was diffeomorphic 
(by $\phi\0{P}$) to the sphere $\bS^m$. We cannot expect to do the same
construction on a general manifold $P$ because $\rC P$ is not in 
general a manifold. But we will canonically construct a complete  Riemannian 
metric on $\rC P-o\0{\rC P}$ that has some of the previous properties.\\

We consider the same data as before: $P^m$, $\cL_P$, $r$, $\xi$, $\s{d}$, $(c,\varsigma)$ satisfying (8.0.1)
but with one change, replace the map $\phi\0{P}$ by a Riemannian metric
$h\0{P}$ on the closed smooth manifold $(P,\cS_P)$. Hence we begin with the 
following data:  $P^m$,  $\cL_P$, $h\0{P}$,
$r$, $\xi$,  $\s{d}$, $(c,\varsigma)$.\\

First note that, by Theorem
7.1, a compatible normal smooth structure on $P$ exists. 
We will assume that $P$ has either dimension $\leq 4$ or $Wh(\pi_1P)=0$,
so that we can apply 7.9.1.
Note also that the sets $\cY(P,\Delta, r)$, $\cY(P,r)$ are defined for general $P$ (no just
for $P=\bS^m$) and satisfy all the properties given in section 6.6. 
Now, since all the links of $P$ are spheres we can define, as in 8.1 and 8.2
the patch metric $\ccP(P,r)=\ccP(P,\cL_P,r, \xi,\s{d},(c,\varsigma))$ on $\rC P-\B_{r\0{m-2}-2}(\rC P)$, and this metric satisfies properties {\bf P2}, {\bf P3}, {\bf P4} given in section
8.2.\\

Recall that in 8.2 this construction is completed by applying the continuation
method to the $\odot$-family of metrics $\Big\{\ccP\big(P,r(r\0{m-2})\big)\Big\}_{r\0{m-2}-\frac{1}{2}}$.
This method consists of two parts: warp forcing and then hyperbolic forcing.
In our more general setting here we can still apply warp forcing, but we cannot
directly apply hyperbolic forcing (at least not in the way given in section 5)
because we do not have $P=\bS^m$. In our case, 
to finish our construction we apply first warp forcing and the a version of
hyperbolic forcing for general $P$; this new version will use the metric 
$h\0{P}$ instead of the canonical metric $\sigma\0{\bS^m}$ on the sphere $\bS^m$.\\

Consider now the  $\odot$-family of metrics $\Big\{\ccP\big(P,r(r\0{m-2})\big)\Big\}
_{r\0{m-2}-\frac{1}{2}}$ and apply warp forcing to obtain
$$
g\0{r\0{m-2}}\,\,=\,\,\cW_{r\0{m-2}-\frac{1}{2}}\Big(\,\cG\big(P,\,r(r\0{m-2})  \big)\,\Big)
$$

\noindent and we have that $g\0{r\0{m-2}}$ is warped on 
$\B_{r\0{m-2}-\frac{1}{2}}(\rC P)-o\0{\rC P}$, specifically we have $g\0{r\0{m-2}}=sinh^2(t)\, g+dt^2$, where $g$
is a Riemannian metric on $P$ (it is the spherical cut of
$\ccP\big(P,r(r\0{m-2}))$ at $r\0{m-2}-\frac{1}{2}$) and $t$ is the distance-to-the-vertex
function on $\rC P$. Let $\rho\0{a,d}$ be the function in 3.4.
Now define the metric 
$\cG(P,h,r)=\cG(P,\cL_P,h,r,\xi,\s{d},(c,\varsigma))$ by
{\small$$
\cG(P,h,r)\,\,=\,\, \left\{
\begin{array}{ll}
h\,+\, \big(\rho\0{r\0{m-2}-d\0{m+1},d\0{m+1}-\frac{1}{2}}\circ t\big)\,\big(g\,-\,h\big) +dt^2 & {\mbox{on}}\,\,\, 
\B_{r\0{m-2}-\frac{1}{2}}(\rC P)\,-\, \B_{r\0{m-2}-d\0{m+1}}(\rC P)\\ \\
\big(\mu \circ t\big)^2 h\,+\, dt^2 & {\mbox{on}}\,\,\, 
 \B_{r\0{m-2}-d\0{m+1}}(\rC P)
\end{array}\right.
$$}

\noindent where $\mu(t)=\frac{e^{t}-e^{\lambda(t)}}{2}$, and
$\lambda=\rho\0{r\0{m-2}-2d\0{m+1},d\0{m+1}}$. 
Also we are assuming $r\0{m-2}-2d\0{m+1}>0$. Note that the metric
$\cG(P,h,r)=\frac{1}{2}e^t h+dt^2$ on $\B_{r\0{m-2}-2d\0{m+1}}(\rC P)-o\0{\rC P}$,
that is for  $0< t\leq r\0{m-2}-2d\0{m+1}$.
We write $\rC P-o\0{\rC P}=P\times (0,\infty)$
and extend the metric $\cG(P,h,r)$ to $P\times \R$ by 
$\frac{1}{2}e^t h+dt^2$ for $-\infty< t\leq r\0{m-2}-2d\0{m+1}$.\\

\noindent {\bf Corollary 8.5.1.} {\it The metrics $\cG(P,h,r)$ and $\ccP(P,r)$ have the following properties}
\begin{enumerate}
\item[{\bf (i)}] {\it $\cG(P,r)$ is a Riemannian metric on $P\times\R$, and it has
the same ray structure as $\sigma\0{\rC P}$ (on $P\times (0,\infty)$).}
\item[{\bf  (ii)}] {\it Properties  {\bf P3} and  {\bf P4}.}
\item[{\bf (iii)}] {\it We have $\cG(P,h,r)=\frac{1}{2}e^t h+dt^2$ for $-\infty< t\leq r\0{m-2}-2d\0{m+1}$.}
\item[{\bf (iv)}] {\it Given $\epsilon>0$ we have that the sectional curvatures of
$\cG(P,h,r)$ are $\epsilon$-pinched to -1 for $ t\geq r\0{m-2}-2d\0{m+1}$ provided
$r-d_i$, $d_i$, $i=2,...,m+1$, and $r-2d\0{m+1}$ are large enough.}
\end{enumerate}\vspace{.2in}

\noindent {\bf Proof.} Item (i) follows from lemma 1.1 and the same
argument used for {\bf P1}  in the spherical case.
Item (ii) is true by construction (see also lemma 8.2.1). Item (iii)
follows from the discussion above, and (iv) from 8.4.2 and Bishop-O'Neill
warp curvature formula \cite{BisOn}, p.27.

\vspace{1in}

\noindent {\bf \large  Section 9. On Charney-Davis Strict Hyperbolization.}\\

In this section we shall state the results  about the Charney-Davis strict
hyperbolization procedure that will be needed in section 10.
We use some of the notation in \cite{ChD}. In particular the canonical $n$-cube
$[0,1]^n$ will be denoted by $\square^n$.
(This differs with the notation used in section 7, where an $n$-cube was denoted by
$q^n$.)
Also $B_n$ is the isometry group
of $\square^n$.\\
\vspace{.8in}

\noindent {\bf  9.1.  Charney-Davis hyperbolization pieces.}\\
A {\it Charney-Davis hyperbolization piece of dimension $n$} is a compact
connected orientable hyperbolic $n$-manifold with corners satisfying the properties stated in lemma 6.2 of \cite{ChD}.
The group $B_n$ acts by isometries on $X^n$
and there is a smooth map $f:X^n\ra\square^n$ constructed in section 5 of
\cite{ChD} with certain properties. 
We collect some facts from \cite{ChD}.

\begin{enumerate}
\item[(1)] For any $k$-face $\square^k$ of $\square^n$ we have that $f^{-1}
(\square^k)$ is totally geodesic in $X^n$ and it is a Charney-Davis hyperbolization
piece of dimension $k$. The totally geodesic submanifold (with corners)
$f^{-1}(\square^k)$ is a $k$-$face$ of $X^n$. Note that the intersection of
faces is a face and every $k$-face is the intersection of exactly $n-k$ distinct
$(n-1)$-faces.
\item[(2)] The map $f$ is $B_n$-equivariant.
\item[(3)] The faces of $X^n$ intersect orthogonally.
\item[(4)] The map $f$ is transversal to the $k$-faces of $\square^n$, $k<n$.
\end{enumerate}

\noindent The $k$-face $f^{-1}(\square^k)$ of $X$ will be 
denoted by $X_{\square^k}$. The interior $f^{-1}(\dot{\square^k})$ will be denoted by $\dX_{\square^k}$.\\

\noindent {\bf Lemma 9.1.1.} {\it For every $n$ and $r>0$ the is a 
Charney-Davis hyperbolization piece of dimension $n$ such that the 
widths of the normal neighborhoods of every $k$-face, $k=0,...n-1$,
are larger that $r$.}\\

\noindent {\bf Proof.} A piece $X^n$ is constructed in section 6 of 
\cite{ChD} by cutting a closed hyperbolic $n$-manifold $M$ along a
system $\{Y_i\}_{i=1}^n$ of codimension one totally geodesic submanifolds of $M$
that intersect orthogonally. The $Y_i$'s are
orientable and two sided.
The group $B_n$ acts by isometries on $M$, permuting the $Y_i$'s.
In particular each $Y_i$ is contained in the fixed point set of a nontrivial
 isometric involution $r_i$. Therefore $r_i$ interchanges
both sides of $Y_i$.\\

The $(n-1)$-faces of $X^n$ correspond to the $Y_i$'s and a $k$-face of
$X^n$ corresponds to the (transverse) intersection of different $n-k$ $Y_i$'s.
Therefore it is enough to show that the intersections of the $Y_i$'s have
 normal neighborhoods with large width. \\

\noindent {\it Claim.} {\it It is
enough to show that the $Y_i$'s have normal neighborhoods with large width.}\\

\noindent {\it Proof of claim.} Let $Z=Y_i\cap Y_j$, $Y_i\neq Y_j$, and assume both $Y_i$, $Y_j$
have normal neighborhoods of width larger than $r$.
Let $\alpha$ be a path with end points in $Z$ of length $<2r$. Then $\alpha$
lies in the normal neighborhood of $Y_i$.
Using the distance decreasing 
normal geodesic deformation of $Y_i$ we can deform $\alpha$, rel end points,
to a shorter path $\beta$ in $Y_i$. Repeat the same argument now with
$\beta$, and using the fact that $Y_i$ and $Y_j$ intersect orthogonally,
we get that the deformation of $\beta$ lies in $Z$. 
The proof for larger intersections is similar. This proves the claim.\\

We continue with the proof of lemma 9.1.1. The claim above and the existence
of the nontrivial isometric involutions $r_i$ imply that it is enough to have
$M$ with large injectivity radius. To obtain this recall that $M$ is given in
section 6 of \cite{ChD} as $M=\HH^n/\Gamma$, where
$\Gamma=\Gamma({\cal{J}})$ is a congruence subgroup
given by the ideal ${\cal{J}}$ of the ring of integers of the totally real
quadratic extension $K=\Q(\sqrt{d})$ of $\Q$. The only conditions required
for ${\cal{J}}$ are that $\Gamma$ is torsion free and $\Gamma\sbs SO_o(n,1)$.
Hence any deeper congruence subgroup $\Gamma'$ will serve as well.
But, by taking deeper congruence subgroups we can increase (in a well-known way)
the injectivity radius as much as we want: let $\gamma_1,...,\gamma_k\in\Gamma$ correspond the closed geodesics of length $<r$, and now take a deeper 
ideal that contains none of the $\gamma_i$'s. This proves the lemma.\\

We need some extra properties for the map $f$, so we give an explicit construction of it.
Recall from the proof above that $X$ is obtained from the closed hyperbolic manifold
$M$ by cutting along the system $\{Y_i\}$. Similarly the map $f$ is obtained in \cite{ChD}
from a map $\varphi :M\ra  \T^n$. And this map has coordinate maps
$\varphi_i:M\ra\bS^1$, which are constructed by applying the Pontryagin-Thom
construction to the framed (two sided) codimension one submanifolds
$Y_i$. Here we need a bit more details so we give an specific construction for
$\varphi$.\\

Let $Y_i\times (-r,r)\sbs M$ be the normal geodesic neighborhood of
$Y_i$ of width $r>2$. Hence for $p=(y,t)\in Y_i\times(-r,r)$, the smooth map $p\mapsto t(p)=t$ gives the signed distance to $Y_i$.
Let  $\eta:\R \ra [-1,1]$ be a smooth map such that $\eta (t)=t/r$ for $t\in (-r+1,r-1)$,\,
$\eta (t)=1$ for $t\geq r$,\, $\eta (t)=-1$ for $t\leq -r$. By identifying $(\bS^1,1)$ with
$\big([-1,1]\big/{\mbox{{\small $-1=1$}}}\,\, ,\,\, 0\big)$,
the smooth map $\eta\circ t$ induces
the smooth map $\varphi_i:Y_i\times (-r,r)\ra \bS^1$ that can be
extended to the whole of $M$.  Note that $\varphi_i^{-1}(1)=Y_i$ and
$\varphi^{-1}(\T^{i-1}\times\{1\}\times\T^{n-i})=Y_i$.
After cutting along the $Y_i$'s we get the map $f:X\ra\square^n$ and each $Y_i$
corresponds to two $(n-1)$-faces $X_{\square\0{0}^{n-1}}$, $X_{\square\0{1}^{n-1}}$
(one for each side of $Y_i$), where $\square\0{j}^{n-1}=\square^{n}\cap
\{x_i=j\}$, $j=0,1$. Moreover the normal neighborhood
$Y_i\times(-r,r)$ corresponds to the two one-sided normal neighborhoods
$X_{\square\0{0}^{n-1}}\times [0,r)$, $X_{\square\0{1}^{n-1}}\times [0,r)$.
Write $f(p)=(f_1(p),...,f_n(p))\in\square^n\sbs\R^n$. Then if $p\in 
X_{\square\0{0}^{n-1}}\times [0,r)$, we have $f_i(p)=\frac{1}{2}\eta (t_i(p))$, where $t_i(p)$ is the distance
to $X_{\square\0{0}^{n-1}}$. Similarly if $p\in X_{\square\0{1}^{n-1}}\times [0,r)$
we have $f_i(p)=1-\frac{1}{2}\eta (t_i(p))$, where $t_i(p)$ is the distance
to $X_{\square\0{1}^{n-1}}$. 
And if $p\in X_{\square\0{0}^{n-1}} \times [0,r-1)$, we have 
$f_i(p)=\frac{1}{2r}\,t_i(p)$, 
and similarly for $X_{\square\0{1}^{n-1}}$.
In particular $p\in X_{\square\0{0}^{n-1}} \times [0,a)$, $a\leq r-1$, if and only
if $f(p)\in\square\0{0}^{n-1} \times [0,\frac{a}{r})$.
In what follows of this paper we assume
$\varphi$ and $f$ are constructed as above.\\

\noindent {\bf Proposition 9.1.2.} {\it The derivative of $f$ sends normal vectors
to $X_{\square^k}$ to normal vectors to $\square^k$.}\\

\noindent {\bf Proof.} For simplicity we assume that $k=n-2$. The proof for general
$k$ is similar. We can also assume that $\square^{n-2}=\square_1^{n-1}\cap\square_2^{n-2}$,
where $\square^{n-1}_i=\{x_i=0\}\cap\square^n$.
Write $U_i=X_{\square_i^{n-1}}$ and $W=X_{\square^{n-2}}=U_1\cap U_2$.
Let $p\in W$. We certainly have that $p\in U_1\times[0,r)$ and 
 $p\in U_2\times[0,r)$. Assume for simplicity that 
$p\in U_3\times[0,r)$ but $p\notin U_i\times[0,r)$, $i\geq 4$.
Let $u\in T_p W$ be normal to $W$ and $\alpha$ be the geodesic with
$\alpha (0)=p$ and $\alpha'(0)=u$.
We have to prove that $(Df_i)_p\,u=0$ for $i\geq 3$.
It follows from $p\notin U_i\times[0,r)$, $i\geq 4$, that  $(Dt_i)_p=0$, hence  $(Df_i)_p=0$, $i\geq 4$.
It remains to prove that $(Df_3)_p\,u=0$. Since $p\in U_3\times[0,r)$
there is a geodesic $\beta$ in $U_3\times[0,r)$, beginning at $U_3$, normal to $U_3$ 
and ending in $p$. Also $t_3(p)$ is equal to the length of $\beta$. Since
$\alpha$ and $\beta$ are perpendicular at $p$ the function $t\mapsto t_3(\alpha(t))$
has a minimum at 0. Hence $D(t_3)_p\, u=0$, therefore $D(f_3)_p \,u= \eta'(t_3(p))\,
D(t_3)_p\, u=0$. This proves the proposition.\\
\\

For a $k$-face $X_{\square^k}$ and $p\in X_{\square^k}$, the set of inward normal vectors to
$X_{\square^k}$ at $p$ can be identified with the canonical all-right $(n-k-1)$-simplex
$\Delta_{\bS^{n-k-1}}$. In this sense we consider $\Delta_{\bS^{n-k-1}}\sbs T_pX$.
Similarly we can consider $\Delta_{\bS^{n-k-1}}\sbs T_q\square^n$, for 
$q\in \square^k$. We make the convention that
the two identifications above are done with respect to an ordering of
the $(n-1)$-faces $X_{\square^{n-1}}$ of $X$ and the corresponding
ordering for $\square^n$. For instance the vectors in
$\Delta_{\bS^{n-k-1}}\sbs T_pX$ tangent to some $X_{\square^{n-1}}$
correspond to the same $(n-1)$-face of $\Delta_{\bS^{n-k-1}}$
as the vectors in $\Delta_{\bS^{n-k-1}}\sbs T_{f(p)}\square^n$ tangent to  
$\square^{n-1}$. With these identifications we get coordinates in
$\Delta_{\bS^{n-k-1}}$: we write $(u_1,...,u_n)=u\in \Delta_{\bS^{n-k-1}}$
where $u_i$ is the angle between $u$ and $X_{\square_i^{n-1}}$ (or $\square_i^{n-1}$).\\

\noindent {\bf Proposition 9.1.3.} {\it For $p\in \dX_k$,
we have that
\begin{enumerate}
\item[(i)] $Df_p$ sends non-zero normal vectors to non-zero normal vectors.
\item[(ii)] For $u\in \Delta_{\bS^{n-k-1}}$ we have $Df_p(u)=\frac{1}{2r}u$.
\item[(ii)] {\bf n} $\circ \,(Df_p|\0{\Delta_{\bS^{n-k-1}}}):\Delta_{\bS^{n-k-1}}
\ra\Delta_{\bS^{n-k-1}}$
is the identity, where {\bf n}$(x)=\frac{x}{|x|}$ is the normalization map.
\end{enumerate}}

\noindent {\bf Remark.} In (ii) we are using  the coordinates on 
$\Delta_{\bS^{n-k-1}}$ mentioned above to identify the normal tangent spaces of
$X_{\square^k}$ and $\square^k$.\\

\noindent {\bf Proof.} 
For simplicity we assume that $k=n-3$. The proof for general
$k$ is similar. As in the proof of proposition 9.1.2 write
$\square^{n-1}_i=\{x_i=0\}\cap\square^n$ and $U_i=X_{\square^i}$.
Write also $W=X_{\square^{n-3}}=U_1\cap U_2\cap U_3$.
Let $u=(u\0{1},u\0{2},u\0{3})\in \Delta_{\bS^{2}}\sbs T_p\, X$.
Then $u_i$ is the spherical distance from $u$ to $T_p U_i$, $i=1,2,3$
(or angle between $u$ and $U_i$).
Let $\alpha$ be the geodesic with $\alpha(0)=p$ and $\alpha'(0)=u$. Then
$u_i$ is the angle, at $p$, between $\alpha$ and $U_i$.
We have to prove that $D(f_i)_p\, u=\frac{1}{2r}u_i$ for $i=1,2,3$. Since the length of $\alpha|_{[0,t]}$
is $t$ and the distance  from $\alpha (t)$ to $U_i$, $i=1,2,3$, is $t_i(\alpha(t))$,
we get a  right hyperbolic triangle with hypotenuse of length $t$
and side equal to $t_i(\alpha(t))$ with opposite angle $u_i$. Hence,
the hyperbolic law of sines implies

\begin{equation*}t_i(\alpha(t))\,\,=\,\, sinh^{-1}\Big(sin\,(u_i)\, sinh\,(t)\Big)
\tag{1}
\end{equation*}

\noindent for $i=1,2,3$. But for $t$ small we have $f_i(\alpha(t))=\frac{1}{2r}t_i(\alpha(t))$
and a simple differentiation, at 0, shows $D(f_i)_p\, u =\frac{1}{2r}u_i$, $i=1,2,3$.
This proves (i), (ii) and (iii) and completes the proof of the proposition.\\

In what follows we assume the width of the
normal neighborhoods of the $X_{\square}$ to be much larger than the number $r>0$. Lemma 9.1.1 asserts this is always possible.
Fix a point $p\in \dX_{\square^k}$ and consider $\Delta_{\bS^{n-k-1}}\sbs T_p X $. 
The cone $\rC_r \Delta_{\bS^{n-k-1}}$ is the set $\{tu\,,\, 0\leq t<r\, ,\, u\in\Delta_{\bS^{n-k-1}}\}$.
We have the exponential map 
$E:\rC_r\Delta_{\bS^{n-k-1}}\ra X$, given by $E(u,t)=exp_p(t\,u)$. 
Write $\square^n=\square^k\times\square^l\sbs\R^k\times\R^l=\R^n$ and
denote by $p\0{i}:\R^n\ra\R^i$, $i=k,j$, the projections.
Also, as in section 6, we write $\bar{R}^n_+=[0,\infty)^n$.\\

\noindent {\bf Lemma 9.1.4.} {\it We have the following properties.}

\begin{enumerate}
\item[{(i)}] {\it The map $E$ respects faces, that is}\,\, {\small $E\,\Big(\,\, \big(\rC_r\Delta_{\bS^{n-k-1}}\big)
\,\,\,\cap \,\,\,T_p X_{\square^j}\,\,\Big)\,\, \,\sbs\,\,\, X_{\square^j}$}  
\item[{(ii)}]  {\it The map $f\circ E$ respects faces, that is}\,\, {\small $(f\circ E)\,\Big(\,\, \big(\rC_r\Delta_{\bS^{n-k-1}}\big)
\,\,\,\cap \,\,\,T_p X_{\square^j}\,\,\Big)\,\, \,\sbs\,\,\,\square^j$}. {\it Hence
$\big(\rC_r\Delta_{\bS^{n-k-1}}\big)
\,\,\,\cap \,\,\,T_p X_{\square^j}=\big(\rC_r\Delta_{\bS^{n-k-1}}\big)
\,\,\,\cap \,\,\,T_{f(p)}\square^j$.}  
\item[{(iii)}] {\it The map $p\0{l}\circ f\circ E$ does not depend on $p$.}
\item[{(iv)}]  {\it The values of the map $p\0{k}\circ f\circ E$ do not depend on the variable $u\in \Delta_{\bS^{n-k-1}}$ (but they do depend on $t$ and $p$).}
\item[{(v)}] {\it Write $T=(t_1,...,t_n)$.
The map $p\0{l}\circ T\circ E:\rC_s\Delta_{\bS^{n-k-1}}\ra\R^{n-k}_+$ is an embedding, provided $s<r-1$.}
\end{enumerate}

\noindent {\bf Remark.} In the second statement of (ii) we are considering both
$\rC_r\Delta_{\bS^{n-k-1}}\sbs T_pX $ and 
$\rC_r\Delta_{\bS^{n-k-1}}\sbs T_{f(p)}\square^n$.\\

\noindent {\bf Proof.} Statement (i) follows from the fact that each $X_{\square}$ is totally geodesic in $X$.
Statement (ii) also follows because $f$ respects faces. Differentiating,  using 
9.1.3 and the fact that the derivative of $E$ (at 0) is the identity we obtain
$\big(\rC_r\Delta_{\bS^{n-k-1}}\big)
\,\,\,\cap \,\,\,T_p X_{\square^j}\sbs T_{f(p)}\square^j$. To get the other inclusion
use 9.1.3 again and count dimensions.\\

To prove (iii) and (iv) we assume for simplicity,   as in the proof of 9.1.3,  that
$k=n-3$ and $\square^{n-3}=\square_1^{n-1}\cap\square^{n-1}_2\cap\square^{n-1}_3$,
where $\square^{n-1}_i=\{x_i=0\}\cap\square^n$. Also $l=3$ and
$\square^3=\square_{4}^{n-1}\cap...\cap\square_n^{n-3}$. Let $u=(u\0{1}, u\0{2},u\0{3})\in\Delta_{\bS^{2}}$. 
From the proof of 9.1.3, we have 

$$ p\0{l}\circ f\circ E(u,t)\,\,=
\Big( f_1\big( E(u,t)  \big)      \, ,\,       f_2\big( E(u,t)  \big)      \, ,\, 
 f_3\big( E(u,t)  \big)      \, ,\, \Big)
$$

\noindent where
$$  f_i\big( E(u,t)  \big)\, =\,\,\eta\,\bigg(\,\, sinh^{-1}\Big(sin\,(u_i)\, sinh\,(t)\Big)  \,\,\bigg) 
$$

\noindent which depend only on $u$ and $t$. This proves (iii).\\

Note that $p\0{k}\circ f\circ E(u,t)=\Big( f_4\big( E(u,t)  \big),..., f_n\big( E(u,t)  \Big)$.
We prove that  $f_4\big( E(u,t)  \big)$ is independent of $u$. The proof is similar for $i>4$. Write $U=X_{\square_4^{n-1}}$, $V=X_{\square^k}$, $W=U\cap V$ and let $p\in\dot{U}$. Note $U$, $V$, $W$ are
totally geodesic. Let $q\in W$ be such that the segment $[q,p]$ has length equal to the distance $d\0{X}(p, W)$ between
$p$ and $W$. 
If $d\0{X}(p, W)\geq r$ we are done because then $f_4\big( E(u,t)  \big)$ is constant. We assume
$L=d\0{X}(p, W)< r$.\\

Let $B$ be 
the union of the images of all geodesics of length $r$ in $U$ beginning at $q$ and perpendicular to $W$. Then
$B$ is isometric to the $r$-cone of the canonical all-right simplex $\Delta_{\bS^{n-4}}$, and it is totally geodesic.
Let $C$ be the union of the images of all geodesics of length $r$ in $X$ beginning at some point in $B$ and perpendicular to $U$. Then $C$ is isometric to $B\times [0,r)$ with the usual
$cosh$-warped product metric. We write $C=B\times [0,r)$.
Therefore $C$ is also totally geodesic, and $C$, $V$ intersect perpendicularly 
at $C\cap V=[q,p]$. Now, if $\alpha$ is a geodesic of length $<r$ beginning at $p$ and perpendicular to
$V$ then $\alpha$ in contained in $C$. Moreover, $\alpha$ is contained in $\ell\times [0,r)\sbs C$
for some ray $\ell\sbs B$ beginning at $q$. Note that $\ell\times [0,r)\sbs C$ is isometric to a convex set in $\HH^2$. Finally
we get that $d\0{X}\big(\alpha(t),U\big)$ can then be computed in $\HH^2$ as the length of the side  $a$ of
a quadrilateral with consecutive  sides $a,b,c,d$, angles $\angle ab=\angle bc=\angle cd =\pi/2$ and $length(a)=t$, $length(b)=L$.
This calculation only depends on $t$ and $L=length [q,p]$ (hence on the choice of $p$) but not on the ``direction" $u$.\\

To prove (v) note that from equation (1) in the proof of 9.1.3 and the fact that $s<r-1$ we get that for $u=(u_1,...,u_{n-k})$, $0\leq t <s$, we have $t_i(tu)=sinh^{-1}\big(\,sin(u_i)\, sinh(t)\,\big)$.
This equation together with $\Sigma_{i=1}^{n-k}sin^2(u_i)=1$ imply 
$$t=sinh^{-1}\Big( \,\big(\,\Sigma_{i=1}^{n-k} \,sinh^2(t_i)\,\big)^{1/2}\,\Big)$$
\noindent Since we also get $u_i=sin^{-1}\big(\,  \frac{sinh(t_i)}{sinh(t)} \,\big)$.
the map $\pi\0{l}\circ T \circ E= (t_1\circ E,...,t_{n-k}\circ E)$ has a continuous inverse. Moreover this inverse is clearly smooth when $t\neq 0$ and all $t_i<t$. But in these cases the derivative of $p\0{l}\circ T\circ E$ can be shown to be injective.
This proves the lemma.\\

\noindent {\bf Remark 9.1.5.} Using the method in the proof above together
with hyperbolic trigonometry we can find an explicit formula for the coordinate functions of the function in  (iv). It can be checked that
these maps are even on the variable $t$. 

\vspace{.8in}

\noindent {\bf  9.2. The Charney-Davis hyperbolization process.}\\
The strict hyperbolization process of Charney and Davis is done by gluing copies
of $X^n$ using the same pattern as the one used to obtain the cube complex $K$
from its cubes. This space is called $K_X$ in \cite{ChD}. 
We call this space the {\it piece-by-piece strict hyperbolization of $K$}.
Note that we get a map
$F:K_X\ra K$, which restricted to each copy of $X$ is just the map $f:X^n\ra\square^n$
in 9.1. We will write $X_{\square^k}=F^{-1}(\square^n)$, for a $k$-cube $\square^k$
of $K$.\\

But to obtain good differential and tangential properties the process described above
is not enough. Therefore in \cite{DJ} and \cite{ChD} an alternative method is
given. We describe this next. As before 
let $X^n$ be a hyperbolization piece
and $K$ be cube complex. We assume there is projection
$p:K\ra \square^n$ (see 7.2 of \cite{ChD}). Now consider $K_X$ given as the fiber product\\

$$\begin{array}{ccc}
K_X&  \stackrel{q\0{X}}{\longrightarrow}& X\\ \\
q\0{K}\downarrow&&\downarrow f\\ \\
K&  \stackrel{p}{\longrightarrow}& \square^n
\end{array} $$\\

\noindent that is $K_X=\{(y,x)\, :\, p(y)=f(x)\}$. Here $q\0{K}$, $q\0{X}$ are projections.
We call this space the {\it fiber-product strict hyperbolization
of $K$}.  We denoted both hyperbolizations by $K_X$ but we shall write $K_X^{{\mbox{\tiny piece-by-piece}}}$ and $K_X^{{\mbox{\tiny fiber-product}}}$
if we need to. We shall write $X_{\square^k}=q\0{K}^{-1}(\square^k)$, for a $k$-cube
$\square^k$ of $K$.\\

\noindent {\bf Remark.} The space $K_X^{{\mbox{\tiny fiber-product}}}$
does depend on the projection map $p$. For instance, if $p$ is a cube map,
that is $p|\0{\square^n}$ is an isometry for every $\square^n\in K$,  then $K_X^{{\mbox{\tiny piece-by-piece}}}$
and $K_X^{{\mbox{\tiny fiber-product}}}$ coincide. But in general these two
hyperbolizations are homeomorphic but the obvious homeomorphism (see below)
does not preserve the natural piecewise differentiable structures.
If needed we shall write $K_X^{{\mbox{\tiny fiber-product}}}(p)$
to show explicitly the dependence on $p$.\\

We now assume that $K$ has a smooth structure $\cS$ compatible with the
cube structure of $K$ (hence $1_K$ is a smooth cubification of the smooth manifold $K$). We assume further that the 
projection $p:K\ra \square^n$ is smooth. Using this and item (4) above (see the beginning of 9.1) it is argued in \cite{ChD} that 0 is a regular value
of the smooth map $(k,x)\mapsto p(k)-f(x)$. Therefore $K_X$ is a smooth submanifold of $K\times X$
(with trivial normal bundle).
Hence if $K^n$ has a smooth structure (compatible with the cube structure $K$) then $K_X$ has a natural smooth structure.  
This is an important point for us, so we need to analyze this with a bit more detail.
First we remark two facts:

\begin{enumerate}
\item[{\bf (i)}] For the regular value argument to work it is assumed (implicitly) in \cite{DJ}, 1f.5, that
the restriction
$p|_{\dsquare}$ of $p$ to every open cube $\dsquare$ of $K$ is an embedding.
(In \cite{DJ} simplices are used instead of cubes.)
\item[{\bf (ii)}]  Let $\square$ be
a $k$-cube of $K$. Then $p(\square)=\square^k$, for some $k$-face $\square^k$ of $\square^n$.
If $K$ has no boundary the smoothness of $p$ implies that  for every $y\in \dsquare$
the derivative $Dp_y(T_yK)=T_{p(y)}\square^k$. (This may not happen if, for
instance, the complex $K$ is equal to $\square^n$.)
\end{enumerate}

We will consider maps $p$ of the form $p=\bar{\rho}\circ c$, where $c$ is a cube map $c:K\ra \square^n$ (i.e. $c|\0{\square^n}$ is an isometry
for every $\square^n$) and $\bar{\rho}$ is
a {\it slow-down-at-the-boundary} map $\bar{\rho}:\square^n\ra\square^n$ given by $\bar{\rho}(x_1,...,x_n)\,=\, (\rho(x_1),...,\rho(x_n))$, with $\rho:I\ra I$ a smooth homeomorphism that
is a smooth diffeomorphism on $(0,1)$ and $\frac{d^k}{dt^k}\rho(0)=\frac{d^k}{dt^k}\rho(1)=0$, $k>0$. In what follows we write 
$K_X^{{\mbox{\tiny fiber-product}}}=K_X^{{\mbox{\tiny fiber-product}}}(\bar{\rho}\circ c)$.\\

\noindent {\bf Remark.} Let $c:K\ra\square^n$ be a cube map. For any $\square^n\in K$ the map
$(c|\0{\square^n})^{-1}:\square^n\ra\square^n\sbs K$ can be identified with
the inclusion $\square^n\hookrightarrow K$. In particular
if $K$ has a smooth structure compatible with the
cube structure of $K$,  the map $(c|\0{\square^n})^{-1}$
is an embedding. Also, note that cube maps $c:K\ra\square^n$ are not smooth
(unless $K=\square^n$).\\

The {\it copy of $X$  over an $n$-cube $\square^n$ of $K$} is $X_{\square^n}=q\0{K}^{-1}(\square^n)$. It is a ``copy" of $X$ because
the projection $q\0{X}|\0{X\0{\square^n}}:X_{\square^n}\ra X$ is a
homeomorphism. Note that the projection of $X\ra X_{\square^n}\ra\square^n$  is given by
$$
x\,\,\mapsto\,\, \Big(\, (p|\0{\square^n})^{-1}(\,  f(x)\,)\, ,\,x\,\Big)\,\,\mapsto
c|\0{\square^n}\circ (p|\0{\square^n})^{-1}(\,  f(x)\,)\,=\,\Big((\bar{\rho})^{-1}\circ f\Big)(x)
$$

\noindent where we are identifying the inclusion $\square^n\hookrightarrow K$
with $(c|\0{\square^n})^{-1}$ (see the remark above). Therefore
$X\ra X_{\square^n}\ra\square^n$ is not equal to $f$.
Notice that even though $q\0{X}$ is smooth, the map $q\0{X}|\0{X\0{\square^n}}$ is  not a diffeomorphism because
the natural (topological) embedding $\big(q\0{X}|\0{X\0{\square^n}}\big)^{-1}:X
\ra K_X$ is not smooth.\\

\noindent \big[If this map were smooth, the composition 
$(c|\0{\square^n})^{-1}\circ q\0{K}\circ \big(q\0{X}|\0{X\0{\square^n}}\big)^{-1}$
would be smooth (see remark above). But this composition is equal to
$(\bar{\rho})^{-1}$, which is not smooth.\big]\\

The price we paid for slowing down the cubes at the boundary
(via $\bar{\rho}$)
was that we sped up the copies of $X$ at the boundaries.
Therefore the natural piecewise hyperbolic (and piecewise differentiable)
structure of $K_X=K_X^{{\mbox{\tiny piece-by-piece}}}$ does not directly
give one in $K_X^{{\mbox{\tiny fiber-product}}}$. To solve this problem we
need the following proposition that says that $\bar{\rho}:\square^n\ra\square^n$ 
can be covered by a homeomorphism $X\ra X$.\\

\noindent {\bf Proposition 9.2.1.} {\it We can choose $\rho:I\ra I$
so that there is a smooth homeomorphism $P:X^n\ra X^n$ such that $fP=\bar{\rho} f$,
i.e. the following diagram commutes}

$$
\begin{array}{ccc}
X&\stackrel{P}{\longrightarrow}&X\\
f\downarrow&&\downarrow f\\
\square^n&\stackrel{\bar{\rho}}{\longrightarrow}&\square^n
\end{array}
$$

\noindent {\it Moreover we can choose $P$ so that
its restriction to every open face $\dX_{\square^k}$ is an embedding.}\\

\noindent {\bf Remark.} With a bit extra work we can get $P$ to be $B_n$-invariant,
but this fact will not be needed.\\

Note that from the construction of the map $\bar{\rho}$ we have that
$D\bar{\rho}|\0{q}.w=0$ for every $\square \in K$, $q\in\square$ and $w$ perpendicular
to $\square$. We have the following addition to lemma 9.2.1.\\

\noindent {\bf Addendum to Proposition 9.2.1.} {\it We can choose $P$ in
proposition 9.2.1 so that for any $\square\in K$ we have that
$DP|\0{p}.v=0$ for every $p\in\square$ and $v$ perpendicular
to $X_\square$.}\\

The proposition and its addendum are proved in appendix D. \\

Now we get a new embedding $X\ra X_{\square^n}$ given by
$\big(q\0{X}|\0{X_{\square^n}}\big)^{-1}\circ P$. This is the ``correct" embedding,
as the next proposition shows.\\

\noindent {\bf Proposition 9.2.2.} {\it 
The embedding $\big(q\0{X}|\0{X_{\square^n}}\big)^{-1}\circ P:
X\ra X_{\square^n}$ is a smooth embedding. Moreover the following
diagram commutes for every $n$-cube $\square^n$ of $K$.}\\

$$
\begin{array}{ccccc}
X&\stackrel{{\mbox{\tiny $\big(q\0{X}|\0{X_{\square^n}}\big)^{-1}\circ P$}}}{\longrightarrow}&X_{\square^n}&
\stackrel{{\mbox{\tiny inclusion}}}{\longrightarrow}& K_X\\\\
f\downarrow&& q\0{K}|\0{X_{\square^n}}\downarrow \,\,\,\,\,\,\,\,\,\,\,\,\,\,\,\,\,\,&&
\,\,\,\,\,\,\downarrow  q\0{K}\\\\
\square^n&\stackrel{1\0{\square^n}}{\longrightarrow}&\square^n&
\stackrel{{\mbox{\tiny inclusion}}}{\longrightarrow}&K
\end{array}
$$\vspace{.2in}

\noindent {\bf Proof.} 
From the (fiber-product) definition of $K_X$ (see diagram at the beginning of 9.2)
and proposition 9.2.1 we get that the following diagram commutes

\begin{equation*}
\begin{array}{ccccc}
K_X&\stackrel{q\0{X}}{\longrightarrow}&X&
\stackrel{P}{\longleftarrow}& X\\\\
q\0{K}\downarrow&&f\downarrow &&
\,\,\,\,\,\,\downarrow f\\\\
K&\stackrel{p}{\longrightarrow}&\square^n&
\stackrel{\bar{\rho}}{\longleftarrow}&\square^n
\end{array}
\tag{1}
\end{equation*}\vspace{.2in}

\noindent and it follows that the left square of the diagram in the statement
of proposition 9.2.2 commutes. The right square commutes by definition.\\

Write $g=\big(q\0{X}|\0{X_{\square^n}}\big)^{-1}\circ P$.
We have that the map $g:X\ra K_X\sbs K\times X$ is smooth if and only if the 
coordinate maps $q\0{K}\circ g$, $q\0{X}\circ g$ are smooth. 
First we have $q\0{X}\circ g=q\0{X}\circ \big(q\0{X}|\0{X_{\square^n}}\big)^{-1}\circ P=P$, which is smooth. \\

\noindent From diagram (1) and definition $p=\bar{\rho}\circ c$ we get 
$$q\0{K}\circ g=q\0{K}\circ 
\big(q\0{X}|\0{X_{\square^n}}\big)^{-1}\circ P=(p|\0{\square^n})^{-1}\circ \bar{\rho}\circ f=
\big(c|\0{\square^n}\big)^{-1}\circ \bar{\rho}^{-1}\circ\bar{\rho}\circ f=
\big(c|\0{\square^n}\big)^{-1}\circ f$$
\noindent Since $\big(c|\0{\square^n}\big)^{-1}$ is just the inclusion
$\square^n\ra K$ (see remark above) we get that $q\0{K}\circ g$ is smooth. 
It remains to prove that $g$ is a smooth embedding. 
Since $q\0{X}$ is smooth  and $P^{-1}$ is smooth on $\dX$
we get that $g^{-1}$ is smooth on $\dX_{\square^n}$. 
Note that the same argument shows that the restriction of $g$ to any
$\dX_{\square^k}$ is an embedding. Therefore if $u$ is  a non-zero
vector tangent to some $\dX_{\square^k}$ then $Dg.u$ is non-zero
and tangent to the corresponding $k$-face $g(X_{\square^k})$.
If $u$ is a vector normal to $\dX_{\square^k}$ then, by lemma 9.1.2, $Df.u$ is non-zero and normal to $\square^k\sbs \square^n$ and certainly $D\big((c|\0{\square^n})^{-1} \circ f\big).u $ is also
non-zero and normal to $\square^k$. 
But from diagram (1)
we have that $q\0{K}\circ g=(c|\0{\square^n})^{-1}\circ f$,
hence $Dg.u$ is non-zero. Moreover,
$q\0{K}$ sends $g(X_{\square^k})$
to $\square^k$, therefore $Dg.u$ is not tangent to $g(X_{\square^k})$.
This proves that $Dg$ is injective on every point of $\dX_{\square^k}$.
This proves the proposition.\\

We can now use $\big(q\0{X}|\0{X_{\square^n}}\big)^{-1}\circ P$  on each
copy of $X$ in $K_X^{{\mbox{\tiny piece-by-piece}}}$ and get a map
$\Phi:K_X^{{\mbox{\tiny piece-by-piece}}}\ra K_X^{{\mbox{\tiny fiber-product}}}$
that is a smooth embedding on each copy of $X$. Hence we can
consider $K_X$ as $K_X^{{\mbox{\tiny piece-by-piece}}}$ with the pulled
back (by $\Phi$) differentiable structure, or $K_X^{{\mbox{\tiny fiber-product}}}$
with the pushed forward piecewise hyperbolic structure. \\

\noindent {\bf Corollary 9.2.3.} {\it The following diagram commutes.}\\

$$
\begin{array}{lcl}
K_X^{{\mbox{\tiny piece-by-piece}}}&\stackrel{\Phi}{\longrightarrow}&K_X^{{\mbox{\tiny fiber-product}}}\\\\
F\downarrow&&\downarrow q\0{K}\\\\
K&\stackrel{1_K}{\longrightarrow}&K
\end{array}
$$\vspace{.2in}

The next result shows that $p=\bar{\rho}\circ c$ is smooth 
on $(K,\cS')$, where $\cS'$ is a normal smooth structure on $K$ for $K$ (see section 7).\\

\noindent {\bf Proposition 9.2.4.} {\it Let $\cS'$ be a normal smooth structure on $K$ for  $K$. Then $p :(K,\cS')\ra\square^n$ is smooth.}\\

The proof is presented in appendix E.\\

Note that we also have that the restriction $p|\0{\dsquare^i}$ on every open cube is an embedding, because the inclusions $\dsquare^i\ra (K,\cS')$ are also embeddings
(see remarks after 7.1).
Therefore the regular value argument in \cite{ChD} (see items (i) and (ii) at the beginning of 9.2)
goes through and we get the following result.\\

\noindent {\bf Corollary 9.2.5.} {\it We have that $K_X$ is a smooth submanifold of 
$(K,\cS')\times X$, with trivial normal bundle.}\\

We denote by $K_X'$ the submanifold $K_X\sbs (K,\cS')\times X$ with its
induced smooth differentiable structure.\\

\noindent {\bf Remark 9.2.6.}
The proof of proposition 9.2.2 also works if we replace $K_X$ by $K_X'$, but with one change: we have to substitute
$X_{\square^n}$ by $\dX_{\square^n}$, that is, the map $\dX\ra\dX_{\square^n}\sbs K_X'$
is an embedding (the key point in the proof is that the inclusion 
$(c|\0{\square^n})^{-1}:\square^n\ra (K,\cS')$ in not an embedding, but its
restriction to $\dsquare^n$ is). It follows that
$\dX_\square$ is a submanifold of $K_X'$, for every $\square\in K$.

\vspace{.8in}

\noindent {\bf  9.3. Normal neighborhoods.}\\
In section 7 we constructed a normal atlas and a normal differentiable structure for a given smooth
cubification of a smooth manifold.
In this section we will construct a similar  atlas on the Charney-Davis hyperbolization of $K$.
 In what follows of this section we will use the notation $K_X$ for $K_X^{{\mbox{\tiny piece-by-piece}}}$.
Recall that we are denoting by $K_X'$ the submanifold $K_X\sbs (K,\cS')\times X$ with its
induced smooth differentiable structure. Recall we have 
a map $\Phi: K_X\ra K_X'$ (see 9.2).
The normal smooth structure $\cS'$ on $K$ has a normal atlas
$\cA=\cA\big(  \cL   \big)$, where
$\cL=\big\{  h_{\square^k}   \big\}\0{\square^k\in K}$ is a
smoothly compatible set of link smoothings for $K$ (see remarks after 7.1).
We assume that the normal bundle of any face of the hyperbolization piece $X$
has width larger than $s\0{0}>0$.  Choose $r$, with $3r<s\0{0}$. By lemma 9.1.1 the number $s\0{0}$ can be taken as large as we want.\\

By lemma 9.1.2 we can use the derivative of
the map $F:K_X\ra K$ (in a piecewise fashion) to identify $\sL(X_{\square^k},K_X)$ with $\sL(\square^k,K)$, where in both cases we consider the ``direction" definition of
link (see remark 6.1.1), that is, the
link $\sL(X_{\square^k},K_X)$ (at $p\in \dX_{\square^k}$) is the set of normal vectors to $X_{\square^k}$ (at $p$) and  the
link $\sL(\square^k,K)$ (at $q\in\dsquare^k$) is the set of normal vectors to $\square^k$ (at $q$). Hence we write
$\sL(X_{\square^k},K_X)=\sL(\square^k,K)$; thus the set of links for $K$ coincides
with the set of links for $K_X$.\\

Let $X_{\square^k}\sbs X_{\square^{n}}$ be a $k$-face of $K_X$,
contained in the copy $X_{\square^n}$ of $X$ over $\square^n$. For a non-zero vector $u$ normal to $X_{\square^k}$ at $p\in X_{\square^k}$, and pointing inside 
$X_{\square^n}$, we have that $exp_p(tu)$ is defined and contained in
$X_{\square^n}$, for $0\leq t< s\0{0}/|u|$. Recall that $h_{\square^k}:\bS^{n-k-1}\ra\sL(\square^k,K)=\sL(X_{\square^k},K_X)$ is the smoothing of the
link corresponding to $\square^k$.
We define the map

$$H\0{\square^k}\,\,\,:\,\,\D^{n-k}\times \dX_{\square^k}
\,\,\,\,\longrightarrow\,\,\,\, K_X $$
\noindent given by

$$H\0{\square^k}(\,t\,v\,\, ,\,\,p\,)\,\,\,=\,\,\, exp\0{p}\,\Big(\,\,2r\,t\,\,h_{\square^k}(v)\,\,\Big)$$\vspace{.2in}

\noindent where $v\in\bS^{n-k-1}$ and $t\in [0,1)$.
For $k=n$ we have that $H\0{\square^n}$ is the inclusion $\dX_{\square^n}\sbs K_X$
(or we can take this as a definition).
Note that $H\0{\square^k}$ is a topological embedding because we are assuming
the width of the normal neighborhood of $X_{\square}$ to be larger than $s\0{0}>2r$.
A chart of the form of $H\0{\square^k}$ (for some link smoothing
$h_{\square^k}$) is called a {\it normal chart for the $k$-face $X_{\square^k}$}.
A collection $\Big\{ H\0{\square^k} \Big\}\0{\square^k\in K}$ of normal charts
is a {\it normal atlas}, and if this atlas is smooth (or $C^k$) the induced
differentiable structure is called a {\it normal smooth (or $C^k$) structure.}\\

\noindent {\bf Proposition 9.3.1.} {\it The normal  atlas $\Big\{ H\0{\square^k} \Big\}\0{\square^k\in K}$ on $K_X$
is smooth.}\\

\noindent {\bf Proof.} 
Since we are assuming $\cA=\cA\big(\{h\0{\square}\}\big)$ smooth we get from proposition 7.8.5 that
the set  of smoothings $\{h\0{\square}\}$ is smoothly compatible, that is
the maps in (7.8.3) are smooth embeddings. For $\square^k\sbs\square^j\in K$
these maps have domains $\D^{n-j}\times \big(\dsquare^j\cap\sL(\square^k,K))\big)$ (the second factor is denoted
by $\dsigma\0{S}^j$ in 7.8) and target space $\bS^{n-k-1}$. We remark that
in this definition (and in section 7.8) we use the ``geometric" definition
of link, while here in section 9 we are using the ``direction" definition
of link (see 6.1.1). But using  (piecewise euclidean or piecewise hyperbolic ) exponential maps in $K$ or $K_X$ we can identify these definitions.
Therefore we can identify $\square^j\cap\sL(\square^k,K))$ with
$\dX_{\square^j}\cap\sL(X_{\square^k},K))$ (the links here are geometric).
This together with the fact that $\sL(\square,K)=\sL(X_\square,K_X)$
imply that we can obtain maps  in the $K_X$ case similar to the maps
in (7.8.3), and these maps have the same domains and target spaces.
Moreover they coincide modulo a slight smooth change (see remark below). Therefore
the $K_X$ versions of (7.8.3) are also ``smoothly compatible".
Now, the proof that
$\big\{ H\0{\square^k} \big\}$ is smooth is similar to the proof that
$\cA=\big\{ h^\bullet\0{\square^k} \big\}$ is smooth (assuming 
$\{h\0{\square}\}$ is smoothly compatible) given in the proof of proposition 7.8.5. This proves proposition 9.3.1.\\

\noindent {\bf Remark.} There is only one adjustment that has to be made in the
proof given in proposition 7.8.5 to be applied to the case of proposition 9.3.1.
In section 7.7 we identified $\rC\sL(\square^k,K)$ as a subset of 
$\sL(\square^j,K)$, $\square^k\sbs\square^j$, using the radial projection
$\Re$ described in remark 7.7.1. In the hyperbolic case (for proposition 9.3.1)
hyperbolic radial projection give a similar identification  (call it $\Re\0{\HH}$). Moreover, 
since ray structures are preserved these two projections coincide in
directions and just differ on the length. This length in the cube case
is given in remark 7.7.1. Using hyperbolic trigonometry the analogous
formula (using the same setting as in 7.7.1) is given by:
$tan^{-1}\Big(\frac{tanh\,(\,\,d\0{K}(v,p)\,\,)}{sinh\,(2r)}\Big)$, which is a 
smooth function. Therefore if  $\{h\0{\square}\}$ is smoothly compatible
using the identifications $\Re$, then  $\{h\0{\square}\}$ is also smoothly compatible
using the identifications $\Re\0{\HH}$.\\\\

We will denote by $\cS\0{K_X}=\cS\0{K_X}\Big( \big\{ h_{\square}\big\} \Big)$ the smooth structure on $K_X$ induced by 
the smooth atlas $\cA\0{K_X}=\Big\{ H\0{\square^k} \Big\}\0{\square^k\in K}$.
Note that $\cA\0{K_X}$ depends uniquely on the smoothly compatible set of
link smoothings $\cL=\{h\0{\square}\}\0{\square\in K}$ for $K$ (hence for $K_X$),
and to express this dependence we will sometimes write 
 $\cA\0{K_X}=\cA\0{K_X}(\cL)$.\\

\noindent {\bf Proposition 9.3.2.} {\it The map $\Phi:\big(K_X,\cS\0{K_X}\big)\ra K_X'$ is a $C^1$-diffeomorphism.}\\

The proof is given in appendix F.\\

Hence the atlas  $\Big\{ H\0{\square^k}\Big\}$ is a normal $C^1$-atlas for the smooth 
manifold $K_X'$. The following result was the main goal of section 9.\\

\noindent {\bf Proposition 9.3.3.} {\it The smooth manifolds $\big(K_X,\cS\0{K_X}\big)$ and $ K_X'$ are smoothly diffeomorphic.}\\

\noindent {\bf Proof.} Just approximate the $C^1$-diffeomorphism $\Phi$ by a smooth diffeomorphism.

\vspace{.8in}

\noindent {\bf 9.4. Normal structures for manifolds with codimension zero singularities.}\\
In this section we treat the case of manifolds with a one point singularity.
The case of manifolds with many (isolated) point singularities
is similar. The results in this section will be used in section 11.\\

We assume the setting and notation of section 7.9.
Let $K_X$ be the Charney-Davis strict hyperbolization of $K$.
Denote also by $p$ the singularity of $K_X$.
Many of the definitions and results given in section 9 still hold (with minor changes)
in the case of manifolds with a one point singularity:\\

\begin{enumerate}
\item[{\bf (1)}] Given a set of link smoothings for $K$ (hence for $K_X$) 
we also get a set of charts $H_\square$ as in section 9.3.
For the vertex $p$ we mean the cone map 
$H_p=\rC h_p:\rC N\ra \rC L\sbs K_X$. We will also denote the restriction
of $H_p$ to $\rC N-\{o\0{\rC N}\}$ by the same notation $H_p$.
As in item (3) of 7.9 here we are identifying $\rC N-\{o\0{\rC N}\}$ with $N\times (0,1]$ with the product smooth structure obtained from
\s{some} smooth structure ${\tilde{\cS}}_N$ on $N$.
As before $\{H_\square\}_{\square\in K}$ is a {\it normal atlas} for $K_X$
(or $K_X-\{p\}$). A normal atlas for $K-\{p\}$ induces a {\it normal smooth
structure on $K_X-\{p\}$.}
\item[{\bf (2)}]  Again we say that
the smooth atlas  $\{H_\square\}$ (or the induced smooth structure, or the set $\{h_\sigma\}$) 
is {\it correct with respect to $N$} if $\cS_N$ is diffeomorphic to ${\tilde{\cS}}_N$.
\item[{\bf (3)}] Let the set $\cL=\{h_\square\}_{\square\in K}$ 
induce a smooth structure on $K-\{p\}$, hence $\cL$ is smoothly compatible 
(see item (3) of 7.9). As proposition 9.3.1 we get that $\{H_\square\}_{\square\in K}$
is a smooth atlas on $K_X-\{p\}$ that induces a normal smooth
structure $\cS_{K_X}$ on $K_X-\{p\}$. Moreover, from Theorem 7.9.1 we get that $\cS_{K_X}$ is correct with respect to $\cS_N$ when $dim\, N\leq 4$ (always) or 
when $dim\, N>4$, provided $Wh(N)=0$. Note that in this case we can take
the domain $\rC N-\{o\0{\rC N}\}=N\times (0,1]$ of $H_p$ with smooth
product structure $\cS_N\times \cS_{(0,1]}$.
\item[{\bf (4)}] It can be verified that a version of proposition 9.3.3 also
holds in this case: $(K_X-\{p\},\cS_{K_X})$ embeds on $(K-\{p\},\cS')\times X$
with trivial normal bundle.
\end{enumerate}

\vspace{1in}

\noindent {\bf \large  Section 10. Proof of the Main Theorem.}\\

In section 2.2  the concept of hyperbolic extension was introduced and it was given
on a subset of hyperbolic space.
We extend next, in the obvious way, this concept to hyperbolic manifolds.\\

As in section 2, let $(N,h)$ be a complete Riemannian manifold with center
$o=o\0{N}$. Let $(P,\sigma\0{P})$ be a hyperbolic manifold.
The {\it hyperbolic extension of $h$  over $P$} is the Riemannian
metric $g=(cosh^2r)\sigma\0{P}+h$ on $P\times N$, where $r:N\ra[0,\infty)$ is the 
distance to $o$ function on $N$ (see 2.2). We write $g=\cE_{P}(h)$
and $(P\times N, g)=\cE_{P}(N,h)$ (or simply $\cE_P(N)$) and we call
$\cE_P(N)$ the {\it hyperbolic extension of $N$ over $P$}.\\
 
We now begin the proof of the Main Theorem. Let $M^n$ be a closed smooth
manifold. 
Let $K$ be a smooth cubification of $M$ (see appendix G) and 
$K_X$ the Charney-Davis strict hyperbolization of $M$, as in section 9.
We can assume that the Charney-Davis hyperbolization piece $X$ has normal
bundles with large widths (see lemma 9.1.1), all larger than a large number $s\0{0}>0$
and $s\0{0}>>3r$, where $r$ is as in 9.1. 
Let $\cA\0{K\0{X}}=\big\{  H_{\square} \big\}\0{\square\in K}$ be a
smooth normal atlas for $K_X$ and $\cS\0{K\0{X}}$ the induced
normal smooth structure on $K_X$ (see 9.3). Recall that the
charts $H_\square$ are constructed from a smoothly compatible set
of link smoothings $\{h\0{\square}\}\0{\square\in K}$ for the links of $K$ (hence also
for $K_X$, see sections 7.8 and 9.3). 
We write $\cL\0{K}=\{h\0{\square}\}\0{\square\in K}$.\\

\noindent {\bf Remark.} The domains of the charts $H_{\square^k}$ are 
the sets $\D^{n-k}\times\dX_{\square^k}$. But in this section, for notational purposes,
we will consider the rescaling of $H_{\square^k}$ given by the composition
$(u,p)\mapsto (\frac{1}{2r}u,p)\mapsto H_{\square^k}(\frac{1}{2r}u,p)$,
defined on $\D^{n-k}(s\0{0}/2r)\times\dX_{\square^k}$. We shall
denote this chart also by $H_{\square^k}$. That is, in this section $H_{\square^k}$
is the chart given by $H_{\square^k}(tv,p)=exp_p(t\,h_{\square^k}(v))$.\\

In what follows, to simplify our notation, we write $\sL(X_\square)=\sL(X_\square,K_X)$.
Recall that the set $\cL_K=\{h\0{\square}\}\0{\square\in K}$ of link smoothings for
the links $\sL(X_\square)=\sL(\square,K)$
of $K_X$ (and of $K$) induce a set of link smoothings 
$\{h\0{\square'}\in\cL_K\, ,\, \square'\subsetneq\square\}$ for the links of
$\sL(X_\square)$ just by restriction (see 7.7). We denote this induced set
of smoothings by $\cL_{\sL(X_\square)}$ or just $\cL\0{\square}$.\\

The space $K_X$ has a natural piecewise hyperbolic metric which we denote by
$\sigma\0{K_X}$. The piecewise hyperbolic metric on the cones
$\rC\sL(X_{\square})$ of the all-right spherical simplices $\sL(X_{\square})$
will be denoted by $\sigma\0{\rC\sL(X_{\square})}$.
Also the restriction of $\sigma\0{K_X}$ to the totally geodesic
space $X_\square$ shall be denoted by $\sigma\0{X_\square}$.\\

For $\square^k\in K$, the {\it (closed) normal neighborhood of $\dX_{\square^k}$ in $K_X$ of width $s<s\0{0}$} is the set $\s{N}\0{s}(\dX_{\square^k},K_X)=H_{\square^k}
\big(\D^{n-k}(s)\times \dX_{\square^k}\big)$. That is, it is the union of the images of
all geodesics of length $\leq s$ in each copy of $X$ containing $X_{\square^k}$,
that begin in (and are normal to) $\dX_{\square^k}$.
Similarly the {\it open
normal neighborhood of $\dX_{\square^k}$ 
of width $s<s\0{0}$}
is the set $\stackrel{\circ}{\s{N}}\0{s}(\dX_{\square^k},K_X)=H_{\square^k}
\big(\,int\,\D^{n-k}(s)\times \dX_{\square^k}\big)$.
Sometimes we will just write  $\s{N}\0{s}(\dX_{\square^k})=\s{N}\0{s}(\dX_{\square^k},K_X)$ and 
$\stackrel{\circ}{\s{N}}\0{s}(\dX_{\square^k})=\,
\stackrel{\circ}{\s{N}}\0{s}(\dX_{\square^k},K_X)$. Since the
normal bundles of the $X_\square$ are canonically trivial (see construction of $X$  in 9.1) we can canonically identify the neighborhood $\s{N}\0{s}(\dX_{\square})$
with $\dX_{\square^k}\times\rC_s\sL(X_{\square^k})$,
where  $\rC_s\sL(X_{\square^k})=\B_s\big( \rC\sL(X_{\square^k})  \big)$
is the closed $s$-cone  of length $s$, that is, it is the ball of radius $s$ on
the (infinite) cone $\rC\sL(X_{\square^k})$ centered at the vertex, see 6.5.
Similarly we have the identification $\stackrel{\circ}{\s{N}}\0{s}(\dX_{\square})= \dX_{\square^k}\times\stackrel{\circ}{\rC}_s\sL(X_{\square^k})$,
where  $\stackrel{\circ}{\rC}_s\sL(X_{\square^k})$
is the open $s$-cone  of length $s$.
Moreover these
identifications are also metric identifications, where we consider 
$\s{N}\0{s}(\dX_{\square},K_X)\sbs K_X$ with the (restricted) piecewise hyperbolic metric 
$\sigma\0{K_X}$ and $ \dX_{\square^k}\times\rC_s\sL(X_{\square^k})$
with the hyperbolic extension metric 
$$\cE_{\dX_{\square}}(\sigma\0{\rC\sL(X_\square)})=cosh^2(t)\,\sigma\0{\dX_{\square}}+
\sigma\0{\rC\sL(X_\square)}$$
\noindent where $t$ is the distance-to-the-vertex function on the cone
$\rC\sL(X_\square)$.\\

\noindent {\bf Remarks.} 

\noindent{\bf 1.} The metric $\sigma\0{\rC\sL(X_\square)}$  is not smooth
but the formula above makes sense, giving a well defined piecewise hyperbolic
metric.

\noindent {\bf 2.} Since we are identifying $\s{N}\0{s}(\dX_{\square})$
with $\dX_{\square^k}\times\rC_s\sL(X_{\square^k})$ we will consider
 $\s{N}\0{s}(\dX_{\square})$ also as a subset of
$\dX_{\square^k}\times\rC\sL(X_{\square^k})$, where
$\rC\sL(X_{\square^k})$ is  the (infinite) cone over $\sL(X_{\square^k})$.
Note that the metric $\cE_{\dX_{\square}}(\sigma\0{\rC\sL(X_\square)})$
is defined on the whole of $\dX_{\square^k}\times\rC\sL(X_{\square^k})$.\\

\noindent {\bf Lemma 10.1.} {\it Let 
 $\square^k=\square^i\cap\square^j$,  $k\geq 0$.
Let $s\0{1}, s\0{2}, s < s\0{0}$ be positive real numbers such that
$\frac{sinh\, s\0{1}}{sinh\, s},\frac{sinh\, s\0{2}}{sinh\, s} \leq\frac{\sqrt{2}}{2}$. Then} 
$$
\sN\0{s\0{1}}\big(\dX_{\square^i}  \big)\,\, \cap\,\,\sN\0{s\0{2}}\big(\dX_{\square^j}  \big)
\,\,\, \sbs\,\,\, \sN\0{s}\big(\dX_{\square^k}  \big)
$$

\noindent {\bf Proof.}
Let $p\in \sN\0{s\0{1}}\big(\dX_{\square^i}  \big)\, \cap\,\sN\0{s\0{2}}\big(\dX_{\square^j}  \big)$. Let $q$ be the closest point to $p$ in the totally geodesic
subspace $X_{\square^k}$. Denote by $t$ the distance between $p$ and $q$. Consider the cone $C=\rC\0{s\0{0}}\big(\sL(X_{\square^k})\big)$
at $q$, that is $C$ is the union of the images of all geodesics of length $s\0{0}$ normal
to $X_{\square^k}$ at $q$. Then $C$ is convex and $p\in C$. Write $A_l=C\cap X_{\square^l}$, $l=i,j$, which are also convex. Let $q\0{l}\in A_l$ be the closest
point to $p$ in $A_l$ and let $\gamma\0{l}$ be the geodesic segment between
$p$ and $q\0{l}$ with $length\, (\gamma\0{l})$ equal to the distance between
$p$ and $q\0{l}$. Hence $a\0{l}=length\, (\gamma\0{l})<s\0{l}$.
 Since $A_l$ is convex we have that $\gamma\0{l}$ is in $A_l$
and $\gamma\0{l}$ is perpendicular to $A_l$ at $q\0{l}$. Hence $q\0{l}$ is also the 
closest point to $p$ in $X_{\square^l}$. We get right triangles with vertices
$p$, $q$, $q\0{l}$ (right at $q\0{l}$), and hypotenuse equal to $t$. 
Let $\theta\0{l}$ be the angle at $p$.
Thus $\theta\0{l}$ is opposite to the side with length $a\0{l}$.
By the hyperbolic law of sines we have $sin\, \theta\0{l}=\frac{sinh\,a\0{l}}{sinh\, s'}$,
and by hypothesis we get
$$
sin\,\theta\0{l}\,\,=\,\, \frac{sinh\,a\0{l}}{sinh\, t}\,\,\leq\,\,
\frac{sinh\,s\0{l}}{sinh\, t}\,\,=\,\, \frac{sinh\,s\0{l}}{sinh\, s}\,
\frac{sinh\, s}{sinh\, t}    \,\,\leq \,\,\frac{\sqrt{2}}{2}\,\frac{sinh\,s}{sinh\, t}
$$

We want to prove that $t\leq s$. Suppose $t>s$. It follows then from the inequality
above that $sin\,\theta\0{l}<\sqrt{2}/2$, thus $\theta\0{l}<\pi/4$.\\

Let $S$ be the link of $X_{\square^k}$ at $q$. The segments
$[q,q\0{l}]$ intersect $S$ in two different vertices $v\0{l}$. Since 
$S$ is an all-right spherical complex and the sets $S\cap A_l$ are disjoint the
(angle) distance $d\0{S}(v\0{i},v\0{j})$ between $v\0{i}$ and $v\0{j}$
is a least $\pi/2$. Also the segment
$[p,q]$ intersects $S$ in a vertex $u$ and we have $\theta\0{l}=d\0{S}(u,v\0{l})$.
Consequently
$$
\frac{\pi}{2}\,\,\leq\,\,  d\0{S}(v\0{i},v\0{j})\,\,\leq\,\, 
d\0{S}(v\0{i},u)\,+\, d\0{S}(u,v\0{j})\,\, =\,\, \theta\0{i}\,+\,\theta\0{j}\,\,
< \frac{\pi}{4}\,+\,\frac{\pi}{4}\,\,=\,\, \frac{\pi}{2}
$$
\noindent which is a contradiction. This proves the lemma.\\

Suppose $\square^j\sbs\square^k\in K$. Then $\square^k$ determines
the all-right spherical simplex $\Delta\0{\sL(\square^j,K)}(\square^k)
=\square^k\cap\sL(\square^j,K)$ in $\sL(\square^j,K)=\sL(X_{\square^j})$.
We will just write $\Delta(\square^k)$ if there is no ambiguity.
(Other definition previously used: $\Delta\0{\sL(\square^j,K)}(\square^k)=
\sL(\square^j,\square^k)$.)
\\

\noindent {\bf Lemma 10.2.} 
{\it Let $\square^j\sbs \square^k$ and $s\0{1},\, s\0{2}<s\0{0}$.
Then}
$$
\sN\0{s\0{1}}(\dX_{\square^j})\,\,\cap\,\,\sN\0{s\0{2}}(\dX_{\square^k})\,\,=\,\,
\dX_{\square^j}\,\,\times\,\,\sN\0{s\0{2}}\Big(\,  \rC\Delta(\square^k), \,\rC\0{s\0{1}} \sL\big(  X_{\square^j} \big)   \,\Big)$$

\noindent {\bf Proof.} Let $p\in \sN\0{s\0{1}}(\dX_{\square^j})$. 
By the identification $\sN\0{s}(\dX_\square)=\dX_\square\times
\rC_s\big(\sL(X_{\square})\big)$ we can write 
$p=(u,x)\in \dX_\square^j\times\rC_{s\0{1}}\big(\sL(X_{\square^j})\big)$.
Since $\{x\}\times\rC_s\big(\sL(X_{\square^j})\big)$ is convex in
$\dX_\square^j\times\rC_s\big(\sL(X_{\square^j})\big)$,
we have that 
$$
d\0{\rC\0{s\0{1}}(\dX_\square^j\times\sL(X_{\square^j}))}\Big(\,p\,,\,
X_{\square^k} \, \Big)\,\,=\,\,
d\0{\rC\0{s\0{1}}(\{x\}\times\sL(X_{\square^j}))}\Big(\,u\,,\,
\rC\Delta(\square^j)\,  \Big)
$$

\noindent where, as usual, $d\0{S}$ denotes distance on a space $S$.
Consequently the first term in the equality above is $<s\0{2}$ if and only
if the second term  is $<s\0{2}$. This proves the lemma.\\

\noindent {\bf Remark.} Clearly the open  version of 
lemma 10.2 also holds:
$$
\stackrel{\circ}{\sN}\0{s\0{1}}(\dX_{\square^j})\,\,\cap\,\,\stackrel{\circ}{\sN}\0{s\0{2}}(\dX_{\square^k})\,\,=\,\,
\dX_{\square^j}\,\,\times\,\,
\stackrel{\circ}{\sN}\0{s\0{2}}\Big(\,  \rC\Delta(\square^j), \,
\stackrel{\circ}{\rC}\0{s\0{1}} \sL\big(  X_{\square^k} \big)   \,\Big)$$
\vspace{.2in}

Now, let $\s{d}$, $r$, $\xi$, $c$ and $\varsigma$ be as in items 1, 2, 3, 4 at the beginning of section 8, and let the numbers $s\0{m,k}=s\0{m,k}(r)$, $r\0{m,k}=r\0{m,k}(r)$
be as in section 6.6. For each  $\square^k\in K$ define the sets

$$\begin{array}{rcl}
\cZ(X_{\square^k})&=& \stackrel{\circ}{\s{N}}_{s\0{n,k}}\big( X_{\square^k} \big)\,\, -\,\, \bigcup\0{i<k}\sN_{r\0{n,i}}
\big(  X_{\square^i}\big)\\\\
\cZ&=& K_X -\,\, \bigcup\0{i<n-1}\sN_{r\0{n,i}}
\big(  X_{\square^i}\big)
\end{array}$$

\noindent By 9.1.1 we can take $s\0{0}$ as large as needed, hence we can assume that $\cZ(X_{\square^k})\sbs\,\,\stackrel{\circ}{\s{N}}_{s\0{0}}(X_{\square^k})$.\\ 

We next use the sets $\cX(P,\Delta,r)$ and $\cX(P,r)$ of section 6.6.
The sets $\cX(\sL(X_{\square^k}),\Delta(\square^j),
r)$ and  $\cX(\sL(X_{\square^k}),r)$
are a subsets of the (infinite) cone $\rC\sL(X_{\square^k})$.\\

\noindent {\bf Lemma 10.3.} {\it We have the following properties}
\begin{enumerate}
\item[(i)] {\it If $\square^i\cap\square^j=\emptyset$ then
$\cZ(X_{\square^i})\cap\cZ(X_{\square^j})=\emptyset$.}
\item[(ii)]  {\it If $\square^k=\square^i\cap\square^j$, $0\leq k<i,j$,  then
$\sN\0{s\0{n,i}}(X_{\square^i})\cap\sN\0{s\0{n,j}}(X_{\square^j})
\sbs \sN\0{r\0{n,k}}(X_{\square^k})$.}
\item[(iii)] {\it If $\square^k=\square^i\cap\square^j$, $0\leq k<i,j$,  then
$\cZ(X_{\square^i})\cap\cZ(X_{\square^j})=\emptyset$.}
\item[(iv)] {\it If $\square^j\sbs\square^k$ then we have (see remark 2 before 10.1)}
$$\cZ(X_{\square^j})\cap\cZ(X_{\square^k})\,\,\sbs\,\,\dX_{\square^j}\,\,\times\,\,
\cX\Big(\rC\sL(X_{\square^j}) ,\, \Delta(\square^k),\,r\Big) $$
\item[(v)] {\it For $k<n-1$ we have}
$$\cZ\cap\cZ(X_{\square^k})\,\,\sbs\,\,\dX_{\square^k} \,\,\times\,\,
\cX\Big(\rC\sL(X_{\square^k}) ,\,r\Big)
$$\end{enumerate}

\noindent {\bf Proof.} Let  $\square^i\cap\square^j=\emptyset$.
Then the distance in $K_X$ from $X_{\square^i}$ to $X_{\square^j}$ is 
at least $2s\0{0}$. This proves (1).\\

Statement (ii) follows from lemma 10.1, item (4) at the beginning
of section 8 and the following calculation
for $l=i,j$
(see 6.6 for the definition of $s\0{n,l}$ and $r\0{n,l}$)
$$
\frac{sinh\, s\0{n,l}}{sinh\, r\0{n,k}}\,\,=\,\, \frac{\Big(\frac{sinh\, r\,\, sin\,\beta\0{l}}
{sin\,\alpha\0{n-2}}\Big)}{\Big(\frac{sinh\, r}{sin\,\alpha\0{n-k-3}} \Big)  }
\,\,=\,\, c\,\varsigma^{l-k}\,\,\leq c\,\varsigma\,\,<\,\,e^{-(4+\xi)}
\,\,<\,\,\frac{\sqrt{2}}{2}
$$

Statement (iii) follows from (ii) and the definition of the sets $\cZ$.
We next prove (iv). Write $Z=\cZ(X_{\square^j})\cap\cZ(X_{\square^k})$.
By the definition of the sets $\cZ$ we have

$$
\begin{array}{ccl}Z&=&
\stackrel{\circ}{\sN}\0{s\0{n,j}}(X_{\square^j})\,\,\cap\,\stackrel{\circ}{\sN}\0{s\0{n,k}}(X_{\square^k})\,\,-\,\, \bigcup\0{l<k}\sN\0{r\0{n,l}}(X_{\square^l})\\
&\sbs&\stackrel{\circ}{\sN}\0{s\0{n,j}}(X_{\square^j})\,\,\cap\,\stackrel{\circ}{\sN}\0{s\0{n,k}}(X_{\square^k})\,\,-\,\, \bigcup\0{j\leq l<k}\sN\0{r\0{n,l}}(X_{\square^l})\\
&\sbs&\stackrel{\circ}{\sN}\0{s\0{n,j}}(X_{\square^j})\,\,\cap\,\stackrel{\circ}{\sN}\0{s\0{n,k}}(X_{\square^k})\,\,-\,\, \bigcup\0{j<l<k}\sN\0{r\0{n,l}}(X_{\square^l})
\,\,-\,\,\sN\0{r\0{n,j}}(X_{\square^j})\\
&\sbs&\stackrel{\circ}{\sN}\0{s\0{n,j}}(X_{\square^j})\,\,\cap\,\stackrel{\circ}{\sN}\0{s\0{n,k}}(X_{\square^k})\,\,-\,\, \bigcup\0{j<l<k}\Big(\sN\0{r\0{n,l}}(X_{\square^l})\,\cap\,  \sN\0{s\0{0}}(X_{\square^j})   \Big)
\,\,-\,\,\sN\0{r\0{n,j}}(X_{\square^j})
\end{array}$$

\noindent This together with lemma 10.2 imply
$Z\sbs \dX_{\square^j}\times A$ where
{\scriptsize $$
A\,\,=\,\,\stackrel{\circ}{\sN}\0{s\0{n,k}}\Big( \rC \Delta(\square^k),\,
\stackrel{\circ}{\rC}\0{s\0{n,j}}\big(\sL(X_{\square^j})\big)\Big)
\,\,-\,\,\bigcup\0{j<l<k}\sN\0{r\0{n,l}}\Big( \rC \Delta(\square^l),\,
\rC\0{s\0{0}}\big(\sL(X_{\square^j})\big)\Big)\,\,-\,\,\B_{r\0{n,j}}\big(\rC\sL(X_{\square^j})\big)
$$}
\noindent hence
{\small $$
A\,\,\sbs \,\,\stackrel{\circ}{\sN}\0{s\0{n,k}}\Big( \rC \Delta(\square^k),\,
\rC\sL(X_{\square^j})\Big)
\,\,-\,\,\bigcup\0{j<l<k}\sN\0{r\0{n,l}}\Big( \rC \Delta(\square^l),\,
\rC\sL(X_{\square^j})\Big)\,\,-\,\,\B_{r\0{n,j}}\big(\rC\sL(X_{\square^j})\big)
$$}
But for $i>j$ we have $s\0{n,i}=s\0{n-j,i-j}$, $r\0{n,i}=r\0{n-j,i-j}$
 and $r\0{n,j}=r\0{n-j-3}$ (see definitions in 6.6).
Therefore $A\sbs \cX\big(\rC\sL(X_{\square^j}) ,\, \Delta(\square^k),\,r\big)$. 
This proves (iv). The proof of (v) is similar to the proof of (iv) with
minor changes.  This proves lemma 10.3.\\ \\

We now smooth the metric $\sigma\0{K_X}$. 
For each $\square\in K$ using the construction in section 8 we get a Riemannian metric
$\cG\Big(\sL(X_\square), \cL_\square, h\0{\square}, r, \xi, \s{d}, (c, \varsigma)\Big)$
on $\sL(X_\square)$,
which we shall simply denote by $\cG\big(\sL(X_\square)\big)$. 
Define  the Riemannian metric $\cG(X_\square)$ on $\stackrel{\circ}{\s{N}}_{s\0{0}}\big(\dX_\square\big)$ by 

$$
\cG\big( X_\square\big)\,\,=\,\,\cE\0{\dX_\square}\Big(\, \cG\big(\sL(X_\square\big)  \, \Big) 
$$

\noindent {\bf Remark.} Recall that we can consider $\stackrel{\circ}{\s{N}}_{s\0{0}}\big(\dX_\square\big)$ contained in the infinite cone
$ \dX_\square\times\rC\sL(X_\square)$ (see remark 2 above), and note that
the definition of $\cG(X_\square)$ makes sense in the whole of
$ \dX_\square\times\rC\sL(X_\square)$.\\

\noindent {\bf Proposition 10.4.} {\it The Riemannian metrics 
$\cG(X_{\square^j})$ and $\cG(X_{\square^k})$ coincide on the intersection
$\cZ(X_{\square^j})\cap\cZ(X_{\square^k})$, $i,j<n-1$.
Also the Riemannian metric $\cG(X_{\square^k})$ coincides with $\sigma\0{K\0{X}}$ on the intersection $\cZ\cap\cZ(X_{\square^k})$.}\\

\noindent {\bf Proof.} For the first statement items (i) and (iii) of lemma 10.3 imply that we only need to
consider the case $\square^j\sbs\square^k$, $j<k<n-1$. 
By item (iv) of lemma 10.3 it is enough to prove that 
$\cG(X_{\square^j})$ and $\cG(X_{\square^k})$ coincide on
$\dX_{\square^j}\,\times\,\cX\Big(\rC\sL(X_{\square^j}) ,\, \Delta(\square^k),\,r\Big)$
(see remark above). Property {\bf P6} in 8.2 implies that the metric $\cG(X_{\square^j})$
coincides with the metric  $$\cE_{\dX_{\square^j}}\bigg[\cE_{\rC\Delta(\square^k)}\bigg( \,\,\cG\Big[\,\,\,\sL\Big(\,\Delta(\square^k)\,,\,\sL(X_{\square^j})\,\Big)\,\,\,\Big]\, \,  \bigg)  \bigg]$$
\noindent on
$\dX_{\square^j}\,\times\,\cX\Big(\rC\sL(X_{\square^j}) ,\, \Delta(\square^k),\,r\Big)$.
But\, $$\sL\big(\,\Delta(\square^k),\sL(X_{\square^j})\,\big)=
\sL\big(\,\Delta(\square^k),\sL(\square^j,K\,)\,\big)=\sL(\square^k,K)=\sL(X_{\square^k})$$
hence we have to prove that 
$$\cE_{\dX_{\square^j}}\big(\cE_{\rC\Delta(\square^k)}(  g)\,   \big) =\cE_{\dX_{\square^k}}(  g)$$
\noindent where $g= \cG\big(\sL(X_{\square^k}))$. And this follows 
from applying proposition 2.3.3 locally.\\

To prove the second statement in proposition 10.4,
using a similar argument as above (with {\bf P7} instead of {\bf P6})
we reduce the problem to showing that on $\dX_{\square^k}\,\times\,\rC\sL(X_{\square^k})$
we have
$\cE_{\dX\0{\square^k}}(\sigma\0{\rC\sL(X\0{\square^k})})=\sigma\0{K\0{X}}$.
And this follows from applying 6.5.2 (7) locally. This proves 
the proposition.\\

Finally define the metric $\cG(K_X)=\cG\big(K_X,\cL,r,\xi,\s{d}, (c,\varsigma)\big)$
to be equal to $\cG(X_{\square^k})$ on $\cZ(X_{\square^k})$, for $\square^k\in K$,
$k<n-1$. And equal to $\sigma\0{K\0{X}}$ on $\cZ$.
By lemmas 10.1 and 10.4 the metric $\cG(K_X)$ is a well defined Riemannian metric
on the smooth manifold $(K_X, \cS\0{K_X})$.\\

\noindent {\bf Corollary 10.5.} {\it Let $\epsilon>0$ and $M^n$  closed.
Choose $\xi$, $c$, $\varsigma$ satisfying (i) and (ii) in 8.4.2, and $\xi\geq n$.
Then the metric
$\cG(K_X)$ has all sectional curvatures $\epsilon$-pinched to -1,
provided  $d_i$, $r-d_i$, $i=2,...,n$, are sufficiently large.}\\

\noindent {\bf Proof.} 
Choose $\epsilon'$, as in 3.2, so that
an $\epsilon'$-hyperbolic manifold with charts of excess $\xi$ has
sectional curvatures $\epsilon$-pinched to -1.
Take $A$ and $a$ so that $A\geq C_3'(n,k,\xi)$ (see 3.5.8),
and  $a\geq a\0{k}(A\epsilon',\xi)$ (see 3.3.10), for $k\leq n$.
Since $M$ is compact we only have finitely many 
cubes in a cubification $K$ of $M$. Hence the set of links of $K$ (hence of $K_X$)
is finite. This together with proposition 8.4.2 imply that all 
$\cG(X_{\square^k})$, $\square^k\in K$, $k<n-1$, are $(\frac{\epsilon'}{A})$-hyperbolic
outside the balls of radius $a$, and are hyperbolic on the
ball of radius $2a$. All this provided $d_i$, $r-d_i$, $i=2,...,n$, are sufficiently large.
We can apply Theorem 3.5.8 (locally, see remark below) to get that
the metrics $\cG(X_{\square^k})$ are $\epsilon'$-hyperbolic outside 
$\cZ(X_{\square^k})\cap N_{a}(X_{\square^k})$, and hyperbolic on
$\cZ(X_{\square^k})\cap N_{2a}(X_{\square^k})$. This proves the corollary.\\

The corollary proves (i) of the main Theorem. Items (ii), (iii) follow from \cite{ChD}.
Item (iv) follows from  proposition 9.3.3 (recall $K_X'$ is a submanifold of 
$((K,\cS')\times X)$, see 9.2, 9.3).
This proves the Main Theorem.\\

\noindent {\bf Remark.} Note that it does not make sense to say that
$\cG(X_{\square^k})$ is $\epsilon'$-hyperbolic because neither $\dX_{\square^k}$
nor $\dX_{\square^k}\times \rC\sL(X_{\square^k})$ have a center. What we mean by
the ``local application of Theorem 3.5.8" mentioned in the proof above is the following.
Take $p\in \cZ(X_{\square^k})$ and let $B\sbs\dX_{\square^k}$ be an open ball
centered at $p$. Note that we can also consider $B\times\rC\sL(X_{\square^k})\sbs
\HH^{n-k}\times\rC\sL(X_{\square^k})=\cE_k(\rC\sL(X_{\square^k}))$ and we can now apply 3.5.8 to $\cE_k(\rC\sL(X_{\square^k}))$, where we are considering $p$
as the center.\\

\vspace{1in}

\noindent {\bf \large  Section 11. Proof of Theorem A.}\\

Let $N$ be a closed smooth manifold that bounds a compact smooth manifold
$M^m$. Denote the given smooth structure of $N$ by $\cS_N$.
Let $Q$ be the smooth $m$-manifold with one point singularity
formed by gluing the cone $\rC_1 N$ to $M$ along $N\sbs M$.
Let $q$ be the singularity of $Q$ and note that it is modeled on $\rC N$
(see 7.9). A triangulation of $Q$ is obtained by coning a smooth triangulation
of the manifold with boundary $M$, and let $f:K\ra Q$ be the induced
cubification (see appendix G). Write $f^{-1}(q)=p$. Note that $(K,f)$ is a smooth cubification of $Q$ in the sense of section 7.9. By item (2) of 7.9 we have that
$Q-\{q\}$ has a a normal smooth structure $\cS'$ for $K$, induced
by a set of links smoothings $\cL$.\\ 

Let $K_X$ be the Charney-Davis strict hyperbolization of $K$. Also denote by
$p$ the singularity of $K_X$. By item (1) of 9.4, the space $K_X-\{p\}$
has a normal smooth atlas $\{H_\square \}_{\square\in K}$ and normal 
smooth structure $\cS_{K_X}$. Moreover, since
we are assuming $Wh(N)=0$ (if $dim\, N>4$) we have that
we can take the domain $\rC N-\{o\0{\rC N}\}=N\times (0,1]$ of $H_p$
with product smooth structure $\cS_N\times\cS_{(0,1]}$ (see item (3) of 9.4).\\

We can now proceed exactly as in section 10 and define the sets
$\cZ(X_\square)$, $\cZ$, and the metrics $\cG(X_\square)$ depending on
$\cL,r,\xi,\s{d}, (c,\varsigma)$. For the special case $\square^0=p$
we use the results in section 8.5. We obtain in this way a
Riemannian metric $\cG(K_X)=\cG(K_X,\cL,r,\xi,\s{d}, (c,\varsigma))$ on $K_X-\{p\}$.
Theorem A and its addendum now follow from 8.5.1 (iii), (iv) and the result of
Belegradek and Kapovich \cite{BK} mentioned in the introduction (before the addendum
to Theorem A). To be able to apply 8.5.1 we need to satisfy the hypothesis
made at the beginning of 8.5: that the Whitehead group $Wh(\pi_1N)$ vanishes.
But this follows from \cite{FH}. This proves Theorem A.

\vspace{1in}

\noindent {\bf \large  Appendix A. Proof of lemma 3.5.5.}\\

Recall that we are considering $\HH^{k+1}$ with two sets of coordinates: the polar coordinates
$(x,t)$ and the $\cE_k(\R)$-coordinates $(y,r)$. We can consider then $x,t,y,r$ as functions defined on $\HH^{k+1}$, specifically:
$x:\HH^{k+1}-\{ o\}\ra \bS^k$, $y:\HH^{k+1}\ra\HH^k$,   $r:\HH^{k+1}\ra\R$,    $t:\HH^{k+1}\ra\R$.
Let $\p_r$ be the gradient vector field of $r$. Then the vectors $\p_r$ are the velocity vectors of the speed one geodesics 
emanating orthogonally from
$\HH^k\sbs \HH^{k+1}$. Also let $\p_t$ be the gradient vector field of $t$ and let $\alpha:\HH^{k+1}-\{ o\}\ra \R$ be the angle between 
$\p_t$ and $\p_r$. Then $\alpha(z)$ is the interior angle, at $z=(y,r)$, of the right triangle with vertices $o$, $y$, $z$. We call 
$\beta(z)$ the interior angle of this triangle at $o$, that is $\beta(z)=\beta(x)$ is the  (signed) spherical distance 
between $x\in \bS^k$ and the equator $\bS^{k-1}\sbs\bS^k$, where $(x,t)$ are the polar coordinates of $z$.
Note that the triangle mentioned above has sides of length $r=r(z)$, $t=t(z)$ and $a=a(y)$, where we are denoting by $a$ the distance
function in $\HH^k$ to $o$.
Using the hyperbolic law of cosines we get:

\begin{equation*}
sin\,\alpha\, =\, \frac{cos\, \beta}{cosh\, r}
\tag{A1}
\end{equation*}

\noindent Therefore

\begin{equation*}
\begin{array}{ccccc}
|sin\,\alpha |\leq \frac{1}{cosh\, r} &\,\,\,\,\,\, &&{\mbox{and}}&\,\,\,\,\,\,\,\,\,\,\,\,\,\,\,\,
|cos\, \alpha |\geq \frac{sinh\, r}{cosh\, r}\,\, =\,\, tanh\, r
\end{array}
\tag{A2}
\end{equation*} \

\noindent Note that the map $sin\, \beta$ is just the height function, i.e. $sin\,\beta (x)$ is the (signed) euclidean distance from 
$x\in\bS^k$ to $\bS^{k-1}$, which is  the last coordinate $x_{k+1}$ of $x=(x_1,...,x_{k+1})$. Therefore the term $sin\, \beta(x')$ that appears in
the definition of $\br$ (see 3.5.4) is the composition 

$$\B^k\stackrel{ {\tiny \frac{1}{sinh\, t_0}}}{\longrightarrow}\,\,\B^k\,\,
\stackrel{exp}{\longrightarrow}\,\,\bS^k\,\,\stackrel{proj}{\longrightarrow}\,\,\R$$

\noindent where the first arrow is multiplication by the constant $1/sinh\, t_0$ and the last arrow is the projection ``take the last coordinate",  and $exp=exp_{x_0}$\\

Write \,$\p_t=\pt$\,
and \,$\p_i=\frac{\p}{\p x_i}$, $i=1,...,k$. More generally, for $v\in\bS^k$, we write $\p_v=\frac{\p}{\p v}$.
Since $exp_{x_0}$ and $proj$ are smooth and the sphere is compact there is a constant $c$ (independent to $x_{0}$)
with $|\, proj\,\circ exp\, |_{C^2}\, <\, c$.\\

\noindent {\bf Remarks.}

\noindent {\bf 1.} The map $\beta$ is continuous but not smooth at the north pole, but  $sin\, \beta$ is smooth.

\noindent {\bf 2.} In what follows we will use the fact that we can take $c=8$. 
Moreover we can take $|\, proj\,\circ exp\, |_{C^1}\, \leq 1$.
A straightforward calculation 
(not given here) can show this.\\

\noindent We have then

\begin{equation*}
\begin{array}{ccccc}
\bigg|\,\, sin\, \big( \beta \, ( x')\big)\,\,\bigg|_{C^1}\,\,\leq \frac{1}{sinh\, t_0}&& {\mbox{and}}&&
\bigg|\,\, sin\, \big( \beta \, ( x')\big)\,\,\bigg|_{C^2}\,\,\leq \frac{8}{sinh^2\, t_0}
\end{array}
\tag{A3}
\end{equation*} \

\noindent Differentiating equation (3.5.4) we get (we write $\bt=t_0+t$)\vspace{.2in}

{\small  {\begin{center} $ \bigg| \pt \br(x,t)\bigg|\, =\, \bigg| \frac{cosh (\bt)\, sin\big( \beta(x')  \big)}{cosh\,\br}\bigg|\,\,
=\, \bigg| \frac{cosh (\br)\, cosh (a)\, sin\big( \beta(x')  \big)}{cosh\,\br}\bigg|\,\,
=\,\, cos\,\alpha \,\,\geq\,\, tanh\, \br  $\end{center}}}

\noindent where the second equality  is obtained from the first hyperbolic law of cosines and the last from the second hyperbolic law of
cosines, and the last inequality comes from  (A2). Note also that we get $|\p_t \br|\leq 1$.
Similarly, further differentiation and (A2)  (and use of the two laws of cosines plus the law of sines) shows

{\small {\begin{center}   $ \bigg| \frac{\p^2}{\p t^2} \br(x,t)\bigg|\, =
\, \bigg| \big( tanh\, \br\, \big)\, \big( sin^2\, \alpha\,\big)\bigg| \,\,\leq\,\, \frac{1}{cosh^2\, \br}  $\end{center}}}

\noindent Also, using (A3) we get

{\small {\begin{center}   $\bigg| \p_v \br(x,t)\bigg|\, =\, \bigg| \frac{sinh (\bt)\,\p_v sin\big( \beta(x')  \big)}{cosh\,\br}\bigg|\,\,
\leq \, \bigg(\frac{sinh\,(t+t_0)}{sinh\, (t_0)}\bigg)\,\,\frac{1}{cosh\, \br}\,\,
\leq \,\,\frac{2\,e^{^{1+\xi}}}{cosh\, \br}\,\,  $\end{center}}}

\noindent provided $t_0\geq ln\, 2$. Differentiating again we obtain

{\small {\begin{center}   $\bigg| \frac{\p^2}{\p t\p x_i} \br(x,t)\bigg|\,
\leq \, \bigg[\bigg(\frac{cosh\,(t+t_0)}{sinh\, (t_0)}\bigg)
\,\,+\,\,\bigg(\frac{sinh\,(t+t_0)}{sinh\, (t_0)}\bigg)\,\bigg]\,\,\frac{1}{cosh\, \br}\,\,
\leq \,\,\frac{4e^{^{1+\xi}}}{cosh\, \br}\,\,  $\end{center}}}

\noindent provided $t_0\geq ln\, 2$. Finally

{\small {\begin{center}   $ \bigg| \frac{\p^2}{\p x_j\p x_i} \br(x,t)\bigg|\,
\leq \, \bigg(\frac{sinh\,(t+t_0)}{sinh\, (t_0)}\bigg)\bigg[ \frac{8}{sinh\, t_0}\,\frac{1}{cosh\, \br}\,\,+\,\,
\frac{(1.02)e^{^{1+\xi}}}{cosh^2\,\br }\,\bigg]\,\,
\leq \,\,(2.05)\frac{e^{^{1+\xi}}}{cosh\, \br}+
(1.05)\Big(\frac{e^{^{1+\xi}}}{cosh\, \br}\Big)^2\,  $\end{center}}}

\noindent provided $t_0\geq\, 2$. Note that all three terms on the right
of the last three equations are less than $\frac{4\,e^{^{2(1+\xi)}}}{cosh\,\br}$.\\

Now, write $F(x,t)=\br (x,t)-(t+r\0{0})$. If we assume $r\0{0}>2$ we get $t_0>2>ln\,2$.
Since all but one of derivatives
of order 1 and 2 of $F$ concide with the ones of $\br$ we get that all such
derivatives are less than $\frac{4\,e^{^{2(1+\xi)}}}{cosh\,\br}$. The remaining derivative is
$\p_t F$. But we have
$$|\p_t F|\leq |1-tanh\,\br|=\frac{e^{-\br}}{cosh\,\br}\leq \frac{1}{cosh\,\br}<
\frac{4\,e^{^{2(1+\xi)}}}{cosh\,\br}$$
It remains to estimate $|F|_{C^0}$.
But for $x\in\B^k$, since $F(0,0)=0$, we have $F(x,0)=\int_0^1 x\,.\, \p _x F(tx,0)\, dt\leq \frac{2e^2}{cosh\, (\br)}$.
Hence $$\Big|\,F(x,t)\,\Big|\,=\,\Big|\,F(x,0)+\int_0^1 \p _t F(x,t)\, dt\,\Big|\,\leq\, \frac{2e^{^{1+\xi}}}{cosh\, (\br)}+\frac{1}{cosh\, (\br)}\,<\, \frac{4\,
e^{^{2(1+\xi)}}}{cosh\,\br}$$

\noindent To finish the proof we need the folowing claim.\\

\noindent {\bf Claim.} {\it We have  $\br>1+ln\,4 +\xi$,
provided $r\0{0}>5+2\,\xi$.}\\

\noindent {\bf Proof of claim.} Note that, since $r\0{0}=\br(0,0)$, we have
{\small $$\begin{array}{lll}
\Big| \, \br(x,t)\,-r\0{0}  \, \Big|&\leq&\Big| \, \br(x,t)\,-\,\br(x,0)  \, \Big|
\,+\,\Big| \, \br(x,0)\,-\,\br(0,0)  \, \Big|\\\\
&\leq&\int^t_0\big| \p\0{t}\br(x,t) \big|\, dt\,\,+\,\,
\int^1_0|x|\,\big| \p\0{x}\br(tx,0)\big|\,dt\\\\
&<& 1+\xi\,+\,\frac{2\,e^{^{1+\xi}}}{cosh\,\br}\,
\end{array}$$}
\noindent  This together with the fact that $\frac{1}{cosh\, \br}<
\frac{2}{e^{\br}}$ imply
\begin{equation*}
\br\,>\,
 r\0{0}\,-\, (1+\xi+\frac{4\, e^{^{1+\xi}}}{e^{\br}})\,
\tag{A4}
\end{equation*}

\noindent Now take 
\begin{equation*}
r\0{0}\,>\,5\,+\,2\,\xi
\tag{A5}
\end{equation*}

It follows from (A4) and (A5) that $\br\neq 1+\xi+ln\,4$.
\big[{\it Proof.} Plug $\br= 1+\xi+ln\,4$ and (A5) in (A4)
to obtain $1+ln\,4>3$ which is a contradiction.\big]
Therefore, since the domain of $\br$ is connected we get that
either $\br>1+ln\,4 +\xi$ or $\br<1+ln\,4+\xi$. But $r\0{0}$ is
in the image of $\br$, and from (A5) we get that $r\0{0}$
does no satisfy  $\br<1+ln\,4+\xi$ (recall $\xi>0$). Consecuently $\br>1+ln\,4 +\xi$.
This proves the claim.\\

\noindent From the claim it follows that $\frac{4\,e^{^{1+\xi}}}{
e^\br}<1$. This together with (A4) imply $\br>r\0{0}-(2+\xi)$,
provided $r\0{0}>5+2\xi$. Thus $\frac{1}{cosh\,\br}<\frac{1}{cosh\, (r\0{0}-
(2+\xi))}$.
This completes the proof of the lemma.
\vspace{.4in}

\noindent {\bf \large  Appendix B. Calculations for the proof of proposition 3.5.1.}\\

We will use the following abbreviations for the partial derivatives: $\p_t=\pt$,
 $\p_i=\frac{\p}{\p u_i}$, $\barp_i=\frac{\p}{\p v_i}$,
where $x_1=(u_1,...,u_k)$ and $x_2=(v_1,...,v_{n-1})$.\\

Write $\kappa=\kappa(r\0{0})=
\frac{4\,e^{^{2(1+\xi)}}}{cosh \big(\,r\0{0}-(2+\xi)\,\big)}$\,\, and\, \, $\zeta(x,t)=\br(x,t)-r\0{0}$.
A calculation shows that \\

\begin{equation*}  \kappa(r\0{0})<e^{-\big(r\0{0}-(7+3\xi)\big)}=
e^{7+3\xi}e^{-r\0{0}},\,\,\,\,
{\mbox{and}} \,\,\,\, {\mbox{for}}\,\,\, r\0{0}>7+3\xi \,\,\,\,
{\mbox{we get}}\,\,\,\kappa<1.
\tag{B1}
\end{equation*}
Write also 

$$c=c_{ij}(x_1,x_2,t)= a_{ij}(x_1, x_2,t)-e^t\delta_{ij}=b_{ij}(x_2,\zeta)-e^t\delta_{ij}$$

\noindent  and

$$d=d(x,r)\, =\, b_{ij}(x,r)\,-\, e^r\delta_{ij}$$\

 \noindent From lemma 3.5.3 and the fact that $\phi$ is $\epsilon$-hyperbolic we have

\begin{equation*}
\begin{array}{ccccc}
|\,\zeta\,\, -\,\, t\,|_{C^2}\,\leq \,\kappa &&{\mbox{and}}\,\,\,&& |d|_{C^2}\, <\, \epsilon\,\,\,\,
\end{array}
\tag{B2}
\end{equation*} \

\noindent From (B2)
we get that

\begin{equation*}
\big| \p_r b_{ij}   \big|_{C^1}\,\leq\, \epsilon\,+\, \big| e^r  \big|_{C^2 }\, <\, \epsilon\,  +\, e^{^{1+\xi}}
\tag{B3}
\end{equation*} \

\noindent Note that (B2) also implies

\begin{equation*}
\begin{array}{cccccccc}
\big| \p_i\zeta \big|_{C^0}\leq\kappa\,\,\,\,\,\,&&
\big| \p_t\zeta \big|_{C^0}\leq 1+\kappa\,\,\,\,\,&&
\big| \p^2_t\zeta \big|_{C^0}\leq\kappa\,\,\,\,\,&&
\barp_j\zeta=0
\end{array}
\tag{B4}
\end{equation*}\\

\noindent and (B1) implies

\begin{equation*}
\zeta<2+\xi
\tag{B5}
\end{equation*}

\noindent provided $r\0{0}>7+3\xi$.\\

\noindent {\bf The $C^0$-norm of $c$.}
Using (B2), (B3) and the Mean Value Theorem  we can write

$$\big| b_{ij}(x_2,\zeta)\,-\,b_{ij}(x_2, t)\big|\,\leq\, \big| \p _r b_{ij}\big|_{C^1}\, \big|\,  \zeta \,-\, t\,\big|  \, \leq\, \,\kappa\,\epsilon\,+\,\kappa \, e^{^{1+\xi}}\,
$$  \\

\noindent And this together with (B1) and (B2) imply 

$$\big| c  \big|_{C^0}\,\,\, \leq\, \,\,\kappa\,\epsilon\,+\,
 \kappa\, e^{^{1+\xi}}\, +\,   \big| b_{ij}(x,t)-e^t\delta_{ij}  \big|_{C^0} \,\,\,\leq\,\,\,\,\kappa\,\epsilon\,+\, \kappa\, e^{^{1+\xi}}\, +\,\epsilon \,\,<\,\,
e^{^{11+4\xi}}\, \Big( \epsilon\,+\, e^{-r\0{0}} \Big)$$ \\

\noindent {\bf The $C^1$-norm of $c$.} We have three types of first derivatives. First:

$$\p_t c\,=\,(\p_r d)\,(\p_t \zeta)\, +\, \big(   (\p_t \zeta)\,-\, 1    \big)\, e^\zeta\,I\, +\, (e^\zeta\,-\, e^t)\,I
$$

\noindent This last equation together with (B1), (B2), (B3), (B4)
and (B5) imply

$$|\,\p_t c\,|\, \leq\, \epsilon \, (1+\kappa)\,+\,  \kappa\, e^{^{2+\xi}}\, +\, \kappa\, e^{^{2+\xi}}
\,\,<\,\,
e^{^{11+4\xi}}\, \Big( \epsilon\,+\, e^{-r\0{0}} \Big)$$

\noindent where we are using the Mean Value Theorem, (B2)
and (B5) to estimate $e^\zeta-e^t$. Analogously

$$|\,\p_i c\,|\,=\,|(\,\p_r d)\,(\p_i \zeta)\, +\,    (\p_i \zeta)\, e^\zeta\, |\,\leq\, \epsilon\, \kappa\,+\, \kappa\,e^{^{2+\xi}}
\,\,<\,\,
e^{^{11+4\xi}}\, \Big( \epsilon\,+\, e^{-r\0{0}} \Big)$$

\noindent and

$$|\,\barp_i c\,|\,=\,|\,\barp_i d |\,<\, \epsilon
\,\,<\,\,
e^{^{11+4\xi}}\, \Big( \epsilon\,+\, e^{-r\0{0}} \Big)$$\\

\noindent {\bf The $C^2$-norm of $c$.} We have six types of first derivatives. As above using (B1),  (B2) and (B3)
we can obtain estimates for them. Here are the first three that do not involve the variable $x_2$:

{\small $$ \begin{array}{lll}
|\, \p^2_t c  \,|\,& =&\, \bigg|\,  (\p_r^2d)(\p_t\zeta)^2\,+\, (\p_rd)(\p_t^2\zeta)\, +\bigg[ (\p_t^2\zeta)\,
+\,(\p_t\zeta)^2\,-\,1  \bigg]\,e^\zeta\,
 +\, (e^\zeta\,-\,e^t  )\, \bigg|\\ \\
& \leq&
\epsilon\,(1+\kappa)^2\, +\, \epsilon\,\kappa\,+\big[  \kappa+(1+\kappa)^2-1 \big]\,e^{^{2+\xi}}\,+\,\kappa\,e^{^{2+\xi}}\\ \\
&<& e^{^{11+4\xi}}\, \Big( \epsilon\,+\, e^{-r\0{0}} \Big)
\\ \\ \\
|\,\p_i\p_t c  \,|\,& =&\, \bigg|\,  (\p_i\p_rd)(\p_i\zeta)(\p_t\zeta)\,+\, (\p_rd)(\p_i\p_t\zeta)\, +
\bigg[ (\p_i\p_t\zeta)\,+\,(\p_i\zeta)(\p_t\zeta)  \bigg]\,e^\zeta  \, \bigg|\\ \\
& \leq\,&
\epsilon\,(1+\kappa)\,\kappa\, +\, \epsilon\,\kappa\,+\big[  \kappa+(1+\kappa)\kappa \big]\,e^{^{2+\xi}}\\ \\
&<& e^{^{11+4\xi}}\, \Big( \epsilon\,+\, e^{-r\0{0}} \Big)

\\ \\ \\
|\,\p_j\p_i c  \,|\,& =&\, \bigg|\,  (\p_j\p_id)(\p_j\zeta)(\p_i\zeta)\,+\, (\p_id)(\p_j\p_i\zeta)\, +
\bigg[ (\p_j\p_i\zeta)\,+\,(\p_j\zeta)(\p_i\zeta)  \bigg]\,e^\zeta  \, \bigg|\\ \\
& \leq\,&
\epsilon\,\kappa^2\, +\, \epsilon\,\kappa\,+\big[  \kappa+\kappa^2 \big]\,e^{^{2+\xi}}\\ \\
&<& e^{^{11+4\xi}}\, \Big( \epsilon\,+\, e^{-r\0{0}} \Big)

\end{array}$$}

\noindent And the ones involving the $x_2$ variable:

{\small$$
\begin{array}{lllll}
|\, \barp_j\barp_i c\,|\,& =&\,| \, \barp_j\barp_i d\,|\,\,\,\,<\,\,\,\,\epsilon
\,\,<\,\, e^{^{11+4\xi}}\, \Big( \epsilon\,+\, e^{-r\0{0}} \Big)&&\\ \\
|\, \p_j\barp_i c\,|\,& =&\,| \, (\p_r\barp_i d)\,(\p_j\zeta)\,|\,\,<\,\,\epsilon\,\kappa
\,\,<\,\, e^{^{11+4\xi}}\, \Big( \epsilon\,+\, e^{-r\0{0}} \Big)\\ \\
|\, \p_t\barp_i c\,|\,& =&\,| \,( \p_r\barp_i d)\,(\p_t\zeta)\,|\,\,<\,\,\epsilon\, (1+\kappa)
\,\,<\,\, e^{^{11+4\xi}}\, \Big( \epsilon\,+\, e^{-r\0{0}} \Big)\end{array}$$}
\vspace{1in}

\noindent {\bf \large  Appendix C. Proof of Lemma 7.3 and 7.4.1(1)}\\

In the first part of this appendix we prove Lemma 7.3,  and in the final part
we prove 7.4.1 (1).
Recall that we are using the following definition: a non-degenerate $PD$ map $K\ra M$, $K$ a complex and $M$ a smooth manifold, is such that its restriction to each simplex (or cube) is smooth and
its derivative at every point is injective (\cite{MunkresLectures} p. 75).
Before we give the proof of lemma 7.3 we need to introduce a variation of the concept
of link.\\

\noindent {\bf C.1. Cubic Links.}\\ 
Let $\epsilon >0$ and write 
\begin{enumerate}
\item[] $I^k(\epsilon)=[0,\epsilon]^k\sbs\R^k$
\item[] $S^k(\epsilon)=\bS^{k-1}(\epsilon)\cap I^k(\epsilon)$
\item[] $C^k(\epsilon)=\D^{k}(\epsilon)\cap I^k(\epsilon)$
\item[] $F^k(\epsilon)=\{(x_1,...x_{k+1})\in\, ,\, x_i=1\, \,\,{\mbox{for some } i}\}$

\end{enumerate}

\noindent where $\bS^{k-1}(\epsilon)$ and $\D^k(\epsilon)$ are the $(k-1)$-sphere and the $k$-disc of radius $\epsilon$ respectively. 
Then $S^k_+(\epsilon)$ is a canonical all-right spherical
simplex of dimension $k$ and radius $\epsilon$ and $I^k(\epsilon)$ is the canonical $k$-cube of side-length $\epsilon$.
Note that, for $\epsilon\leq \epsilon' << \epsilon''$ we have  $C^k(\epsilon)\sbs I^k(\epsilon')\sbs C^k(\epsilon'')$.
In this situation by deforming along rays emanating from the origin we can construct a radial $PD$ self-homeomorphism
$\Theta=\Theta(\epsilon, \epsilon' , \epsilon'')$ on $C^k(\epsilon'')$ such that\, {\bf (1)} it is the identity on 
$C^k(\epsilon/2)$, and near $S^k(\epsilon'')$ \, {\bf (2)} sends $I^k(\epsilon')$ to $C^k(\epsilon)$. 
The inverse of $\Theta$ is also $PD$.\\

Let $K^k$ be a cubic complex, $\sigma^i\in K$ and $x\in\dsigma^i$. Write $j=k-i$.
Then the $\epsilon$-link $\sL_\epsilon(\dsigma^i,K)$
and $\epsilon$-star  $\rC\sL_\epsilon(\dsigma,K)$ of $\sigma^i$ at $x$ are built by gluing copies of
$S^j(\epsilon)$ and $C^j(\epsilon)$ respectively. The $\epsilon$-{\it cubic link}\,\, $\Box\sL_\epsilon(\dsigma,K)$ of
$\sigma^i$ at $x$ is obtained by gluing copies of $F^j(\epsilon)$ in the same way as the $S^j(\epsilon)$.
Analogously, the $\epsilon$-{\it cubic cone link}\,\, $\Box\rC\sL_\epsilon(\dsigma,K)$ of
$\sigma^i$ at $x$ is obtained by gluing copies of $I^j(\epsilon)$ in the same way as the $C^j(\epsilon)$.
Then we can write   $$\rC\sL_\epsilon(\dsigma,K)\,\,\sbs\,\, \Box\rC\sL_{\epsilon'}(\dsigma,K)\,\,\sbs\,\,\rC\sL_{\epsilon''}(\dsigma,K)$$
\noindent for $\epsilon\leq\epsilon' << \epsilon''$. By gluing the maps $\Theta$ mentioned above simplexwise we
obtain a radial $PD$ self-homeomorphism $\Theta$ (we use the same letter) on $\rC\sL_{\epsilon''}(\dsigma,K)$ 
with similar properties: {\bf (1)} it is the identity on 
$\rC\sL_{\epsilon/2}(\dsigma,K)$, and  near $\rC\sL_{\epsilon''}(\dsigma,K)$ 
\, {\bf (2)} sends $\Box\rC\sL _{\epsilon'}(\dsigma,K)$ to $\rC\sL_\epsilon(\dsigma,K)$. 
Again the inverse of $\Theta$ is also $PD$.\\

\noindent {\bf The spherical case.}\\
Consider $\bS^k\sbs\R^{k+1}=\R\times\R^{k-1}\times\R$.
We write $\R^{k-1}=\{0\}\times\R^{k-1}\times\{0\}$
and $\R^k=\{0\}\times \R^k$.
For $\theta\0{1}\in (-\pi/2,\pi/2)$ let $P_{\theta\0{1}}$ be the hyperplane
in $\R^{k+1}$ containing $\R^{k-1}$ and an making ``oriented" angle $\theta\0{1}$
with $\R^k$. The orientation is such that $\theta\0{1}$ is positive if and only if
the intersection of $P_{\theta\0{1}}$ with $\{x_{n+1}>0\}\sbs\R^{k+1}$ lies on
$\{x_1>0\}$. Let $v(\theta\0{1})$ be the unique upward (i.e. its $(n+1)$-coordinate
is non-negative) vector in $P_{\theta\0{1}}$ perpendicular to $\R^{k-1}$.
Then $v(\theta\0{1})=cos(\theta\0{1}) e_{n+1}+sin(\theta\0{1}) e_{1}$, where
the $e_i$ form the canonical bases of $\R^{k+1}$. We get a function
$\theta\0{1}:\{x_{n+1}>0\}\ra  (-\pi/2,\pi/2)$ given by
$\theta\0{1}(p)=\theta\0{1}$ if and only if $p\in P_{\theta\0{1}}$.
This function is the composition of the orthogonal projection $\R^{k+1}\ra \R\times\{0\}\times\R=\{x_i=0,\, i=2,...,n\}$, with angle function of polar
coordinates. Therefore $\theta\0{1}$ is smooth. We can repeat this construction
with other coordinates lines and obtain similar functions $\theta\0{i}$, $i=2,...,n$,
and we get a smooth map $\theta=(\theta\0{1},...,\theta\0{n})$.
We denote the restriction of $\theta$ to the upper hemisphere $\bS^k_+$ of
$\bS^k$ by the same letter $\theta$. Let $p=e_{n+1}\in\bS^k_+$ (the north pole).
Then $\theta(p)=0\in\R^k$. Let $\alpha\0{1}(s)=(sin(s),0,...,0,cos(s))\in\bS^k_+$
($s$ small). Then $\theta\0{1}(\alpha(s))=s$. It follows that $\theta$ is a 
diffeomorphism near $p$ and it has an inverse $\vartheta:\B_\mu(\R^k)\ra
\bS^k$, where $\B_\mu(\R^k)$ is the ball of (small) radius
$\mu$ in $\R^k$ centered at the origin.\\

If $K$ is an all-right spherical complex instead of a cubic one, similar concepts
can be defined as above, with
the $\epsilon$ objects replaced by their images in $\bS^k$ by $\vartheta$,
were we take $\epsilon$ small enough.

\vspace{.5in}

\noindent {\bf C.2. Proof of Lemma 7.3.}\\
We write $L=\sL(q^i,K)$. 
Then $L$ is an all-right spherical complex $PL$ equivalent to the four sphere $\bS^4$.
From section 7.6 we have that
the smooth structure $\cS_L$ on $f(L)$ given by lemma 7.2 has an atlas of the form
$$\cA_{f(L)}\, =\,\Big\{ \, \big(\,  h^\bullet_{\sigma^i}  \, ,\, \D^{4-i}\times\dsigma \,\big) \,\Big\}_{\sigma\in L}$$

Denote the smooth four manifold $(f(L),\cS_L)$ by $N$.
The identity map $L\ra N$ is not $PD$ (see remark 3 after Theorem 7.1).
To prove that $N$ is diffeomorphic to $\bS^4$ we will modify the identity
$L\ra N$ to a non-degenerate $PD$ homeomorphism $\phi:L\ra N$, which implies that
$\phi:L\ra N$ is a smooth triangulation of $N$.
 We will need the following Lemma.\\

\noindent {\bf Lemma C.2.1.} {\it  Let $\varphi:J\ra\bS^k$ be a smooth triangulation of
$\bS^k$ and $\psi:\bS^k\ra\bS^k$ homeomorphism which is non-degenerate $PD$ with respect to
$J$. Then $\varphi$ extends to a smooth triangulation $\varphi':J'\ra\D^{k+1}$ and $\psi$
extends to a homeomorphism $\psi':\D^{k+1}\ra\D^{k+1}$ which is 
non-degenerate $PD$ with respect to $J'$.}\\

\noindent {\bf Addendum to Lemma C.2.1.} {\it We can assume that $J'$ and $\psi'$
are radial outside a small ball centered at the origin}.\\

\noindent {\bf Proof.} We use lemma 9.8 of \cite{MunkresLectures}. Take
$K=J\times [\eta,1]$ and $f(x,t)=t\varphi(x)$, which is non-degenerate
$PD$ onto its image (hence a smooth triangulation of its image \cite{MunkresLectures}, p.77). Here $\eta>0$ is small.
Also take $K_1=J\times [\eta,2\eta]$ and $g: K_1\ra \R^{k+1}$ the secant approximation of $f|_{K'_1}$ (\cite{MunkresLectures}, p. 87) where
$K_1'$ is a subdivision on $K_1$ so that $f|_{K'_1}$ is an embedding
(\cite{MunkresLectures}, Th. 8.8).
Taking $\epsilon$ in 9.8 of
\cite{MunkresLectures} sufficiently small we get a triangulation $K'$ of
$J\times [\eta,1]$ and a non-degenerate $PD$ map $h$ on $K'$ which coincides
with $f$ outside (say) $J\times [0,3\eta]$ and is simplicial on $J\times [\eta,2\eta]$.
Take $J'=(\rC J) \,\coprod _{J=J\times\{\eta\}}K$ and extend $h$ simplicially on
$\rC J\sbs J'$ to obtain $h':J'\ra \D^{k+1}$. If $\epsilon$ is sufficiently small
$h'$ is non-degenerate. It can be checked that $J'$ depends only on the 
bounds of the derivatives of $\varphi$. Therefore we can assume that when we
apply the same argument to $\psi\circ\varphi:J\ra \bS^k$ we obtain a similar
map $h''$ defined also on the same $J'$, that is $h'':J'\ra\D^{k+1}$ is non-degenerate,
is radial outside a small neighborhood of the cone point of $J'$ and coincides
with $\psi\circ \varphi$ on $J\sbs J'$.
Finally take $\psi'=h''\circ (h')^{-1}$. This proves the lemma.\\

Now, the identity map $L\ra N$ is already $PD$ outside the two skeleton $L^2$ of
$L$. For $\sigma\in L$ and $\delta>0$ denote by $\sigma\0{\delta}$ the set of points
in $\sigma$ that lie at a distance $>\delta$ from $\p\sigma$. Hence 
$\sigma\0{\delta}\sbs\dsigma$.\\

Using the charts
 $\big(\,  h^\bullet_{\sigma^i}  \, ,\, \D^{4-i}\times\dsigma \,\big)$ we identify
$ \D^{4-i}\times\dsigma $ with its image by $ h^\bullet_{\sigma^i}$.
In particular, for $x\in\dsigma$, we are identifying $\D^{n-i}\times\{x\}$ 
with the $\rC\sL(\sigma,K)=\rC\sL_\epsilon(\dsigma,K)$
based at $x$, for some $\epsilon>0$. \\

We choose $\delta>0$ small and $\epsilon >0$ even smaller to make sure that
the open sets $\D^{2}\times\sigma\0{\delta}$, for all 2-simplices $\sigma^2\in L$, have
disjoint closures.
Let $\sigma^2\in L$ be a 2-simplex in $L$ and $x\in\sigma\0{\delta}^2$. Using lemma C.2.1
we can modify the identity $\D^2\times\{x\}=\rC\sL (\dsigma^2,L)\ra\D^2\times\{x\}$ 
near the origin $(0,x)\in \D^2\times\{x\}$ to make it
a non-degenerate $PD$ homeomorphism. Crossing with the identity $1_\dsigma$ we obtain a non-degenerate
$PD$ self-homeomorphism of $\D^2\times\sigma\0{\delta}^2$,
which is the identity outside a small neighborhood of $\sigma^2\0{\delta}$. Patching all these maps for all $\sigma^2\in L$
together with the identity we
obtain a non-degenerate $PD$ homeomorphism $\psi$  on $L(2)=(L-L^2)\cup\big(\bigcup _{\sigma^2\in L}\D^2\times\sigma^2\0{\delta}\big)$.
Note that $\psi$ is still a product map on each $\D^2\times\sigma^2\0{\delta}$.\\

Now let $\sigma^1\in L$ be a 1-simplex and $x\in\sigma^1\0{\delta'}$, where $\delta'$ is small but with
 $\delta<<\delta'$ (we may have to take $\delta$ and $\epsilon$ even smaller). Also let $\epsilon'$ be such that $\epsilon <<\epsilon'$ so that 
$\sL_{\epsilon'}(\dsigma^1,L)\sbs L(2)$ hence we can apply $\psi$ to   $\sL_{\epsilon'}(\dsigma^1,L)$.
We would like to apply lemma C.2.1 to extend $\psi$ near the $\sigma^1\0{\delta}$ but the problem is that
$\psi$ does not map $\sL_{\epsilon'}(\dsigma^1,L)$ to itself. To correct this we work with cubic links.
Note that we also get that $\Box\sL_{\epsilon'}(\dsigma^1,L)\sbs L(2)$. And, since 
$\psi$ is still a product map on each $\D^2\times\sigma^2\0{\delta}$ we have that
now $\psi$ maps $\Box\sL_{\epsilon'}(\dsigma^1,L)$ to itself. We can now apply lemma C.2.1 to the restriction of
$\Theta\,\psi\,\Theta^{-1}$ to $\sL_{\epsilon'}(\dsigma^1,L)$ to obtain a 
non-degenerate $PD$ extension $\psi'$. Take now
$\Theta^{-1}\, \psi'\,\Theta$ and cross this map with $1_{\sigma^1\0{\delta'}}$ to extend $\psi$ near $\sigma^1\0{\delta'}$.
Doing this for all $\sigma^1\in L$ and patching these maps with the previous $\psi$  we
obtain a non-degenerate $PD$ homeomorphism $\psi$  on $L(1)=L(2)\cup\big(\bigcup _{\sigma^1\in L}\D^2\times\sigma^1\0{\delta}\big)$.
Note that $\psi$ is still a product map on each $\D^2\times\sigma^2\0{\delta''}$ and each $\D^3\times\sigma^1\0{\delta'}$ (where $\delta''$ is slightly smaller than
$\delta$).
To extend $\psi$ to the whole of $L$ we proceed one step further (now for 0-simplices $\sigma^0$) in a similar way.
This proves lemma 7.3.\\

Variations of the argument used in the proof above give the following corollaries. 
For the first corollary we use the notation in Theorem 7.1 and its proof.\\

\noindent {\bf Corollary C.2.2.} {\it Let $f:K\ra (M^n,\cS)$ be as in Theorem 7.1.
Let $\cS'=\cS_5$ (as in the proof of 7.1). If \,$n=4$ then $(M,\cS')$ and $K$
are $PL$ homeomorphic (hence $\cS'$ is diffeomorphic
to $\cS$).}\\

\noindent {\bf Proof.} Replace spherical simplexes by cubes in the proof
of lemma 7.3 given above to obtain a non-degenerate $PD$ homeomorphism
$K\ra (M,\cS')$. Since $f:K\ra (M,\cS)$ is also a non-degenerate $PD$ homeomorphism
we get that  $(M,\cS')$ and  $(M,\cS)$ are $PL$ equivalent. But we
are in dimension four hence we can conclude that  $(M,\cS')$
$(M,\cS')$ and  $(M,\cS)$ are diffeomorphic. This proves the corollary.\\

We can also generalize the arguments in the proof of lemma 7.3 given
above to higher dimensions. 
Suppose $L$ is a (cube of all-right spherical) complex of dimension $n$. We can choose the $\epsilon$'s and 
$\delta$'s above properly and define 
$L(n-2)=(L-L^{n-2})\cup\big(\bigcup _{\sigma^{n-2}\in L}\D^2\times\sigma^{n-2}\0{\delta}\big)$ and
$L(k)=L(k-1)\cup\big(\bigcup _{\sigma^k\in L}\D^2\times\sigma^k\0{\delta}\big)$.\\

\noindent {\bf Corollary C.2.3.} {\it Let $f:K\ra (M^n,\cS)$ be as in Theorem 7.1.
Let $\cS_k$ be as in the proof of 7.1 satisfying condition
} {\bf S($k$)}. {\it Then there is a non-degenerate $PD$ homeomorphism
$\psi:L(k)\ra f(L(k))$.}\\

\noindent {\bf Proof.} Use exactly the same arguments and extend the induction 
in the proof of lemma 7.3. This proves the corollary.\\

We now prove statement (1) of proposition 7.4.\\\\

\noindent {\bf Proof of C($4$).}\\
The smooth structure $\cS_4$ is defined on $M-M_{n-4}$ and we have to prove
that it extends to a smooth structure $\cS_4'$ on $M$, and $\cS_4'$ is
diffeomorphic to $\cS$.\\

Write $A=M-M_{n-4}$.
Let $L(k)$ be as in corollary C.2.3 and write $B=L(n-3)\sbs K$, $D=f(B)\sbs M$. 
Note that $D\sbs A$. Here is an important
remark:\\

\noindent {\bf (C.2.4)}\hspace{.2in} {\it the inclusion $D\hookrightarrow A$ is a
homotopy equivalence.}\\

By corollary C.2.3 there is a non-degenerate $PD$ homeomorphism
$\psi:B\ra D$. We have that $\psi$ extends to a homeomorphism
(we use the same letter) $\psi:K\ra M$.\\

\noindent {\bf Remarks.}\\
\noindent {\bf 1.} The last statement follows from the inductive proof of corollary C.2.3, taking $k=0$. But we can also use a weaker version that does not use
C.2.1 and uses Alexander's trick beyond codimension 3.\\
\noindent {\bf 2.} We have to take a slightly ``smaller" set $B$ to obtain this
extension.\\

Write $f'=f|_B$ and note that, by hypothesis, 
$f':B\ra D$ is also a non-degenerate $PD$ homeomorphism.
Consider the smooth structures $\cS_\psi=\psi^*\cS_4$ and
$\cS_f=(f')^*\cS=(f^*\cS)|_B$ on $B$. The $PL$ structure on $B$
(induced by $K$) is thus Whitehead compatible with both differentiable
structures $\cS_\psi$ and $\cS_f$. But $A=M-M_{n-4}$ has the homotopy type
of a 3-complex hence, by C.2.3, so do $D$ and $B$. Since $PL/O$ is 6-connected
it follows from the theory of smoothings of $PL$-manifolds (see Theorem
4.2 in the second essay in \cite{HirschLectures}) that $\cS_\psi$ and $\cS_f$ are concordant.
Consequently, pushing forward everything to $D$ by $\psi$ we have that the smooth structures $\cS_4|_D$ and $\cS''=\psi_*\cS_f$ are concordant.   Note that
$\cS''=\Big((f\circ \psi^{-1})^*\cS\Big)_{D}$ hence we get that\\

\noindent {\bf (C.2.5)}\,\,\,\,\,\,\,
{\it  the smooth structures $\cS_4|_D$ and $\Big((f\circ \psi^{-1})^*\cS\Big)_{D}$ on $D$ are concordant}\\

\noindent It follows from C.2.4, C.2.5 and
the theory of smoothings of topological manifolds (see the Classification Theorem 10.1,
p .194, in \cite{KiSi} and its naturality for restrictions) that the
smooth structures $\cS_4$ and $\Big((f\circ \psi^{-1})^*\cS\Big)_{A}$ on $A$
are concordant. Therefore we can find a self-homeomorphism $g$ on $A$
such that $\Big((f\circ \psi^{-1})^*\cS\Big)_{A}=g^*\cS_4$.
Moreover we can assume (see Theorem 4.1, p. 25, in \cite{KiSi}) that
that $d\0{M}(x,g(x))\ra 0$, as $x\ra\p A$ (here $d\0{M}$ is any metric on $M$
inducing the topology on $M$). Therefore we can extend $g$ to a self-homeomorphism
on $M$ by defining $g(x)=x$, for $x\notin A$. To finish the proof just take
$\cS'_4=(f\circ \psi^{-1}\circ g^{-1})^*\cS$, which
is defined on the whole $M$ and extends $\cS_4$. This proves statement {\bf C}(4).

\vspace{1in}

\noindent {\bf \large  Appendix D. Proof of Proposition 9.2.1 and its Addendum.}\\

As always we write $I=[0,1]$.
Recall that the function $\bar{\rho}$ is defined as $\bar{\rho}(x_1,...,x_n)=
(\rho(x_1),...,\rho(x_n))$, where $\rho:I\ra I$ is as in section 9.2. We will
assume the following extra condition on $\rho$:\\

\noindent {\bf (D.1.)} \hspace{1.5in} $\rho(x)=x$ for $\delta\leq x\leq 1-\delta$\\

\noindent for some small $\delta>0$. 
Let $\square^{n-1}$ be an 
$(n-1)$- face of $\square^n=\{(x_1,...,x_n)\, ,\, 0\leq x_i\leq 1\}$. For simplicity  write $\square^n=\square^{n-1}\times I$, and
consider the vector field on $\square^n$, depending on $\square^{n-1}$, given by $V\0{\square^{n-1}}(x)=e_n=(0,...,0,1)$. This vector field is perpendicular
to $\square^{n-1}$ and generates the collar $\eta\0{\square^{n-1}}:\square^{n-1}\times I\ra\square^n$,  of $\square^{n-1}$ in $\square^n$
(which for the decomposition $\square^n=\square^{n-1}\times I$ is just the identity). \\

Let $\hat{\rho}$ be the smooth self-homeomorphism on $\square^{n-1}\times [0,\delta]$
given by $\hat{\rho}(x,t)=(x,\rho(t))$. 
Let $\Lambda\0{\square^{n-1}}$ to be the smooth self-homeomorphism on
$\square^n\ra\square^n$ that is the identity outside $\eta\0{\square^{n-1}}\big( \square^{n-1}\times [0,\delta) \big)$  and on 
the image of $\eta\0{\square^{n-1}}$ it is equal
to $\eta\0{\square^{n-1}}\circ\hat{\rho}
\circ\eta\0{\square^{n-1}}^{-1}$. Hence we can write\\

\noindent {\bf (D.2.)} \hspace{1.8in} 
$\bar{\rho}\,\,=\,\, \Lambda\0{\square_1^{n-1}}\circ...\circ \Lambda\0{\square^{n-1}_{2n}}$\\

\noindent for any ordering $\square_1^{n-1},...,\square^{n-1}_{2n}$ of all the
$(n-1)$-faces of $\square^n$.\\

We will assume that the width of the normal neighborhoods of the
$X_{\square}$ in $X$ are larger than $3r$ (see 9.1). Recall that by lemma 9.1.1
we can assume $r$ as large as we want.\\

\noindent {\bf Lemma D.3.} {\it For each $\square^{n-1}$ the vector field
$V\0{\square^{n-1}}$ has a lifting $W_{\square^{n-1}}$ to $X$ near $X_{\square^{n-1}}$.
Moreover $W_{\square^{n-1}}$
is perpendicular to $X_{\square^{n-1}}$.}\\

\noindent {\bf Remark.}
By $W_{\square^{n-1}}$ being a lifting of $V_{\square^{n-1}}$
 near $X_{\square^{n-1}}$ we mean that $W_{\square^{n-1}}$ is
defined on a normal neighborhood of $X_{\square^{n-1}}$ of width $\leq r$,
and $Df\,.\, W_{\square^{n-1}}=V_{\square^{n-1}}$.\\


Before we present the proof of lemma D.3 we show how it implies
proposition 9.2.1. The addendum to 9.2.1 will be proved later, at the end of
this appendix. There is an $s'$ such that all $W_{\square^{n-1}}$
are defined on the normal neighborhood of $X_{\square^{n-1}}$ of width $s'$.
Using the vector fields $W_{\square^{n-1}}$ we get collars
$\tau\0{\square^{n-1}}:X_{\square^{n-1}}\times[0,a]\ra X$, for some fixed $a>0$.
Since $W_{\square^{n-1}}$ is a lifting of $V_{\square^{n-1}}$
we get\\

\noindent {\bf (D.4.)}  \hspace{1.4in} 
$f\,\Big( \tau\0{\square^{n-1}}\big( x,t  \big) \Big)\,\,=\,\,
\eta\0{\square^{n-1}}\big( f(x),t  \big)$\\

\noindent For instance, in the special case of the trivial decomposition
$\square^n=\square^{n-1}\times I$ we get 
$f\,\big( \tau\0{\square^{n-1}}\big( x,t  \big) \big)=
\big( f(x),t  \big)$ because, in this case $\eta\0{\square^{n-1}}$
is just the identity.
Let  now $\theta\0{\square^{n-1}}$ be the smooth self-homeomorphism on
$X_{\square^{n-1}}\times[0,a]$ given by\\

\noindent {\bf (D.5.)}  \hspace{1.7in} 
$\theta\0{\square^{n-1}}\big( x,t  \big) \,\,=\,\,\big( x,\rho(t)  \big)$\\

Assuming $\delta>0$ in (D.1) such that $\delta<a$, we get that  $\theta\0{\square^{n-1}} $ is the identity outside 
$X_{\square^{n-1}}\times[0,\delta]\sbs X_{\square^{n-1}}\times[0,a)$.
Finally define $\Theta\0{\square^{n-1}}$ to be the the smooth self-homeomorphism on $X$ that is the identity outside $\tau\0{\square^{n-1}}\big( X_{\square^{n-1}}\times [0,\delta) \big)$  and on 
the image of $\tau\0{\square^{n-1}}$ is equal
to $\tau\0{\square^{n-1}}\circ\theta\0{\square^{n-1}}
\circ\tau\0{\square^{n-1}}^{-1}$.\\

\noindent {\bf Claim D.6.} {\it For every \,$\square^{n-1}$ we have that}\,\,\,
$ f\circ \Theta_{\square^{n-1}}\,\,=\,\,\Lambda_{\square^{n-1}}\circ f.$\\

\noindent {\bf Proof of claim.} By (D.4) we have that
$f\big( \tau\0{\square^{n-1}}\big(X_{\square^{n-1}}\times [0,\delta)\big) \big)=
\eta\0{\square^{n-1}}\big(\square^{n-1}\times [0,\delta)\big)$.
Hence a point $p\in X$ is in $\tau\0{\square^{n-1}}\big(X_{\square^{n-1}}\times [0,\delta)\big)$ if and only if its image $f(p)$ is in
$\eta\0{\square^{n-1}}\big(\square^{n-1}\times [0,\delta)\big)$
(see paragraph before 9.1.2).
If $p$ is not in $\tau\0{\square^{n-1}}\big(X_{\square^{n-1}}\times [0,\delta)\big)$ we get that $\Theta_{\square^{n-1}}(p)=p$ and $\Lambda_{\square^{n-1}}(f(p))=f(p)$ and the claim is true in this case. Assume now that $p$ is in
$\tau\0{\square^{n-1}}\big(X_{\square^{n-1}}\times [0,\delta)\big)$.
Write $\tau\0{\square^{n-1}}(x,t)=p$, thus $f(p)=\eta\0{\square^{n-1}}( f(x),t)$.
By applying (D.4) and (D.5) several times we get

$$\begin{array}{ccl}
f\circ \Theta\0{\square^{n-1}}(p)&=&f\circ \tau\0{\square^{n-1}}\circ\theta\0{\square^{n-1}}\circ\tau\0{\square^{n-1}}^{-1}(p)\\
&=&
f\circ \tau\0{\square^{n-1}}\circ\theta\0{\square^{n-1}}\big(x,t\big)\\ &=&
f\circ \tau\0{\square^{n-1}}\big(x,\rho(t)\big)\\ &=& \eta\0{\square^{n-1}}\big(f(x),\rho(t)\big)
\\ &=& \eta\0{\square^{n-1}}\circ\hat{\rho}\big(f(x),t\big)
\\ &=& \eta\0{\square^{n-1}}\circ\hat{\rho}\circ\eta\0{\square^{n-1}}^{-1}\circ f(p)\\
&=&\Lambda\0{\square^{n-1}}\circ f (p)
\end{array}
$$

\noindent This proves the claim.\\

To finish the proof of proposition 9.2.1 just define $P=\Theta_{\square^{n-1}_1}\circ...\circ\Theta_{\square^{n-1}_{2n}}$.
The fact that $f\circ P=P\circ f$ follows from (D.2) and claim D.6. This proves
proposition 9.2.1. \\

It remains to prove lemma D.3.\\

\noindent {\bf Proof of lemma D.3.}
Fix $\square^{n-1}$. Without loss of generality we assume $\square^{n-1}=
\square^{n-1}_1$, where $\square^{n-1}_i=\{\,x_i=0\,\}\cap\square^n$.
We have $\square^n=I\times\square^{n-1}_1$. Write
$V=V_{\square^{n-1}_1}$ and $W=W_{\square^{n-1}_1}$. 
Now, since the condition $Df.W=V$ is linear we have, using a partition of
unity
and taking $\delta$ small in (D.1), that it is enough to find
a lift of $V$ just locally, that is:\\

\begin{enumerate}\item[{\bf (D.7.)}] {\it  for every $p\in\square^{n-1}_1$
there is a neighborhood $U$ of $p$ in $X$ and vector field
$W$ on $U$ such that $Df.W=V$}\end{enumerate}


Let $\square^k\sbs\square^{n}$. Write
$D^j(s)=\rC_s\Delta_{\bS^j}=\{tu\in\R^{j+1},\,t\in[0,s],
\,u\in\Delta_{\bS^j}\}$.
We identify the closed normal neighborhood 
$N_s(\square^k)$ of $\square^k$ 
of width $s$
with $\square^k\times D^{n-k}(s)$ (here $s<1$).
Similarly we identify the closed normal neighborhood 
$N_s(X_{\square^k})$ of $X_{\square^k}$ 
of width $s$
(via the exponential map) with $X_{\square^k}\times D^{n-k}(s)$.
Note that for $\square^k\sbs\square^n$ 
we can write $\square^k=\bigcap\0{\square^k\sbs\square^{n-1}}\square^{n-1}$.
Define

$$A_s(\square^k)=\bigcap\0{\square^k\sbs\square^{n-1}}N_{s}\Big(\square^{n-1}  \Big)$$

 $$A_s(X_{\square^k})=\bigcap\0{\square^k\sbs\square^{n-1}}N_{s}\Big(X_{\square^{n-1} } \Big)$$
\noindent and  for $k<n$

$$
\begin{array}{ccl}L(X_{\square^k})&=& A_{3r}\big(\, X_{\square^k}   \,  \big)
\,\,\,-\,\,\, \bigcup\0{\square^k\not\subset\square^{n-1}}N_{2r}\Big(X_{\square^{n-1}}  \Big)\\
\\
&=&
\bigcap\0{\square^k\sbs\square^{n-1}}N_{3r}\Big(X_{\square^{n-1}}  \Big)\,\,\,-\,\,\, \bigcup\0{\square^k\not\subset\square^{n-1}}N_{2r}\Big(X_{\square^{n-1}}  \Big)
\end{array}$$\\

Note that $A_{s}(X_{\square^k})\sbs N_{s'}(X_{\square^k})$ for large $s'$
(how large $s'$ should be with respect to $s$ can be calculated using hyperbolic
trigonometry). Hence $L(X_{\square^k})\sbs N_{s}(X_{\square^k})$ for large $s$.\\


\noindent {\bf Claim D.8.} {\it We have 
 $X_{\square^{n-1}_1}\sbs\bigcup\0{\square^k\sbs \square^{n-1}_1}
L(X_{\square^k})$.}\\

\noindent If $p\notin L(X_{\square^{n-1}_1})$ then $p\in N_{2r}(X_{\square^{n-1}})$ for some
$\square^{n-1}$. Hence $p\in A_{2r}\big( X_{\square^{n-1}\cap\,\square^{n-1}_1}\big)$. Therefore we either have $p\in   L_{2r}\big( X_{\square^{n-1}\cap\square^{n-1}_1}\big)$, or $p\in N_{2r}(X_{\square_2^{n-1}})$,
for some $\square^{n-1}_2$ different from $\square^{n-1}_1$ and $\square^{n-1}$.
Arguing in the same way by induction we get that if $p\notin L(X_{\square^k})$
for all $\square^k\sbs\square^{n-1}_1$ with $k>0$ then $p\in A_{2r}(X_{\square^0})\sbs L(X_{\square^0})$, for some vertex $\square^0$. This proves claim D.8.\\

We now prove statement (D.7).
We use the construction of the map $f$ given in section 9.1.
Let $p\in \square^{n-1}_1$.
From claim D.8 we can assume that $p\in L(X_{\square^k})$, for some $\square^k\sbs\square^{n-1}_1$. Write $l=n-k$.
Note that  $L(X_{\square^k})\sbs N_{s}(X_{\square^k})=X_{\square^k}
\times D^l(s)$ (for large $s$), hence we will sometimes write $p=(p,0)\in X_{\square^k}
\times D^l(s)=N_{s}(X_{\square^k})$.\\

For simplicity we assume $\square^k=\square^{n-1}_{1}\cap...\cap\square^{n-1}_{l}$, $l=n-k$. Hence, using the notation in lemma 9.1.4, we have that 
$p\0{l}\circ f=(f_1,...,f_l)$ and $p\0{l}\circ T=(t_1,...,t_l)$.\\

\noindent {\bf Claim D.9.} 
{\it We have }
\begin{enumerate}
\item[{\it (a)}] {\it if $i >l$ then $f_i(p)=1/2$,}
\item[{\it (b)}] {\it if $(q,u)\in L(X_{\square^k})$ is
close to $p=(p,0)$, then $f_i(q,u)=1/2$, $i>l$,}
\item[{\it (c)}] {\it let $U=U'\times D\sbs X_{\square^k}\times D^l(s)$
be a product neighborhood where (b) holds for every $(q,u)\in U$.
Then $p\0{l}\circ T $ is an embedding on $ \{q\}\times D$, for every $q\in U'$.}
\end{enumerate}

\noindent Since $p\in L(\square^k)$ we have that $p\notin N_{2r}(X_{\square^{n-1}_i})$,
for $i>l$. Therefore $t_i(p)=d\0{X}(p,X_{\square^{n-1}_i})>2r>r$, $i>l$, and (a)
follows. Item (b) follows from (a), continuity and the fact that
the sets $N_s$ are closed. Item (c) follows from
lemma 9.1.4 (v). This proves claim D.9.\\

To finish the proof of (D.7) on $U$ just take $W(q,u)=\frac{1}{2r}\Big(\big( p\0{l}\circ T \big)|\0{\{      q\}\times D}\Big)^*(e_1)$, where $e_1$ is the constant vector field
$(1,0,...,0)$ on $\R^l$. (Note that $W$ is different from the gradient, with respect to
the hyperbolic metric on $X$, of the distance to $X_{\square^{n-1}_1}$ function $t_1$.)
It follows now from (b) of claim D.9 and the fact that $f_1(x)=\frac{1}{2r}\rho(t_1(x))=\frac{t_1(x)}{2r}$, if $x$ is close to $X_{\square^{n-1}_1}$, 
that $Df.W=e_1=V$. This proves (D.7). It can be verified from the
construction that the second statement of D.3 holds.  This proves lemma D.3 and completes the proof of proposition 9.2.1. \\

\noindent {\bf Proof of the Addendum to Proposition 9.2.1.}
Note that from the second statement in lemma D.3, (D.5) and the definition of
$\Theta\0{\square^{n-1}}$ we get\\

\noindent {\bf (D.10.)}\hspace{1.5in}
$D\Theta\0{\square^{n-1}}.W\0{\square^{n-1}}\,=\,0$\\

We need a lemma, which is is essentially an initial value version of claim D.3.\\

\noindent {\bf Lemma D.11.} {\it Let $U$ be a (not necessarily tangent) vector field on $X\0{\square^{n-1}}$. Suppose that
$Df.U=V\0{\square^{n-1}}$. Then there is a self-diffeomorphism
$g$ on $X$ covering the identity $1\0{\square^{n-1}}:\square^{n-1}\ra\square^{n-1}$
(see diagram) with $Dg.U=W\0{\square^{n-1}}$.}
$$
\begin{array}{ccc}
X&\stackrel{g}{\longrightarrow}&X\\
f\downarrow&&\downarrow f\\
\square^{n}&\stackrel{1\0{\square^n}}{\longrightarrow}&\square^n
\end{array}
$$

\noindent {\bf Proof.} Using collars and integral curves the problem is reduced to
finding and extension $U$ of $U$ to a neighborhood of $X\0{\square^{n-1}}$,
with $Df.U=V\0{\square^{n-1}}$ and $U=W\0{\square^{n-1}}$
outside an even smaller neighborhood of $X_{\square^{n-1}}$
(the argument uses the integral curves of $-U$).
The proof that such an extension exists is similar to that of claim D.3.
(without the perpendicularity condition). The only change needed is at the very end of the proof of (D.7) (after the proof of (D.9)). In our present case 
we have that $U(q)=U(q,0)=W(q,0)+T(q)$, where $T(q)$ is tangent
to $X\0{\square^{n-1}}$ (this is because $Df.U=V\0{\square^{n-1}}$ and $q\in L(X\0{\square^{k}})$).
Now take
$U(q,v)=W(q,u)+\rho(|v|)T(q)$, where $\rho(t)$ is equal to 1 near $t=0$ and equal to
0 for $t\geq \mu$, for some small $\mu>0$. This proves the lemma.\\

We now prove the addendum. 
Recall that at the beginning of appendix D we ordered the $(n-1)$-cubes: $\square^{n-1}_1,...
,\square^{n-1}_{2n}$, and we constructed the corresponding
$\Lambda_{\square^{n-1}_i}$, $\Theta_i=\Theta_{\square^{n-1}_i}$.
Write $V_i=V\0{\square^{n-1}_i}$
and  $W_i=W\0{\square^{n-1}_i}$.
We will need the following statement
which follows from the definition of the $\Lambda_{\square^{n-1}_i}$.\\

\noindent {\bf (D.12.)}\hspace{1.5in} $D\Lambda_{\square^{n-1}_i}.V_j=V_j$,
{\it \,\,\,for $i\neq j$}.\\

Take now $P=\Theta_{2n}\circ g\0{2n-1}\circ...\circ g\0{1}\circ\Theta_1$, where
the $g\0{i}$ are obtained in the following way. 
From claim D.6. (D.12) and the fact that $Df.W_i=V_i$ we get that $Df(D\Theta_1. W_2)=V_2$, hence we can apply lemma
D.11 to get a self-diffeomorphism $g\0{1}:X\ra X$ lifting the identity and 
satisfying $Dg\0{1}.(D\Theta_1.W_2)=W_2$. Next note that from
 D.6, (D.12), lemma D.11 and the fact that $Df.W_i=V_i$ we get that 
$Df(D(g\0{1}\circ \Theta_1). W_3)=V_3$ and 
we can apply lemma
D.11 to get a self-diffeomorphism $g\0{2}:X\ra X$ lifting the identity and 
satisfying $Dg\0{2}.(D(g\0{1}\circ\Theta_1).W_3)=W_3$, and so on.
From the choice of the $g\0{i}$ and (D.10) we get that $DP. W_i=0$.
Also from claim D.6 and the fact that $g\0{i}$ lifts the identity we
get $f\circ P=\bar{\rho}\circ f$. This proves the addendum to proposition 9.2.1.

\vspace{1in}

\noindent {\bf \large  Appendix E. Proof of Proposition 9.3.1.}\\

We shall demand the following condition on $\rho$: that the derivatives of $\rho$ {\it approach zero exponentially fast, at 0 and 1}. That is\\

\begin{enumerate}
\item[{\bf (E.1.)}]\hspace{.3in} {\it for every $k$ there are positive $a$ and $b$ such that $|\frac{d^k}{dt^k}\rho(t)|\leq a e^{-\frac{b}{(1-t)t}}$.}
\end{enumerate}\vspace{.2in}

  Let $\cA=\Big\{ \big(   h_{\square^i}^\bullet,\D^{n-i}
\times\dsquare^i    \big)  \Big\}$ be a normal atlas inducing $\cS'$.
We write $W_{\square^i}$ for the image of $h^\bullet_{\square^i}$. Note that
$W_{\square^i}$ is a normal neighborhood $\rC\sL (\square, K)\times\dsquare^i$
of $\dsquare^i$.
Write $c=(c_1,...,c_n)$.
We will prove that $\mu=\rho\circ c_1$ is smooth. The prove for $\rho\circ  c_i$ is the same.\\

We state three facts about the map $c_1$, which can be verified by inspecting
each of them cube by cube.

\begin{enumerate} 
\item[{\bf (1)}] There are three possibilities for a cube $\square\in K$:
First  $c_1(\square)=\{0\}$ and we say $\square$ is a {\it 0-valued-cube},
second $c_1(\square)=\{1\}$ and we say $\square$ is a {\it 1-valued-cube}
and finally $c_1|_{\square}$ is onto $I=[0,1]$ and we say in this last
case that $\square$ is an {\it I-cube}.
In what follows everything we do for 0-valued-cubes can be done for 1-valued-cubes, so we will
just ignore 1-valued-cubes.
\item[{\bf (2)}] For a 0-valued-cube $\square^i$ the map $c_1$ (and hence $\mu$ and all
its derivatives)
is a product map on a neighborhood of $\square^i$. Specifically $c_1$ factors
through a composition
$$
W_{\square^i}\,\,=\,\, \rC\sL (\dsquare^i,K)\times \square^i\,\,
\stackrel{{\mbox{{\tiny projection}}}}{\longrightarrow} \,\,\rC\sL (\dsquare^i,K)
\,\,\longrightarrow \,\, I
$$
\item[{\bf (3)}]
For a I-cube $\square^i$ the map $c_1$ (and hence $\mu$ and all its derivatives)
is a product map on a neighborhood of $\square^i$. Specifically $c_1$ factors
through a composition
$$
W_{\square^i}\,\,=\,\, \rC\sL (\dsquare^i,K)\times \dsquare^i\,\,
\stackrel{{\mbox{{\tiny projection}}}}{\longrightarrow} \,\,\square^i
\,\,\longrightarrow \,\, I
$$
\noindent where the last arrow is also a projection: $(x_1,...,x_n)\mapsto x_1$.
\end{enumerate}\vspace{.2in}

We prove that $\mu=\rho\circ c_1$ is smooth by showing that its representative
$\mu_{\square}=\mu\circ (h_{\square}^\bullet)^{-1}$ on each chart is smooth. We prove this by induction
on the decreasing dimension of the cubes. Consider first the following two statements
that depend on the $i$-cube $\square^i$:\\

\begin{enumerate}
\item[{\bf A($\square^i$):}]  We have that $\mu_{\square^i}$ is smooth on $\D^{n-i}\times\dsquare^i$. (Hence $\mu$ is smooth on $W_{\square^i}$.)
\item[{\bf B($\square^i$):}]  For every $\square^j<\square^i$, $\square^j$ a 0-valued-cube,
the map $\mu_{\square^i}$ and all its derivatives approach zero
exponentially fast with respect to the distance to $\square^j$. That is,
 for every $k$ there are positive $a$ and $b$ such than $|\frac{d^k}{dt^k}
\mu_{\square^i}(p)|\leq a e^{-\frac{b}{(1-t)t}}$, where $t=d\0{\D^{n-i}\times\square^i}(p,\square^j)$
\end{enumerate}

Recall that the chart maps $h_{\square^i}^\bullet$ respect the product $\D^{n-i}
\times\dsquare^i$ and the inclusion maps of cubes $\dsquare\ra (M,\cS')$ are
embeddings.
Therefore item (3) above implies that for an $I$-cube $\square^i$ the map 
$\mu_{\square^i}$ is just projection given by the composition
$$
\D^{n-i}\times\dsquare^i\stackrel{{\mbox{{\tiny projection}}}}{\longrightarrow}
\dsquare^i\stackrel{{\mbox{{\tiny projection}}}}{\longrightarrow}\,\,\, I\,\,\,
\stackrel{\rho}{\longrightarrow}\,\,\, I
$$
\noindent where the last projection is projection to the $x_1$ coordinate.
Since $\mu_{\square^i}$ is linear {\bf A($\square$)} is true for every $I$-cube $\square$.
Also, if $\square^i$ is an $I$-cube and $\square^j<\square^i$ is a 0-valued-cube,
we can write $\square^i=\square^{i-1}\times \square^1$, $\square^j<\square^{i-1}$,
where the projection on to the $x_1$ coordinate is $\square^{i-1}\times \square^1
\ra \square^1=I$. But for $p\in \D^{n-i}\times\dsquare^i$ we have
$$\mu_{\square^i}(p)\,\,\,\leq\,\,\,d\0{\D^{n-i}\times\square^i}(p,\square^{i-1})\,\,\,\leq\,\,\, d\0{\D^{n-i}\times\square^i}(p,\square^j)$$
\noindent and, since $\mu_{\square^i}$ is linear, it follows that {\bf B($\square^i$)} is also true for every $I$-cube $\square^i$.\\

For a 0-valued $(n-1)$-cube $\square^{n-1}$ it is straightforward to verify that
{\bf A($\square^{n-1}$)} and {\bf B($\square^{n-1}$)} hold true.
Assume now that {\bf A($\square^i$)} and {\bf B($\square^i$)} hold true for every
0-valued-cube $\square^i$, $i> k$. We prove the same is true for 0-valued-cubes $\square^k$.\\

Let $\square^k$ be a cube of dimension $k$.
Since $W_{\square^k}-\dsquare^k\sbs \bigcup_{i>k}W_{\square^i}$
it follows from the inductive hypothesis
{\bf A($\square^i$)}, $i>k$, that $\mu_{\square^k}$ is smooth on $\D^{n-k}\times\dsquare^k- \dsquare^k$, where we are writing $\dsquare^k=\{0\}\times \dsquare^k$.
By item (2) above $\mu_{\square^k}$ is a product on $\D^{n-k}\times\dsquare^k$,
hence,  it is enough to prove that the restriction $\nu=\mu\0{\square^k}|_{\D^{n-k}}:\D^{n-k}\ra I$
is smooth at  $0\in\D^{n-i}$. And
by the Mean Value Theorem we only need to prove that all partial
derivatives of $\nu:\D^{n-k}\ra I$ tent to zero as
a point $p$  tends to $0\in\D^{n-i}$.\\  

Let $v_m=t_mu_m\in\D^{n-k}$,  $u\in\bS^{n-k-1}=\p \D^{n-k}$, $t_m\in(0,1)$, $t_m\ra 0$.
We want to prove that all partial derivatives of $\nu$ at $v_m$ tend to zero
as $m\ra \infty$.
We can assume (arguing by contradiction) that  $u_m\ra u\in\bS^{n-k-1}$.\\

It was shown in the proof of Theorem 7.1 (see lemma 7.2) that the link $S=\sL(\square^k,K)$ is
a submanifold of $(M,\cS')$. The open sets 
$U_{\square^i}=S\cap W_{\square^i}$, $\square^i>\square^k$, form an open cover of $S$. 
Note that $U_{\square^i}$ is a normal neighborhood of $\square^i\cap S$
in $S$.
Write $x_m=h_{\square^k}^\bullet (v_m)$ and $y_m=h_{\square^k}^\bullet (u_m)\in S$. (Rigorously $h_{\square^k}^\bullet$ is not defined on $\p\D^{n-k}$ but,
after rescaling, we can assume this does happen).
Since $h_{\square^k}^\bullet$ restricted to $\D^{n-k}$ is, by definition, a cone map 
we can write $x_m=t_my_m$, where this last product is realized on the cone
link $\rC\sL (\square^k,K)$ of $\square^k$. And we also get $y_m\ra z=h_{\square^k}^\bullet
(u)$. \\

We
have that $z\in\dsquare^i\cap S$, for some $\square^i>\square^k$. Let $V$ be a small neighborhood of $z$ in $S$ with
$\bar{V}\sbs U_{\square^i}$ and we assume $y_m\in V$ for all $m$.
Write $$\nu=\mu\circ \big(h_{\square^k}^\bullet|_{\D^{n-i}}\big)= \Big(\mu\circ h_{\square^i}^\bullet\Big)\circ\Big( \big(h_{\square^i}^\bullet\big)^{-1}
\big(h_{\square^k}^\bullet|_{\D^{n-i}}\big)\Big)$$

By {\bf B($\square^i$)} all partial derivatives of the first term
$\mu\circ h_{\square^i}^\bullet$ approach zero exponentially fast as a point
get close to $\square^k$. Likewise, by corollary 7.8.2 the derivatives of the second term
$ \big(h_{\square^i}^\bullet\big)^{-1}
\big(h_{\square^k}^\bullet|_{\D^{n-i}}\big)$ grow at most polynomially fast.
Therefore, by applying the chain rule to the composition above we get that all partial
derivatives of $\nu$ tend to zero as $v_m\ra 0\in \D^{n-i}$. 
This proves {\bf A($\square^k$)}.\\

Note that the convergence of the derivatives of
$\nu$ to zero shown above is exponentially fast. This together with the fact that
(see item (2) above) the map $\mu_{\square^k}$ is a product on $\D^{n-i}\times\dsquare^k$
imply {\bf B($\square^k$)}. This proves the proposition.\\

\vspace{1in}

\noindent {\bf \large  Appendix F. Proof of Proposition 9.3.2.}\\

Recall that $\cA=\Big\{\big( h^\bullet_{\square^k} \,,\, \D^{n-k}\times\dsquare^k
  \big)   \Big\}$ is a normal atlas on $K$ (for $K$, see 7.1), that generates
the normal smooth structure $\cS'$. Also $\{H\0{\square}\}$ is a normal atlas for
$K_X=K_X^{{\mbox{\tiny piece-by-piece}}}$, generating the smooth structure $\cS\0{K\0{X}}$. We will assume that the charts $H\0{\square^k}:
\D^{n-k}\times\dsquare^k\ra K_X$ are defined on the larger sets
$\D^{n-k}(1+\delta)\times\dsquare^k$  (here $\D(1+\delta)$ is the open disc of
radius $1+\delta$). We can obtain this using 9.1.1.\\

Write $H'_{\square^k}=\Phi\circ H\0{\square^k}$.
It is enough to prove that the maps  $H'_{\square^k}:\D^{n-k}\times\dX\0{\square}
\ra K_X'$ are $C^1$-embeddings.
To prove this we need to prove that the following
coordinate maps are both $C^1$

$$
\begin{array}{llll}
q\0{K}\circ H'\0{\square^k}:&\D^{n-k}\times\dX\0{\square^k}&\longrightarrow&(K,\cS')\\
q\0{X}\circ H'\0{\square^k}:&\D^{n-k}\times\dX\0{\square^k}&\longrightarrow&X
\end{array}
$$

\noindent  We prove this by induction down the dimension of the skeleta.
First for $k=n$ recall that $H\0{\square^n}:\dX\0{\square^n}\ra K_X$ is just the inclusion. Hence proposition 9.2.2 implies that $q\0{K}\circ H'\0{\square^n}:\dX_{\square^n}\ra \dot{\square}^n$ is 
the map $\iota\circ f$, where $\iota:\dot{\square}^n\ra (K,\cS')$ is the inclusion, which is smooth.
(Recall that the inclusion $\square^n\ra (K,\cS')$ is not necessarily differentiable but its restriction $\iota$ to $\dsquare^n$ is smooth, see 7.1.)
Therefore $q\0{K}\circ H'\0{\square^n}$ is smooth.
Also, by the definition of the map $\Phi$, we have  $q\0{X}\circ H'\0{\square^n}=P$, which is also smooth. Moreover, by proposition 9.2.1,  $P|\0{\dot{\square}^n}$ is an embedding. Therefore $H'\0{\square^n}$ is a smooth embedding for every $n$-cube $\square^n\in K$.\\

Assume we have proved that  $H'\0{\square^j}$ is a $C^1$-embedding for  every $j$-cube $\square^j\in K$, $j>k$.
We have to prove that the same is true for all $k$-cubes.
We prove this in three parts. In the first part we prove that
$q\0{X}\circ H'\0{\square^k}$ is $C^1$. In the second part we prove that
$q\0{K}\circ H'\0{\square^k}$ is $C^1$. This two parts imply that
$H'\0{\square^k}$ is $C^1$. Finally, in the third part we prove that
$H'\0{\square^k}$ is an embedding.  Fix a $k$-cube $\square^k$.\\

\noindent {\bf FIRST PART.} {\it The map  $q\0{X}\circ H'\0{\square^k}$ is $C^1$.}\\

\noindent {\bf Proof.}
Denote by $V_{\square}$ the image of $H\0{\square}$. For each $\square^i$ with $\square^k<\square^i$ we have that on
$U_{\square^i}=(H\0{\square^k})^{-1}(V_{\square^i})$
we can write $$q\0{X}\circ H_{\square^k}'\,\,=\,\,\Big(q\0{X}\circ\Phi\circ H\0{\square^i}\Big)\circ\Big(H\0{\square^i}^{-1}\circ H\0{\square^k} \Big)
\,\,=\,\,\Big(q\0{X}\circ H'\0{\square^i}\Big)\circ\Big(H\0{\square^i}^{-1}\circ H\0{\square^k} \Big)
$$
\noindent  which is $C^1$ by inductive hypothesis and proposition
9.3.1. Since $\Big(\D^{n-k}-\{0\}\Big)\times\dX_{\square^k}$ is contained in the union of the $U_{\square^i}$, $i>k$, we have that
$q\0{X}\circ H'\0{\square^k}$ is $C^1$ outside $\dX_{\square^k}=\{0\}\times\dX_{\square^k}$.\\

Since the map $q\0{X}\circ \Phi|\0{\square^n}$ can be identified with the map $P:X\ra X$ for a $n$-cube $\square^n$
(recall $X\0{\square^n}$ is a copy of $X$), proposition 9.2.1 implies that the derivatives at a point $(0,p)\in \{0\}\times\dX_{\square^k}$
in the $X$ directions $(0,v)$ exist because $q\0{X}\circ H'\0{\square^i}$ 
on $\dX_{\square^k}=\{0\}\times\dX_{\square^k}$ is an embedding.
We next show that the derivatives in the radial directions also exist and vanish.
For this take a ray
$\alpha(t)= (tu,p)\in \D^{n-k}\times\dX_{\square^k}$ and write
$\beta(t)= H\0{\square^k}(\alpha(t))$. 
Note that $\beta'(0)=DH\0{\square^k}.u$ is normal to $X\0{\square^k}$.
We have that
the image of $\beta$ is contained in some $X_{\square^n}$.
As mentioned above the map $q\0{X}\circ \Phi$ on $X_{\square^n}$
can be identified with the map $P:X\ra X$. The fact that the radial derivative
in the direction $u$ exits and vanishes now follows from the addendum to 9.2.1.\\

Finally we need to prove that the first derivatives are continuous. But
this follows from a result analogous to lemma 7.8.6 with
$X_{\square^i}$ replacing $i$-cubes, which can easily be verified.
This concludes the proof of the first part.\\

\noindent {\bf SECOND PART.} {\it The map  $q\0{K}\circ H'\0{\square^k}:
\D^{n-k}\times\dX_{\square^k}\ra (K,\cS')$ is $C^1$.}\\

\noindent {\bf Proof.} This proof will take the next five pages.
First note that, by corollary 9.2.3 and the definition
of $H\0{\square}$ we have

\begin{equation*} q\0{K}\circ H'\0{\square^i}(\,t\,v\,,\, p\,)=F\,\Big( exp\0{p}\,\big(\,2\,r\,t\,\,h_{\square^i}(v)\,\big)\,\Big)
\tag{1}
\end{equation*}

\noindent Write $G_{\square^i}=\big(h^\bullet_{\square^i}\big)^{-1}\circ q\0{K}\circ H'_{\square^i}:\D^{n-i}\times\dX_{\square^i}\longrightarrow\D^{n-i}\times\dsquare^i$. Since $\{h^\bullet\0{\square}\}$ is an atlas for $(K,\cS')$,
by inductive hypothesis we have that $G_{\square^i}$ is $C^1$,
for $i>k$, and \\

\noindent {\bf (F.1.)}  \,{\it to prove that $q\0{K}\circ H'\0{\square^k}$ is $C^1$
 it is enough to prove that $G_{\square^k}$ is \, $C^1$.}\\

\noindent Write also $G_{\square}=(R_\square\, ,\, T_\square)$. For
$u=tv\in \D^{n-i}$, $t=|u|$, we have

\begin{equation*}
G_{\square^i}(\,t\,v\,,\, p\,)\,\,=\,\,\Big(R_{\square^i}(\,t\,v\,,\, p\,)\, ,\, T_{\square^i}(\,t\,v\,,\, p\,)  \Big)\,\,=\,\, \big(h^\bullet_{\square^i}\big)^{-1}
\circ F\,\Big( exp\0{p}\,\big(\,2\,r\,t\,\,h_{\square^i}(v)\,\big)\,\Big)
\tag{2}
\end{equation*}

\noindent 
It follows from (2) and lemma 9.1.4 (iii), (iv), that we can write

\begin{equation*}
G_{\square^i}(\,u\,,\, p\,)\,\,=\,\,\Big(R_{\square^i}(\,u\,)\, ,\, T_{\square^i}(\,|u|\,,\, p\,)  \Big)\tag{3}
\end{equation*}

\noindent  that is, $R$ does not depend on $p$ and $T$ depends on $p$ and the
length $|u|$ of $u$ (not on the direction of $u$). Also it can be checked from
9.1.4 (iv) and 9.1.5 that $(u, p)\mapsto T_{\square}(|u|, p)$ is smooth.
This together with (F.1) and (3) imply that\\

\noindent {\bf (F.2.)} {\it to prove that $q\0{K}\circ H'\0{\square^k}$ is $C^1$ it is enough to prove that
$R_{\square^k}:\D^{n-k}\ra\D^{n-k}$ is \, $C^1$.}\\

\noindent {\bf Claim F.3.} {\it The map\, $R_{\square^k}$ is $C^1$ on $\D^{n-k}-\{ 0\}$.}\\

\noindent {\bf Proof of claim F.3.} Recall that by inductive hypothesis we have that $G_{\square^i}$, $R_{\square^i}$ and $q\0{K}\circ H\0{\square^i}'$ are $C^1$, for all $\square^i$, $i>k$.
Denote by $V_{\square}$ the image of $H\0{\square}$. For each $\square^i$ with $\square^k<\square^i$ we have that on
$U_{\square^i}=(H\0{\square^k})^{-1}(V_{\square^i})$
we can write $$q\0{K}\circ H_{\square^k}'\,\,=\,\,\Big(q\0{k}\circ\Phi\circ H\0{\square^i}\Big)\circ\Big(H\0{\square^i}^{-1}\circ H\0{\square^k} \Big)
\,\,=\,\,\Big(q\0{k}\circ H'\0{\square^i}\Big)\circ\Big(H\0{\square^i}^{-1}\circ H\0{\square^k} \Big)
$$
\noindent  which is $C^1$ by inductive hypothesis and proposition
9.3.1. Since $\D^{n-k}-\{0\}=(\D^{n-k}-\{0\})\times\{p\}$ (for any $p\in
\dX\0{\square^k}$) is contained in the union of the $U_{\square^i}$, $i>k$, we have that
$q\0{K}\circ H'\0{\square^k}$ (hence $G_{\square^k}$ and $R_{\square^k}$) is $C^1$ outside $0$.  This proves claim F.3.\\

From  lemma 9.1.3 (ii) and the fact that the derivative of the exponential (at 0)
is the identity we get

\begin{equation*}
\frac{\p}{\p v}R_{\square^k}(0)=v
\tag{4}
\end{equation*}

\noindent That is, all directional derivatives at 0 of $R_{\square^k}$ exist and
if $R_{\square^k}$ were differentiable its derivative at 0 would be the identity matrix 1.
It follows from (4), (F.2) and claim F.3 that\\

\noindent {\bf (F.4.)} {\it to prove that $q\0{K}\circ H'\0{\square^k}$ is $C^1$ it suffices to prove
$DR_{\square^k}|\0{q}\ra 1$ (the identity matrix)  

\hspace{.21in} as $q\ra 0$.}\\

Write $S=\sL(\square^k,K)=\sL(X_{\square^k}, K_X)$ (at some point $F(p)\in\dsquare^k$ and
$p\in\dX\0{\square^k}$, respectively, and recall we are using ``direction" links).
Also write $\D^{n-k}=\D^{n-k}\times \{p\}\sbs \D^{n-k}\times \dX_{\square^k}$.
For $\square^j$, $\square^k\sbs\square^j$ set $\sigma\0{\square^j}=
\D^{n-k}\cap(H_{\square^K})^{-1}(X\0{\square^j})$ and
$\dsigma\0{\square^j}=
\D^{n-k}\cap(H_{\square^K})^{-1}(\dX\0{\square^j})$.
Note that the sets $\sigma\0{\square^j}$ and $\dsigma\0{\square^j}$ are cone sets.
That is, if $u\in\sigma\0{\square^j}$ then $tu\in\sigma\0{\square^j}$, $t\in [0,1]$.
Similarly for $\dsigma\0{\square^j}$ (see 9.1.4).\\

\noindent {\bf Claim F.5.} {\it We have 
 $R_{\square^k}(\sigma\0{\square^j})=\sigma\0{\square^j}$
\,\,and\,\, $R_{\square^k}(\dsigma\0{\square^j})=\dsigma\0{\square^j}$.}\\

\noindent {\bf Proof of claim F.5.} 
We prove the first identity, the second one is similar.
Let $u=tv\in \dsigma\0{\square^j}$, $|v|=1$. We assume $t>0$.
Then $exp\0{p}(2\,r\,t\,h\0{\square^k}(v))=H\0{\square^k}(u,p)\in X\0{\square^j}$. 
By (1) and lemma 9.1.4 (ii) we have that $h^\bullet\0{\square^k}\circ G\0{\square^k}(u,p)\in \square^j$. By the definition of  $h^\bullet\0{\square^k}$ we get
$\big(a\,h\0{\square^k}(\frac{1}{a}\,R_{\square^k}(u))\,,\,T_{\square^k}(|u|,p)\,\big)\in \square^j$, where $a$ is the length of $R_{\square^k}(u)$.
By 9.1.4 (ii) we get  $h\0{\square^k}(\frac{1}{a}\,R_{\square^k}(u))\in T_pX\0{\square^j}$, which implies $\frac{1}{a}\,R_{\square^k}(u)\in\sigma\0{\square^j}$. Since this set is a cone set it follows that
$R_{\square^k}(u)\in\sigma\0{\square^j}$. This proves claim F.5.\\

Now, let $q_n\ra0$ in $\D^{n-k}$. We can assume (arguing by contradiction) that 
$q_n= t_n u_n$, with $(t_n,u_n)\in\R^+ \times\bS^{n-k-1}$, $t_n\ra 0$, $u_n\ra u\in \dsigma\0{\square^j}$
for some $\square^j\in K$ containing $\square^k$. Hence:\\

\noindent{\bf  (F.6.)} {\it to prove that $q\0{K}\circ H'\0{\square^k}$ is $C^1$ it suffices to prove that
$DR_{\square^k}|\0{(t_nu_n)}\ra 1$, as $n\ra\infty$, 

\hspace{.26in}where
$u_n\ra u\in\dsigma\0{\square^j}$ and $t_n\ra 0$.}\\

\noindent {\bf Claim F.7.} {\it Statement (F.6) holds for  $j=n$.}\\
 
\noindent {\bf Proof of claim F.7.} We have that $u\in\dsigma\0{\square^n}$ for some $\square^n\in K$. Therefore there is a small compact neighborhood $V$
of $u\in\bS^{n-k-1}\sbs\D^{n-k}$ such that (we can assume that) all $q_n$ and $tu$, $t\in (0,1]$, 
lie in the interior of the cone $\rC V$. Denote by $h:\D^{n-k}\ra\rC S$ the map $\rC h_{\square^k}$, where
$h_{\square^k}$ is the link smoothing of $S=\sL(\square^k,K)$. Since $F|\0{\square^n}=f$ on $\rC V$ we can write
$R_{\square^k}=h^{-1}\circ (\pi\circ f\circ e)\circ h$, where $e$ is the exponential map given by
$e(v)=exp\0{p}(2r v)$, and $\pi:\rC S\times \square^k\ra \rC S$ is the
projection. (Note that $e(tv)=E(2rt,v)$, where $E$ is as in 9.1.4.) Hence $(DR_{\square^k})|\0{q_n}=(Dh|\0{y_n})^{-1}\, D(\pi\circ f\circ e)|\0{h(q_n)}\, Dh|\0{q_n}$, where 
$y_n=h^{-1}(\pi\circ f\circ e)(h(q_n))$. By lemma 9.1.3 (ii) and the fact that
the derivative of the exponential at 0 is the identity we have that $ D(\pi\circ f\circ e)|\0{h(q_n)}\ra 1$, as $q_n\ra 0$.
On the other hand, since $h$ is a cone map, by lemma 7.8.1 we get that  $Dh$ and $Dh^{-1}$ are both bounded on $\rC V$. Moreover
$Dh|\0{q_n}=Dh|\0{u_n}$ (see remark after 7.8.1) and $Dh|\0{y_n}=Dh|\0{\frac{y_n}{|y_n|}}$.
But, since $\pi\circ f\circ e$ is smooth and $D(\pi\circ f\circ e)|\0{p}=1$ we have that for any $v$ we get

{\small \begin{equation*}
\frac{\pi\circ f\circ e\,\big(\, t_n\, v\big)}{t_n}\,\,\longrightarrow\,\, v
\tag{5}
\end{equation*}}

\noindent From the fact that $h^{-1}$ is a cone map and (5) we have

{\small $$\frac{t_n}{|\,\,\,h^{-1}\Big((\pi\circ f\circ e)\,\big(t_n\,h\, (  u_n) \big)  \Big)\,\,\,  |}=
\Big[\,\,|\,h^{-1}\Big(\frac{(\pi\circ f\circ e)\,\big(t_n\,h\, (  u_n)\big) }{t_n} \,\, \Big)\,\,|\,\,\Big]^{-1}\,\,\longrightarrow\,
lim\0{n\ra\infty}\,|u_n|^{-1}\,\,=\,\,1 $$}

\noindent  This together with the fact the $h$ and $h^{-1}$ are cone maps imply
{\small $$\begin{array}{lll}
\frac{y_n}{|y_n|}&=&\frac{h^{-1}\Big((\pi\circ f\circ e)\,h\, \big( t_n\, u_n \big)  \Big)}
{|\,\,\,h^{-1}\Big((\pi\circ f\circ e)\,h\, \big( t_n\, u_n \big)  \Big) \,\,\, |}\,\,=\,\,
h^{-1}\bigg(\,\,
\frac{\Big((\pi\circ f\circ e)\,h\, \big( t_n\, u_n \big)  \Big)}
{|\,\,\,h^{-1}\Big((\pi\circ f\circ e)\,h\, \big( t_n\, u_n \big)  \Big)\,\,\,  |} \,\, \bigg)\\\\
&=&h^{-1}\bigg(\,\,
\frac{\Big((\pi\circ f\circ e)\,h\, \big( t_n\, u_n \big)  \Big)}
{t_n} \frac{t_n}{|\,\,\,h^{-1}\Big((\pi\circ f\circ e)\,\big(t_n\,h\, (  u_n) \big)  \Big)\,\,\,  |} \,\,\bigg)
\end{array}$$}

\noindent  Consequently
$\frac{y_n}{|y_n|}\ra u$. Therefore $DR_{\square^k}|\0{q_n}\ra 1$. This proves claim F.7.\\

We will prove statement (F.6) by decreasing induction on $j$. Claim F.7 was
the first step of this induction. We assume statement F.6 holds for all
$\square^i$, $j<i$.\\

\noindent {\bf Claim F.8.} {\it Statement (F.6) holds for  $j$.}\\
  
\noindent {\bf Proof of claim F.8.} 
The proof has three steps.\\

\noindent {\it Step 1. It is enough to assume that $u_n=u$. }\\
As in the proof of claim F.7 let $V$ be 
a small compact neighborhood
of $u\in\bS^{n-k-1}\sbs\D^{n-k}$ such that all $q_n$ and $tu$, $t\in (0,1]$, 
lie in the interior of the cone $\rC V$. From the definition of $G_\square$
we have that on $\rC V$ we can  write

$$
G_{\square^k}\,\,=\,\, \Big( \,\big( h^\bullet_{\square^j}\big)^{-1}\circ\,\, h^\bullet_{\square^k}   \Big)^{-1}\,\,
\circ\,\, G_{\square^j}\,\,\circ\,\, \Big( \,\big( H_{\square^j}\big)^{-1}\circ H_{\square^k}   \Big)
$$

To simplify the notation write $h= \big( h^\bullet_{\square^j}\big)^{-1}\circ\, h^\bullet_{\square^k}$ and $H=\big( H_{\square^j}\big)^{-1}\circ H_{\square^k}$. Hence 

\begin{equation*}DG_{\square^k}=Dh^{-1}.DG_{\square^j}.DH
\tag{6}
\end{equation*}

We next compare $DG_{\square^k}|\0{(t\0{n} u\0{n},p)}$ and $DG_{\square^k}|\0{(t\0{n} u,p)}$. We analyze the three terms $DH$, $DG_{\square^j}$, $Dh$ in (6).\\

\noindent {\bf First term:  $DH$.} Since $H$ is a cone map we get that
$DH|\0{t\0{n}u_{n}}-DH|\0{t\0{n}u}\ra 0$.\\

\noindent {\bf Remark.} The map $H$ is an euclidean-to-hyperbolic cone map, and it
is not an euclidean cone map but it is an euclidean
cone map up to a smooth change of coordinates on a compact set.\\

\noindent {\bf Second term:  $DG_{\square^j}.$}
Differentiating (3) we get

\begin{equation*} DG_{\square^j}|\0{(u,y)} (v,w)=\Big(\,  DR_{\square^j}|\0{u}\,.\,v\,\,,\,\, 
\frac{\p}{\p t}T_{\square^j}|\0{(t,y)}\frac{u.v}{|u|}\,+\, \frac{\p}{\p y}T_{\square^j}|\0{(t,y)} .w  \,\Big)
\tag{7}
\end{equation*}

\noindent where $t=|u|$. Since $T_{\square^j}$ can be extended to a smooth map
on $\bar{\D}^{n-k}\times X_{\square^j}$ (which is compact) 
the second term (i.e the $T_\square$ term) is  lipschitz on the variables
$t$ and $y$. Since the distance between 
$H(t_nu_n,p)$ and $H(t_n u,p)$   goes to zero it follows that the $T_\square$ terms
in the right hand side of the equation above evaluated at $H(t_nu_n,p)$ and $H(t_n u,p)$
get close as $n\ra\infty$. Also, by inductive hypothesis, the first terms tend both to 1.
Therefore we get
$$
DG_{\square}|\0{H(t\0{n}u\0{n})}-DG_{\square}|\0{H(t\0{n}u)}\longrightarrow\,\,
\,\,0
$$

\noindent as $n\ra \infty$. \\

\noindent {\bf Third term: $Dh$.}
Since $h$ is a cone map to prove that $Dh\0{G_{\square^j}(H(t\0{n}u\0{n}))}$
and $Dh_{G\0{\square^j}(H(t\0{n}u))}$ are close we need to prove that
the directions of
$G_{\square^j}(H(t\0{n}u\0{n}))$ and $G_{\square^j}(H(t\0{n}u))$
 are close. This is equivalent to proving that the directions of their images by $h$ are close. Since $G_{\square^k}=h\circ G_{\square^j}\circ H$ this means
proving that the directions of $G_{\square^k}(t\0{n}u\0{n})$ and $G_{\square^k}(t\0{n}u)$ are close.  
Let $E$ and $p\0{l}$ be the maps in lemma 9.1.4.
We can assume (arguing by contradiction) that $t_nu_n$, $t_n u$ lie on
$ X_{\square^n}$, for some $\square^n$.
Since $h^\bullet_{\square^k}$ is also
a cone map (on the first variable) it is enough to prove that the directions of $p\0{l}\circ f\circ E(2r t_n, h_{\square^k}(u_n))$
and $p\0{l}\circ f\circ E(2rt_n, h_{\square^k}(u))$ are close.
But this is true because $p\0{l}\circ f\circ E$ is smooth.
This concludes step 1.\\

\noindent{\it Step 2. We prove that
 $DR_{\square^k}|\0{(t_nu)}\ra 1$, as $n\ra\infty$, where
$ u\in\dsigma\0{\square^j}$ and $t_n\ra 0$.}\\

Note that
every inclusion $\dsigma\0{\square^j}\hookrightarrow \D^{n-k}$ is a smooth embedding (see section 7). We have two cases.\\

\noindent {\bf First case.} {\it We have that
 $DR_{\square^k}|\0{(t_nu)}v\ra v$, as $n\ra\infty$, when $v$ is tangent to $c\sigma\0{\square^j}$.}\\

\noindent This follows from an argument similar to the one given
in the proof of F.7 (recall that from F.5 we have  $R_{\square^k}(\sigma\0{\square^j})=\sigma\0{\square^j}$). This proves the first case.\\\\

Recall that $H$ is the change of variables
$H=\big( H_{\square^j}\big)^{-1}\circ H_{\square^k}$.
Let $u\in\dsigma^j$ and $v\in\R^{n-k}=T_u \D^{n-k}$.
We say that $v$ is a $X$-{\it fiber vector at $u$} if
$DH|\0{u}v=(z,0)\in
\R^{n-j}\times T_{H(u)}\dX_{\square^j}=T_{H(u)}(\D^{n-j}\times\dX_{\square^j})$, for some $z\in\R^{n-j}$. We write $v=v\0{z}^X$
(the reason for the upper index $X$ will be clear in a moment). Fixing $z$ we
obtain a constant vector field $(z,0)$, hence we obtain the corresponding
vector field $v\0{z}^X$ of $X$-fiber vectors on $\dsigma\0{\square^j}$. Thus 
$v\0{z}^X$ is characterized by $DH|\0{u}.v\0{z}^X(u)=(z,0)$.\\

Similarly, we can work on $K$ instead of $K_X$, and $h$ instead of $H$
and obtain vector fields of $\square$-{\it fiber vectors} $v\0{z}^\square$
on $\dsigma\0{\square^j}$ with the property that $Dh|\0{u}.v\0{z}^\square(u)=(z,0)
\in\R^{n-j}\times T_{h(u)}\square^j=T_{h(u)}(\D^{n-j}\times\square^j)$.\\

\noindent {\bf Claim F.9.} {\it We have $v\0{z}^X=v\0{z}^\square$.}\\

\noindent {\bf Proof of claim F9.} Fix $u$. Then the (hyperbolic) geodesic
$t\mapsto H\0{\square^k}(tu)$, and the straight segment 
$t\mapsto h\0{\square^k}^\bullet(tu)$ are both contained in $X\0{\square^n}$
and $\square^n$, respectively, for some $\square^n$. This together with
the fact that both $h\0{\square^k}^\bullet$ and $H\0{\square^k}$ use
the same link smoothing $h\0{\square^k}$ in their definition imply that
we can reduce our problem to the following setting. Consider $\R^n=\R^{n-k}\times
\R^k$, $\R^k\sbs\R^j$, the metrics $\sigma\0{\R^n}$ and $\sigma\0{\HH^n}=
\sigma\0{\HH^{n-k}}+cosh^2(r)\sigma\0{\HH^{k}}$, and $z$ (a constant vector field) perpendicular to $\R^j$. Here $r$ is the distance to $\HH^k$. In this case $h\0{\square^k}^\bullet$
corresponds to the perpendicular (to $\R^k$) exponential map from a point $p\in \R^k$ 
and $H\0{\square^k}$ correspond to
the perpendicular (to $\HH^k=(\R^k,\sigma\0{\HH^k})$) exponential map from $p$. 
The former exponential is just the inclusion and the latter exponential is done with respect to $\sigma\0{\HH^n}$. But in this setting
these two exponentials coincide, hence the preimage of $z$ by them also
coincide. This proves the claim.\\

Given $z$ as above we write $v\0{z}$  to denote $v\0{z}^X$ and $v\0{z}^\square$
and we say that $v=v\0{z}$ is  a {\it fiber vector}.
The following statement can be easily verified in the euclidean case (i.e for
$v\0{z}^\square$).

\begin{equation*}
v\0{z}(tu)\,\,=\,\,v\0{z}(u)
\tag{8}
\end{equation*}\vspace{.2in}

\noindent {\bf Second case.} {\it 
If $v$ is fiber vector at  $u$, then 
$DR_{\square^k}|\0{(t\0{n}u)}v\ra v$, as $n\ra\infty$.}\\

\noindent We have $v=v\0{z}(u)$, hence the image of $v$ by $DH|\0{(t\0{n} u)}$ is  $(z,0)$. Using formula
(7), the inductive hypothesis and the fact that $(0,q\0{t})=H(tu)\in \{ 0\}\times\dX_{\square^j}$
we see that  for any $z'\in\R^{n-j}$ we have

$$
DG_{\square^j}|\0{(H(t_n u))}(z',0)=(z',0)
$$

\noindent This together with the equation (6) imply

$$
DG_{\square^k}|\0{(t_nu)}(v,0)=Dh^{-1}\0{G\0{\square^j}(H(t_nu))}.DH\0{t\0{n} u}.v
$$

\noindent Therefore, from (8) and claim 9 we get

$$
DG_{\square^k}|\0{(t_nu)}(v,0)=DG_{\square^k}|\0{(t_nu)}(v\0{z}(u),0)=
v\0{z}\big( G_{\square^k}(t_nu)\big)
$$


\noindent  This together the fact that $h$ is a cone map imply
that to prove case 2 it is enough to prove that the directions of $G_{\square^k}(t_nu)=h^{-1}\circ G_{\square^j}\circ H(t_nu)$ tend to $u$, as $n\ra\infty$. 
But this last statement is implied by $$lim_{n\ra\infty}\frac{G_{\square^k}(t_nu)}{t_n}=DG_{\square^k}|\0{(0)}.u=\Big( DR_{\square^k}|\0{0}.u\,,\, 0 \Big)=u$$

\noindent which follows from equation (4). (The fact that the second coordinate
of the third term is 0 follows from (7) and the fact that the function $T_\square$
is even on the variable $t$.)
This proves the second case, step 2, claim F.8 and concludes the second part.\\

\vspace{.2in}

\noindent {\bf THIRD PART.} {\it The maps $H'_{\square}$ are $C^1$-embeddings.}\\

\noindent{\bf Proof.} Again by induction. This is true for
$H'_{\square^n}:\dX_{\square^n}\hookrightarrow K'_X$, which can be
identified with $P|\0{\dsquare^n}$ (see 9.2.1).
Assume that  the $H'_{\square^j}$ are $C^1$-embeddings for $j>k$. Fix a $\square^k$.
Using the argument used in the first and second parts we get that 
$H'_{\square^k}$ is a $C^1$-embedding outside $\dX_{\square^k}
=\{0\}\times\dX_{\square^k}$.
From the second part we see that the derivative $DH'_{\square^k}$
maps non-zero vectors $v$ at $\dX_{\square^k}$ in the $\D^{n-k}$ direction to non-zero vectors (see equation (4)).
As mentioned in the first part the map $q\0{X}\circ H'_{\square^k}|\0{\dX_{\square^k}}$ can be identified with $P|\0{\dX\0{\square^k}}$, which
is a diffeomorphism (see 9.2.1). Hence $DH'_{\square^k}$ maps
non-zero vectors $w$ in the $X_{\square^k}$ direction to non-zero vectors.
Moreover $DH'_{\square^k}.v$ is perpendicular to $DH'_{\square^k}.w$.
It follows that $H'_{\square^k}$ is an embedding on $\{0\}\times\dX_{\square^k}$.
This concludes the third part and completes the proof of proposition 9.3.2.

\vspace{1in}

\noindent {\bf \large  Appendix G. Cubification of Manifolds.}\\

We need the following lemma in section 10.\\

\noindent {\bf Lemma G.1.} {\it Every smooth manifold admits a smooth cubification.}\\

\noindent{\bf Proof.} Since every smooth manifold admits a smooth triangulation
it is enough to give a cubification of the canonical simplex $\Delta^{n}=\{(x_1,...x_{n+1})\,,\, \Sigma x_i=1\,,\, x_i\geq 0\}$ that is equivariant by its isometry group $S_{n+1}$,
which is the permutation group of $\{1,...,n+1\}$. To do this let $F^n=\{(x_1,...x_{n+1})\in\square^{n+1}\,,\,
 x_i=1$ for some $i\,\}$. Then $F^n$ is a cubical complex formed by the proper faces
of $\square^{n+1}$ that do not contain the origin. Let $P$ the the $n$-plane containing
$\Delta^n$ and $T: \bar{\R^{n+1}}_+-\{0\}\ra P$ the projection using rays from the
origin. Then $T(F^n)=\Delta^n$
and $T$ is $S_{n+1}$-equivariant. Moreover we have that restriction of $T$
to a cube in $F^n$ is a smooth embedding. This proves the lemma.\\

\noindent {\bf Remark.} Of course the canonical cubification of the simplex
given in the proof above gives a simple way of ``cubifying" any simplicial
complex. We use this observation in section 11.


Pedro Ontaneda

SUNY, Binghamton, N.Y., 13902, U.S.A.

\end{document}